\documentclass{article}
\usepackage{amssymb,amsmath,amsthm,tikz,multirow,nccrules,float,pst-solides3d,graphicx,makecell,caption,yhmath,bm}
\usepackage[hidelinks]{hyperref}
\usepackage{marvosym}
\usepackage{authblk}
\usepackage{abstract}
\renewcommand{\abstractname}{}
\usetikzlibrary{arrows,calc}

\usepackage[twoside]{geometry}
\title{Tilings of the sphere by congruent regular triangles and congruent rhombi$^{\ast}$}

\author{Qi Yuan$^{1}$, Erxiao Wang$^{2}$}

\allowdisplaybreaks[4]  

\newcommand{\sub}{\subset}
\newcommand{\pa}{\partial}
\newcommand{\mb}{\mathbf}
\newcommand{\mc}{\mathcal}
\newcommand{\bb}{\mathbb}
\newcommand\aaa{\alpha}
\newcommand\bbb{\beta}
\newcommand\ccc{\gamma}
\newcommand\ddd{\delta}
\newcommand{\arccot}{\mathrm{arccot}\,}
\newcommand{\tabincell}[2]{\begin{tabular}{@{}#1@{}}#2\end{tabular}}

\newcommand{\dash}{\hspace{0.1em}\dashrule{0.7}{2.4 1 2.4 1 2.4}\hspace{0.1em}} 
\newcommand{\thin}{\hspace{0.1em}\rule{0.7pt}{0.8em}\hspace{0.1em}}
\newcommand{\thick}{\hspace{0.1em}\rule{1.5pt}{0.8em}\hspace{0.1em}}

\newcommand{\tri}{\vartriangle}
\newcommand{\dia}{\lozenge}

\newcommand{\bo}{\boldsymbol}

\newtheorem{theorem}{Theorem}
\newtheorem{lemma}[theorem]{Lemma}

\newtheorem{proposition}[theorem]{Proposition}
\newtheorem*{theorem*}{Theorem}

\theoremstyle{definition}
\newtheorem*{definition*}{Definition}
\newtheorem*{case*}{Case}
\newtheorem*{subcase*}{Subcase}

\theoremstyle{remark}

\numberwithin{equation}{section}

\def\nd{\noindent}

\begin{document}
	
	\date{}
	\maketitle

	\renewcommand{\thefootnote}{\fnsymbol{footnote}}
	
	\footnotetext{\hspace*{-5mm} \begin{tabular}{@{}r@{}p{13.4cm}}
			$^1$ & School of Mathematical Sciences, Zhejiang Normal University, Jinhua, Zhejiang Province, 321004, China. E-mail: qiyuan@zjnu.edu.cn\\
			$^2$ & Corresponding author. School of Mathematical Sciences, Zhejiang Normal University, Jinhua, Zhejiang Province, 321004, China. E-mail: wang.eric@zjnu.edu.cn\\
			$^{\ast}$ & Research was supported by Key projects of Zhejiang Natural Science Foundation No. LZ22A010003 and ZJNU Shuang-Long Distinguished Professorship Fund No. YS304319159.
	\end{tabular}}
	
	\renewcommand{\thefootnote}{\arabic{footnote}}

	\begin{onecolabstract}
\nd{}{\bf Abstract} All edge-to-edge tilings of the sphere by congruent regular triangles and congruent rhombi are classified as: (1) a $1$-parameter family of protosets each admitting a unique $(2a^3,3a^4)$-tiling like a triangular prism; (2) a $1$-parameter family of protosets each admitting 2 different $(8a^3,6a^4)$-tilings like a cuboctahedron and a triangular orthobicupola respectively; (3) a sequence of protosets each admitting a unique $(2a^3,(6n-3)a^4)$-tiling like a generalized anti-triangular prism for each $n\ge3$; (4) 26 sporadic protosets, among which nineteen admit a unique tiling, one admits 3 different tilings, one admits 5 different tilings, three admit 2 different tilings, two admit too many tilings to count. The moduli of parameterized tilings and all geometric data are provided.
		
		\vskip 4.5mm
		
\nd{}{\bf Keywords} spherical tiling, prototile, protoset, regular triangle, rhombus, classification.\\
\nd{}{\bf 2000 MR Subject Classification } 52C20, 05B45
\end{onecolabstract}

	\section{Introduction}
	 In recent years, all edge-to-edge monohedral tilings of the sphere by congruent simple polygons (assuming the degree of all vertices $\ge 3$) have been fully classified after many authors' efforts, see \cite{so,ua,sa,wy1,wy2,wy3,awy,lqwx1,lw,lqwx2,cly}. Multihedral tilings of the sphere by regular polygons have just been classified in both edge-to-edge and non-edge-to-edge cases, see \cite{aehj}. However, multihedral tilings of the sphere by general polygons are rarely studied.
	
	In this paper, we will start the study of edge-to-edge dihedral tilings of the sphere by congruent triangles and congruent quadrilaterals, and fully classify the equilateral case. In other words, we classify edge-to-edge tilings of the sphere by congruent regular triangles and congruent rhombi (see Fig.\,\ref{t,q}) with edge lengths $a$, and all vertices have degree
$\ge3$. We will simply call such tilings $(a^3,a^4)$-tilings. We also denote the angle of the regular triangle by $\aaa$ and denote the angles of the rhombus by $\bbb,\ccc$. Without loss of generality, we assume that $\beta\ge\gamma$. 

	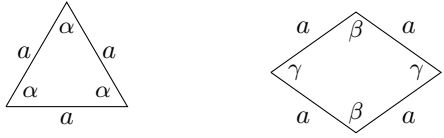
\begin{figure}[htp]
		\centering
		\begin{tikzpicture}[>=latex,scale=0.8]       
		\draw (0,0) -- (2,0) -- (1,{sqrt(3)}) -- (0,0); 
		\node at (0.4,0.25) {\small $\aaa$};
		\node at (1.6,0.25) {\small $\aaa$};
		\node at (1,1.3) {\small $\aaa$};
        \node at (0.3,0.9) {$a$};
        \node at (1.7,0.9) {$a$};
        \node at (1,-0.2) {$a$};
		\end{tikzpicture}\hspace{50pt}
		\begin{tikzpicture}[>=latex,scale=0.8]       
		\draw (-1.4,0) -- (0,1) -- (1.4,0)  -- (0,-1) -- (-1.4,0);
		\node at (0,0.68) {\small $\bbb$};
		\node at (0,-0.68) {\small $\bbb$};
		\node at (1,0) {\small $\ccc$};
		\node at (-1,0) {\small $\ccc$};
        \node at (0.87,0.74) {$a$};
        \node at (0.87,-0.74) {$a$};
        \node at (-0.87,0.74) {$a$};
        \node at (-0.87,-0.74) {$a$};
		\end{tikzpicture}
		\caption{The regular triangle and the rhombus.}
		\label{t,q}
	\end{figure}

We use $\aaa^k\bbb^l\ccc^m$ to mean a vertex having $k$ copies of $\aaa$, $l$ copies of $\bbb$, etc. The {\em anglewise vertex combination}, abbreviatled as \textbf{AVC}, is the collection of all vertices in a tiling. Then
the notation $T(2a^3,3a^4;6\aaa\bbb\ccc)$ means the tiling has exactly 2 regular triangles, 3 rhombi and 6 vertices $\aaa\bbb\ccc$, and is uniquely determined by them. In general there may exist several different tilings with
the same set of vertices. Then we use $\{8a^3,6a^4;12\aaa^2\bbb\ccc:2\}$ to mean that there are 2 different tilings.

The icosahedron naturally gives a monohedral tiling $T(20a^3;12\aaa^5)$, in which any two adjacent regular triangles can be merged into a rhombus $i.e.$ $\bbb=2\aaa=2\ccc$ that produces some $(a^3,a^4)$-tilings. We will call all such $(a^3,a^4)$-tilings to be of the \textbf{icosahedral type}. We also have the \textbf{octahedral type} (the tetrahedral type is dismissed since it has  vertices of degree $2$). 

\begin{figure}[htp]
		\centering
		\begin{tikzpicture}[>=latex,scale=0.8]       
		\draw (0,0) -- (1,1) -- (2,0) -- (2,-1)--(0,-1)--(0,0)--(2,0) (1,1)--(1,-1); 
		\node at (0.8,0.15) {\small $\aaa$};
		\node at (0.8,0.55) {\small $\aaa$};
		\node at (0.35,0.15) {\small $\aaa$};
		\node at (1.2,0.15) {\small $\aaa$};
		\node at (1.2,0.55) {\small $\aaa$};
		\node at (1.65,0.15) {\small $\aaa$};
		\node at (0.8,-0.25) {\small $\bbb$};
		\node at (1.2,-0.25) {\small $\bbb$};
		\node at (0.2,-0.8) {\small $\bbb$};
		\node at (1.8,-0.8) {\small $\bbb$};
		\node at (0.2,-0.2) {\small $\ccc$};
		\node at (1.8,-0.2) {\small $\ccc$};
		\node at (0.8,-0.8) {\small $\ccc$};
		\node at (1.2,-0.8) {\small $\ccc$};
		\end{tikzpicture}\hspace{50pt}
		\begin{tikzpicture}[>=latex,scale=0.8]       
		\draw (0,0) -- (1,1) -- (2,1)  -- (2,0) -- (1,-1)--(0,-1)--(0,0)--(2,0) (1,1)--(1,-1);
		\node at (0.8,0.15) {\small $\aaa$};
		\node at (0.8,0.55) {\small $\aaa$};
		\node at (0.35,0.15) {\small $\aaa$};
		\node at (1.2,-0.15) {\small $\aaa$};
		\node at (1.2,-0.55) {\small $\aaa$};
		\node at (1.65,-0.15) {\small $\aaa$};
		\node at (0.8,-0.25) {\small $\bbb$};
		\node at (1.2,0.25) {\small $\bbb$};
		\node at (0.2,-0.8) {\small $\bbb$};
		\node at (1.8,0.8) {\small $\bbb$};
		\node at (0.2,-0.2) {\small $\ccc$};
		\node at (1.8,0.2) {\small $\ccc$};
		\node at (0.8,-0.8) {\small $\ccc$};
		\node at (1.2,0.8) {\small $\ccc$};
		\end{tikzpicture}
		\caption{ Different Adjacent Angle Deductions of $\aaa^2\bbb^2$.}
		\label{aad}
	\end{figure}
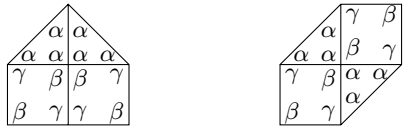

The very useful tool \textit{adjacent angle deduction} (abbreviated as \textbf{AAD}) has been introduced in \cite[Section 2.5]{wy1}. We give a quick review here using Figure \ref{aad}. Let ``$\thin$" denote an $a$-edge. Then we indicate the arrangements of angles and edges by denoting the vertices as $\thin\aaa\thin\aaa\thin\bbb\thin\bbb\thin$ or $\thin\aaa\thin\bbb\thin\aaa\thin\bbb\thin$. The notation can be reversed, such as $\thin\aaa\thin\aaa\thin\bbb\thin\bbb\thin=\thin\bbb\thin\bbb\thin\aaa\thin\aaa\thin$; and it
can be rotated, such as $\thin\aaa\thin\aaa\thin\bbb\thin\bbb\thin=\thin\aaa\thin\bbb\thin\bbb\thin\aaa\thin$. We also denote the first vertex in Figure \ref{aad} as $\bbb\thin\bbb\cdots,\thin\bbb\thin\bbb\thin\cdots,\aaa\thin\bbb\cdots,\thin\aaa\thin\aaa\thin\bbb\thin\cdots$, and denote the consecutive angle segments as $\bbb\thin\bbb,\thin\bbb\thin\bbb\thin,\aaa\thin\bbb,\thin\aaa\thin\aaa\thin\bbb\thin$.
	
We write $^\lambda\theta^\mu$ to mean $\lambda,\mu$ are the two angles adjacent to $\theta$ in a tile. The first picture has the AAD $\thin^\aaa\aaa^\aaa\thin^\aaa\aaa^\aaa\thin^\ccc\bbb^\ccc\thin^\ccc\bbb^\ccc\thin$, which gives $\thin\aaa\thin\aaa\thin\cdots,\thin\aaa\thin\ccc\thin\cdots,$ $\thin\ccc\thin\ccc\thin\cdots$; the second has the AAD $\thin^\aaa\aaa^\aaa\thin^\ccc\bbb^\ccc\thin^\aaa\aaa^\aaa\thin^\ccc\bbb^\ccc\thin$, which gives $\thin\aaa\thin\ccc\thin\cdots$.

Let $f_\tri,f_\dia$ be the numbers of triangle tiles and quadrilateral tiles; $S_\tri,S_\dia$ be the areas of triangle tiles and quadrilateral tiles. We have $f_\tri S_\tri + f_\dia S_\dia =4\pi$ (the radius of the sphere is always 1).

	\begin{theorem*}
	All $(a^3,a^4)$-tilings are classified as:
	\begin{enumerate}
		\item A $1$-parameter family of protosets each admitting a unique tiling $T(2a^3,3a^4;6\aaa\bbb\ccc)$ like a triangular prism. $($See the 1st of Fig.\,\ref{Fig2}$)$;
			
		\item A $1$-parameter family of protosets each admitting 2 different tilings $\{8a^3,6a^4;12\aaa^2\bbb\ccc:2\}$ like a cuboctahedron and a triangular orthobicupola respectively. All $\aaa^2\bbb\ccc$ vertices are $\thin\aaa\thin\bbb\thin\aaa\thin\ccc\thin$ in the first tiling. By flipping half of the first, we get the second tiling which has six $\thin\aaa\thin\aaa\thin\bbb\thin\ccc\thin$ near the equator. $($See the 2nd and 3rd of Fig.\,\ref{Fig2}$)$;
 
		\item A sequence of protosets each admitting a unique tiling $T(2a^3,(6n-3)a^4;6\aaa\bbb\ccc^n,(6n-6)\bbb^2\ccc)$ like a generalized anti-triangular prism for each $n\ge3$. $($See the 4th, 5th and 6th of Fig.\,\ref{Fig2} when $n=3,4,5)$;

	\begin{figure}[htp]
		\centering
		\includegraphics[scale=0.195]{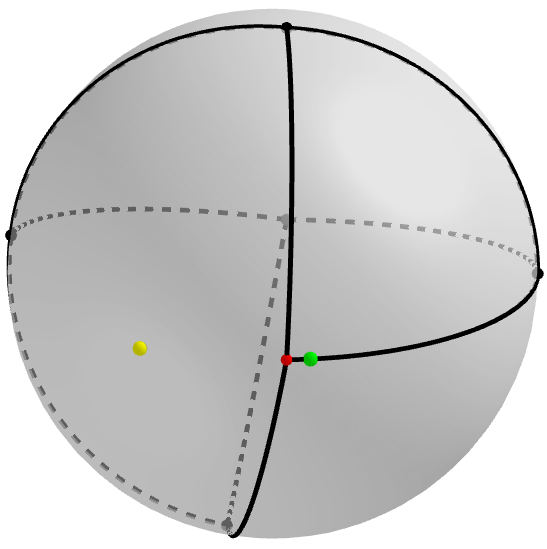}
		\includegraphics[scale=0.24]{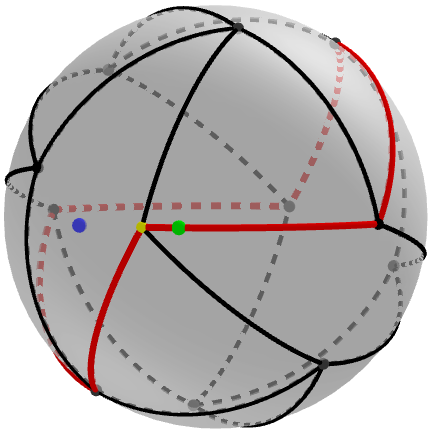}
		\includegraphics[scale=0.29]{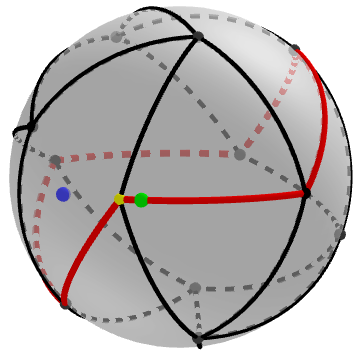}
		\includegraphics[scale=0.23]{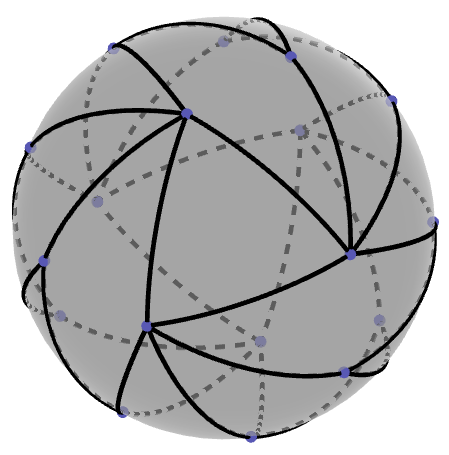}
		\includegraphics[scale=0.24]{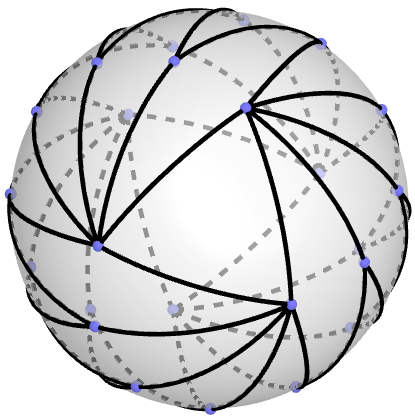}
		\includegraphics[scale=0.234]{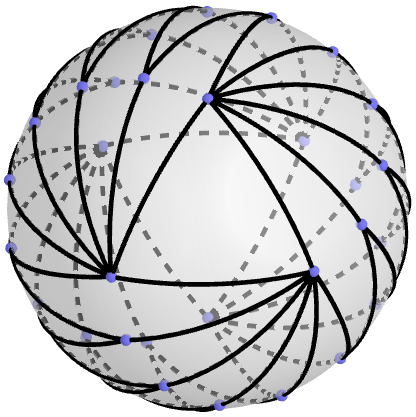}
		\caption{Some 3D pictures in the first, second and third classes.}
		\label{Fig2}
	\end{figure}

		\item  Twenty-six sporadic protosets as listed in Table \ref{t}, among which nineteen admit a unique tiling, one admits 3 different tilings, one admits 5 different tilings, three admit 2 different tilings, two admit too many tilings to count.

\begin{table}[H]
\centering
\begin{tabular}{c|c|c|c}
$(f_\tri,f_\dia)$ & Page & $\aaa,\bbb,\ccc;\,a$ & all vertices \& tilings  \\
\hline    
(4,1) & \multirow{7}{*}{\pageref{b=pi}} & $\frac{1}{2},1,\ccc=\bbb;\,\frac{1}{2}$ & $4\aaa^2\bbb,1\aaa^4$\\
\cline{1-1}\cline{3-4}
(6,3) & \multirow{7}{*}{} & $\frac{1}{2},1,\frac{1}{6};\,\frac{1}{2}$ & $6\aaa^2\bbb,2\aaa^3\ccc^3$\\
\cline{1-1}\cline{3-4}       
\tabincell{c}{(4,4)\\ (6,2)} & \multirow{7}{*}{} & $\frac{1}{2},1,\frac{1}{4};\,\frac{1}{2}$ & \tabincell{c}{$4\aaa\bbb\ccc^2,4\aaa^2\bbb$\\ $4\aaa^2\bbb,2\aaa^3\ccc^2,1\aaa^4$}\\
\cline{1-1}\cline{3-4}
\tabincell{c}{(4,2)\\ (4,2)\\ (6,1)} & \multirow{7}{*}{} & $\frac{1}{2},1,\frac{1}{2};\,\frac{1}{2}$ & \tabincell{c}{$4\aaa^2\bbb,2\aaa^2\ccc^2$\\ $2\aaa\bbb\ccc,2\aaa^2\bbb,2\aaa^3\ccc$\\ $2\aaa^2\bbb,2\aaa^3\ccc,2\aaa^4$}\\
\hline    \hline 
(8,2) & \multirow{3}{*}{\pageref{ab3}} & 0.4335,0.6992,$\ccc=\bbb$;\,0.4158 & $8\aaa^3\bbb$ \\
\cline{1-1}\cline{3-4}
(32,6) & \multirow{3}{*}{} & 0.3621,0.5513,$\ccc=\bbb$;\,0.2427 & $24\aaa^4\bbb$ \\
\cline{1-1}\cline{3-4}
(8,18) & \multirow{3}{*}{} & 0.3596,0.5467,$\ccc=\bbb$;\,0.2326 & $24\aaa\bbb^3$\ :\ 2 \\
\hline \hline     
(4,3) & \pageref{rt} & $\frac{4}{9},\frac{7}{9},\frac{2}{3};\,0.4326$ & $3\aaa\bbb^2,3\aaa^3\ccc,1\ccc^3$\\
\hline  
(4,4) & \pageref{4,4} & 0.4296,0.7851,0.5703;\,0.4094 & $4\aaa\bbb^2,4\aaa^2\ccc^2$\\
\hline  
(8,3) & \pageref{8,3} & 0.4195,0.7412,0.5804;\,0.3918 & $6\aaa^3\bbb,3\aaa^2\ccc^2$\\
\hline  
(4,12) & \pageref{4,12} & $0.3754,\frac{2}{3},0.4789$;\,0.2884 & $12\aaa\bbb\ccc^2,4\bbb^3$\\
\hline
(8,12) & \pageref{8,12} & $0.3701,0.6298,\frac{1}{2}$;\,0.2716 & $6\aaa^2\bbb^2,12\aaa\bbb\ccc^2$\\
\hline 
(20,6) & \pageref{20,6} & 0.3807,0.8578,0.2385;\,0.3040 & $12\aaa^3\bbb,6\aaa^4\ccc^2$\\
\hline   
(8,24) & \pageref{8,24} & $0.3541,0.5729,\frac{1}{2}$;\,0.2082 & $24\aaa\bbb^2\ccc,6\ccc^4$\\
\hline
(20,12) & \pageref{20,12.1} & 0.3614,0.6385,0.4577;\,0.2401 & $12\aaa^2\bbb^2,12\aaa^3\ccc^2$\\
\hline 
(20,12) & \pageref{20,12.2} & 0.3733,0.8798,0.1866;\,0.2820 & $12\aaa^3\bbb,12\aaa^2\bbb\ccc^2$\\
\hline 
(20,24) & \pageref{20,24.1} & 0.3510,0.5877,0.4734;\,0.1927 & $24\aaa\bbb^2\ccc,12\aaa^3\ccc^2$\\
\hline   
(44,12) & \pageref{44,12} & 0.3590,0.9229,0.1024;\,0.2301 & $24\aaa^3\bbb,12\aaa^5\ccc^2$\\
\hline  
(20,60) & \pageref{20,60} & $0.3421,0.6289,\frac{2}{5}$;\,0.1379 & $60\aaa\bbb^2\ccc,12\ccc^5$\\
\hline 
(16,6) & \multirow{3}{*}{\pageref{16,6}} & 0.3861,0.8415,0.2805;\,0.3188 & $12\aaa^3\bbb,4\aaa^3\ccc^3$\\
\cline{1-1}\cline{3-4}
(32,12) & \multirow{3}{*}{} & 0.3650,0.9049,0.1349;\,0.2536 & $24\aaa^3\bbb,6\aaa^4\ccc^4$\\
\cline{1-1}\cline{3-4}
(80,30) & \multirow{3}{*}{} & 0.3480,0.9558,0.0519;\,0.1766 & $60\aaa^3\bbb,12\aaa^5\ccc^5$\\
\hline  \hline
(20,36) & \pageref{20,36} & 0.3465,0.6089,0.4356;\,0.1675 & $36\aaa\bbb^2\ccc,12\aaa^2\ccc^3$\ :\ 2\\
\hline   
(32,6) & \pageref{32,6} & 0.3686,0.8939,0.1566;\,0.2665 & $12\aaa^3\bbb,12\aaa^5\ccc$\ :\ 5\\
\hline   
\hline 
\tabincell{c}{$(20-2m,m)$\\ $1\le m \le 9$} & \pageref{ico} & $\frac{2}{5},\frac{4}{5},\frac{2}{5}$;\,0.3524 & \tabincell{c}{$\aaa\bbb^2,\bbb^2\ccc,\aaa^3\bbb$,\\ $\aaa^2\bbb\ccc,\aaa\bbb\ccc^2,\bbb\ccc^3$,\\ $\aaa^5,\aaa^4\ccc,\aaa^3\ccc^2$,\\$\aaa^2\ccc^3,\aaa\ccc^4,\ccc^5$\ :\ ?}\\
\hline 
(20,24) & \pageref{20,24.2} & 0.3579,0.8210,0.2315;\,0.2257 & \tabincell{c}{$(12+k)\aaa\bbb^2,k\aaa^3\ccc^4$,\\ $(24-2k)\aaa^2\bbb\ccc^2$\ :\ ?\\ $0\le k <12$}\\
\hline  
\end{tabular}
\caption{26 sporadic protosets.}
\label{t}
\end{table}
	\end{enumerate}
\end{theorem*} 
	
\begin{figure}[H]
		\centering
		\includegraphics[scale=0.3]{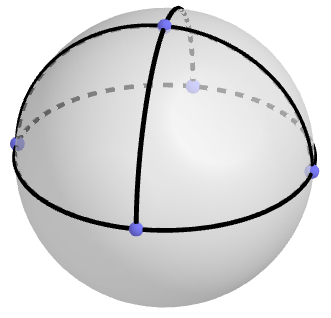}
		\includegraphics[scale=0.24]{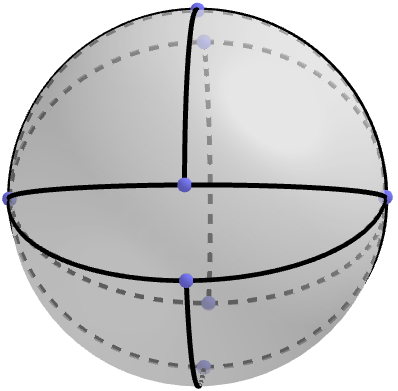}
		\includegraphics[scale=0.26]{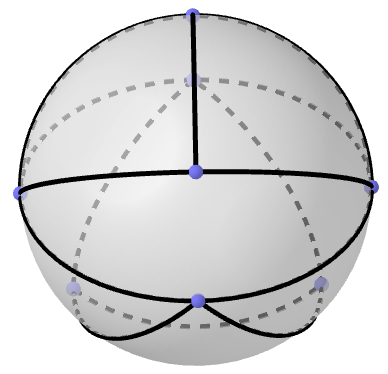}
		\includegraphics[scale=0.25]{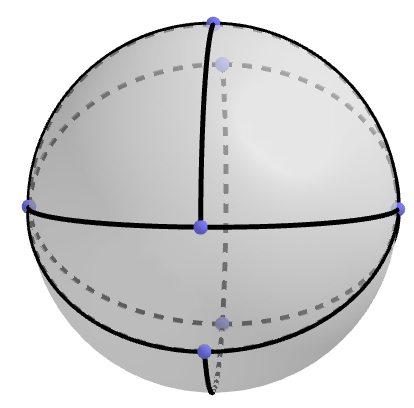}
		\includegraphics[scale=0.3]{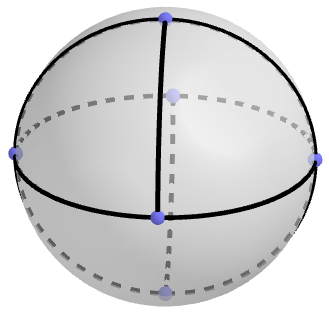}
		\includegraphics[scale=0.275]{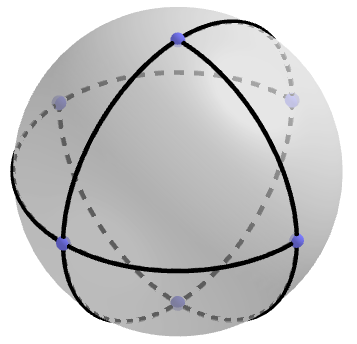}

		\includegraphics[scale=0.272]{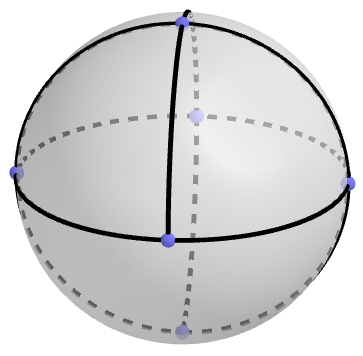}\hspace{2pt}
		\includegraphics[scale=0.34]{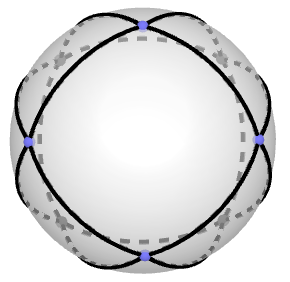}
		\includegraphics[scale=0.245]{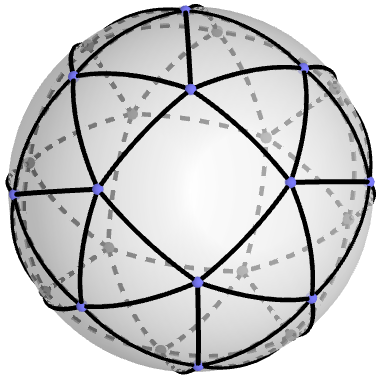}\hspace{4pt}
		\includegraphics[scale=0.25]{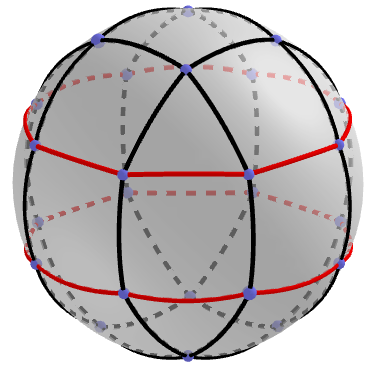}\hspace{4pt}
		\includegraphics[scale=0.14]{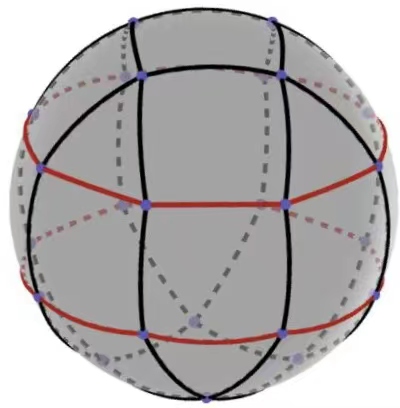}\hspace{4pt}
		\includegraphics[scale=0.29]{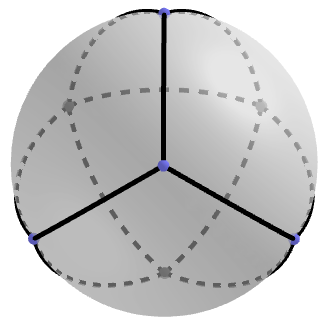}

		\includegraphics[scale=0.255]{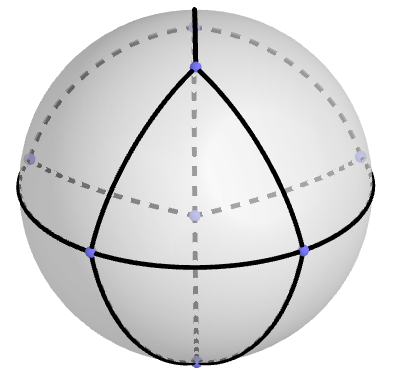}
		\includegraphics[scale=0.28]{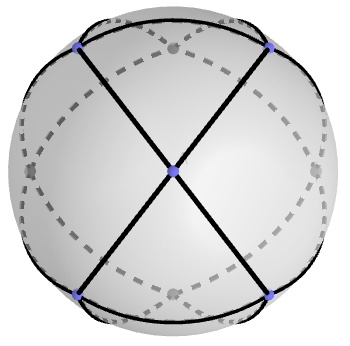}
		\includegraphics[scale=0.28]{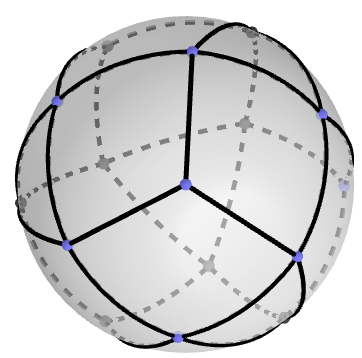}
		\includegraphics[scale=0.25]{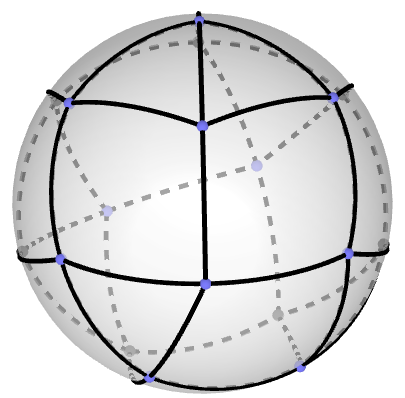}
		\includegraphics[scale=0.25]{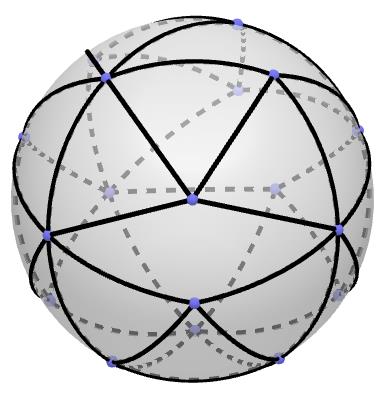}
		\includegraphics[scale=0.235]{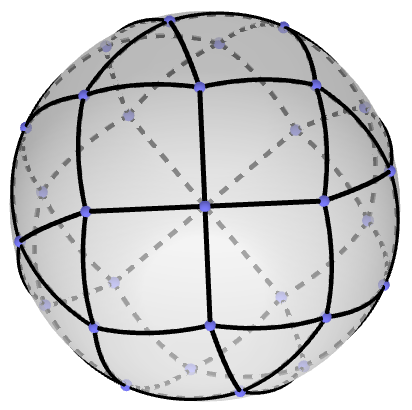}

		\includegraphics[scale=0.265]{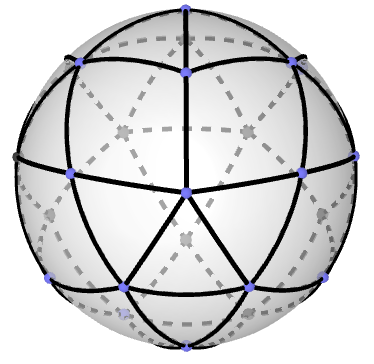}\hspace{2pt}
		\includegraphics[scale=0.225]{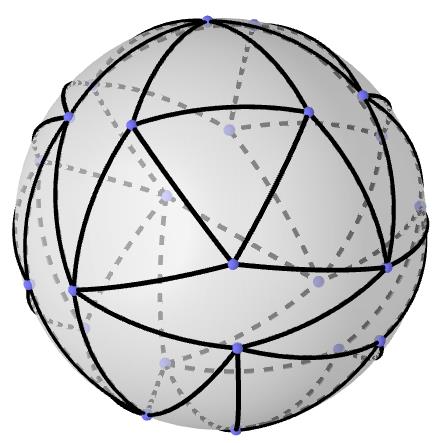}\hspace{2pt}
		\includegraphics[scale=0.17]{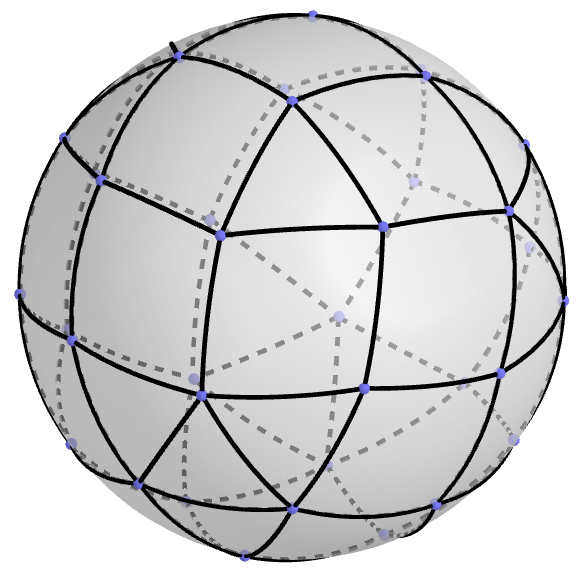}
		\includegraphics[scale=0.25]{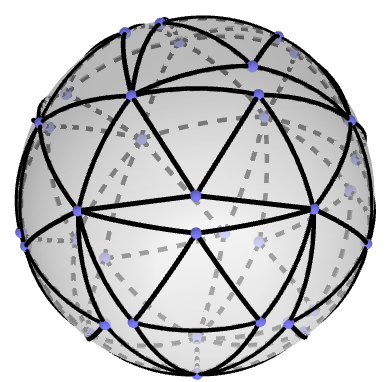}
		\includegraphics[scale=0.15]{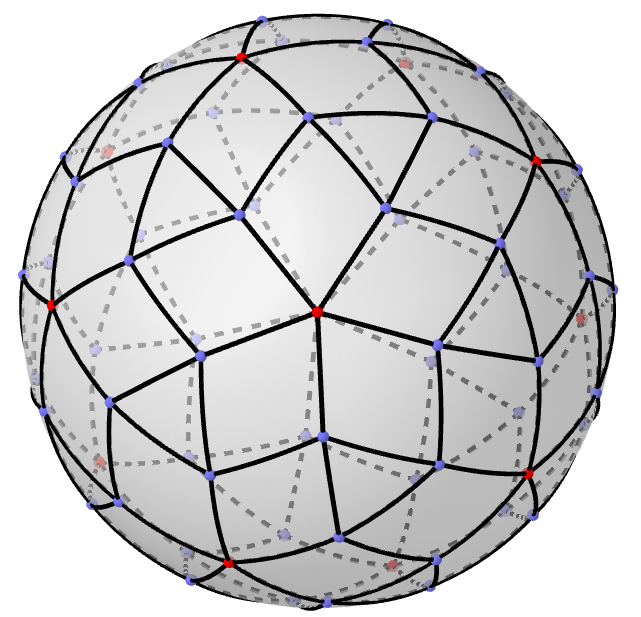}
		\includegraphics[scale=0.3]{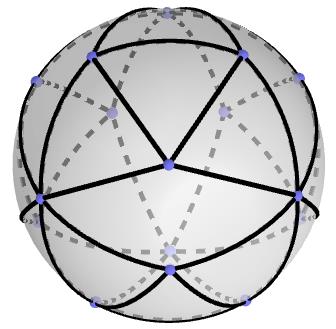}

		\includegraphics[scale=0.25]{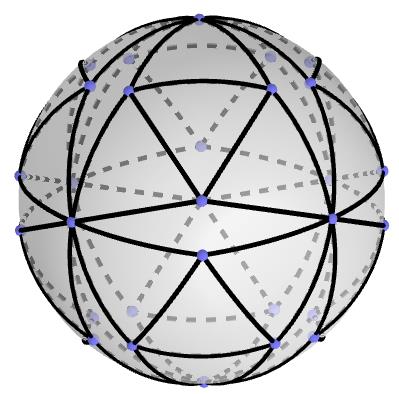}\hspace{3pt}
		\includegraphics[scale=0.17]{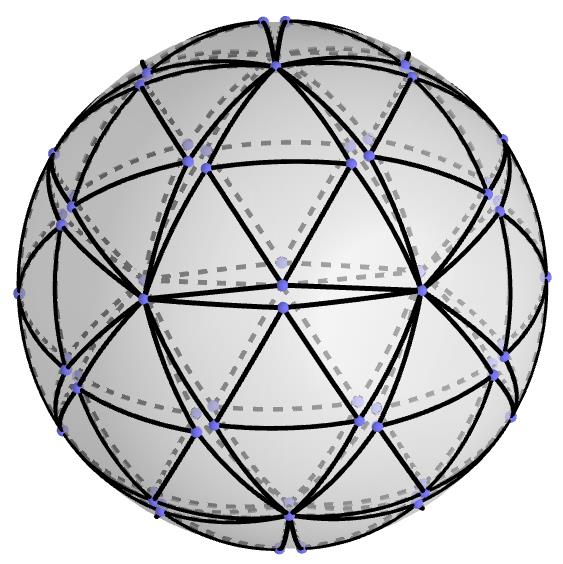}\hspace{2pt}
		\includegraphics[scale=0.175]{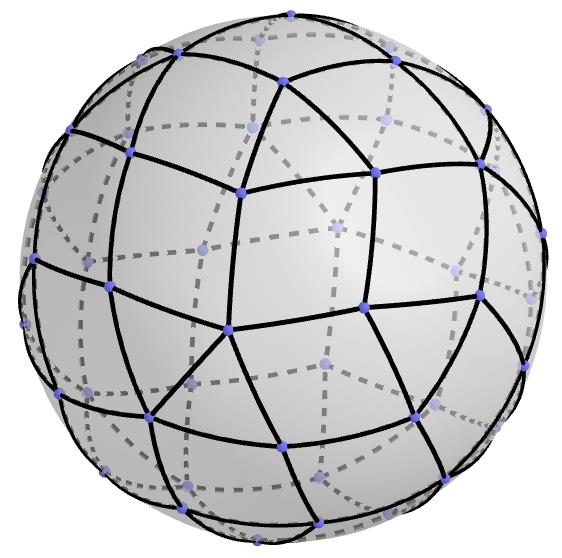}\hspace{2pt}
		\includegraphics[scale=0.235]{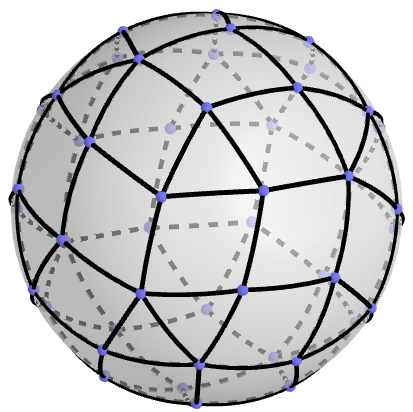}
		\includegraphics[scale=0.24]{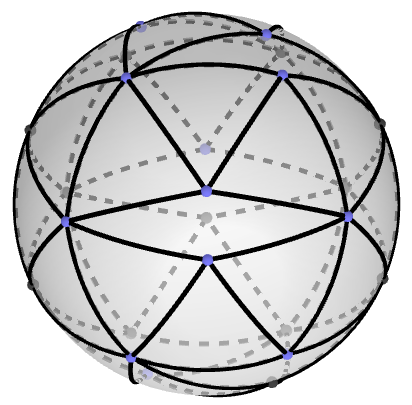}
		\includegraphics[scale=0.23]{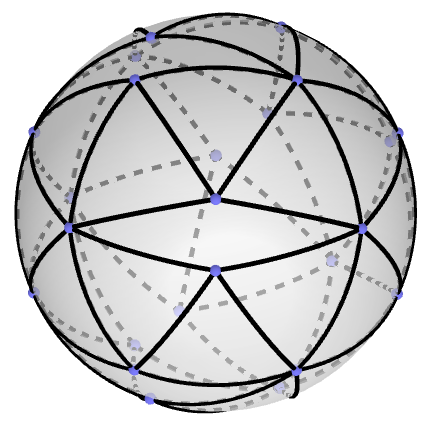}

		\includegraphics[scale=0.28]{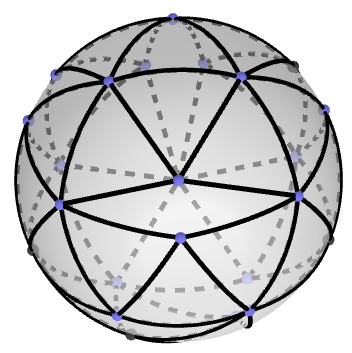}
		\includegraphics[scale=0.25]{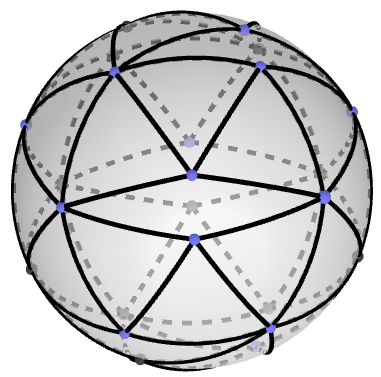}
		\includegraphics[scale=0.23]{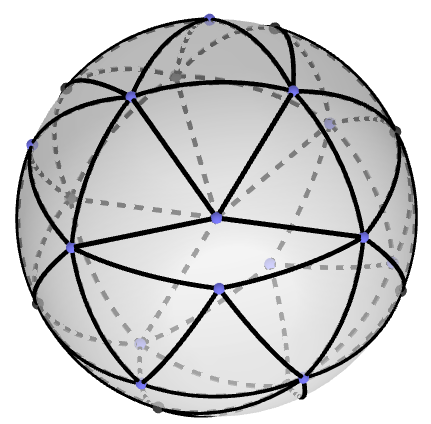}
		\caption{All 3D pictures in the fourth class except for two polymorphic ones.}
		\label{Fig3}
	\end{figure}

The notation ? means the protoset admits a large number of different tilings, and we have not figured them out yet. The other 3D pictures are shown in Figure \ref{Fig3}.

The numerical geometric data is listed in Table \ref{t}, where the angles and edge lengths are expressed in units of $\pi$, and the last column counts all vertices and also all tilings when they are not uniquely determined by the vertices. A rational fraction means the precise value. A decimal expression, such as $a=0.4158$, means an approximate value  $0.4158\pi < a < 0.4159\pi$. Moreover, the exact formulas are provided in the appendix.

{\bf Acknowledgement}\ \  We would like to thank junior student Zhihong Lin for very helpful discussions in Case $\bbb^2\ccc$. Thank two junior students Fangbin Chen and Nan Zhang for showing us how to draw Figure \ref{Fig2} and \ref{Fig3} using GeoGebra.

	\section{Basic Facts}\label{basic_facts}

	Let $v,e$ be the numbers of vertices and  edges. Let $v_k$ be the number of vertices of degree $k$. We have Euler's formula and basic counting equalities: 
	\begin{align*}
	2&=v-e+f_\tri +f_\dia , \\
	2e=3f_\tri +4f_\dia 
	&=\sum_{k=3}^{\infty}kv_k=3v_3+4v_4+5v_5+\cdots, \\
	v
	&=\sum_{k=3}^{\infty}v_k=v_3+v_4+v_5+\cdots.
	\end{align*}
	
    Then it is easy to derive   
	\begin{align}
    v &=2+\frac{1}{2}f_\tri +f_\dia , \label{eq1} \\
	3f_\tri +2f_\dia  &=12+2\sum_{k=4}^{\infty}(k-3)v_k=12+2v_4+4v_5+6v_6+\cdots, \label{eq2} \\
	v_3+f_\tri  &=8+\sum_{k=5}^{\infty}(k-4)v_k=8+v_5+2v_6+\cdots, \label{eq3} \\
    3v_3+2v_4+v_5 &=12+2f_\dia +\sum_{k=7}^{\infty}(k-6)v_k=12+2f_\dia +v_7+2v_8+\cdots. \label{eq4}
	\end{align}
	 
     We note that $f_\tri $ is even and $\ge2$ by (\ref{eq1}). The equality (\ref{eq2}) show that $f_\tri =2,f_\dia =3$ for $v_k=0\,(k\ge4)$, which is the simplest one. There must be some vertex of degree 3, 4 or 5 by (\ref{eq4}).

    We also have some geometric equations:

    \begin{align}
    \cos a &=\cos^2 a+\sin^2 a \cos \alpha  \nonumber \\
           &=\cot \alpha \cot \frac{\alpha}{2} \nonumber \\
           &=\cot \frac{\beta}{2}\cot \frac{\gamma}{2}. \label{tr3} 
    \end{align}
Then we deduce that
    \begin{align}
    (\cot \frac{\beta}{2}\cot \frac{\gamma}{2})^2 +[1-(\cot \frac{\beta}{2}\cot \frac{\gamma}{2})^2] \cos \alpha- \cot \frac{\beta}{2}\cot \frac{\gamma}{2}=0,   \label{tr4}\\
    \cot \frac{\beta}{2}\cot \frac{\gamma}{2}-\cot \alpha \cot \frac{\alpha}{2}=0.\label{tr5}
    \end{align}

		\begin{lemma}\label{AAD}
          The simple application of the AAD:
			\begin{itemize}
		\item There must be $\thin\aaa\thin\bbb\thin\cdots$ and $\thin\aaa\thin\ccc\thin\cdots$ in an $(a^3,a^4)$-tiling.
		\item If $\bbb^2\cdots$ or $\thin\bbb\thin\bbb\thin\cdots$ is not a vertex, then $\thin\ccc\thin\ccc\thin\cdots$ can never be a vertex.
	       \end{itemize}
	   \end{lemma}

		\begin{lemma}\label{angle}
		In an $(a^3,a^4)$-tiling, we have $\frac{1}{3}\pi<\aaa<\pi,\bbb>\frac{1}{2}\pi,\bbb+\ccc>\pi$ and $0<a<\frac{2}{3}\pi$.
	\end{lemma}

	\begin{proof}
Since $3\aaa=\pi+S_\tri$ and $2(\bbb+\ccc)=2\pi+S_\dia$, we have $\aaa>\frac{1}{3}\pi$ and $\bbb+\ccc>\pi$. By $\bbb\ge\ccc$, we get $\bbb>\frac{1}{2}\pi$.
    
		If $\aaa\ge \pi$, then $S_\tri\ge2\pi$ and $f_\tri S_\tri\ge4\pi$. This implies that $f_\dia S_\dia\le0$, a contradiction. So $\aaa<\pi$ $i.e.$ $\frac{1}{3}\pi<\aaa<\pi$.

If $\aaa\to \pi$, then $S_\tri\to 2\pi$, which imply the three vertices of the regular triangle are on the same equator and $a\to \frac{2}{3}\pi$. So $0<a<\frac{2}{3}\pi$.
	\end{proof}

We will call an $(a^3,a^4)$-tiling or its protoset \textbf{convex} when $\bbb<\pi$, \textbf{concave} when $\bbb>\pi$, and \textbf{degenerate} when $\bbb=\pi$.

		\begin{lemma}\label{type}
			In any $(a^3,a^4)$-tiling, we have
	\begin{itemize}
		\item $0<a<\frac{1}{2}\pi \Longleftrightarrow \frac{1}{3}\pi<\aaa<\frac{1}{2}\pi<\bbb<\pi$;
		\item $a=\frac{1}{2}\pi \Longleftrightarrow \aaa=\frac{1}{2}\pi,\bbb=\pi$;
		\item $\frac{1}{2}\pi<a<\frac{2}{3}\pi \Longleftrightarrow \frac{1}{2}\pi<\aaa<\pi<\bbb<2\pi$.
	\end{itemize}
	\end{lemma}

	\begin{figure}[htp]
		\centering
		\begin{tikzpicture}[>=latex,scale=0.8]       
		\draw (0,0) circle (2);   
	    \fill (0,0) circle (1pt);
		\node at (0.3,0) {\small $O$};
		\draw[dashed,color=red] (2,0) arc (0:180:2 and 0.4); 
		\draw[color=red] (2,0) arc (0:-180:2 and 0.4); 
		\draw (0,2) arc (142.66:182.42:3.65); 
		\draw (0,2) arc (37.34:-2.42:3.65); 
        \node at (0,2) {$\bullet$};
		\node at (0,2.3) {\small $N$};
        \node at (0.75,-0.4) {$\bullet$};
		\node at (1,-0.6) {\small $B$};
        \node at (-0.75,-0.4) {$\bullet$};
		\node at (-1,-0.6) {\small $A$};
		\node at (0,-0.65) {\small $l$};
		\end{tikzpicture}\hspace{30pt}
		\begin{tikzpicture}[>=latex,scale=0.8]       
		\draw (0,0) circle (2);   
	    \fill (0,0) circle (1pt);
		\node at (0.3,0) {\small $O$};
		\draw[dashed,color=red] (2,0) arc (0:180:2 and 0.4); 
		\draw[color=red] (2,0) arc (0:-180:2 and 0.4); 
		\draw[dashed] (1.9,0.66) arc (0:180:1.9 and 0.39); 
		\draw (1.9,0.66) arc (0:-180:1.9 and 0.39); 
		\draw (0,2) arc (138.11:177.18:2.77); 
		\draw (0,2) arc (41.89:2.82:2.77); 
		\draw[color=green] (0.7,0.3) arc (35.46:144.54:0.86); 
        \node at (0,2) {$\bullet$};
		\node at (0,2.3) {\small $N$};
		\node at (0,-0.65) {\small $l$};
		\node at (-1.2,0.65) {\small $l_1$};
        \node at (0.7,0.3) {$\bullet$};
		\node at (1,0.1) {\small $B$};
        \node at (-0.7,0.3) {$\bullet$};
		\node at (-1,0.1) {\small $A$};
		\end{tikzpicture}\hspace{30pt}
		\begin{tikzpicture}[>=latex,scale=0.8]       
		\draw (0,0) circle (2);   
	    \fill (0,0) circle (1pt);
		\node at (0.3,0) {\small $O$};
		\draw[dashed,color=red] (2,0) arc (0:180:2 and 0.4); 
		\draw[color=red] (2,0) arc (0:-180:2 and 0.4); 
		\draw[dashed] (1.9,-0.66) arc (0:180:1.9 and 0.39); 
		\draw (1.9,-0.66) arc (0:-180:1.9 and 0.39); 
		\draw (0,2) arc (136.02:202.26:2.84); 
		\draw (0,2) arc (43.98:-22.26:2.84); 
		\draw[color=green] (0.58,-1.03) arc (-37.24:-142.76:0.73); 
        \node at (0,2) {$\bullet$};
		\node at (0,2.3) {\small $N$};
		\node at (0,-0.65) {\small $l$};
		\node at (-1.2,-1.2) {\small $l_2$};
        \node at (0.58,-1.03) {$\bullet$};
		\node at (1,-0.75) {\small $B$};
        \node at (-0.58,-1.03) {$\bullet$};
		\node at (-1,-0.75) {\small $A$};
		\end{tikzpicture}
		\caption{Proof of Lemma \ref{type}.}
		\label{angle type}
	\end{figure}

	\begin{proof}
        Let $\bigtriangleup NAB$ be the spherical isosceles triangle with $\wideparen{NA}=\wideparen{NB}$, $N$ be the pole, $O$ be the centre of sphere, $l$ be the equator and $l_1,l_2$ be the latitude lines.

If $\wideparen{NA}=\frac{1}{2}\pi$, then $\wideparen{AB}$ is on $l$ and $\angle{NAB}=\frac{1}{2}\pi$. If $\wideparen{NA}<\frac{1}{2}\pi$, then $A,B$ are on $l_1$ and $\wideparen{AB}$ is above the $l_1$. Since $\wideparen{NA}\perp l_1$, we get $\angle{NAB}<\frac{1}{2}\pi$. Similarly if $\wideparen{NA}>\frac{1}{2}\pi$, then $\angle{NAB}>\frac{1}{2}\pi$.

By considering the base angle, we have the same things happening for $\wideparen{NA}$.
		Therefore we get this conclusion.
	\end{proof}

In fact, we see that $2\bbb+2\ccc-2\pi<4\ccc$ $i.e.$ $\bbb<\pi+\ccc$ in the third of Figure \ref{angle type}.

	\begin{lemma}\label{a3}
		In an $(a^3,a^4)$-tiling, $\aaa^3$ can never be a vertex.
	\end{lemma}

	\begin{proof}
		The vertex $\aaa^3$ means $S_\tri=\pi$ and $3\le f_\tri \le4$. Since $f_\tri $ is even, $f_\tri =4$ but $f_\dia S_\dia=0$, a contradiction. Therefore, $\aaa^3$ is not a vertex.
	\end{proof}

		\begin{lemma}\label{deg345}
		In a convex $(a^3,a^4)$-tiling with $\ccc<\bbb$, all vertices of degree 3, 4 or 5 are shown in Table \ref{v>}.
	\end{lemma}

\begin{table}[htp]
\centering
\begin{tabular}{|c|c|}
\hline 
degree & vertex \\
\hline \hline 
$3$
& $\aaa\bbb^2,\bbb^3,\aaa\bbb\ccc,\bbb^2\ccc,\aaa\ccc^2,\bbb\ccc^2,\ccc^3$ \\
\hline 
$4$
& $\aaa^3\bbb,\aaa^2\bbb^2,\aaa\bbb^3,\aaa^3\ccc,\aaa^2\bbb\ccc,\aaa\bbb^2\ccc,\aaa^2\ccc^2,\aaa\bbb\ccc^2,\aaa\ccc^3,\bbb\ccc^3,\ccc^4$ \\
\hline 
$5$
& $\aaa^5,\aaa^4\bbb,\aaa^4\ccc,\aaa^3\ccc^2,\aaa^2\bbb\ccc^2,\aaa^2\ccc^3,\aaa\bbb\ccc^3,\aaa\ccc^4,\bbb\ccc^4,\ccc^5$\\
\hline
\end{tabular}
\caption{Vertices of degree 3, 4, 5 for $\ccc<\bbb<\pi$.}
\label{v>}
\end{table}

	\begin{proof}
By Lemma \ref{angle} and \ref{type}, we have $\bbb+\ccc>\pi$, $\frac{\pi}{3}<\aaa<\frac{\pi}{2}<\bbb<\pi$, and $0<\ccc<\bbb$. If $\aaa^n\bbb^m\ccc^l$ is a vertex, then we get $0\le n\le5 ,0\le m\le3,l\ge0, n+m+l\ge3$. 

        By Lemma \ref{a3} and $2\aaa+\ccc<2\aaa+\bbb<2\pi$, we know $\aaa^3,\aaa^2\bbb$ or $\aaa^2\ccc$ can never be a vertex.

        By $3\bbb+\ccc>2\bbb+2\ccc>2\pi$ and $4\aaa<2\pi$, we know $\bbb^3\ccc,\bbb^2\ccc^2$ or $\aaa^4$ can never be a vertex.
 
       When the degree is 5, we know $\bbb^3\cdots$ can never be a vertex by $3\bbb+\ccc>2\pi$ and $2\aaa+3\bbb>2\pi$. By $2\bbb+2\ccc>2\pi$, $2\aaa+2\bbb+\ccc>2\pi$ and $3\aaa+2\bbb>2\pi$, we know $\bbb^2\cdots$ can never be a vertex. By $3\aaa+\bbb+\ccc>2\pi$, we know $\aaa^3\bbb\ccc$ can never be a vertex.

Therefore, we get all vertices of degree 3, 4 or 5 as listed in Table \ref{v>}.
	\end{proof}

		\begin{lemma}\label{b(r)}
		In a convex $(a^3,a^4)$-tiling, if $\bbb=\ccc$, then $\aaa<\bbb<2\aaa$.
	\end{lemma}

	\begin{figure}[H]
		\centering
		\includegraphics[scale=0.34]{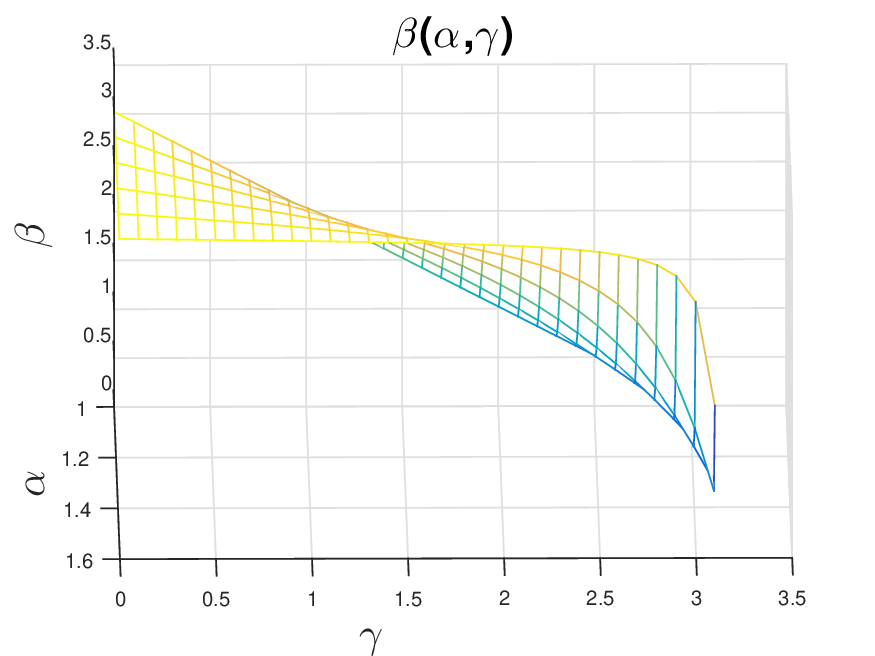}
		\includegraphics[scale=0.34]{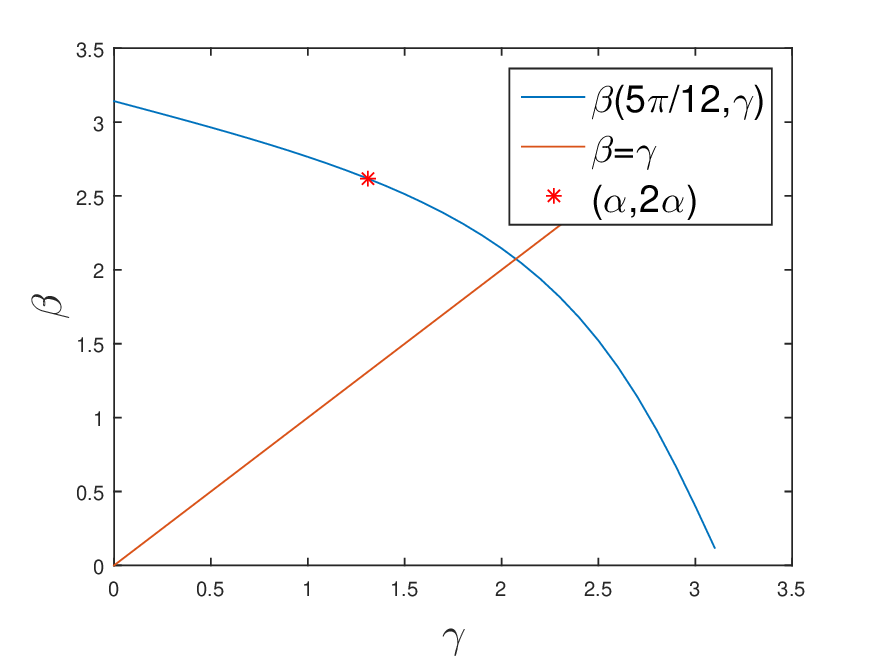}
		\caption{Proof of Lemma \ref{b(r)}.}
		\label{b(r)1}
	\end{figure}

	\begin{proof}
        By Lemma \ref{type}, we have $0<a<\frac{1}{2}\pi,\frac{1}{3}\pi<\aaa<\frac{1}{2}\pi$, and $0<\ccc\le\bbb<\pi$. By (\ref{tr3}), we get
\[
\beta=2\arccot\left(\cos a\tan \frac{\gamma}{2}\right)=2\arccot\left(\cot \alpha \cot \frac{\alpha}{2}\tan \frac{\gamma}{2}\right).
\]
If $a$ or $\aaa$ is fixed, then we derive that
\[
\beta'=-\frac{\cos a}{\cos^2\frac{\gamma}{2}+(\cos a\sin\frac{\gamma}{2})^2}<0.
\]
So $\beta$ is strictly decreasing with respect to $\gamma$ (see the second of Figure \ref{b(r)1}). Moreover, there is a special case that two regular triangles can be merged into a rhombus $i.e.$ $\bbb=2\aaa>\aaa=\ccc$. When $\bbb=\ccc$, we have $\aaa<\ccc$ and $\bbb<2\aaa$ $i.e.$ $\aaa<\bbb<2\aaa$. 
	\end{proof}

	\begin{lemma}\label{ico}
		In a convex $(a^3,a^4)$-tiling, if $\bbb=2\aaa=2\ccc$, then the tiling must be of the icosahedral type.
	\end{lemma}

	\begin{proof}
		The $\bbb=2\aaa=2\ccc$ means the rhombus is made of two regular triangles. So the $(a^3,a^4)$-tiling induces an $a^3$-tiling. By Lemma \ref{type}, we have $\frac{1}{3}\pi<\aaa<\frac{1}{2}\pi$. Hence $\aaa^5$ is the only vertex and we get the icosahedron. 
	\end{proof}

    \begin{lemma}\label{ico1}
		In a convex $(a^3,a^4)$-tiling, if $\ccc=\aaa$ or $\bbb=k\ccc+(2-k)\aaa$ for some constant $k\ge0$, then the tiling must be of the icosahedral type.
	\end{lemma}

	\begin{proof}
        When $\ccc=\aaa$, we get $\bbb=2\aaa$ by (\ref{tr5}). In the proof of Lemma \ref{b(r)}, we know $\beta=2\arccot\left(\cos a\tan \frac{\gamma}{2}\right)$ is strictly decreasing with respect to $\gamma$ for fixed $a$ or $\aaa$. When $k\ge0$, $\bbb=k\ccc+(2-k)\aaa$ is increasing with respect to $\gamma$. So the two equations together have at most one solution. But $\ccc=\aaa,\bbb=2\aaa$ is an obvious solution to both equations. By Lemma \ref{ico}, the convex $(a^3,a^4)$-tiling must be of the icosahedral type in either case.
	\end{proof}

We will call an $(a^3,a^4)$-tiling or its protoset \textbf{rational} when all angles $(\aaa,\bbb,\ccc)$ are rational multiples of $\pi$.

	\begin{lemma}\label{rational}
		If the protoset is convex and rational, then the angles $(\aaa,\bbb,\ccc)$ must be
		$(\frac{4}{9},\frac{7}{9},\frac{2}{3})\pi$, $(\frac{3}{7},\frac{17}{21},\frac{11}{21})\pi$ or $(\aaa,2\aaa,\aaa)$ for any $\aaa \in  (\mathbb{Q}\cap(\frac{1}{3},\frac{1}{2}))\pi$.
	\end{lemma}

	\begin{proof}	
By Lemma \ref{type}, we have $\frac{1}{3}\pi<\aaa<\frac{1}{2}\pi<\bbb<\pi$, and $\bbb\ge\ccc$, $\bbb+\ccc>\pi$. Let $x=e^{i\aaa},y=e^{i\bbb},z=e^{i\ccc}$. Then we derive $x^2y + x^2z - xyz - xy - xz - x + y + z=0$ by (\ref{tr5}).

Now we use the programs of \cite{cj,kkpr} in \href{https://github.com/kedlaya/tetrahedra/}{https://github.com/kedlaya/tetrahedra/} to compute the cyclotomic points for this equation. It turns out that the only solutions within our range are the ones stated.	 
	\end{proof}

\begin{proposition}\label{rt}	
	Any convex rational $(a^3,a^4)$-tiling is either of the icosahedral type or $T(4a^3,3a^4;$ $3\aaa^3\ccc,3\aaa\bbb^2,\ccc^3)$.
\end{proposition}

\begin{proof}
	When $(\aaa,\bbb,\ccc)=(\aaa,2\aaa,\aaa)$, the $(a^3,a^4)$-tiling must be of the icosahedral type by Lemma \ref{ico}. When $(\aaa,\bbb,\ccc)=(\frac{3}{7},\frac{17}{21},\frac{11}{21})\pi$, $\aaa\ccc^3$ is the only vertex, contradicting Lemma \ref{AAD}. 

	\begin{figure}[htp]
		\centering
			\begin{tikzpicture}[>=latex,scale=0.7]
				\foreach \a in {0,1,2}
				\draw[rotate=120*\a]
				(0,0)--(0,1.2)
				(-1.45,0.84)--(0,1.2)--(1.45,0.84)
                (1.45,0.84) arc (10.65:169.35:1.47);
		\node at (0.3,0) {\small $\ccc$};
		\node at (-0.3,0) {\small $\ccc$};
		\node at (0,-0.3) {\small $\ccc$};
	\node[draw,shape=circle, inner sep=0.2] at (0.8,0.5) {\small $1$};
	\node[draw,shape=circle, inner sep=0.2] at (-0.8,0.5) {\small $2$};
	\node[draw,shape=circle, inner sep=0.2] at (0,-0.9) {\small $3$};
	\node[draw,shape=circle, inner sep=0.2] at (0,1.5) {\small $4$};
	\node[draw,shape=circle, inner sep=0.2] at (1.43,-0.79) {\small $5$};
	\node[draw,shape=circle, inner sep=0.2] at (-1.43,-0.79) {\small $6$};
	\node[draw,shape=circle, inner sep=0.2] at (2,1) {\small $7$};
			\end{tikzpicture}
		\caption{$T(4a^3,3a^4;3\aaa^3\ccc,3\aaa\bbb^2,\ccc^3)$.}
		\label{r3}
	\end{figure}
When $(\aaa,\bbb,\ccc)=(\frac{4}{9},\frac{7}{9},\frac{2}{3})\pi$, we have the AVC $\subset\{\aaa^3\ccc,\aaa\bbb^2,\ccc^3\}$. The AAD $\thin\ccc^\bbb\thin^\aaa\aaa\thin\cdots$ of $\aaa^3\ccc$ gives $\thin\aaa\thin\bbb\thin\cdots=\aaa\bbb^2$ and the AAD $\thin\bbb^\ccc\thin^\ccc\bbb\thin\cdots$ of $\aaa\bbb^2$ gives $\thin\ccc\thin\ccc\thin\cdots=\ccc^3$. So $\ccc^3$ must appear, which determines $T_1,T_2,T_3$ in Figure \ref{r3}. Then $\bbb^2\cdots=\aaa\bbb^2$ determines $T_4,T_5,T_6$; $\aaa^2\ccc\cdots=\aaa^3\ccc$ determines $T_7$. The tiling is completed and the 3D picture is the 12th of Figure \ref{Fig3}.
\end{proof}

A vertex $\aaa^{n_1}\bbb^{n_2}\ccc^{n_3}$ can be efficiently represented by its vector type $\boldsymbol{n}=(n_1\, n_2 \, n_3)\in \mathbb{N}^3$. We will use both representations interchangeably afterwards for convenience.

\begin{lemma}[Irrational Angle Lemma]\label{irrational}
	
	Given two different vertices $\boldsymbol{l},\boldsymbol{m}\in \mathbb{N}^3$ in an $(a^3,a^4)$-tiling with some irrational angle, if there is another vertex $\boldsymbol{n}$, then $\boldsymbol{l},\boldsymbol{m},\boldsymbol{n}$ are linearly dependent. In other words, the square matrix of $\boldsymbol{l},\boldsymbol{m},\boldsymbol{n}$ has determinant $0$. 
\end{lemma}
\begin{proof} 
	By the vertices $\boldsymbol{l},\boldsymbol{m}, \boldsymbol{n}$, the angles satisfy a linear system of equations
	\[\left(\begin{matrix} 
		\boldsymbol{l}\\
        \boldsymbol{m}\\
        \boldsymbol{n}
	\end{matrix}\right)\left(\begin{matrix} 
		\aaa\\
        \bbb\\
        \ccc
	\end{matrix}\right)=\left(\begin{matrix} 
		2\\
        2\\
        2
	\end{matrix}\right)\pi.
	\]
	If $\boldsymbol{l},\boldsymbol{m},\boldsymbol{n}$ are linearly independent, then the system has a unique rational solution, a contradiction. 
\end{proof}

In Section \ref{irr1} and \ref{irr2}, we will use Lemma \ref{irrational} to handle the convex $(a^3,a^4)$-tiling with some irrational angle, which imposes strong constraints on all possible vertices. But before that, we will discuss the short and easier concave and degenerate cases first in the next section.

	\section{Concave and degenerate cases $\bbb\ge\pi$}

	Since all vertices have degree
$\ge3$, $\bbb^2\cdots$ is not a vertex. By Lemma \ref{AAD}, $\thin\ccc\thin\ccc\thin\cdots$ can never be a vertex and $\thin\aaa\thin\bbb\thin\cdots$ must appear.

     \subsection*{$\bo{\bbb>\pi}$}
      By Lemma \ref{type}, we have $\frac{1}{2}\pi<\aaa<\pi$ and $\thin\aaa\thin\bbb\thin\cdots=\aaa\bbb\ccc$, which determines $T_1,T_2,T_3$ in Figure \ref{abrb>}. Then $\thin\bbb_3\thin\ccc_2\thin\cdots=\thin\bbb_3\thin\ccc_2\thin\aaa\thin$ determines $T_4$, $\thin\aaa_1\thin\bbb_3\thin\cdots=\aaa\bbb\ccc$ determines $T_5$ and we get a tiling like a triangular prism.
	\begin{figure}[H]
		\centering
		\begin{tikzpicture}[>=latex,scale=0.75]       
		\draw (0,0) -- (0,2) -- (1,1) -- (0,0);
        \draw (0,0) -- (3,0) -- (2,1) -- (1,1)
              (2,1) -- (3,2)
              (0,2) -- (3,2) -- (3,0);
		\node at (0.7,1) {\small $\aaa$};
		\node at (1,0.7) {\small $\bbb$};
		\node at (1,1.3) {\small $\ccc$};
		\node at (-0.2,2.2) {\small $\ccc$};
	\node[draw,shape=circle, inner sep=0.2] at (0.3,1) {\small $1$};
	\node[draw,shape=circle, inner sep=0.2] at (1.5,0.5) {\small $2$};
	\node[draw,shape=circle, inner sep=0.2] at (1.5,1.5) {\small $3$};
	\node[draw,shape=circle, inner sep=0.2] at (2.5,1) {\small $4$};
	\node[draw,shape=circle, inner sep=0.2] at (1.5,2.3) {\small $5$};
		\end{tikzpicture}
		\caption{$T(2a^3,3a^4;6\aaa\bbb\ccc)$.}
		\label{abrb>}
	\end{figure}

     \subsection*{$\bo{\bbb=\pi}$}\label{b=pi}
     By Lemma \ref{type}, we have $\aaa=\frac{1}{2}\pi$ and $\thin\aaa\thin\bbb\thin\cdots=\aaa^2\bbb,\aaa\bbb\ccc$. 

     \subsubsection*{Case $\aaa\bbb\ccc$}
We have $\ccc=\frac{1}{2}\pi$. In Figure \ref{abr2b=p}, $\aaa\bbb\ccc$ determines $T_1,T_2,T_3$; $\thin\bbb\thin\ccc\thin\cdots=\thin\bbb\thin\ccc\thin\aaa\thin$ determines $T_4$. When one of $\thin\aaa_1\thin\bbb_3\thin\cdots$ and $\thin\aaa_4\thin\bbb_2\thin\cdots$ is $\aaa\bbb\ccc$, the tiling is the same as Figure \ref{abrb>}. When they are $\aaa^2\bbb$, $T_5,T_6$ are determined and the 3D picture is the 6th of Figure \ref{Fig3}.

	\begin{figure}[H]
		\centering
		\begin{tikzpicture}[>=latex,scale=0.75]       
		\draw (0,0) -- (0,2) -- (1,1) -- (0,0);
        \draw (0,0) -- (3,0) -- (2,1) -- (1,1)
              (2,1) -- (3,2)
              (0,2) -- (3,2) -- (3,0) (3,2)--(3.5,2.3) (0,0)--(-0.5,-0.3);
		\node at (0.7,1) {\small $\aaa$};
		\node at (1,0.7) {\small $\bbb$};
		\node at (1,1.3) {\small $\ccc$};
	\node[draw,shape=circle, inner sep=0.2] at (0.3,1) {\small $1$};
	\node[draw,shape=circle, inner sep=0.2] at (1.5,0.5) {\small $2$};
	\node[draw,shape=circle, inner sep=0.2] at (1.5,1.5) {\small $3$};
	\node[draw,shape=circle, inner sep=0.2] at (2.5,1) {\small $4$};
	\node[draw,shape=circle, inner sep=0.2] at (1.5,2.3) {\small $5$};
	\node[draw,shape=circle, inner sep=0.2] at (1.5,-0.3) {\small $6$};
		\end{tikzpicture}
		\caption{$T(4a^3,2a^4;2\aaa\bbb\ccc,2\aaa^2\bbb,2\aaa^3\ccc)$.}
		\label{abr2b=p}
	\end{figure}

     \subsubsection*{Case $\aaa^2\bbb$}
By $\aaa^2\bbb=\thin\bbb^\ccc\thin^\aaa\aaa\thin\aaa\thin$ and similar proof of Lemma \ref{deg345}, we get $\thin\aaa\thin\ccc\thin\cdots=\aaa\bbb\ccc^2,\aaa^2\ccc,\aaa^2\ccc^2,\aaa^3\ccc,$ $\aaa^3\ccc^2,\aaa^3\ccc^3$, and $\thin\aaa\thin\bbb\thin\cdots$ can only be $\aaa^2\bbb$.

	\begin{figure}[htp]
		\centering
		\begin{tikzpicture}[>=latex,scale=0.8]       
		\draw (0,1.2)--(0,-1.2) (1,1.2)--(1,-1.2) (2,1.2)--(2,-1.2) (4,1.2)--(4,-1.2) (5,1.2)--(5,-1.2);
        \draw (1,0) -- (0,0) (2,0) -- (3,0.5)--(4,0)--(3,-0.5)--(2,0) (3,0.5)--(3,-0.5);
		\node at (0.5,1.2) {\small $\aaa$};
		\node at (0.5,-1.2) {\small $\aaa$};
		\node at (2.35,0) {\small $\aaa$};
		\node at (3.65,0) {\small $\aaa$};
		\node at (3,1.3) {\small $\bbb$};
		\node at (3,-1.3) {\small $\bbb$};
		\node at (1.5,1.2) {\small $\ccc$};
		\node at (1.5,-1.2) {\small $\ccc$};
		\node at (4.5,1.2) {\small $\ccc$};
		\node at (4.5,-1.2) {\small $\ccc$};
	\node[draw,shape=circle, inner sep=0.2] at (0.5,0.5) {\small $1$};
	\node[draw,shape=circle, inner sep=0.2] at (1.5,0) {\small $2$};
	\node[draw,shape=circle, inner sep=0.2] at (3,0.8) {\small $3$};
	\node[draw,shape=circle, inner sep=0.2] at (4.5,0) {\small $4$};
	\node[draw,shape=circle, inner sep=0.2] at (2.75,0) {\small $5$};
	\node[draw,shape=circle, inner sep=0.2] at (3.25,0) {\small $6$};
	\node[draw,shape=circle, inner sep=0.2] at (3,-0.8) {\small $7$};
	\node[draw,shape=circle, inner sep=0.2] at (0.5,-0.5) {\small $8$};
		\end{tikzpicture}\hspace{30pt}
		\begin{tikzpicture}[>=latex,scale=0.8]  
    \draw (0,1.2)--(0,-1.2) (1,1.2)--(1,-1.2) (2,1.2)--(2,-1.2) (3,1.2)--(3,-1.2) (0,0)--(2,0);
		\node at (0.5,1.2){\small $\aaa$};
		\node at (1.5,1.2){\small $\aaa$};
		\node at (0.5,-1.2){\small $\aaa$};
		\node at (1.5,-1.2){\small $\aaa$};
		\node at (2.5,1.2){\small $\ccc$};
		\node at (2.5,-1.2){\small $\ccc$};
	\node[draw,shape=circle, inner sep=0.2] at (0.5,0.5) {\small $1$};
	\node[draw,shape=circle, inner sep=0.2] at (1.5,0.5) {\small $2$};
	\node[draw,shape=circle, inner sep=0.2] at (2.5,0) {\small $3$};
	\node[draw,shape=circle, inner sep=0.2] at (0.5,-0.5) {\small $4$};
	\node[draw,shape=circle, inner sep=0.2] at (1.5,-0.5) {\small $5$};
		\end{tikzpicture}\hspace{30pt}
		\begin{tikzpicture}[>=latex,scale=0.8] 
	\foreach \a in {0,1}
	{
		\begin{scope}[xshift=2*\a cm] 
		\draw (0,1.2)--(0,-1.2)
		(1,1.2)--(1,-1.2) (0,0)--(1,0)
		(2,1.2)--(2,-1.2);
		\node at (0.5,1.2){\small $\aaa$};
		\node at (1.5,1.2){\small $\ccc$};
		
		\node at (0.5,-1.2){\small $\aaa$};
		\node at (1.5,-1.2){\small $\ccc$};
		\end{scope}
	}
	\node[draw,shape=circle, inner sep=0.2] at (0.5,0.5) {\small $1$};
	\node[draw,shape=circle, inner sep=0.2] at (1.5,0) {\small $2$};
	\node[draw,shape=circle, inner sep=0.2] at (2.5,0.5) {\small $3$};
	\node[draw,shape=circle, inner sep=0.2] at (3.5,0) {\small $4$};
	\node[draw,shape=circle, inner sep=0.2] at (0.5,-0.5) {\small $5$};
	\node[draw,shape=circle, inner sep=0.2] at (2.5,-0.5) {\small $6$};
		\end{tikzpicture}
		\begin{tikzpicture}[>=latex,scale=0.8] 
	\foreach \a in {0,1}
	{
		\begin{scope}[xshift=2*\a cm] 
		\draw (0,1.2)--(0,-1.2)
		(1,1.2)--(1,-1.2) (0,0)--(1,0)
		(2,1.2)--(2,-1.2);
		\node at (0.5,1.2){\small $\aaa$};
		
		\node at (0.5,-1.2){\small $\aaa$};
		\end{scope}
	}
    \draw (1,0)--(2,0);
		\node at (1.5,1.2){\small $\aaa$};
		\node at (1.5,-1.2){\small $\aaa$};
		\node at (3.5,1.2){\small $\ccc$};
		\node at (3.5,-1.2){\small $\ccc$};
	\node[draw,shape=circle, inner sep=0.2] at (0.5,0.5) {\small $1$};
	\node[draw,shape=circle, inner sep=0.2] at (1.5,0.5) {\small $2$};
	\node[draw,shape=circle, inner sep=0.2] at (2.5,0.5) {\small $3$};
	\node[draw,shape=circle, inner sep=0.2] at (3.5,0) {\small $4$};
	\node[draw,shape=circle, inner sep=0.2] at (0.5,-0.5) {\small $5$};
	\node[draw,shape=circle, inner sep=0.2] at (1.5,-0.5) {\small $7$};
	\node[draw,shape=circle, inner sep=0.2] at (2.5,-0.5) {\small $6$};
		\end{tikzpicture}\hspace{15pt}
		\begin{tikzpicture}[>=latex,scale=0.8] 
	\foreach \a in {0,1}
	{
		\begin{scope}[xshift=2*\a cm] 
		\draw (0,1.2)--(0,-1.2)
		(1,1.2)--(1,-1.2) (0,0)--(1,0)
		(2,1.2)--(2,-1.2);
		\node at (0.5,1.2){\small $\aaa$};
		\node at (1.5,1.2){\small $\ccc$};
		
		\node at (0.5,-1.2){\small $\aaa$};
		\node at (1.5,-1.2){\small $\ccc$};
		\end{scope}
	}
    \draw (-1,1.2)--(-1,-1.2) (-1,0)--(0,0);
		\node at (-0.5,1.2){\small $\aaa$};
		\node at (-0.5,-1.2){\small $\aaa$};
	\node[draw,shape=circle, inner sep=0.2] at (-0.5,0.5) {\small $1$};
	\node[draw,shape=circle, inner sep=0.2] at (0.5,0.5) {\small $2$};
	\node[draw,shape=circle, inner sep=0.2] at (1.5,0) {\small $3$};
	\node[draw,shape=circle, inner sep=0.2] at (2.5,0.5) {\small $4$};
	\node[draw,shape=circle, inner sep=0.2] at (3.5,0) {\small $5$};
	\node[draw,shape=circle, inner sep=0.2] at (0.5,-0.5) {\small $7$};
	\node[draw,shape=circle, inner sep=0.2] at (-0.5,-0.5) {\small $6$};
	\node[draw,shape=circle, inner sep=0.2] at (2.5,-0.5) {\small $8$};
		\end{tikzpicture}\hspace{15pt}
		\begin{tikzpicture}[>=latex,scale=0.8] 
	\foreach \a in {0,1,2}
	{
		\begin{scope}[xshift=2*\a cm] 
		\draw (0,1.2)--(0,-1.2)
		(1,1.2)--(1,-1.2) (0,0)--(1,0)
		(2,1.2)--(2,-1.2);
		\node at (0.5,1.2){\small $\aaa$};
		\node at (1.5,1.2){\small $\ccc$};
		
		\node at (0.5,-1.2){\small $\aaa$};
		\node at (1.5,-1.2){\small $\ccc$};
		\end{scope}
	}
	\node[draw,shape=circle, inner sep=0.2] at (0.5,0.5) {\small $1$};
	\node[draw,shape=circle, inner sep=0.2] at (1.5,0) {\small $2$};
	\node[draw,shape=circle, inner sep=0.2] at (2.5,0.5) {\small $3$};
	\node[draw,shape=circle, inner sep=0.2] at (3.5,0) {\small $4$};
	\node[draw,shape=circle, inner sep=0.2] at (4.5,0.5) {\small $5$};
	\node[draw,shape=circle, inner sep=0.2] at (5.5,0) {\small $6$};
	\node[draw,shape=circle, inner sep=0.2] at (0.5,-0.5) {\small $7$};
	\node[draw,shape=circle, inner sep=0.2] at (2.5,-0.5) {\small $8$};
	\node[draw,shape=circle, inner sep=0.2] at (4.5,-0.5) {\small $9$};
		\end{tikzpicture}
		\caption{$\thin\aaa\thin\ccc\thin\cdots=\aaa\bbb\ccc^2,\aaa^2\ccc,\aaa^2\ccc^2,\aaa^3\ccc,\aaa^3\ccc^2,\aaa^3\ccc^3$.}
		\label{ar}
	\end{figure}

When $\thin\aaa\thin\ccc\thin\cdots=\aaa\bbb\ccc^2$, we have $\ccc=\frac{1}{4}\pi$. The vertex $\aaa\bbb\ccc^2=\thin\aaa\thin\ccc\thin\bbb\thin\ccc\thin$ determines $T_1,\cdots,T_4$ in the first of Figure \ref{ar}; $\thin\bbb\thin\ccc\thin\cdots=\thin\bbb\thin\ccc\thin\aaa\thin\cdots$ $=\thin\bbb\thin\ccc\thin\aaa\thin\ccc\thin$ determines $T_5,T_6,T_7$; $\thin\aaa\thin\bbb\thin\cdots=\aaa^2\bbb$ determines $T_8$. 

We omit the details and claim that the other possibilities of $\thin\aaa\thin\ccc\thin\cdots$ produce exactly the other  tilings in Figure \ref{ar}. The 3D pictures are the 3rd, 1st, 5th, 7th, 4th and 2nd of Figure \ref{Fig3}, respectively.

	\section{Convex case $\ccc=\bbb<\pi$ with some irrational angle}\label{irr1}

     By Lemma \ref{type} and \ref{b(r)}, we have $\frac{\pi}{3}<\aaa<\frac{\pi}{2}<\bbb<2\aaa$. By Lemma \ref{AAD} and similar proof of Lemma \ref{deg345}, we get $\thin\aaa\thin\bbb\thin\cdots=\aaa\bbb^2,\aaa\bbb^3,\aaa^2\bbb^2,\aaa^3\bbb,\aaa^4\bbb$. The Lemma \ref{irrational} implies that there is only one kind of vertex in all such tilings.
 
    \subsection*{Case $\thin\aaa\thin\bbb\thin\cdots=\aaa\bbb^2$, AVC $\subset\{\aaa\bbb^2\}$}

By the AVC, we can get a tiling which is the same as Figure \ref{abrb>}.

     \subsection*{Case $\thin\aaa\thin\bbb\thin\cdots=\aaa\bbb^3$, AVC $\subset\{\aaa\bbb^3\}$}

      The vertex $\aaa\bbb^3$ determines $T_1,T_2,T_3,T_4$ in Figure \ref{ab3v}; $\aaa\bbb\cdots=\aaa\bbb^3$ determines $T_5,T_6,T_7$. Since the symmetry of the partial tiling and $\bbb^2\cdots=\aaa\bbb^3$, we have four cases:

	\begin{figure}[htp]
		\centering
		\begin{tikzpicture}[>=latex,scale=0.6]       
		\draw (0,0) -- (1,0) -- (0.5,0.87) -- (0,0)
              (1,0) -- (1.87,0.5) -- (1.37,1.37) -- (0.5,0.87)
              (0,0) -- (-0.87,0.5) -- (-0.37,1.37) -- (0.5,0.87)
              (-0.37,1.37) -- (0.5,1.87) -- (1.37,1.37)
              (0,0) -- (0,-1) -- (1,-1) -- (1,0)
              (1.87,0.5) -- (1.87,-0.5) -- (1,-1)
              (-0.87,0.5) -- (-0.87,-0.5) -- (0,-1);
        \draw[dashed] (0,-1) -- (0.5,-1.87) -- (1,-1) 
                      (-0.37,1.37) -- (-1.37,1.37) -- (-0.87,0.5)
                      (1.37,1.37) -- (2.37,1.37) -- (1.87,0.5);
	\node[draw,shape=circle, inner sep=0.2] at (0.5,0.3) {\small $1$};
	\node[draw,shape=circle, inner sep=0.2] at (1.2,0.65) {\small $2$};
	\node[draw,shape=circle, inner sep=0.2] at (0.5,1.3) {\small $3$};
	\node[draw,shape=circle, inner sep=0.2] at (-0.2,0.65) {\small $4$};
	\node[draw,shape=circle, inner sep=0.2] at (-0.35,-0.3) {\small $5$};
	\node[draw,shape=circle, inner sep=0.2] at (0.5,-0.5) {\small $6$};
	\node[draw,shape=circle, inner sep=0.2] at (1.35,-0.3) {\small $7$};
		\end{tikzpicture}\hspace{10pt}
		\begin{tikzpicture}[>=latex,scale=0.6]       
		\draw (0,0) -- (1,0) -- (0.5,0.87) -- (0,0)
              (1,0) -- (1.87,0.5) -- (1.37,1.37) -- (0.5,0.87)
              (0,0) -- (-0.87,0.5) -- (-0.37,1.37) -- (0.5,0.87)
              (-0.37,1.37) -- (0.5,1.87) -- (1.37,1.37)
              (0,0) -- (0,-1) -- (1,-1) -- (1,0)
              (1.87,0.5) -- (1.87,-0.5) -- (1,-1)
              (-0.87,0.5) -- (-0.87,-0.5) -- (0,-1);
        \draw[dashed] (0,-1) -- (0.5,-1.87) -- (1,-1) 
                      (-0.37,1.37) -- (-1.37,1.37) -- (-0.87,0.5)
                      (1.87,0.5) -- (2.73,0.95) -- (2.26,1.84) -- (1.37,1.37);
		\end{tikzpicture}\hspace{10pt}
		\begin{tikzpicture}[>=latex,scale=0.6]       
		\draw (0,0) -- (1,0) -- (0.5,0.87) -- (0,0)
              (1,0) -- (1.87,0.5) -- (1.37,1.37) -- (0.5,0.87)
              (0,0) -- (-0.87,0.5) -- (-0.37,1.37) -- (0.5,0.87)
              (-0.37,1.37) -- (0.5,1.87) -- (1.37,1.37)
              (0,0) -- (0,-1) -- (1,-1) -- (1,0)
              (1.87,0.5) -- (1.87,-0.5) -- (1,-1)
              (-0.87,0.5) -- (-0.87,-0.5) -- (0,-1);
        \draw[dashed] (0,-1) -- (0.5,-1.87) -- (1,-1)
                      (-0.87,0.5) -- (-1.73,0.95) -- (-1.26,1.84) -- (-0.37,1.37)
                      (1.87,0.5) -- (2.73,0.95) -- (2.26,1.84) -- (1.37,1.37)
                      (-1.26,1.84)--(0.5,1.87)--(2.26,1.84);
		\node at (0.5,1.87){\small $\bullet$};
		\end{tikzpicture}\hspace{10pt}
		\begin{tikzpicture}[>=latex,scale=0.6]       
		\draw (0,0) -- (1,0) -- (0.5,0.87) -- (0,0)
              (1,0) -- (1.87,0.5) -- (1.37,1.37) -- (0.5,0.87)
              (0,0) -- (-0.87,0.5) -- (-0.37,1.37) -- (0.5,0.87)
              (-0.37,1.37) -- (0.5,1.87) -- (1.37,1.37)
              (0,0) -- (0,-1) -- (1,-1) -- (1,0)
              (1.87,0.5) -- (1.87,-0.5) -- (1,-1)
              (-0.87,0.5) -- (-0.87,-0.5) -- (0,-1);
        \draw[dashed] (0,-1) -- (0,-2) -- (1,-2) -- (1,-1) 
                      (-0.87,0.5) -- (-1.73,0.95) -- (-1.26,1.84) -- (-0.37,1.37)
                      (1.87,0.5) -- (2.73,0.95) -- (2.26,1.84) -- (1.37,1.37)
                      (-1.26,1.84)--(0.5,1.87)--(2.26,1.84);
		\node at (0.5,1.87){\small $\bullet$};
		\end{tikzpicture}\hspace{10pt}
		\caption{$\aaa\bbb^3$ is the vertex.}
		\label{ab3v}
	\end{figure}

   In the third and fourth cases, we get $\aaa^2\bbb\cdots$, contradicting the AVC. In the first case, we can determine $T_8,T_9,\cdots,T_{26}$ as the first of Figure \ref{ab=r} and its 3D picture is the 10th of Figure \ref{Fig3} corresponding to the rhombicuboctahedron. In the second case, we get a different tiling as the second of Figure \ref{ab=r} and its 3D picture is the 11th of Figure \ref{Fig3} corresponding to the elongated square gyrobicupola. By rotating the partial tiling with the red boundary $45^\circ$, the first tiling will become the second. \label{ab3}

     \subsection*{Case $\thin\aaa\thin\bbb\thin\cdots=\aaa^2\bbb^2$, AVC $\subset\{\aaa^2\bbb^2\}$}     

The vertex $\aaa^2\bbb^2$ is $\thin\aaa\thin\aaa\thin\bbb\thin\bbb\thin$ or $\thin\aaa\thin\bbb\thin\aaa\thin\bbb\thin$. In the third of Figure \ref{ab=r}, $\thin\aaa\thin\aaa\thin\bbb\thin\bbb\thin$ determines $T_1,T_2,T_3,T_4$; $\bbb^2\cdots=\aaa^2\bbb^2$ determines $T_5,T_6$; $\aaa^2\cdots=\aaa^2\bbb^2$ determines $T_7,\cdots,T_{10}$. Similarly, $T_{11},T_{12},T_{13},T_{14}$ are determined and its 3D picture is a triangular orthobicupola. Next $\aaa^2\bbb^2$ can only be $\thin\aaa\thin\bbb\thin\aaa\thin\bbb\thin$. Then we get a new tiling in the fourth of Figure \ref{ab=r} and its 3D picture is a cuboctahedron. By rotating the partial tiling with the red boundary $60^\circ$ or flipping the partial tiling with the red boundary, the first tiling will become the second.

     \subsection*{Case $\thin\aaa\thin\bbb\thin\cdots=\aaa^3\bbb$, AVC $\subset\{\aaa^3\bbb\}$}     

By the AVC, we can get a tiling in the fifth of Figure \ref{ab=r} and its 3D picture is the 8th of Figure \ref{Fig3}  corresponding to one of the anti-prisms.

     \subsection*{Case $\thin\aaa\thin\bbb\thin\cdots=\aaa^4\bbb$, AVC $\subset\{\aaa^4\bbb\}$}     
In the sixth of Figure \ref{ab=r}, $\aaa^4\bbb$ determines $T_1,\cdots,T_5$, $\aaa\bbb\cdots=\aaa^4\bbb$ determines $T_6,\cdots,T_{13}$. Since the symmetry of the partial tiling and $\aaa^3\cdots=\aaa^4\bbb$, we might as well take $\thin\aaa_3\thin\aaa_4\thin\aaa_6\thin\cdots=\thin\aaa_3\thin\aaa_4\thin\aaa_6\thin\aaa\thin\bbb\thin$ which determines $T_{14},T_{15}$. Then we can determine $T_{16},\cdots,T_{38}$ and its 3D picture is the 9th of Figure \ref{Fig3} corresponding to a snub cube.

	\begin{figure}[H]
		\centering
			\begin{tikzpicture}[>=latex,scale=0.8]
				\foreach \a in {0,1,2}
				\draw[rotate=120*\a]
				(0,0.58)--(0.87,1.08)--(0,1.58)--(-1.08,2.63)--(1.08,2.63)--(0,1.58)
				(0.87,1.08)--(1.87,1.08)--(1.08,2.63)--(2.91,1.68)--(2.82,-0.38)--(1.87,1.08)--(1.37,0.21)--(0.87,1.08) (1.37,0.21)--(0.5,-0.29)--(0,0.58) (0,1.58)--(-0.87,1.08)
                (2.91,1.68) arc (24.5:155.5:3.19);
               \draw[line width=1,color=red] (-2.91,1.68)--(-1.08,2.63)--(0,1.58)--(0.87,1.08)--(1.37,0.21)--(1.37,-0.79)--(1.74,-2.26)--(0,-3.36);
               \draw[rotate=120][line width=1,color=red] (2.91,1.68) arc (24.5:155.5:3.19);
	\node[draw,shape=circle, inner sep=0.2] at (0,0) {\small $1$};
	\node[draw,shape=circle, inner sep=0.2] at (0.7,0.35) {\small $2$};
	\node[draw,shape=circle, inner sep=0.2] at (0,1) {\small $3$};
	\node[draw,shape=circle, inner sep=0.2] at (-0.7,0.35) {\small $4$};
	\node[draw,shape=circle, inner sep=0.2] at (-0.9,-0.6) {\small $5$};
	\node[draw,shape=circle, inner sep=0.2] at (0,-0.8) {\small $6$};
	\node[draw,shape=circle, inner sep=0.2] at (0.9,-0.6) {\small $7$};
	\node[draw,shape=circle, inner sep=0.2] at (1.4,0.79) {\small $8$};
	\node[draw,shape=circle, inner sep=0.2] at (-1.4,0.79) {\small $9$};
	\node[draw,shape=circle, inner sep=0.2] at (0,-1.6) {\footnotesize $10$};
	\node[draw,shape=circle, inner sep=0.2] at (-1,1.72) {\footnotesize $11$};
	\node[draw,shape=circle, inner sep=0.2] at (-1.99,0.01) {\footnotesize $12$};
	\node[draw,shape=circle, inner sep=0.2] at (-2.35,1.36) {\footnotesize $13$};
	\node[draw,shape=circle, inner sep=0.2] at (-1,-1.59) {\footnotesize $14$};
	\node[draw,shape=circle, inner sep=0.2] at (1,-1.59) {\footnotesize $15$};
	\node[draw,shape=circle, inner sep=0.2] at (0,-2.67) {\footnotesize $16$};
	\node[draw,shape=circle, inner sep=0.2] at (1,1.72) {\footnotesize $17$};
	\node[draw,shape=circle, inner sep=0.2] at (1.99,0.01) {\footnotesize $18$};
	\node[draw,shape=circle, inner sep=0.2] at (2.35,1.36) {\footnotesize $19$};
	\node[draw,shape=circle, inner sep=0.2] at (0,2.22) {\footnotesize $20$};
	\node[draw,shape=circle, inner sep=0.2] at (-1.93,-1.11) {\footnotesize $21$};
	\node[draw,shape=circle, inner sep=0.2] at (1.93,-1.11) {\footnotesize $22$};
	\node[draw,shape=circle, inner sep=0.2] at (0,3.07) {\footnotesize $23$};
	\node[draw,shape=circle, inner sep=0.2] at (-2.66,-1.54) {\footnotesize $24$};
	\node[draw,shape=circle, inner sep=0.2] at (2.66,-1.54) {\footnotesize $25$};
	\node[draw,shape=circle, inner sep=0.2] at (3.28,2.09) {\footnotesize $26$};
			\end{tikzpicture}
			\begin{tikzpicture}[>=latex,scale=0.8]
				\foreach \a in {0,1,2}
				\draw[rotate=120*\a]
				(-0.5,-0.29)--(-1.37,0.21)--(-0.87,1.08)--(0,1.58)--(0.87,1.08)--(0,0.58)--(-0.5,-0.29) 
				(-0.87,1.08)--(0,0.58);
				\foreach \a in {0,1}
				\draw[rotate=120*\a]
                (-1.37,-0.79)--(-2.82,-0.38)--(-2.91,1.68)--(-1.08,2.63)--(0,1.58) (-2.82,-0.38)--(-1.87,1.08)--(-1.08,2.63) (-1.37,0.21)--(-1.87,1.08)--(-0.81,1.08);
				\foreach \a in {1}
				\draw[rotate=120*\a]
                (-1.08,2.63)--(1.08,2.63)
                (2.91,1.68) arc (24.5:155.5:3.19);
          \draw (0.87,1.08)--(1.73,1.58)--(2.23,0.71)--(1.37,0.21) (-1.08,2.63)--(2.19,2.91)--(1.73,1.58) (1.74,-2.26)--(3.61,0.44)--(2.23,0.71) (2.19,2.91)--(3.61,0.44)
                (1.73,1.58) arc (44.89:135.11:1.22) (2.23,0.71) arc (15.11:-75.11:1.22)
                (2.19,2.91) arc (53.24:154.03:3.4) (3.61,0.44) arc (6.76:-94.03:3.4);
               \draw[line width=1,color=red] (-2.91,1.68)--(-1.08,2.63)--(0,1.58)--(0.87,1.08)--(1.37,0.21)--(1.37,-0.79)--(1.74,-2.26)--(0,-3.36);
               \draw[rotate=120][line width=1,color=red] (2.91,1.68) arc (24.5:155.5:3.19);
	\node[draw,shape=circle, inner sep=0.2] at (0,0) {\small $1$};
	\node[draw,shape=circle, inner sep=0.2] at (0.7,0.35) {\small $2$};
	\node[draw,shape=circle, inner sep=0.2] at (0,1) {\small $3$};
	\node[draw,shape=circle, inner sep=0.2] at (-0.7,0.35) {\small $4$};
	\node[draw,shape=circle, inner sep=0.2] at (-0.9,-0.6) {\small $5$};
	\node[draw,shape=circle, inner sep=0.2] at (0,-0.8) {\small $6$};
	\node[draw,shape=circle, inner sep=0.2] at (0.9,-0.6) {\small $7$};
	\node[draw,shape=circle, inner sep=0.2] at (1.4,0.79) {\small $8$};
	\node[draw,shape=circle, inner sep=0.2] at (-1.4,0.79) {\small $9$};
	\node[draw,shape=circle, inner sep=0.2] at (0,-1.6) {\footnotesize $10$};
	\node[draw,shape=circle, inner sep=0.2] at (-1,1.72) {\footnotesize $11$};
	\node[draw,shape=circle, inner sep=0.2] at (-1.99,0.01) {\footnotesize $12$};
	\node[draw,shape=circle, inner sep=0.2] at (-2.35,1.36) {\footnotesize $13$};
	\node[draw,shape=circle, inner sep=0.2] at (-1,-1.59) {\footnotesize $14$};
	\node[draw,shape=circle, inner sep=0.2] at (1,-1.59) {\footnotesize $15$};
	\node[draw,shape=circle, inner sep=0.2] at (0,-2.67) {\footnotesize $16$};
	\node[draw,shape=circle, inner sep=0.2] at (0.9,1.5) {\footnotesize $17$};
	\node[draw,shape=circle, inner sep=0.2] at (1.8,0.01) {\footnotesize $18$};
	\node[draw,shape=circle, inner sep=0.2] at (2.35,1.36) {\footnotesize $22$};
	\node[draw,shape=circle, inner sep=0.2] at (0,2.22) {\footnotesize $20$};
	\node[draw,shape=circle, inner sep=0.2] at (-1.93,-1.11) {\footnotesize $21$};
	\node[draw,shape=circle, inner sep=0.2] at (1.93,-1.11) {\footnotesize $19$};
	\node[draw,shape=circle, inner sep=0.2] at (0,3.07) {\footnotesize $23$};
	\node[draw,shape=circle, inner sep=0.2] at (-2.66,-1.54) {\footnotesize $24$};
	\node[draw,shape=circle, inner sep=0.2] at (2.66,-1.54) {\footnotesize $25$};
	\node[draw,shape=circle, inner sep=0.2] at (3.28,2.09) {\footnotesize $26$};
			\end{tikzpicture}
        \begin{tikzpicture}[>=latex,scale=0.75]       
		\draw (0,0) -- (2,0)  (0,1) -- (2,1) (1,-1) -- (1,2)
              (0,0) -- (1,-1) -- (2,0) -- (2,1) -- (1,2) -- (0,1) -- (0,0)
              (2,0) -- (3,0.5) -- (2,1) (0,0) -- (-1,0.5) -- (0,1) (1,2) -- (1,3.12) (1,-1) -- (1,-2.12)
              (3,0.5) -- (1,3.12) -- (-1,0.5) -- (1,-2.12) -- (3,0.5);
        \draw (1,3.12)--(1,3.8) (1,-2.12)--(1,-2.8);
        \draw[line width=1,color=red] (1,3.8)--(1,-2.8);
	\node[draw,shape=circle, inner sep=0.2] at (0.65,1.3) {\small $1$};
	\node[draw,shape=circle, inner sep=0.2] at (1.35,1.3) {\small $2$};
	\node[draw,shape=circle, inner sep=0.2] at (1.5,0.5) {\small $3$};
	\node[draw,shape=circle, inner sep=0.2] at (0.5,0.5) {\small $4$};
	\node[draw,shape=circle, inner sep=0.2] at (0.65,-0.3) {\small $5$};
	\node[draw,shape=circle, inner sep=0.2] at (1.35,-0.3) {\small $6$};
	\node[draw,shape=circle, inner sep=0.2] at (0.5,1.9) {\small $7$};
	\node[draw,shape=circle, inner sep=0.2] at (1.5,1.9) {\small $8$};
	\node[draw,shape=circle, inner sep=0.2] at (0.5,-0.9) {\small $9$};
	\node[draw,shape=circle, inner sep=0.2] at (1.5,-0.9) {\footnotesize $10$};
	\node[draw,shape=circle, inner sep=0.2] at (-0.3,0.5) {\footnotesize $11$};
	\node[draw,shape=circle, inner sep=0.2] at (2.3,0.5) {\footnotesize $12$};
	\node[draw,shape=circle, inner sep=0.2] at (3.5,0.5) {\footnotesize $13$};
	\node[draw,shape=circle, inner sep=0.2] at (-1.5,0.5) {\footnotesize $14$};
		\end{tikzpicture}\hspace{50pt}
		\begin{tikzpicture}[>=latex,scale=0.8]    
		\draw (0,0) -- (0,1) -- (1,1) -- (1,0) -- (0,0);
        \draw (1,0) -- (2,0.5) -- (1,1) -- (0.5,2) -- (0,1) -- (-1,0.5) -- (0,0) -- (0.5,-1) -- (1,0)
              (0.5,2) -- (1.79,1.79) -- (2,0.5) -- (1.79,-0.79) -- (0.5,-1) -- (-0.79,-0.79) -- (-1,0.5) -- (-0.79,1.79) -- (0.5,2);
        \draw (1.79,1.79) arc (17.71:162.29:1.35);
        \draw (1.79,1.79) arc (72.29:-72.29:1.35);
        \draw (-0.79,-0.79) arc (-162.29:-17.71:1.35);
        \draw (-0.79,-0.79) arc (-107.71:-252.29:1.35);
        \draw[line width=1,color=red] (-0.79,-0.79) arc (-107.71:-252.29:1.35) 
             (-0.79,1.79)--(0.5,2)--(1,1)--(1,0)--(0.5,-1)--(-0.79,-0.79);
	\node[draw,shape=circle, inner sep=0.2] at (0.5,1.3) {\small $1$};
	\node[draw,shape=circle, inner sep=0.2] at (0.5,0.5) {\small $2$};
	\node[draw,shape=circle, inner sep=0.2] at (-0.3,0.5) {\small $3$};
	\node[draw,shape=circle, inner sep=0.2] at (-0.4,1.4) {\small $4$};
	\node[draw,shape=circle, inner sep=0.2] at (1.3,0.5) {\small $5$};
	\node[draw,shape=circle, inner sep=0.2] at (1.4,1.4) {\small $6$};
	\node[draw,shape=circle, inner sep=0.2] at (0.5,-0.3) {\small $7$};
	\node[draw,shape=circle, inner sep=0.2] at (1.4,-0.4) {\small $8$};
	\node[draw,shape=circle, inner sep=0.2] at (-0.4,-0.4) {\small $9$};
	\node[draw,shape=circle, inner sep=0.2] at (0.5,2.35) {\footnotesize $10$};
	\node[draw,shape=circle, inner sep=0.2] at (2.35,0.5) {\footnotesize $11$};
	\node[draw,shape=circle, inner sep=0.2] at (0.5,-1.35) {\footnotesize $12$};
	\node[draw,shape=circle, inner sep=0.2] at (-1.35,0.5) {\footnotesize $13$};
	\node[draw,shape=circle, inner sep=0.2] at (-2.1,0.5) {\footnotesize $14$};
		\end{tikzpicture}\hspace{50pt}
			\begin{tikzpicture}[>=latex,scale=0.7]
            \begin{scope}
				\foreach \a in {0,1,2,3}
                \draw[rotate=90*\a]   
				 (-0.6,0.6)--(0.6,0.6)
				(0,1.57)--(0.6,0.6)--(1.57,0)
                (0,1.57) arc (123.74:-33.74:1.13);
	\node[draw,shape=circle, inner sep=0.2] at (-0.9,0) {\small $1$};
	\node[draw,shape=circle, inner sep=0.2] at (-1,1) {\small $2$};
	\node[draw,shape=circle, inner sep=0.2] at (0,0.9) {\small $3$};
	\node[draw,shape=circle, inner sep=0.2] at (0,0) {\small $4$};
	\node[draw,shape=circle, inner sep=0.2] at (1,1) {\small $5$};
	\node[draw,shape=circle, inner sep=0.2] at (0.9,0) {\small $6$};
	\node[draw,shape=circle, inner sep=0.2] at (0,-0.9) {\small $7$};
	\node[draw,shape=circle, inner sep=0.2] at (1,-1) {\small $8$};
	\node[draw,shape=circle, inner sep=0.2] at (-1,-1) {\small $9$};
	\node[draw,shape=circle, inner sep=0.2] at (2.2,0) {\footnotesize $10$};
			\end{scope}\hspace{50pt}
            \begin{scope}[xshift=6.5cm]
			\foreach \a in {0,1,2,3}
            \draw[rotate=90*\a] (0.6,0.6)--(-0.6,0.6)--(-1.44,1.44)--(-2.35,2.42)--(0.54,3.08)--(2.42,2.35)
				(-0.6,0.6)--(0,1.73)--(0.6,0.6)--(1.44,1.44)--(0,1.73)--(-1.44,1.44)
                (0,1.73)--(0.54,3.08)--(1.44,1.44)
                (-2.35,2.42)--(0.77,4.66)--(2.42,2.35)
                (0.77,4.66)--(0.54,3.08)
                (4.66,-0.77) arc (-4.28:75.57:5.2);
	\node[draw,shape=circle, inner sep=0.2] at (-0.9,0) {\small $1$};
	\node[draw,shape=circle, inner sep=0.2] at (-1.2,0.7) {\small $2$};
	\node[draw,shape=circle, inner sep=0.2] at (-0.7,1.2) {\small $3$};
	\node[draw,shape=circle, inner sep=0.2] at (0,0.9) {\small $4$};
	\node[draw,shape=circle, inner sep=0.2] at (0,0) {\small $5$};
	\node[draw,shape=circle, inner sep=0.2] at (0.7,1.2) {\small $6$};
	\node[draw,shape=circle, inner sep=0.2] at (1.2,0.7) {\small $7$};
	\node[draw,shape=circle, inner sep=0.2] at (0.9,0) {\small $8$};
	\node[draw,shape=circle, inner sep=0.2] at (0,-0.9) {\small $9$};
	\node[draw,shape=circle, inner sep=0.2] at (0.7,-1.2) {\footnotesize $10$};
	\node[draw,shape=circle, inner sep=0.2] at (1.2,-0.7) {\footnotesize $11$};
	\node[draw,shape=circle, inner sep=0.2] at (-1.2,-0.7) {\footnotesize $12$};
	\node[draw,shape=circle, inner sep=0.2] at (-0.7,-1.2) {\footnotesize $13$};
	\node[draw,shape=circle, inner sep=0.2] at (0.7,2) {\footnotesize $14$};
	\node[draw,shape=circle, inner sep=0.2] at (-0.9,2.1) {\footnotesize $15$};
	\node[draw,shape=circle, inner sep=0.2] at (-2.2,1.5) {\footnotesize $16$};
	\node[draw,shape=circle, inner sep=0.2] at (-2,0.6) {\footnotesize $17$};
	\node[draw,shape=circle, inner sep=0.2] at (-2.1,-0.9) {\footnotesize $18$};
	\node[draw,shape=circle, inner sep=0.2] at (-1.5,-2.2) {\footnotesize $19$};
	\node[draw,shape=circle, inner sep=0.2] at (-0.6,-2.1) {\footnotesize $20$};
	\node[draw,shape=circle, inner sep=0.2] at (0.9,-2.1) {\footnotesize $21$};
	\node[draw,shape=circle, inner sep=0.2] at (2.2,-1.5) {\footnotesize $22$};
	\node[draw,shape=circle, inner sep=0.2] at (2,-0.6) {\footnotesize $23$};
	\node[draw,shape=circle, inner sep=0.2] at (2.1,0.9) {\footnotesize $24$};
	\node[draw,shape=circle, inner sep=0.2] at (1.5,2.2) {\footnotesize $25$};
	\node[draw,shape=circle, inner sep=0.2] at (1.08,3.47) {\footnotesize $26$};
	\node[draw,shape=circle, inner sep=0.2] at (-0.16,3.39) {\footnotesize $27$};
	\node[draw,shape=circle, inner sep=0.2] at (-2.3,3.23) {\footnotesize $28$};
	\node[draw,shape=circle, inner sep=0.2] at (-3.47,1.08) {\footnotesize $29$};
	\node[draw,shape=circle, inner sep=0.2] at (-3.39,-0.16) {\footnotesize $30$};
	\node[draw,shape=circle, inner sep=0.2] at (-3.23,-2.3) {\footnotesize $31$};
	\node[draw,shape=circle, inner sep=0.2] at (-1.08,-3.47) {\footnotesize $32$};
	\node[draw,shape=circle, inner sep=0.2] at (0.16,-3.39) {\footnotesize $33$};
	\node[draw,shape=circle, inner sep=0.2] at (2.3,-3.23) {\footnotesize $34$};
	\node[draw,shape=circle, inner sep=0.2] at (3.47,-1.08) {\footnotesize $35$};
	\node[draw,shape=circle, inner sep=0.2] at (3.39,0.16) {\footnotesize $36$};
	\node[draw,shape=circle, inner sep=0.2] at (3.23,2.3) {\footnotesize $37$};
	\node[draw,shape=circle, inner sep=0.2] at (5.05,1.59) {\footnotesize $38$};
             \end{scope} 
			\end{tikzpicture}
		\caption{$\thin\aaa\thin\bbb\thin\cdots=\aaa\bbb^3,\aaa^2\bbb^2,\aaa^3\bbb,\aaa^4\bbb.$}
		\label{ab=r}
	\end{figure}

	\section{Convex case $\ccc<\bbb<\pi$ with some irrational angle}\label{irr2}
	
By Lemma \ref{angle} and \ref{type}, we have $\bbb+\ccc>\pi$, $\frac{\pi}{3}<\aaa<\frac{\pi}{2}<\bbb<\pi$, and $0<\ccc<\bbb$. By Lemma \ref{deg345}, we only need to discuss the vertices in Table \ref{v>} case by case. In each case, we can get the AVC by Lemma \ref{irrational}.

	\subsection{The vertex types of degree 3}\label{deg3}

\subsubsection*{Case $\aaa\ccc^2$}

The AAD $\thin\aaa\thin\ccc^\bbb\thin^\bbb\ccc\thin$ of $\aaa\ccc^2$ gives $\thin\bbb\thin\bbb\thin\cdots$. Similar to the proof of Lemma \ref{deg345}, $\bbb^2\cdots$ can never be a vertex.

\subsubsection*{Case $\bbb\ccc^2$}

	By $3\bbb>\bbb+2\ccc=2\pi$, we have $\pi>\bbb>\frac{2\pi}{3}$, $\aaa<\frac{\pi}{2}<\ccc<\frac{2\pi}{3}<2\aaa$. The vertex $\bbb\ccc^2$ determines $T_1,T_2,T_3$ in the left of Figure \ref{br2abr}; $\bbb^2\cdots=\aaa\bbb^2$ determines $T_4$; $\bbb\ccc\cdots=\bbb\ccc^2$ determines $T_5$. But $\aaa\bbb\ccc\cdots$ appears, contradicting $\bbb\ccc\cdots=\bbb\ccc^2$.

		\begin{figure}[htp]
		\centering
		\begin{tikzpicture}[>=latex,scale=0.9]       
		\draw (0.87,0.5) -- (0,0) -- (-0.87,0.5) (0,-1) -- (0,0);
        \draw (0.87,0.5) -- (0,1) -- (-0.87,0.5) -- (-0.87,-0.5)--(0,-1)--(0.87,-0.5)--(0.87,0.5) (0,1)--(1.59,0.92)--(0.87,-0.5);
        \draw (0.87,-0.5) arc (15.57:-195.57:0.9);
		\node at (0,0.25) {\small $\bbb$};
		\node at (0.24,-0.13) {\small $\ccc$};
		\node at (-0.24,-0.13) {\small $\ccc$};
		\node at (0.9,0.5) {\small $\ccc$};		
        \node at (0.87,-0.5){\small $\bullet$};
	\node[draw,shape=circle, inner sep=0.2] at (0,0.7) {\small $1$};
	\node[draw,shape=circle, inner sep=0.2] at (0.5,-0.26) {\small $2$};
	\node[draw,shape=circle, inner sep=0.2] at (-0.5,-0.26) {\small $3$};
	\node[draw,shape=circle, inner sep=0.2] at (0,-1.25) {\small $4$};
	\node[draw,shape=circle, inner sep=0.2] at (1.2,0.71) {\small $5$};
		\end{tikzpicture}\hspace{50pt}
		\begin{tikzpicture}[>=latex,scale=0.8]       
		\draw (0,0) -- (0,2) -- (1,1) -- (0,0);
        \draw (0,0) -- (3,0) -- (2,1) -- (1,1)
              (2,1) -- (3,2)
              (0,2) -- (3,2) -- (3,0);
        \draw (3,0)--(1.5,-1)--(0,0)--(-1,1)--(0,2); 
		\node at (0.7,1) {\small $\aaa$};
		\node at (1,0.7) {\small $\bbb$};
		\node at (1,1.3) {\small $\ccc$};
		\node at (0,0) {\small $\bullet$};
	\node[draw,shape=circle, inner sep=0.2] at (0.3,1) {\small $1$};
	\node[draw,shape=circle, inner sep=0.2] at (1.5,0.5) {\small $2$};
	\node[draw,shape=circle, inner sep=0.2] at (1.5,1.5) {\small $3$};
	\node[draw,shape=circle, inner sep=0.2] at (2.5,1) {\small $4$};
	\node[draw,shape=circle, inner sep=0.2] at (1.5,-0.5) {\small $5$};
	\node[draw,shape=circle, inner sep=0.2] at (-0.3,1) {\small $6$};
		\end{tikzpicture}
		\caption{$\bbb\ccc^2$ or $\aaa\bbb\ccc$ is a vertex.}
		\label{br2abr}
	\end{figure}

\subsubsection*{Case $\aaa\bbb\ccc$}

Similar to the proof of Lemma \ref{deg345}, $\bbb^2\cdots$ can never be a vertex. So $\thin\ccc\thin\ccc\thin\cdots$ can never be a vertex by Lemma \ref{AAD}. In the right of Figure \ref{br2abr}, $\aaa\bbb\ccc$ determines $T_1,T_2,T_3$, $\thin\bbb\thin\ccc\thin\cdots=\thin\bbb\thin\ccc\thin\aaa\thin$ determines $T_4$, and we get $\aaa\bbb\cdots=\aaa\bbb\ccc,\aaa^2\bbb,\aaa^3\bbb$. 

When one of $\thin\aaa_1\thin\bbb_3\thin\cdots$ and $\thin\aaa_4\thin\bbb_2\thin\cdots$ is $\aaa\bbb\ccc$, the tiling is the same as Figure \ref{abrb>}. When one of them is $\aaa^2\bbb$, we get $\ccc=\aaa$. Then we conclude that either there is no tilings or there is tilings which must be of the icosahedral type by Lemma \ref{ico1}. In a similar discussion later, we will omit the details. When they are $\aaa^3\bbb$, $T_5,T_6$ are determined and we get $\ccc=2\aaa$. But $\aaa^3\ccc\cdots$ can never be a vertex.

\subsubsection*{Case $\bbb^3$}

We have $\bbb=\frac{2}{3}\pi,\ccc>\frac{1}{3}\pi$. By $\bbb^3=\thin\bbb^\ccc\thin^\ccc\bbb\thin\bbb\thin$ and similar proof of Lemma \ref{deg345}, we get $\thin\ccc\thin\ccc\thin\cdots=\aaa\bbb\ccc^2,\bbb\ccc^3,\aaa^3\ccc^2,\aaa^2\ccc^2,\aaa^2\ccc^3,\aaa\ccc^3,\aaa\ccc^4,\ccc^4,\ccc^5$.

	\begin{figure}[htp]
		\centering
			\begin{tikzpicture}[>=latex,scale=0.7]
				\foreach \a in {0,1,2}
				\draw[rotate=120*\a]
				(0,0)--(0,1.5)
				(-1.3,0.75)--(-0.42,2.55)--(2.42,0.91)--(1.3,0.75)
				(-1.3,0.75)--(0,1.5)--(1.3,0.75)
				(0,1.5)--(-0.42,2.55)--(-1.16,3.26)--(-2,1.64)
                (3.41,-0.62) arc (-1.75:101.03:3.84);
		\node at (0,-0.3) {\small $\bbb$};
		\node at (0.2,0.2) {\small $\bbb$};
		\node at (-0.2,0.2) {\small $\bbb$};
		\node at (0.2,1.6) {\small $\bbb$};
		\node at (-1.31,1.1) {\small $\ccc$};
		\node at (-1.6,0.5) {\small $\ccc$};
		\node at (0.3,-1.6) {\small $\ccc$};
		\node at (-0.25,-1.75) {\small $\ccc$};
		\node at (1.5,0.55) {\small $\ccc$};
		\node at (-0.2,2.6) {\small $\ccc$};
		\node at (-2.2,-1.6) {\small $\ccc$};
		\node at (2.4,-1.15) {\small $\ccc$};
	\node[draw,shape=circle, inner sep=0.2] at (0,-0.75) {\small $1$};
	\node[draw,shape=circle, inner sep=0.2] at (0.6,0.375) {\small $2$};
	\node[draw,shape=circle, inner sep=0.2] at (-0.6,0.375) {\small $3$};
	\node[draw,shape=circle, inner sep=0.2] at (-0.5,1.7) {\small $4$};
	\node[draw,shape=circle, inner sep=0.2] at (0.8,1.5) {\small $5$};
	\node[draw,shape=circle, inner sep=0.2] at (-1.65,0) {\small $6$};
	\node[draw,shape=circle, inner sep=0.2] at (-1.3,2) {\small $7$};
	\node[draw,shape=circle, inner sep=0.2] at (-1.2,-1.2) {\small $8$};
	\node[draw,shape=circle, inner sep=0.2] at (0.8,-1.4) {\small $9$};
	\node[draw,shape=circle, inner sep=0.2] at (-1.1,-2) {\footnotesize $10$};
	\node[draw,shape=circle, inner sep=0.2] at (1.65,-0.45) {\footnotesize $11$};
	\node[draw,shape=circle, inner sep=0.2] at (2.5,0) {\footnotesize $12$};
	\node[draw,shape=circle, inner sep=0.2] at (1.5,2) {\footnotesize $13$};
	\node[draw,shape=circle, inner sep=0.2] at (-2.5,0) {\footnotesize $14$};
	\node[draw,shape=circle, inner sep=0.2] at (1.5,-2) {\footnotesize $15$};
	\node[draw,shape=circle, inner sep=0.2] at (2,3) {\footnotesize $16$};
			\end{tikzpicture}
		\caption{$T(4a^3,12a^4;12\aaa\bbb\ccc^2,4\bbb^3)$.}
		\label{b31}
	\end{figure}

When $\thin\ccc\thin\ccc\thin\cdots=\aaa\bbb\ccc^2$, we get the AVC $\subset\{\bbb^3,\aaa\bbb\ccc^2\}$ by Lemma \ref{irrational}. The vertex $\bbb^3$ determines $T_1,T_2,T_3$ in Figure \ref{b31}. By the symmetry of the partial tiling and $\ccc^2\cdots=\aaa\bbb\ccc^2$, we might as well take $\thin\ccc_2\thin\ccc_3\thin\cdots=\thin\ccc_2\thin\ccc_3\thin\aaa\thin\bbb\thin$ which determines $T_4,T_5$; $\aaa\bbb\cdots=\aaa\bbb\ccc^2$ determines $T_6,T_7$. Similarly, $T_8,T_9,\cdots,T_{16}$ are determined and the 3D picture is the 15th of Figure \ref{Fig3}.\label{4,12}

	\begin{table}[htp]                        
		\centering     
		\begin{tabular}{|c|c|c|}\hline	 			
				$\thin\ccc\thin\ccc\thin\cdots$ & AVC $\subset$ a set  & Conclusion \\
				\hline  
				$\bbb\ccc^3$&$\{\bbb^3,\bbb\ccc^3\}$& \multirow{5}{*}{No $\aaa\bbb\cdots$, contradicting Lemma \ref{AAD}} \\
				\cline{1-2}
				$\aaa^3\ccc^2$&$\{\bbb^3,\aaa^3\ccc^2\}$&\multirow{5}{*}{} \\
				\cline{1-2} 
				$\aaa^2\ccc^k,k=2,3$&$\{\bbb^3,\aaa^2\ccc^k\}$&\multirow{5}{*}{} \\			
				\cline{1-2}
				$\aaa\ccc^k,k=3,4$&$\{\bbb^3,\aaa\ccc^k\}$&\multirow{5}{*}{} \\
				\cline{1-2} 
				$\ccc^k,k=4,5$&$\{\bbb^3,\ccc^k\}$&\multirow{5}{*}{} \\
				\hline 
		\end{tabular}
	\caption{The simple cases when $\bbb^3$ appears.}\label{Tab-1}  
	\end{table}

\subsubsection*{Case $\bbb^2\ccc$} \label{b2c}
By $\bbb^2\ccc=\thin\ccc^\bbb\thin^\ccc\bbb\thin\bbb\thin$ and similar proof of Lemma \ref{deg345}, we get $\thin\bbb\thin\ccc\thin\cdots=\bbb^2\ccc,\aaa\bbb\ccc^h(h\ge3),\aaa^2\bbb\ccc^j(j$ $\ge2),\bbb\ccc^k(k\ge4)$.

\vspace{9pt}
\noindent{Subcase $\thin\bbb\thin\ccc\thin\cdots=\bbb\ccc^k(k\ge4)$, AVC $\subset\{\bbb^2\ccc,\bbb\ccc^k,\ccc^{2k-1}\}$}

There is no the vertex $\aaa\bbb\cdots$, contradicting Lemma \ref{AAD}.

\vspace{9pt}
\noindent{Subcase $\thin\bbb\thin\ccc\thin\cdots=\aaa\bbb\ccc^h(h\ge3)$, AVC $\subset\{\bbb^2\ccc,\aaa\bbb\ccc^h,\aaa^2\ccc^{2h-1}\}$}

When $\thin\ccc\thin\aaa\thin\ccc\thin\cdots$ appears, $T_1,T_2,T_3$ are determined in Figure \ref{rar1}. Then one of $\thin\aaa_1\thin\bbb_2\thin\cdots$ and $\thin\aaa_1\thin\bbb_3\thin\cdots$ must be $\aaa\bbb^2\cdots$, contradicting the AVC. So $\aaa\bbb\ccc^h=\thin\aaa\thin\bbb\thin\ccc\thin\cdots\thin\ccc\thin$. 

		\begin{figure}[htp]
			\centering
			\begin{tikzpicture}[>=latex,scale=0.8] 
			\draw (1.5,0.13)--(1,1)--(-1,1)--(-1.5,0.13)--(1.5,0.13) (0.5,0.13)--(0,1)--(-0.5,0.13);
			\node at (-0.37,0.77){\small $\ccc$};
			\node at (0.37,0.77){\small $\ccc$};
	\node[draw,shape=circle, inner sep=0.2] at (0,0.5) {\small $1$};
	\node[draw,shape=circle, inner sep=0.2] at (0.7,0.5) {\small $2$};
	\node[draw,shape=circle, inner sep=0.2] at (-0.7,0.5) {\small $3$};
			\end{tikzpicture}
			\caption{$\thin\ccc\thin\aaa\thin\ccc\thin\cdots$ appears.}\label{rar1}
		\end{figure}

If $h=3$, then $\thin\aaa\thin\bbb\thin\ccc\thin\ccc\thin\ccc\thin$ determines $T_1,\cdots,T_5$ in Figure \ref{abr3}; $\aaa\bbb\cdots=\aaa\bbb\ccc^3$ determines $T_6,\cdots,T_{10}$; $\bbb^2\cdots=\bbb^2\ccc$ determines $T_{11},\cdots,T_{16}$; $\bbb\ccc^3\cdots=\aaa\bbb\ccc^3$ determines $T_{17}$. The 3D picture is the 4th of Figure \ref{Fig2}. If $h>3$, then we have the same tilings $T(2a^3,xa^4;6\aaa\bbb\ccc^h,y\bbb^2\ccc)$. By $y+3=x$ and $3h+\frac{y}{2}=x$, we get $x=6h-3,y=6h-6$. 
	\begin{figure}[htp]
		\centering
		\begin{tikzpicture}[>=latex,scale=0.7]
				\foreach \a in {0,1,2}
				\draw[rotate=120*\a]
				(0.6,-0.35)--(0,0.69)--(0,1.89)--(0.92,2.66)--(0.92,1.46)--(2.1,1.26)--(1.18,0.48)--(1.64,-0.95)--(0.6,-0.35) (1.18,0.48)--(0,0.69)--(0.92,1.46)
                (1.85,-2.13)arc(-19.03:10.47:6.67) (1.85,-2.13)arc(-68.82:90.47:2.48);
		\node at (-0.2,0.7) {\small $\bbb$};
		\node at (0.25,1.2) {\small $\ccc$};
		\node at (0.45,0.8) {\small $\ccc$};
		\node at (0.4,0.3) {\small $\ccc$};
		\node at (0.9,-0.82) {\small $\ccc$};
		\node at (0.43,-0.95) {\small $\ccc$};
		\node at (0.1,-0.6) {\small $\ccc$};
		\node at (-1.18,-0.4) {\small $\ccc$};
		\node at (-1,0.1) {\small $\ccc$};
		\node at (1.2,1.6) {\small $\ccc$};
		\node at (1.5,0.45) {\small $\ccc$};
		\node at (-0.47,-1.51) {\small $\ccc$};
		\node at (0.88,-2) {\small $\ccc$};
		\node at (-1.1,1.1) {\small $\ccc$};
		\node at (-2.1,0.15) {\small $\ccc$};
	\node[draw,shape=circle, inner sep=0.2] at (0,0) {\small $1$};
	\node[draw,shape=circle, inner sep=0.2] at (-0.6,0.75) {\small $2$};
	\node[draw,shape=circle, inner sep=0.2] at (0.48,1.77) {\small $3$};
	\node[draw,shape=circle, inner sep=0.2] at (1.2,1) {\small $4$};
	\node[draw,shape=circle, inner sep=0.2] at (1,-0.15) {\small $5$};
	\node[draw,shape=circle, inner sep=0.2] at (1.35,-1.4) {\small $6$};
	\node[draw,shape=circle, inner sep=0.2] at (0.27,-1.58) {\small $7$};
	\node[draw,shape=circle, inner sep=0.2] at (-0.6,-0.85) {\small $8$};
	\node[draw,shape=circle, inner sep=0.2] at (-1.76,-0.48) {\small $9$};
	\node[draw,shape=circle, inner sep=0.2] at (-1.49,0.55) {\footnotesize $10$};
	\node[draw,shape=circle, inner sep=0.2] at (1.83,-0.24) {\footnotesize $11$};
	\node[draw,shape=circle, inner sep=0.2] at (2.79,0.63) {\footnotesize $12$};
	\node[draw,shape=circle, inner sep=0.2] at (-1.05,-1.5) {\footnotesize $13$};
	\node[draw,shape=circle, inner sep=0.2] at (-0.55,-2.78) {\footnotesize $14$};
	\node[draw,shape=circle, inner sep=0.2] at (-0.75,1.6) {\footnotesize $15$};
	\node[draw,shape=circle, inner sep=0.2] at (-2,2) {\footnotesize $16$};
	\node[draw,shape=circle, inner sep=0.2] at (1.35,3.25) {\footnotesize $17$};
		\end{tikzpicture}
		\caption{$T(2a^3,15a^4;6\aaa\bbb\ccc^3,12\bbb^2\ccc).$}
		\label{abr3}
	\end{figure}

By $\ccc=2\pi-2\bbb,\aaa=(1-h)2\pi+(2h-1)\bbb$ and (\ref{tr4}), we obtain $2\cos \aaa\cos \bbb-\cos \bbb-\cos \aaa=\cos \aaa(2\cos \bbb-1)-\cos \bbb=0$. Let $f(h,\bbb)=\cos \aaa(2\cos \bbb-1)-\cos \bbb$, we have $f'_\bbb=-(2h-1)\sin \aaa(2\cos \bbb-1)+\sin \bbb(1-2\cos \aaa)$.

Since $\aaa\in(\frac{\pi}{3},\frac{\pi}{2})$, we get $\bbb\in(\bbb_{min},\bbb_{max})$, where $\bbb_{min}=\frac{\frac{\pi}{3}+(h-1)2\pi}{2h-1}=\pi-\frac{\frac{2}{3}\pi}{2h-1}>\frac{2}{3}\pi,\bbb_{max}=\frac{\frac{\pi}{2}+(h-1)2\pi}{2h-1}=\pi-\frac{\frac{\pi}{2}}{2h-1}<\pi$. Obviously, $f'_\bbb>0,f(h,\bbb_{min})<0,f(h,\bbb_{max})>0$. This implies that $\forall h\ge3,f(h,\bbb)=0$ has a unique solution in $(\bbb_{min},$ $\bbb_{max})$. Moreover,
\[
\frac{d\bbb}{dh}=-\frac{f'_h}{f'_\bbb}=-\frac{-(2\cos \bbb-1)\sin \aaa(2\bbb-2\pi)}{f'_\bbb}>0.
\]
This means $\bbb$ is monotonically increasing with respect to $h$. 

By $\lim\limits_{h\to+\infty}\beta_{min}=\lim\limits_{h\to+\infty}\beta_{max}=\pi$, we get $\lim\limits_{h\to+\infty}\beta=\pi$, $\lim\limits_{h\to+\infty}\ccc=0$. We observe that there are three $a$-edges between the poles, so $\lim\limits_{h\to+\infty}a=\frac{\pi}{3}$ and $\lim\limits_{h\to+\infty}\aaa=\arccos\frac{1}{3}$.

\vspace{9pt}
\noindent{Subcase $\thin\bbb\thin\ccc\thin\cdots=\aaa^2\bbb\ccc^j(j\ge2)$, AVC $\subset\{\bbb^2\ccc,\aaa^2\bbb\ccc^j,\aaa^4\ccc^{2j-1}\}$}

As shown in figure \ref{rar1}, when $\thin\ccc\thin\aaa\thin\ccc\thin\cdots$ appears, we get two $\thin\aaa\thin\bbb\thin\cdots=\thin\ccc\thin\aaa\thin\aaa\thin\bbb\thin\cdots$. But $\aaa\bbb^2\cdots$ appears, contradicting the AVC. So $\thin\bbb\thin\ccc\thin\cdots=\aaa^2\bbb\ccc^j=\thin\aaa\thin\aaa\thin\ccc\thin$ $\cdots$ determines $T_1,T_2,T_3$ in Figure \ref{ba}. Then $\thin\bbb_3\thin\aaa_2\thin\cdots$ is either $\thin\bbb_3\thin\aaa_2\thin\aaa\thin\cdots$ or $\thin\bbb_3\thin\aaa_2\thin\ccc\thin\cdots$.

	\begin{figure}[htp]
		\centering
		\begin{tikzpicture}[>=latex,scale=0.8]       
		\draw (-1,0) -- (-1.50,-0.87) -- (2.50,-0.87) -- (2,-1.73)  (-1,0) -- (1,0);
        \draw (0,-1.73) -- (2,-1.73)  (0,0) -- (1.5,-2.6) -- (0.5,-2.6) (1,-1.73) -- (1.5,-0.87)
              (-0.5,-0.87) -- (0.5,-2.6) (-1,-1.73) -- (0,0)
              (-1,-1.73) -- (-0.5,-2.6) -- (1,0);
		\node at (-0.3,-0.2) {\small $\ccc$};
		\node at (-0.5,-1.2) {\small $\ccc$};
		\node at (0.2,-2) {\small $\ccc$};
		\node at (1.3,-1.5) {\small $\ccc$};
		\node at (0.5,-0.87) {$\bullet$};
	\node[draw,shape=circle, inner sep=0.2] at (0.5,-0.3) {\small $1$};
	\node[draw,shape=circle, inner sep=0.2] at (0,-0.5) {\small $2$};
	\node[draw,shape=circle, inner sep=0.2] at (-0.7,-0.5) {\small $3$};
	\node[draw,shape=circle, inner sep=0.2] at (0,-1.2) {\small $4$};
	\node[draw,shape=circle, inner sep=0.2] at (-0.5,-1.8) {\small $5$};
	\node[draw,shape=circle, inner sep=0.2] at (0.5,-1.4) {\small $6$};
	\node[draw,shape=circle, inner sep=0.2] at (0.7,-2.2) {\small $7$};
	\node[draw,shape=circle, inner sep=0.2] at (1,-1.2) {\small $8$};
	\node[draw,shape=circle, inner sep=0.2] at (1.8,-1.2) {\small $9$};
		\end{tikzpicture}\hspace{50pt}
		\begin{tikzpicture}[>=latex,scale=0.8]       
		\draw (-1,0) -- (2,0) -- (1.52,-0.73) --(0.5,-0.87)-- (-1.5,-0.87) -- (-1,0) (0,0) -- (0.5,-0.87);
        \draw (0.5,0.87) -- (-0.5,0.87) -- (0,0) -- (-0.5,-0.87) -- (0,-1.73) -- (1,-1.73) -- (0.5,-0.87) -- (1,0) -- (0.5,0.87) (3,0) -- (2,0) -- (1.5,0.87) -- (1,0) (1.5,0.87) -- (1.5,1.87) (2,0) -- (2.5,0.87) -- (1.5,0.87)
              (0.5,-0.87)--(1.41,-1.47)--(2.39,-1.21)--(1.52,-0.73) (3,0)--(2.39,-1.21)
              (1,0)--(0.72,1.17)--(1.5,1.87)
              (1.5,0.87)--(2.14,1.77)--(3.2,1.77)--(2.5,0.87)--(3.68,1.11)--(3.12,0.21)--(2,0);
		\node at (-0.3,-0.2) {\small $\ccc$};
		\node at (-0.1,-1.1) {\small $\ccc$};
		\node at (0.9,-0.6) {\small $\ccc$};
		\node at (0.7,0.2) {\small $\ccc$};
		\node at (1.1,-1.05) {\small $\ccc$};
		\node at (1.07,0.4) {\small $\ccc$};
		\node at (1.9,-0.7) {\small $\ccc$};
		\node at (1.93,1.1) {\small $\ccc$};
		\node at (2.41,0.29) {\small $\ccc$};
		\node at (2.5,0.87) {$\bullet$};
	\node[draw,shape=circle, inner sep=0.2] at (0.5,-0.3) {\small $1$};
	\node[draw,shape=circle, inner sep=0.2] at (0,-0.5) {\small $2$};
	\node[draw,shape=circle, inner sep=0.2] at (-0.7,-0.5) {\small $3$};
	\node[draw,shape=circle, inner sep=0.2] at (0.2,-1.35) {\small $4$};
	\node[draw,shape=circle, inner sep=0.2] at (1.28,-0.4) {\small $5$};
	\node[draw,shape=circle, inner sep=0.2] at (1.6,-1.13) {\small $6$};
	\node[draw,shape=circle, inner sep=0.2] at (2.3,-0.5) {\small $7$};
	\node[draw,shape=circle, inner sep=0.2] at (0.2,0.5) {\small $8$};
	\node[draw,shape=circle, inner sep=0.2] at (1.5,0.3) {\small $9$};
	\node[draw,shape=circle, inner sep=0.2] at (1.15,1) {\footnotesize $10$};
	\node[draw,shape=circle, inner sep=0.2] at (2,0.5) {\footnotesize $11$};
	\node[draw,shape=circle, inner sep=0.2] at (2.35,1.32) {\footnotesize $12$};
	\node[draw,shape=circle, inner sep=0.2] at (2.91,0.55) {\footnotesize $13$};
		\end{tikzpicture}
		\caption{$\thin\aaa\thin\aaa\thin\ccc\thin\cdots$ appears.}
		\label{ba}
	\end{figure}

When $\thin\bbb_3\thin\aaa_2\thin\aaa\thin\cdots$ appears, $\thin\bbb_3\thin\aaa_2\thin\aaa\thin\cdots=\thin\bbb_3\thin\aaa_2\thin\aaa\thin\ccc\thin\cdots$ determines $T_4,T_5$ in the first; $\thin\bbb_5\thin\aaa_4\thin\cdots=\thin\bbb_5\thin\aaa_4\thin\aaa\thin\ccc\thin\cdots$ determines $T_6,T_7$. Similarly, $T_8,T_9$ are determined but $\aaa^5\cdots$ appears, contradicting the AVC. 

When $\thin\bbb_3\thin\aaa_2\thin\ccc\thin\cdots$ appears, $T_4$ is determined in the second. Then $\thin\bbb_4\thin\aaa_2\thin\aaa_1\thin\cdots=\aaa^2\bbb\ccc^j$ determines $T_5,T_6$; $\bbb^2\cdots=\bbb^2\ccc$ determuines $T_7$. Since $\thin\aaa_2\thin\aaa_1\thin\ccc_5\thin\cdots$ appears, $\thin\bbb_5\thin\aaa_1\thin\cdots=\thin\bbb_5\thin\aaa_1\thin\ccc\thin\cdots=\thin\ccc\thin\aaa\thin\bbb_5\thin\aaa_1\thin\ccc\thin\cdots$ determines $T_8,T_9,T_{10}$. Similarly, we can determine $T_{11},T_{12},T_{13}$ but $\aaa\bbb^2\cdots$ appears, contradicting the AVC.

\vspace{9pt}
\noindent{Subcase $\thin\bbb\thin\ccc\thin\cdots$ can only be $\bbb^2\ccc$}

Then we have the monohedral tiling in Figure \ref{bcdpic}, a contradiction.

		\begin{figure}[htp]
			\centering
			\begin{tikzpicture}[>=latex,scale=0.6] 
			\foreach \a in {0,1,2}
			{
				\begin{scope}[xshift=2*\a cm] 
				\draw (0,0)--(0,-2)
				(2,0)--(2,-2)
				(0.5,-3)--(0.5,-5);
				\draw  (0,-2)--(0.5,-3);
				\draw (0.5,-3)--(2,-2);
				\node at (0.6,-2.55){\small $\ccc$};
				\node at (0.2,-3.2){\small $\bbb$};
				\node at (0.8,-3.2){\small $\bbb$};
				\node at (1,0){\small $\ccc$};
				\node at (1.5,-5){\small $\ccc$};
				
				\node at (1.9,-2.45){\small $\ccc$};
				\node at (1.7,-1.8){\small $\bbb$};
				\node at (2.3,-1.8){\small $\bbb$};
				\end{scope}
			}
			
			\draw (6,-2)--(6.5,-3);
			\draw  (6.5,-3)--(6.5,-5);

			\fill (8,-2) circle (0.05); 
			\fill (8.2,-2) circle (0.05);
			\fill (8.4,-2) circle (0.05);
			
			\node[draw,shape=circle, inner sep=0.5] at (1,-1) {\small $1$};
			\node[draw,shape=circle, inner sep=0.5] at (3,-1) {\small $2$};
			\node[draw,shape=circle, inner sep=0.5] at (5,-1) {\small $5$};
			\node[draw,shape=circle, inner sep=0.5] at (1.5,-4) {\small $3$};
			\node[draw,shape=circle, inner sep=0.5] at (3.5,-4) {\small $4$};
			\node[draw,shape=circle, inner sep=0.5] at (5.5,-4) {\small $6$};
			
			\node at (6.2,-3.2){\small $\bbb$};
			\node at (0.3,-1.8){\small $\bbb$};
			
			\end{tikzpicture}
			\caption{$T(f_\dia  \, \bbb^2\ccc, 2\ccc^{f_\dia /2})$.} \label{bcdpic}
		\end{figure}

\subsubsection*{Case $\aaa\bbb^2$}
By $2\bbb+2\ccc>2\pi=\aaa+2\bbb$, we have $2\ccc>\aaa>\frac{1}{3}\pi$ $i.e.$ $\ccc>\frac{1}{6}\pi$. The AAD $\thin\aaa\thin\bbb^\ccc\thin^\ccc\bbb\thin$ of $\aaa\bbb^2$ gives $\thin\ccc\thin\ccc\thin\cdots$. Similar to the proof of Lemma \ref{deg345}, all possibilities of $\thin\ccc\thin\ccc\thin\cdots$ are listed in Table \ref{rr}.

\begin{table}[H]\small
\centering
\begin{tabular}{|c|}
\hline 
$\bbb\ccc^4,\bbb\ccc^5,\bbb\ccc^6,\aaa\bbb\ccc^3,\aaa\bbb\ccc^4,\aaa^2\bbb\ccc^2$ \\
\hline 
$\ccc^3,\ccc^4,\ccc^6,\ccc^7,\ccc^8,\ccc^9,\ccc^{10},\ccc^{11}$ \\
\hline   
$\aaa\ccc^3,\aaa\ccc^5,\aaa\ccc^6,\aaa\ccc^7,\aaa\ccc^8,\aaa\ccc^9$\\
\hline
$\aaa^2\ccc^2,\aaa^2\ccc^4,\aaa^2\ccc^5,\aaa^2\ccc^6,\aaa^2\ccc^7$\\
\hline
$\aaa^3\ccc^3,\aaa^3\ccc^4,\aaa^3\ccc^5$\\
\hline
$\aaa^4\ccc^2,\aaa^4\ccc^3$\\
\hline
\end{tabular}
\caption{All possibilities of $\thin\ccc\thin\ccc\thin\cdots$.}
\label{rr}
\end{table}

	\begin{table}[htp]\small                   
		\centering     
		\begin{tabular}{|c|c|c|}\hline	 			
				$\thin\ccc\thin\ccc\thin\cdots$ & AVC $\subset$ a set  & Conclusion \\
				\hline  
				$\bbb\ccc^k,k=4,5,6$&$\{\aaa\bbb^2,\bbb\ccc^k\}$&\multirow{2}{*}{No $\aaa\ccc\cdots$, contradicting Lemma \ref{AAD}} \\
				\cline{1-2} 
				$\ccc^k,k=3,4,6,\cdots,11$&$\{\aaa\bbb^2,\ccc^k\}$&\multirow{2}{*}{} \\			
				\hline
		\end{tabular}
	\caption{The simple cases when $\aaa\bbb^2$ appears.}\label{Tab-2}  
	\end{table}

\vspace{9pt}
\noindent{Subcase $\thin\ccc\thin\ccc\thin\cdots=\aaa\bbb\ccc^k (k=3,4)$, AVC $\subset \{\aaa\bbb^2,\aaa\bbb\ccc^k,\aaa\ccc^{2k}\}$}

As shown in Figure \ref{rar1}, when $\thin\ccc\thin\aaa\thin\ccc\thin\cdots$ appears, we get two $\thin\aaa\thin\bbb\thin\cdots$ and one of them must be $\thin\aaa\thin\bbb\thin\ccc\thin\cdots\thin\ccc\thin$. Since $\thin\ccc\thin\ccc\thin\cdots$ is either $\thin\ccc\thin\aaa\thin\ccc\thin\cdots$ or $\thin\aaa\thin\bbb\thin\ccc\thin\cdots\thin\ccc\thin$, there must be $\thin\ccc\thin\ccc\thin\ccc\thin$ $\cdots$ which determines $T_1,T_2,T_3$ in Figure \ref{rar}. Then $\bbb^2\cdots=\aaa\bbb^2$ determines $T_4,T_5$. But $\aaa^2\ccc\cdots$ appears, contradicting the AVC.

		\begin{figure}[htp]
			\centering
			\begin{tikzpicture}[>=latex,scale=0.8] 
				\foreach \a in {0,1,2}
				\draw[rotate=45.43*\a]
                 (0,0)--(1.14,0.46)--(1.62,1.6)--(0.48,1.14)--(0,0);
            \draw (-1.62,1.6)--(0,2.28)--(1.62,1.6);
			\node at (-0.42,0.4){\small $\ccc$};
			\node at (0,0.53){\small $\ccc$};
			\node at (0.42,0.4){\small $\ccc$};
			\node at (0,2.28){\small $\bullet$};
	\node[draw,shape=circle, inner sep=0.2] at (-0.81,0.8) {\small $1$};
	\node[draw,shape=circle, inner sep=0.2] at (0,1.14) {\small $2$};
	\node[draw,shape=circle, inner sep=0.2] at (0.81,0.8) {\small $3$};
	\node[draw,shape=circle, inner sep=0.2] at (-0.63,1.61) {\small $4$};
	\node[draw,shape=circle, inner sep=0.2] at (0.63,1.61) {\small $5$};
			\end{tikzpicture}
			\caption{$\thin\ccc\thin\ccc\thin\ccc\thin\cdots$ appears.}\label{rar}
		\end{figure}

\vspace{9pt}
\noindent{Subcase $\thin\ccc\thin\ccc\thin\cdots=\aaa^2\bbb\ccc^2$, AVC $\subset\{\aaa\bbb^2,\aaa^2\bbb\ccc^2,\aaa^3\ccc^4\}$} 

Similarly, when $\thin\ccc\thin\aaa\thin\ccc\thin\cdots$ appears, one of the two $\thin\aaa\thin\bbb\thin\cdots$ must be $\aaa^2\bbb\ccc^2$. Moreover, the vertex $\aaa^3\ccc^4$ has four arrangements in Figure \ref{a3r4}. In the first, if $\thin\aaa_{10}\thin\ccc_5\thin\aaa_9\thin\cdots=\thin\aaa_{10}\thin\ccc_5\thin\aaa_9\thin\ccc\thin\cdots$, then $\aaa^2\bbb\ccc^2$ appears; if $\thin\aaa_{10}\thin\ccc_5\thin\aaa_9\thin$ $\cdots=\thin\aaa_{10}\thin\ccc_5\thin\aaa_9\thin\bbb\thin\cdots$, then it is $\aaa^2\bbb\ccc^2$; if $\thin\aaa_{10}\thin\ccc_5\thin\aaa_9\thin\cdots=\thin\aaa_{10}\thin\ccc_5\thin\aaa_9\thin\aaa\thin\cdots=\thin\ccc\thin\aaa_{10}\thin\ccc_5\thin\aaa_9\thin\aaa\thin\ccc\thin\ccc\thin$, then $\aaa^2\bbb\ccc^2$ appears. Since the others have $\thin\ccc\thin\aaa\thin\ccc\thin\cdots$, there must be the vertex $\aaa^2\bbb\ccc^2$.
	\begin{figure}[htp]
		\centering
		\begin{tikzpicture}[>=latex,scale=0.8]       
		\draw (-1,0) -- (1,0)
              (-0.71,0.71) -- (0,0) -- (-0.5,-0.87) 
              (0.71,0.71) -- (0,0) -- (0.5,-0.87)
              (0,1) -- (-0.71,1.71) -- (-0.71,0.71) -- (-1.71,0.71) -- (-1,0) -- (-0.5,-0.87) -- (0.5,-0.87) -- (1,0) -- (1.71,0.71) -- (0.71,0.71) -- (0.71,1.71) -- (0,1) --(0,0)
              (-1.71,0.71)--(-0.71,1.71)--(0.71,1.71)--(1.71,0.71);
		\node at (-0.45,0.2) {\small $\ccc$};
		\node at (-0.2,0.4) {\small $\ccc$};
		\node at (0.45,0.2) {\small $\ccc$};
		\node at (0.2,0.4) {\small $\ccc$};
	\node[draw,shape=circle, inner sep=0.2] at (0.5,-0.3) {\small $1$};
	\node[draw,shape=circle, inner sep=0.2] at (0,-0.5) {\small $2$};
	\node[draw,shape=circle, inner sep=0.2] at (-0.5,-0.3) {\small $3$};
	\node[draw,shape=circle, inner sep=0.2] at (-0.85,0.3) {\small $4$};
	\node[draw,shape=circle, inner sep=0.2] at (-0.3,0.8) {\small $5$};
	\node[draw,shape=circle, inner sep=0.2] at (0.3,0.8) {\small $6$};
	\node[draw,shape=circle, inner sep=0.2] at (0.85,0.3) {\small $7$};
	\node[draw,shape=circle, inner sep=0.2] at (1.1,1) {\small $8$};
	\node[draw,shape=circle, inner sep=0.2] at (0,1.45) {\small $9$};
	\node[draw,shape=circle, inner sep=0.2] at (-1,1) {\footnotesize $10$};
		\end{tikzpicture}\hspace{15pt}
		\begin{tikzpicture}[>=latex,scale=0.8]       
		\draw (-1,0) -- (0,0) -- (0.97,-0.26)
              (-0.71,0.71) -- (0,0) -- (-0.5,-0.87) 
              (0.71,0.71) -- (0,0) -- (0.5,-0.87)
              (0,1) -- (-0.71,1.71) -- (-0.71,0.71) -- (-1.71,0.71) -- (-1,0) -- (-0.5,-0.87) -- (0.5,-0.87) -- (1.47,-1.12) -- (0.97,-0.26) -- (0.71,0.71) -- (0.71,1.71) -- (0,1) --(0,0)
              (-1.71,0.71)--(-0.71,1.71)--(0.71,1.71);
		\node at (-0.45,0.2) {\small $\ccc$};
		\node at (-0.25,0.4) {\small $\ccc$};
		\node at (0.25,-0.25) {\small $\ccc$};
		\node at (0.2,0.4) {\small $\ccc$};
	\node[draw,shape=circle, inner sep=0.2] at (0.7,-0.5) {\small $1$};
	\node[draw,shape=circle, inner sep=0.2] at (0,-0.5) {\small $2$};
	\node[draw,shape=circle, inner sep=0.2] at (-0.5,-0.3) {\small $3$};
	\node[draw,shape=circle, inner sep=0.2] at (-0.85,0.3) {\small $4$};
	\node[draw,shape=circle, inner sep=0.2] at (-0.3,0.8) {\small $5$};
	\node[draw,shape=circle, inner sep=0.2] at (0.3,0.8) {\small $6$};
	\node[draw,shape=circle, inner sep=0.2] at (0.53,0.2) {\small $7$};
	\node[draw,shape=circle, inner sep=0.2] at (-1.1,1) {\small $8$};
	\node[draw,shape=circle, inner sep=0.2] at (0,1.45) {\small $9$};
		\end{tikzpicture}\hspace{15pt}
		\begin{tikzpicture}[>=latex,scale=0.8]       
		\draw (-1,0) -- (0,0) -- (0.97,-0.26)
              (-0.71,0.71) -- (0,0) -- (-0.5,-0.87) 
              (0.87,0.5) -- (0,0) -- (0.5,-0.87)
              (0,1) -- (-0.71,1.71) -- (-0.71,0.71) -- (-1.71,0.71) -- (-1,0) -- (-0.5,-0.87) -- (0.5,-0.87) -- (1.47,-1.12) -- (0.97,-0.26) -- (1.83,0.24) -- (0.87,0.5) -- (0,1) --(0,0)
              (-1.71,0.71)--(-0.71,1.71) (1.83,0.24)--(1.47,-1.12);
		\node at (-0.45,0.2) {\small $\ccc$};
		\node at (-0.2,0.4) {\small $\ccc$};
		\node at (0.25,-0.25) {\small $\ccc$};
		\node at (0.4,0) {\small $\ccc$};
	\node[draw,shape=circle, inner sep=0.2] at (0.7,-0.5) {\small $1$};
	\node[draw,shape=circle, inner sep=0.2] at (0,-0.5) {\small $2$};
	\node[draw,shape=circle, inner sep=0.2] at (-0.5,-0.3) {\small $3$};
	\node[draw,shape=circle, inner sep=0.2] at (-0.85,0.3) {\small $4$};
	\node[draw,shape=circle, inner sep=0.2] at (-0.3,0.8) {\small $5$};
	\node[draw,shape=circle, inner sep=0.2] at (0.3,0.5) {\small $6$};
	\node[draw,shape=circle, inner sep=0.2] at (0.85,0.15) {\small $7$};
	\node[draw,shape=circle, inner sep=0.2] at (-1.1,1) {\small $8$};
	\node[draw,shape=circle, inner sep=0.2] at (1.4,-0.4) {\small $9$};
		\end{tikzpicture}\hspace{15pt}
		\begin{tikzpicture}[>=latex,scale=0.8]       
		\draw (-1,0) -- (0,0) -- (0.97,-0.26)
              (-0.71,0.71) -- (0,0) -- (-0.5,-0.87) 
              (0.87,0.5) -- (0,0) -- (0.26,-0.97)
              (0,1) -- (-0.71,1.71) -- (-0.71,0.71) -- (-1.71,0.71) -- (-1,0) -- (-0.5,-0.87) -- (-0.24,-1.83) -- (0.26,-0.97) -- (0.97,-0.26) -- (1.83,0.24) -- (0.87,0.5) -- (0,1) --(0,0)
              (-1.71,0.71)--(-0.71,1.71);
		\node at (-0.45,0.2) {\small $\ccc$};
		\node at (-0.2,0.4) {\small $\ccc$};
		\node at (0.4,0) {\small $\ccc$};
		\node at (-0.1,-0.45) {\small $\ccc$};
	\node[draw,shape=circle, inner sep=0.2] at (0.35,-0.4) {\small $1$};
	\node[draw,shape=circle, inner sep=0.2] at (-0.1,-0.9) {\small $2$};
	\node[draw,shape=circle, inner sep=0.2] at (-0.5,-0.3) {\small $3$};
	\node[draw,shape=circle, inner sep=0.2] at (-0.85,0.3) {\small $4$};
	\node[draw,shape=circle, inner sep=0.2] at (-0.3,0.8) {\small $5$};
	\node[draw,shape=circle, inner sep=0.2] at (0.3,0.5) {\small $6$};
	\node[draw,shape=circle, inner sep=0.2] at (0.85,0.15) {\small $7$};
	\node[draw,shape=circle, inner sep=0.2] at (-1.1,1) {\small $8$};
		\end{tikzpicture}
		\caption{The vertex arrangements for $\aaa^3\ccc^4$.}
		\label{a3r4}
	\end{figure}

\begin{figure}[H]
		\centering
		\includegraphics[scale=0.1675]{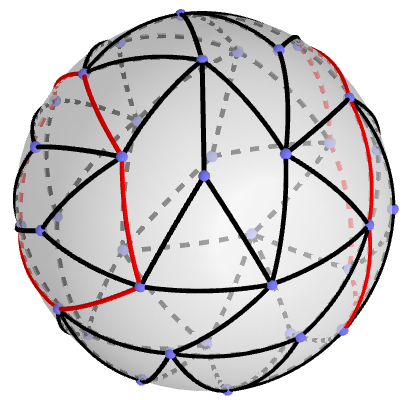}\hspace{0.5pt}
		\includegraphics[scale=0.182]{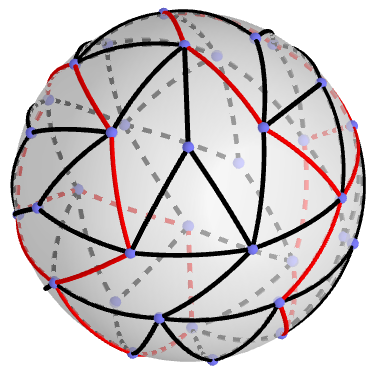}
		\includegraphics[scale=0.182]{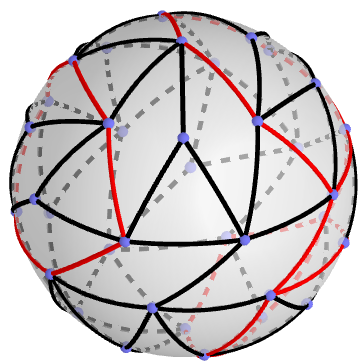}
		\includegraphics[scale=0.206]{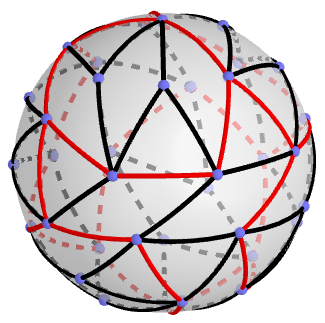}
		\includegraphics[scale=0.2]{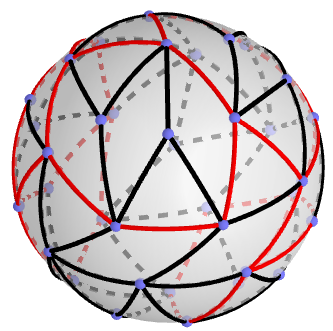}
		\includegraphics[scale=0.195]{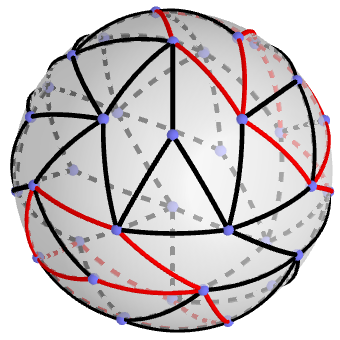}
		\includegraphics[scale=0.18]{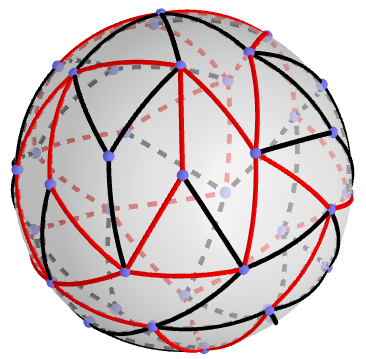}
		\includegraphics[scale=0.176]{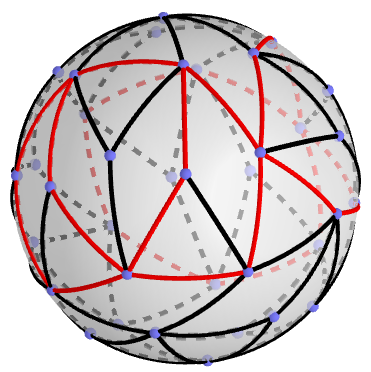}\hspace{3pt}

		\includegraphics[scale=0.203]{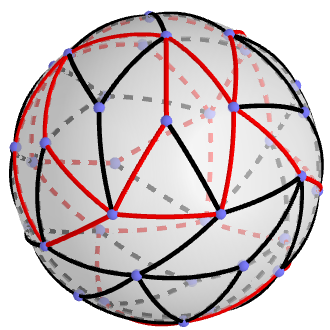}\hspace{0.6pt}
		\includegraphics[scale=0.184]{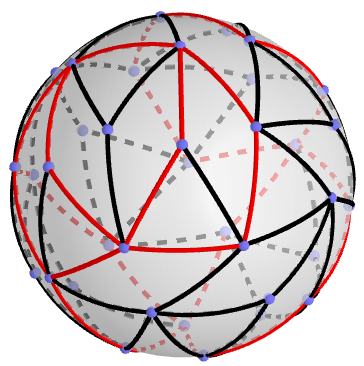}\hspace{0.6pt}
		\includegraphics[scale=0.184]{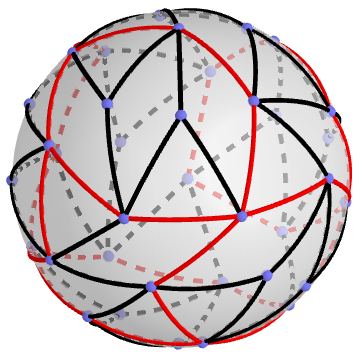}
		\includegraphics[scale=0.184]{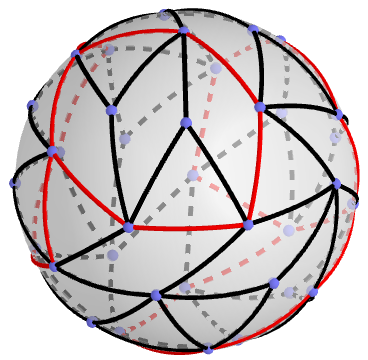}
		\includegraphics[scale=0.185]{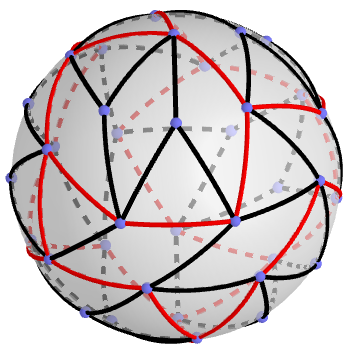}
		\includegraphics[scale=0.191]{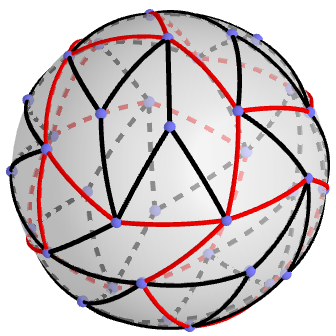}
		\includegraphics[scale=0.185]{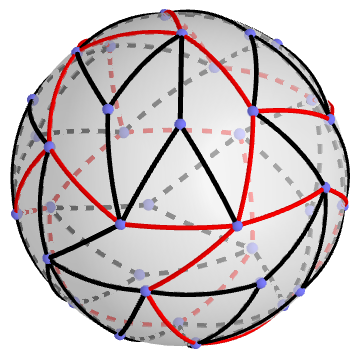}
		\includegraphics[scale=0.189]{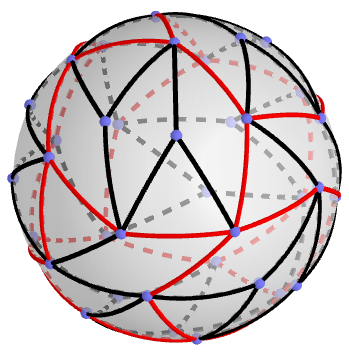}\hspace{3pt}

		\includegraphics[scale=0.188]{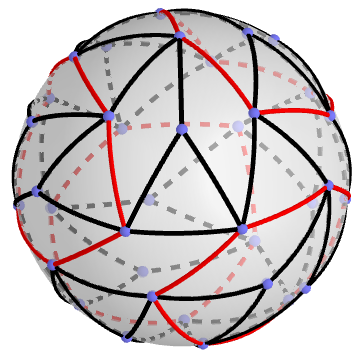}
		\includegraphics[scale=0.182]{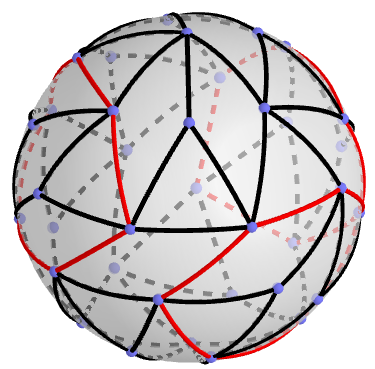}
		\includegraphics[scale=0.1925]{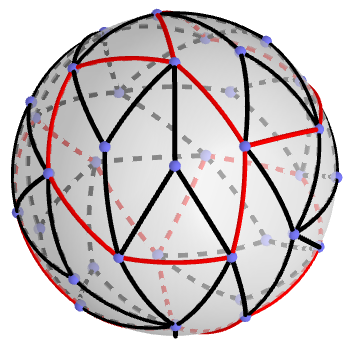}
		\includegraphics[scale=0.1955]{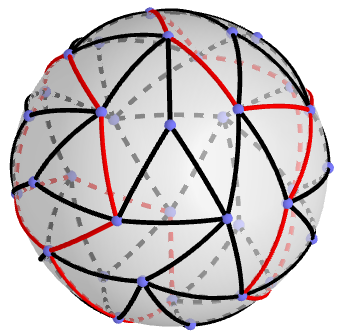}
		\includegraphics[scale=0.191]{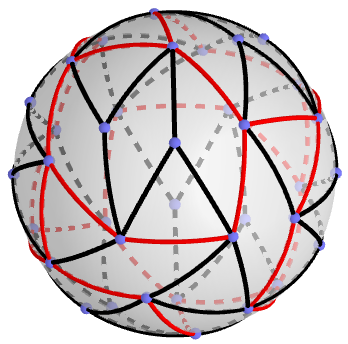}
		\includegraphics[scale=0.1825]{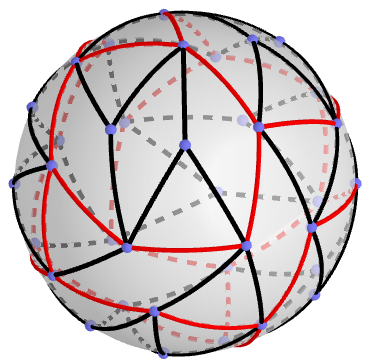}
		\includegraphics[scale=0.2]{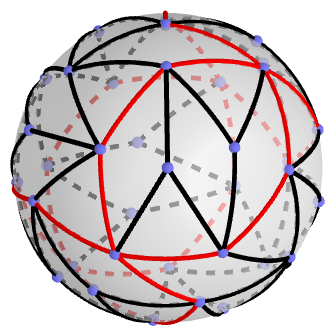}
		\includegraphics[scale=0.1877]{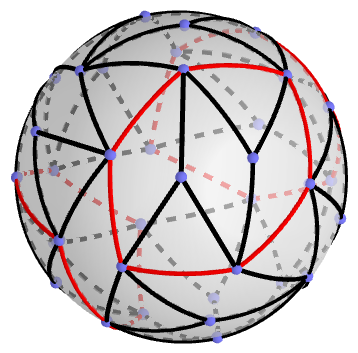}
		\caption{$T(20a^3,24a^4;12\aaa\bbb^2,24\aaa^2\bbb\ccc^2)$.}
		\label{18}
	\end{figure}

Without $\aaa^3\ccc^4$, we have the tiling $T(20a^3,24a^4;12\aaa\bbb^2,24\aaa^2\bbb\ccc^2)$ with 24 different tilings in Figure \ref{18}. \label{20,24.2}

\begin{figure}[H]
		\centering
		\begin{minipage}[c]{0.4\linewidth}
		\begin{tikzpicture}[>=latex,scale=0.6]  
        \draw (0,3.49) -- (1.24,2.54) -- (1.24,0.95) -- (0,0) -- (-1.24,0.95) -- (-1.24,2.54) -- (0,3.49);
	    \fill (1.24,2.54) circle (2.5pt);
	    \fill (-1.24,2.54) circle (2.5pt);
        \draw[dashed] (0,4)--(0,-1);
		\node at (0.3,1.8) {\small $L$}; 
		\node at (-0.8,2.54) {\small $A$};
		\node at (0.8,2.54) {\small $B$};
		\node at (0,3.7) {\small $\ccc\aaa$};
		\node at (1.9,2.54) {\small $\aaa\ccc\ccc$};
		\node at (1.7,0.95) {\small $\bbb$};
		\node at (0,-0.3) {\small $\aaa\ccc$};
		\node at (-1.9,0.95) {\small $\aaa\ccc\ccc$};
		\node at (-1.7,2.54) {\small $\bbb$};
		\end{tikzpicture}
		\end{minipage}\hspace{-50pt}%
		\begin{minipage}[c]{0.4\linewidth}\hspace*{3em}
		\includegraphics[scale=0.3]{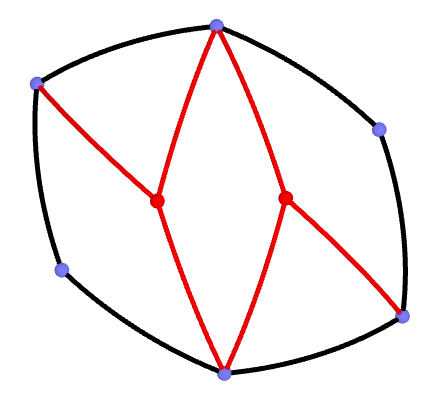}
		\end{minipage}
		\caption{The basic module.}
		\label{ab24}
	\end{figure}

In these tilings, we find a module in Figure \ref{ab24}. The flip of the module with respect to line $L$ keeps the angle sums of all vertices and transforms the tiling to a new tiling $T(20a^3,24a^4;(12+k)\aaa\bbb^2,(24-2k)\aaa^2\bbb\ccc^2,k\aaa^3\ccc^4)$, where $0\le k\le11$. For example, $A,B$ are vertices of degree 5, by flipping, $A$ become a vertex of degree 7 and $B$ become a vertex of degree 3. This indicates that for every vertex of degree 7 added, one vertex of degree 3 will be added and two vertices of degree 5 will be reduced. 

Since the diversity of the number, arrangement and position of $\aaa^3\ccc^4$ and $\aaa^2\bbb\ccc^2$, the complexity of the tilings is greatly increased.

\vspace{9pt}
\noindent{Subcase $\thin\ccc\thin\ccc\thin\cdots=\aaa\ccc^k(k=3,5,7,9)$, AVC $\subset\{\aaa\bbb^2,\aaa\ccc^k\}$}

There must be $\thin\ccc\thin\ccc\thin\ccc\thin\cdots$. As shown in figure \ref{rar}, we get $\aaa^2\ccc\cdots$, contradicting the AVC.

\vspace{9pt}
\noindent{Subcase $\thin\ccc\thin\ccc\thin\cdots=\aaa^2\ccc^2$, AVC $\subset\{\aaa\bbb^2,\aaa^2\ccc^2\}$}

\begin{figure}[htp]
		\centering
		\begin{tikzpicture}[>=latex,scale=0.6]  
        \draw (0,3.2) -- (0,2.3) -- (1.14,1.56) -- (1.11,-1.16) -- (-1.11,-1.16) -- (-1.14,1.56) -- (0,2.3)--(0,0);
        \draw (-1.85,-1.64)--(-1.11,-1.16)--(0,0)--(1.11,-1.16)--(1.85,-1.64) (1.14,1.56)--(0,0.93)--(-1.14,1.56);
		\node at (0.24,0.15) {\small $\bbb$};
		\node at (-0.24,0.15) {\small $\bbb$};
		\node at (1.35,1.56) {\small $\bbb$};
		\node at (-1.35,1.56) {\small $\bbb$};
	\node[draw,shape=circle, inner sep=0.2] at (0,-0.7) {\small $1$};
	\node[draw,shape=circle, inner sep=0.2] at (0.75,0.2) {\small $2$};
	\node[draw,shape=circle, inner sep=0.2] at (-0.75,0.2) {\small $3$};
	\node[draw,shape=circle, inner sep=0.2] at (0.35,1.5) {\small $4$};
	\node[draw,shape=circle, inner sep=0.2] at (-0.35,1.5) {\small $5$};
	\node[draw,shape=circle, inner sep=0.2] at (1.5,0) {\small $6$};
	\node[draw,shape=circle, inner sep=0.2] at (-1.5,0) {\small $7$};
	\node[draw,shape=circle, inner sep=0.2] at (0,-1.5) {\small $8$};
		\end{tikzpicture}
		\caption{$T(4a^3,4a^4;4\aaa\bbb^2,4\aaa^2\ccc^2)$.}
		\label{abb1}
	\end{figure}

In Figure \ref{abb1}, $\aaa\bbb^2$ determines $T_1,T_2,T_3$; $\ccc^2\cdots=\aaa^2\ccc^2$ determines $T_4,T_5$; $\aaa^2\cdots=\aaa^2\ccc^2$ determines $T_6,T_7$; $\aaa\ccc^2\cdots=\aaa^2\ccc^2$ determines $T_8$. The 3D picture is the 13th of Figure \ref{Fig3}. \label{4,4}

\vspace{9pt}
\noindent{Subcase $\thin\ccc\thin\ccc\thin\cdots=\aaa^2\ccc^k(k=4,\cdots,7)$, AVC $\subset\{\aaa\bbb^2,\aaa^2\ccc^k\}$}

When $\thin\ccc\thin\aaa\thin\ccc\thin\cdots$ appears, one of the two $\thin\aaa\thin\bbb\thin\cdots$ must be $\thin\ccc\thin\aaa\thin\bbb\thin\cdots$, contradicting the AVC. So $\aaa^2\ccc^k=\thin\aaa\thin\aaa\thin\ccc\thin\cdots\thin\ccc\thin$. As shown in figure \ref{rar}, $\thin\aaa\thin\ccc\thin\aaa\thin\cdots$ appears, a contradiction. Then $\thin\ccc\thin\ccc\thin\cdots=\aaa^3\ccc^k(k=3,5)$ is also impossible.

\vspace{9pt}
\noindent{Subcase $\thin\ccc\thin\ccc\thin\cdots=\aaa^4\ccc^k(k=2,3)$, AVC $\subset\{\aaa\bbb^2,\aaa^4\ccc^k\}$}

Similarly, $\thin\ccc\thin\aaa\thin\ccc\thin\cdots$ is not a vertex. When $k=2$, $\aaa\bbb^2$ determines $T_1,T_2,T_3$ as the first of Figure \ref{ab212}; $\ccc^2\cdots=\aaa^4\ccc^2$ determines $T_4,\cdots,T_7$; $\aaa^2\cdots=\aaa^4\ccc^2$ determines $T_8,T_9,T_{10},T_{11}$; $\aaa\bbb\cdots=\aaa\bbb^2$ determines $T_{12},T_{13}$. But $\thin\ccc_{12}\thin\aaa_4\thin\aaa_5\thin\aaa_8\thin\ccc_{13}\thin\cdots=\thin\ccc_{12}\thin\aaa_4\thin\aaa_5\thin\aaa_8\thin\ccc_{13}\thin\aaa\thin$, a contradiction.
	\begin{figure}[H]
		\centering
		\begin{tikzpicture}[>=latex,scale=0.85]
				\draw(-2,0)--(-2,1.73)--(2,1.73) (2,0)--(-2,0) (-1.5,0.87)--(1.5,0.87) 
                     (-2,1.73)--(-1.5,0.87)--(-2,0) (-1,1.73)--(-0.5,0.87)--(0,0) (0,1.73)--(-0.5,0.87)--(-1,0)
                     (2,1.73)--(1.5,0.87)--(2,0) (1,1.73)--(0.5,0.87)--(0,0) (0,1.73)--(0.5,0.87)--(1,0)
				(-0.42,2.69)--(0,1.73)--(0.42,2.69)
				(-1.66,2.69)--(-2,1.73) (1.66,2.69)--(2,1.73);
		\node at (-1.4,0.65) {\small $\bbb$};
		\node at (-1.4,1.1) {\small $\bbb$};
		\node at (-1,2) {\small $\bbb$};
		\node at (0.8,0.65) {\small $\ccc$};
		\node at (0.8,1.1) {\small $\ccc$};
		\node at (1,2) {\small $\bbb$};
		\node at (0,1.73) {\small $\bullet$};
	\node[draw,shape=circle, inner sep=0.2] at (-1.8,0.8) {\small $1$};
	\node[draw,shape=circle, inner sep=0.2] at (-1.2,0.35) {\small $2$};
	\node[draw,shape=circle, inner sep=0.2] at (-1.2,1.5) {\small $3$};
	\node[draw,shape=circle, inner sep=0.2] at (-0.5,1.5) {\small $4$};
	\node[draw,shape=circle, inner sep=0.2] at (0,1.1) {\small $5$};
	\node[draw,shape=circle, inner sep=0.2] at (0,0.6) {\small $6$};
	\node[draw,shape=circle, inner sep=0.2] at (-0.5,0.3) {\small $7$};
	\node[draw,shape=circle, inner sep=0.2] at (-1,2.5) {\footnotesize $12$};
	\node[draw,shape=circle, inner sep=0.2] at (0.5,1.5) {\small $8$};
	\node[draw,shape=circle, inner sep=0.2] at (0.5,0.3) {\small $9$};
	\node[draw,shape=circle, inner sep=0.2] at (1.2,0.4) {\footnotesize $11$};
	\node[draw,shape=circle, inner sep=0.2] at (1.2,1.4) {\footnotesize $10$};
	\node[draw,shape=circle, inner sep=0.2] at (1,2.5) {\footnotesize $13$};
		\end{tikzpicture}\hspace{40pt}
			\begin{tikzpicture}[>=latex,scale=0.8]
            \draw (-1.14,0.46)--(-0.48,1.14)--(0,2.28)--(0.48,1.14)--(1.14,0.46) (0.48,1.14)--(0,0)--(-0.48,1.14)
             (-1.62,1.6)--(-1.14,0.46)--(0,0)--(1.14,0.46)--(1.62,1.6)--(0,2.28)--(-1.62,1.6);
			\node at (-0.42,1.83){\small $\ccc$};
			\node at (0,0.53){\small $\ccc$};
			\node at (0.42,1.83){\small $\ccc$};
	\node[draw,shape=circle, inner sep=0.2] at (-0.58,0.53) {\small $1$};
	\node[draw,shape=circle, inner sep=0.2] at (0,1.14) {\small $2$};
	\node[draw,shape=circle, inner sep=0.2] at (0.58,0.53) {\small $3$};
	\node[draw,shape=circle, inner sep=0.2] at (-0.96,1.33) {\small $4$};
	\node[draw,shape=circle, inner sep=0.2] at (0.96,1.33) {\small $5$};
			\end{tikzpicture}
		\caption{$\thin\ccc\thin\ccc\thin\cdots=\aaa^4\ccc^k(k=2,3)$.}
		\label{ab212}
	\end{figure}
When $k=3$, we get $\thin\ccc\thin\ccc\thin\cdots=\thin\aaa\thin\aaa\thin\aaa\thin\aaa\thin\ccc\thin\ccc\thin\ccc\thin$, $\thin\ccc\thin\ccc\thin\aaa\thin\aaa\thin\ccc\thin\aaa\thin\aaa\thin$. The latter determines $T_1,T_2,T_3$ as the second of Figure \ref{ab212}; $\aaa\bbb\cdots=\aaa\bbb^2$ determines $T_4,T_5$. Then $\thin\ccc_4\thin\ccc_2\thin\ccc_5\thin\cdots=\thin\aaa\thin\aaa\thin\aaa\thin\aaa\thin\ccc_4\thin\ccc_2\thin\ccc_5\thin$ and $\thin\aaa\thin\aaa\thin\aaa\thin\aaa\thin\ccc\thin\ccc\thin\ccc\thin$ must appear. Similar to the first of Figure \ref{ab212}, $\thin\ccc\thin\aaa\thin\ccc\thin\cdots$ appears, a contradiction. 

\subsubsection*{Case $\ccc^3$}
We have $\ccc=\frac{2\pi}{3}<2\aaa$. The AAD $\thin\ccc^\bbb\thin^\bbb\ccc\thin\ccc\thin$ of $\ccc^3$ gives $\thin\bbb\thin\bbb\thin\cdots$. Similar to the proof of Lemma \ref{deg345}, $\bbb^2\cdots$ is not a vertex.

	\subsection{The vertex types of degree 4}

If there is no vertices of degree 3, then $f_\tri \ge8$ by (\ref{eq3}).

\subsubsection*{Case $\ccc^4$}
We have $\ccc=\frac{\pi}{2}>\aaa$. By $\ccc^4=\thin\ccc^\bbb\thin^\bbb\ccc\thin\ccc\thin\ccc\thin$ and similar proof of Lemma \ref{deg345}, we get $\thin\bbb\thin\bbb\thin\cdots=\aaa\bbb^3,\aaa\bbb^2\ccc,\aaa^2\bbb^2$.

\vspace{9pt}
\noindent{Subcase $\thin\bbb\thin\bbb\thin\cdots=\aaa\bbb^3$, AVC $\subset \{\ccc^4,\aaa\bbb^3\}$}

There is no the vertex $\aaa\ccc\cdots$, contradicting Lemma \ref{AAD}.

\vspace{9pt}
\noindent{Subcase $\thin\bbb\thin\bbb\thin\cdots=\aaa^2\bbb^2$, AVC $\subset \{\ccc^4,\aaa^2\bbb^2,\aaa\bbb\ccc^2\}$}

As shown in Figure \ref{rar}, we get $\aaa^2\ccc\cdots$, contradicting the AVC.

\vspace{9pt}
\noindent{Subcase $\thin\bbb\thin\bbb\thin\cdots=\aaa\bbb^2\ccc$, AVC $\subset \{\ccc^4,\aaa\bbb^2\ccc\}$}

In Figure \ref{r4}, $\ccc^4$ determines $T_1,T_2,T_3,T_4$. By the symmetry of the partial tiling and $\bbb^2\cdots=\aaa\bbb^2\ccc$, we might as well take $\thin\bbb_1\thin\bbb_2\thin\cdots=\thin\bbb_1\thin\bbb_2\thin\ccc\thin\aaa\thin$ which determines $T_5,T_6$; $\thin\bbb_1\thin\bbb_4\thin\cdots=\thin\bbb_1\thin\bbb_4\thin\aaa\thin\ccc\thin$ determines $T_7,T_8$. Similarly, we can determine $T_9,\cdots,T_{12}$. Then $T_{13},T_{14},\cdots,T_{32}$ are determined and the 3D picture is the 18th of Figure \ref{Fig3}. \label{8,24}
	\begin{figure}[htp]
		\centering
		\begin{tikzpicture}[>=latex,scale=0.65] 
				\foreach \a in {0,1,2,3}
				\draw[rotate=90*\a]
        (-5.16,3.73)--(-4.18,2.59)--(-2.01,2.96)--(-2.75,1.35)--(-1.47,1.47)--(0,1.06)--(1.47,1.47)
        (2.59,4.18)--(0.61,4.15)--(1.35,2.75)--(-0.08,2.41)--(0,1.06)
        (-1.47,1.47)--(-0.08,2.41)--(-2.01,2.96)--(0.61,4.15) (0,0)--(0,1.06)
        (0.61,4.15) arc (71.28:144.83:4.21);
		\node at (0.2,0.2) {\small $\ccc$};
		\node at (-0.2,0.2) {\small $\ccc$};
		\node at (-0.2,-0.25) {\small $\ccc$};
		\node at (0.2,-0.25) {\small $\ccc$};
		\node at (0.2,1.3) {\small $\ccc$};
		\node at (-1.35,0.2) {\small $\ccc$};
		\node at (-0.2,-1.35) {\small $\ccc$};
		\node at (1.35,-0.2) {\small $\ccc$};
		\node at (-1.5,1.7) {\small $\bbb$};
		\node at (-1.7,-1.5) {\small $\bbb$};
		\node at (1.5,-1.7) {\small $\bbb$};
		\node at (1.7,1.5) {\small $\bbb$};
		\node at (0,2.65) {\small $\bbb$};
		\node at (-2.65,0) {\small $\bbb$};
		\node at (0,-2.65) {\small $\bbb$};
		\node at (2.65,0) {\small $\bbb$};
		\node at (1.45,2.95) {\small $\ccc$};
		\node at (-2.95,1.45) {\small $\ccc$};
		\node at (-1.45,-2.95) {\small $\ccc$};
		\node at (2.95,-1.45) {\small $\ccc$};
		\node at (0.6,4.4) {\small $\ccc$};
		\node at (-4.4,0.6) {\small $\ccc$};
		\node at (-0.6,-4.4) {\small $\ccc$};
		\node at (4.4,-0.6) {\small $\ccc$};
	\node[draw,shape=circle, inner sep=0.2] at (-0.6,0.6) {\small $1$};
	\node[draw,shape=circle, inner sep=0.2] at (0.6,0.6) {\small $2$};
	\node[draw,shape=circle, inner sep=0.2] at (0.6,-0.6) {\small $3$};
	\node[draw,shape=circle, inner sep=0.2] at (-0.6,-0.6) {\small $4$};
	\node[draw,shape=circle, inner sep=0.2] at (0.6,1.8) {\small $5$};
	\node[draw,shape=circle, inner sep=0.2] at (-0.5,1.65) {\small $6$};
	\node[draw,shape=circle, inner sep=0.2] at (-1.8,0.6) {\small $7$};
	\node[draw,shape=circle, inner sep=0.2] at (-1.65,-0.5) {\small $8$};
	\node[draw,shape=circle, inner sep=0.2] at (-0.6,-1.8) {\small $9$};
	\node[draw,shape=circle, inner sep=0.2] at (0.5,-1.65) {\footnotesize $10$};
	\node[draw,shape=circle, inner sep=0.2] at (1.8,-0.6) {\footnotesize $11$};
	\node[draw,shape=circle, inner sep=0.2] at (1.65,0.5) {\footnotesize $12$};
	\node[draw,shape=circle, inner sep=0.2] at (-1.6,2.35) {\footnotesize $13$};
	\node[draw,shape=circle, inner sep=0.2] at (-2.35,-1.6) {\footnotesize $14$};
	\node[draw,shape=circle, inner sep=0.2] at (1.6,-2.35) {\footnotesize $15$};
	\node[draw,shape=circle, inner sep=0.2] at (2.35,1.6) {\footnotesize $16$};
	\node[draw,shape=circle, inner sep=0.2] at (0.3,3.3) {\footnotesize $17$};
	\node[draw,shape=circle, inner sep=0.2] at (-3.3,0.3) {\footnotesize $18$};
	\node[draw,shape=circle, inner sep=0.2] at (-0.3,-3.3) {\footnotesize $19$};
	\node[draw,shape=circle, inner sep=0.2] at (3.3,-0.3) {\footnotesize $20$};
	\node[draw,shape=circle, inner sep=0.2] at (2,3.5) {\footnotesize $21$};
	\node[draw,shape=circle, inner sep=0.2] at (-3.5,2) {\footnotesize $22$};
	\node[draw,shape=circle, inner sep=0.2] at (-2,-3.5) {\footnotesize $23$};
	\node[draw,shape=circle, inner sep=0.2] at (3.5,-2) {\footnotesize $24$};
	\node[draw,shape=circle, inner sep=0.2] at (-2,3.5) {\footnotesize $25$};
	\node[draw,shape=circle, inner sep=0.2] at (-3.5,-2) {\footnotesize $26$};
	\node[draw,shape=circle, inner sep=0.2] at (2,-3.5) {\footnotesize $27$};
	\node[draw,shape=circle, inner sep=0.2] at (3.5,2) {\footnotesize $28$};
	\node[draw,shape=circle, inner sep=0.2] at (0,5) {\footnotesize $29$};
	\node[draw,shape=circle, inner sep=0.2] at (-5,0) {\footnotesize $30$};
	\node[draw,shape=circle, inner sep=0.2] at (0,-5) {\footnotesize $31$};
	\node[draw,shape=circle, inner sep=0.2] at (5,0) {\footnotesize $32$};
		\end{tikzpicture}
		\caption{$T(8a^3,24a^4;24\aaa\bbb^2\ccc,6\ccc^4)$.}
		\label{r4}
	\end{figure}

\subsubsection*{Case $\aaa^2\bbb^2$}
By $2\bbb+2\ccc>2\pi=2\aaa+2\bbb$, we have $\ccc>\aaa$. If there is a vertex $\aaa^2\bbb^2=\thin\aaa\thin\aaa\thin\bbb^\ccc\thin^\ccc\bbb\thin$, then we get $\thin\ccc\thin\ccc\thin\cdots=\aaa\bbb\ccc^2,\bbb\ccc^3,\aaa^3\ccc^2,\aaa^2\ccc^3,\aaa\ccc^3,\aaa\ccc^4,\ccc^5$ by similar proof of Lemma \ref{deg345}. 
	\begin{table}[htp]                  
		\centering     
		\begin{tabular}{|c|c|c|}\hline	 			
				$\thin\ccc\thin\ccc\thin\cdots$ & AVC $\subset$ a set  & Conclusion \\
				\hline  
				$\bbb\ccc^3$&$\{\aaa^2\bbb^2,\bbb\ccc^3\}$&\multirow{2}{*}{No $\aaa\ccc\cdots$, contradicting Lemma \ref{AAD}} \\
				\cline{1-2} 
				$\ccc^5$&$\{\aaa^2\bbb^2,\ccc^5\}$&\multirow{2}{*}{} \\			
				\hline 
		\end{tabular}
	\caption{The simple cases when $\thin\aaa\thin\aaa\thin\bbb\thin\bbb\thin$ appears.}\label{Tab-4}  
	\end{table}

\vspace{9pt}
\noindent{Subcase $\thin\ccc\thin\ccc\thin\cdots=\aaa\ccc^k(k=3,4)$, AVC $\subset \{\aaa^2\bbb^2,\aaa\ccc^k\}$}

By $\aaa\ccc^k=\thin\ccc\thin\ccc\thin\ccc\thin\cdots$, we get $\aaa^2\ccc\cdots$ as shown in Figure \ref{rar}, contradicting the AVC. 

\vspace{9pt}
\noindent{Subcase $\thin\ccc\thin\ccc\thin\cdots=\aaa^2\ccc^3$, AVC $\subset\{\aaa^2\bbb^2,\aaa^2\ccc^3\}$}

	\begin{figure}[htp]
		\centering
		\begin{tikzpicture}[>=latex,scale=0.8]
				\draw(0,0)--(2,0)--(2,1)--(1,2)--(0,1)--(0,0)--(0.5,-0.87)--(2.5,-0.87)--(2,0) 
                     (0.5,-0.87)--(1,-1.73)--(2,-1.73)--(2.5,-0.87)--(3.18,0.42)--(2,1)
                     (1.5,-0.87)--(2,-1.73) (3.18,0.42)--(2,0) 
                     (1,2)--(1,0) --(0.5,-0.87) (1,0)--(1.5,-0.87) (0,1)--(2,1)--(2.62,1.8)--(1.69,2.88)--(1,2)
                     (2.62,1.8)--(4.47,0.9)--(3.18,0.42) (2.5,-0.87)--(3.6,-1.08)--(4.47,0.9)
                     (2,-1.74)--(3.7,-2.2)--(3.6,-1.08);
		\node at (1.2,0.75) {\small $\bbb$};
		\node at (0.8,0.75) {\small $\bbb$};
		\node at (1.2,-0.2) {\small $\ccc$};
		\node at (1.2,0.75) {\small $\bbb$};
		\node at (0.75,-1.1) {\small $\ccc$};
		\node at (1.98,1.35) {\small $\ccc$};
		\node at (2.43,1.06) {\small $\ccc$};
		\node at (3.3,0.15) {\small $\bbb$};
		\node at (2.62,-1.22) {\small $\ccc$};
		\node at (2,-1.74) {\small $\bullet$};
	\node[draw,shape=circle, inner sep=0.2] at (0.7,1.3) {\small $1$};
	\node[draw,shape=circle, inner sep=0.2] at (1.3,1.3) {\small $2$};
	\node[draw,shape=circle, inner sep=0.2] at (1.5,0.35) {\small $3$};
	\node[draw,shape=circle, inner sep=0.2] at (0.5,0.35) {\small $4$};
	\node[draw,shape=circle, inner sep=0.2] at (1.65,-0.4) {\small $5$};
	\node[draw,shape=circle, inner sep=0.2] at (1,-0.6) {\small $6$};
	\node[draw,shape=circle, inner sep=0.2] at (0.5,-0.3) {\small $7$};
	\node[draw,shape=circle, inner sep=0.2] at (2.35,0.5) {\small $8$};
	\node[draw,shape=circle, inner sep=0.2] at (2.5,-0.3) {\small $9$};
	\node[draw,shape=circle, inner sep=0.2] at (3.12,1.1) {\footnotesize $10$};
	\node[draw,shape=circle, inner sep=0.2] at (1.83,2.09) {\footnotesize $11$};
	\node[draw,shape=circle, inner sep=0.2] at (3.28,-0.45) {\footnotesize $12$};
	\node[draw,shape=circle, inner sep=0.2] at (1.25,-1.35) {\footnotesize $13$};
	\node[draw,shape=circle, inner sep=0.2] at (2,-1.15) {\footnotesize $14$};
	\node[draw,shape=circle, inner sep=0.2] at (3.02,-1.59) {\footnotesize $15$};
		\end{tikzpicture}\hspace{50pt}
		\begin{tikzpicture}[>=latex,scale=0.8]
				\draw(0,3.17)--(2.11,3.57)--(3.41,2.2)--(2.44,0.79)--(1.75,-1.03)--(1.35,-2.3)--(0,-1.75)--(0.5,-0.87)--(0,0)--(1,0)--(1,1)--(0,2)
                (0,3.17)--(1.75,2.22)--(3.41,2.2) (0,1)--(1,1)--(1.75,2.22)--(2.44,0.79)--(1,1) 
                (0.5,-0.87)--(1,0)--(1.75,-1.03)--(0.5,-0.87);
				\draw(0,3.17)--(-2.11,3.57)--(-3.41,2.2)--(-2.44,0.79)--(-1.75,-1.03)--(-1.35,-2.3)--(0,-1.75)--(-0.5,-0.87)--(0,0)--(-1,0)--(-1,1)--(0,2)
                (0,3.17)--(-1.75,2.22)--(-3.41,2.2) (0,1)--(-1,1)--(-1.75,2.22)--(-2.44,0.79)--(-1,1) 
                (-0.5,-0.87)--(-1,0)--(-1.75,-1.03)--(-0.5,-0.87) (0,0)--(0,3.17);
		\node at (0.2,0.75) {\small $\bbb$};
		\node at (-0.2,0.75) {\small $\bbb$};
		\node at (0,-0.4) {\small $\ccc$};
		\node at (0.6,-1.2) {\small $\bbb$};
		\node at (1.2,0.09) {\small $\bbb$};
		\node at (0.3,2.09) {\small $\bbb$};
		\node at (1.88,2.5) {\small $\bbb$};
		\node at (-0.6,-1.2) {\small $\bbb$};
		\node at (-1.2,0.09) {\small $\bbb$};
		\node at (-0.3,2.09) {\small $\bbb$};
		\node at (-1.88,2.5) {\small $\bbb$};
		\node at (0,3.17) {\small $\bullet$};
	\node[draw,shape=circle, inner sep=0.2] at (-0.3,1.3) {\small $1$};
	\node[draw,shape=circle, inner sep=0.2] at (0.3,1.3) {\small $2$};
	\node[draw,shape=circle, inner sep=0.2] at (-0.5,0.35) {\small $4$};
	\node[draw,shape=circle, inner sep=0.2] at (0.5,0.35) {\small $3$};
	\node[draw,shape=circle, inner sep=0.2] at (-0.5,-0.3) {\small $5$};
	\node[draw,shape=circle, inner sep=0.2] at (0,-1) {\small $6$};
	\node[draw,shape=circle, inner sep=0.2] at (0.5,-0.3) {\small $7$};
	\node[draw,shape=circle, inner sep=0.2] at (1.04,-0.58) {\small $8$};
	\node[draw,shape=circle, inner sep=0.2] at (0.93,-1.61) {\small $9$};
	\node[draw,shape=circle, inner sep=0.2] at (1.67,0.24) {\footnotesize $10$};
	\node[draw,shape=circle, inner sep=0.2] at (1.7,1.38) {\footnotesize $11$};
	\node[draw,shape=circle, inner sep=0.2] at (0.8,2.12) {\footnotesize $12$};
	\node[draw,shape=circle, inner sep=0.2] at (1.62,3.04) {\footnotesize $13$};
	\node[draw,shape=circle, inner sep=0.2] at (2.47,1.67) {\footnotesize $14$};
	\node[draw,shape=circle, inner sep=0.2] at (-1.04,-0.58) {\footnotesize $15$};
	\node[draw,shape=circle, inner sep=0.2] at (-0.93,-1.61) {\footnotesize $16$};
	\node[draw,shape=circle, inner sep=0.2] at (-1.67,0.24) {\footnotesize $17$};
	\node[draw,shape=circle, inner sep=0.2] at (-1.7,1.38) {\footnotesize $18$};
	\node[draw,shape=circle, inner sep=0.2] at (-0.8,2.12) {\footnotesize $19$};
	\node[draw,shape=circle, inner sep=0.2] at (-1.62,3.04) {\footnotesize $20$};
	\node[draw,shape=circle, inner sep=0.2] at (-2.47,1.67) {\footnotesize $21$};
		\end{tikzpicture}
		\caption{$\thin\ccc\thin\ccc\thin\cdots=\aaa^2\ccc^3$.}
		\label{a2b24}
	\end{figure}

In Figure \ref{a2b24}, $\thin\aaa\thin\aaa\thin\bbb\thin\bbb\thin$ determines $T_1,T_2,T_3,T_4$. When $\thin\ccc_4\thin\ccc_3\thin\cdots=\thin\ccc_4\thin\ccc_3\thin\ccc\thin\aaa\thin\aaa\thin$, $T_5,T_6,T_7$ are determined in the first. Then we can determine $T_8,\cdots,T_{15}$ but $\aaa\bbb\ccc\cdots$ appears, contradicting the AVC. Therefore, $\thin\ccc_4\thin\ccc_3\thin\cdots=\thin\ccc_4\thin\ccc_3\thin\aaa\thin\ccc\thin\aaa\thin$ determines $T_5,T_6,T_7$ in the second. Then we can determine $T_8,\cdots,T_{21}$ but $\ccc^4\cdots$ appears, contradicting the AVC.

\vspace{9pt}
\noindent{Subcase $\thin\ccc\thin\ccc\thin\cdots=\aaa\bbb\ccc^2$, AVC $\subset\{\aaa^2\bbb^2,\aaa\bbb\ccc^2\}$}

In the first of Figure \ref{a2b21}, $\thin\aaa\thin\aaa\thin\bbb\thin\bbb\thin$ determines $T_1,T_2,T_3,T_4$. By the symmetry of the partial tiling and $\ccc^2\cdots=\aaa\bbb\ccc^2$, we might as well take $\thin\ccc_3\thin\ccc_4\thin\cdots=\thin\ccc_3\thin\ccc_4\thin\aaa\thin\bbb\thin$ which determines $T_5,T_6$. Then we can determine $T_7,T_8,\cdots,T_{20}$ and the 3D picture is the 16th of Figure \ref{Fig3}. \label{8,12}

	\begin{figure}[htp]
		\centering
		\begin{tikzpicture}[>=latex,scale=0.65] 
				\foreach \a in {0,1,2}
				\draw[rotate=120*\a]      
	   (0.54,2.65)--(2.76,1.4)--(3.06,-2.39)--(2.02,-1.79)--(0.91,-1.17)--(0,-0.69)--(0.6,0.35)--(1.86,0.28)--(2.02,-1.79)  (0.56,1.37)--(2.76,1.4)--(1.86,0.28)--(0.91,-1.17) 
       (0.55,3.84) arc (85.2:198.63:4.02);
		\node at (0.8,1.1) {\small $\bbb$};
		\node at (0.3,1.1) {\small $\bbb$};
		\node at (0.6,0) {\small $\bbb$};
		\node at (2.1,0.15) {\small $\ccc$};
		\node at (2.9,1.45) {\small $\ccc$};
		\node at (0.3,2.8) {\small $\bbb$};
		\node at (-1.39,1.19) {\small $\bbb$};
		\node at (-0.7,-0.1) {\small $\ccc$};
		\node at (0.05,-1) {\small $\ccc$};
		\node at (-0.92,-2) {\small $\ccc$};
		\node at (-0.16,-3.3) {\small $\ccc$};
		\node at (-2.8,1.67) {\small $\ccc$};
	\node[draw,shape=circle, inner sep=0.2] at (0.1,1.8) {\small $1$};
	\node[draw,shape=circle, inner sep=0.2] at (1,1.8) {\small $2$};
	\node[draw,shape=circle, inner sep=0.2] at (1.5,0.8) {\small $3$};
	\node[draw,shape=circle, inner sep=0.2] at (-0.2,0.8) {\small $4$};
	\node[draw,shape=circle, inner sep=0.2] at (0,0) {\small $5$};
	\node[draw,shape=circle, inner sep=0.2] at (0.8,-0.5) {\small $6$};
	\node[draw,shape=circle, inner sep=0.2] at (-0.8,2.3) {\small $7$};
	\node[draw,shape=circle, inner sep=0.2] at (2.23,2.36) {\small $8$};
	\node[draw,shape=circle, inner sep=0.2] at (2.47,-0.5) {\small $9$};
	\node[draw,shape=circle, inner sep=0.2] at (1.6,-1) {\footnotesize $10$};
	\node[draw,shape=circle, inner sep=0.2] at (-1.3,0.5) {\footnotesize $11$};
	\node[draw,shape=circle, inner sep=0.2] at (-0.6,-0.6) {\footnotesize $12$};
	\node[draw,shape=circle, inner sep=0.2] at (0,-1.7) {\footnotesize $13$};
	\node[draw,shape=circle, inner sep=0.2] at (1.09,-1.77) {\footnotesize $14$};
	\node[draw,shape=circle, inner sep=0.2] at (-2.1,0) {\footnotesize $15$};
	\node[draw,shape=circle, inner sep=0.2] at (-1.6,-0.75) {\footnotesize $16$};
	\node[draw,shape=circle, inner sep=0.2] at (-1.74,-1.8) {\footnotesize $17$};
	\node[draw,shape=circle, inner sep=0.2] at (1.5,-2.8) {\footnotesize $18$};
	\node[draw,shape=circle, inner sep=0.2] at (-3.1,0) {\footnotesize $19$};
	\node[draw,shape=circle, inner sep=0.2] at (1.5,4) {\footnotesize $20$};
		\end{tikzpicture}
		\begin{tikzpicture}[>=latex,scale=0.7] 
				\foreach \a in {0,1,2}
				\draw[rotate=120*\a]            
		(0,3.38)--(1.65,2.34)--(3.8,2.18)--(2.85,0.26)--(2.93,-1.69) (1.65,2.34)--(2.85,0.26)
        (0,1.89)--(1.65,2.34)--(1.2,0.69)--(2.85,0.26)--(1.64,-0.95)
        (0,4.66)--(3.8,2.18)--(4.04,-2.33)--(2.93,-1.69)--(1.64,-0.95)--(0.6,-0.35)--(0,0.69)--(1.2,0.69)--(0.6,-0.35)
        (0,4.66) arc (103.38:-43.38:4.21);
		\node at (-0.2,1.65) {\small $\bbb$};
		\node at (0.2,1.65) {\small $\bbb$};
		\node at (-1.68,2.62) {\small $\ccc$};
		\node at (1.68,2.62) {\small $\ccc$};
		\node at (-1.3,0.35) {\small $\bbb$};
		\node at (1.25,0.35) {\small $\bbb$};
		\node at (-0.65,-0.7) {\small $\ccc$};
		\node at (0.65,-0.7) {\small $\ccc$};
		\node at (-3.1,0.2) {\small $\ccc$};
		\node at (3.1,0.2) {\small $\ccc$};
		\node at (-2.88,-2.12) {\small $\bbb$};
		\node at (2.88,-2.12) {\small $\bbb$};
	\node[draw,shape=circle, inner sep=0.2] at (-0.5,2.5) {\small $1$};
	\node[draw,shape=circle, inner sep=0.2] at (0.5,2.5) {\small $2$};
	\node[draw,shape=circle, inner sep=0.2] at (0.8,1.5) {\small $3$};
	\node[draw,shape=circle, inner sep=0.2] at (-0.8,1.5) {\small $4$};
	\node[draw,shape=circle, inner sep=0.2] at (0.6,0.25) {\small $5$};
	\node[draw,shape=circle, inner sep=0.2] at (0,0) {\small $6$};
	\node[draw,shape=circle, inner sep=0.2] at (-0.6,0.25) {\small $7$};
	\node[draw,shape=circle, inner sep=0.2] at (-1.2,3.2) {\small $8$};
	\node[draw,shape=circle, inner sep=0.2] at (-2.7,1.5) {\small $9$};
	\node[draw,shape=circle, inner sep=0.2] at (-1.8,1) {\footnotesize $10$};
	\node[draw,shape=circle, inner sep=0.2] at (1.2,3.2) {\footnotesize $11$};
	\node[draw,shape=circle, inner sep=0.2] at (2.7,1.5) {\footnotesize $12$};
	\node[draw,shape=circle, inner sep=0.2] at (1.8,1) {\footnotesize $13$};
	\node[draw,shape=circle, inner sep=0.2] at (-1.6,-0.2) {\footnotesize $14$};
	\node[draw,shape=circle, inner sep=0.2] at (1.6,-0.2) {\footnotesize $15$};
	\node[draw,shape=circle, inner sep=0.2] at (0,-0.7) {\footnotesize $16$};
	\node[draw,shape=circle, inner sep=0.2] at (-0.8,-1.5) {\footnotesize $17$};
	\node[draw,shape=circle, inner sep=0.2] at (0.8,-1.5) {\footnotesize $18$};
	\node[draw,shape=circle, inner sep=0.2] at (0,-2.15) {\footnotesize $19$};
	\node[draw,shape=circle, inner sep=0.2] at (-2.4,-0.82) {\footnotesize $20$};
	\node[draw,shape=circle, inner sep=0.2] at (-1.92,-1.67) {\footnotesize $21$};
	\node[draw,shape=circle, inner sep=0.2] at (2.4,-0.82) {\footnotesize $22$};
	\node[draw,shape=circle, inner sep=0.2] at (1.92,-1.67) {\footnotesize $23$};
	\node[draw,shape=circle, inner sep=0.2] at (-3.43,-0.56) {\footnotesize $24$};
	\node[draw,shape=circle, inner sep=0.2] at (3.43,-0.56) {\footnotesize $25$};
	\node[draw,shape=circle, inner sep=0.2] at (-2.2,-2.65) {\footnotesize $26$};
	\node[draw,shape=circle, inner sep=0.2] at (2.2,-2.65) {\footnotesize $27$};
	\node[draw,shape=circle, inner sep=0.2] at (0,-3.3) {\footnotesize $28$};
	\node[draw,shape=circle, inner sep=0.2] at (-4.5,1.65) {\footnotesize $29$};
	\node[draw,shape=circle, inner sep=0.2] at (1,-4.5) {\footnotesize $30$};
	\node[draw,shape=circle, inner sep=0.2] at (3,3.5) {\footnotesize $31$};
	\node[draw,shape=circle, inner sep=0.2] at (0,5.3) {\footnotesize $32$};
		\end{tikzpicture}
		\caption{$\thin\ccc\thin\ccc\thin\cdots=\aaa\bbb\ccc^2,\aaa^3\ccc^2$.}
		\label{a2b21}
	\end{figure}

\vspace{9pt}
\noindent{Subcase $\thin\ccc\thin\ccc\thin\cdots=\aaa^3\ccc^2$, AVC $\subset\{\aaa^2\bbb^2,\aaa^3\ccc^2\}$}

In the second of Figure \ref{a2b21}, $\thin\aaa\thin\aaa\thin\bbb\thin\bbb\thin$ determines $T_1,T_2,T_3,T_4$; $\ccc^2\cdots=\aaa^3\ccc^2$ determines $T_5,T_6,T_7$. When $\thin\ccc_4\thin\aaa_1\thin\cdots=\thin\ccc_4\thin\aaa_1\thin\aaa\thin\cdots$, we get $\thin\aaa_2\thin\aaa_1\thin\aaa\thin$ $\cdots=\thin\aaa_2\thin\aaa_1\thin\aaa\thin\ccc\thin\ccc\thin$. But $\thin\bbb\thin\aaa_2\thin\ccc_3\thin$ $\cdots$ appears, contradicting the AVC. So $\thin\ccc_4\thin\aaa_1\thin\cdots=\thin\ccc_4\thin\aaa_1\thin\ccc\thin\aaa\thin\aaa\thin$ determines $T_8,T_9,T_{10}$. Then $T_{11},\cdots,T_{32}$ are determined and the 3D picture is the 19th of Figure \ref{Fig3}.\label{20,12.1}

	\begin{figure}[htp]
		\centering
		\begin{tikzpicture}[>=latex,scale=0.75]
				\draw(0,1)--(1.14,0.96)--(1.11,-0.43)--(0,-1)--(-1.14,-0.96)--(-1.11,0.43)--(0,1)--(0,-1) (-1.11,0.43)--(1.11,-0.43);
		\node at (0.29,0.19) {\small $\bbb$};
		\node at (-0.29,-0.19) {\small $\bbb$};
	\node[draw,shape=circle, inner sep=0.2] at (-0.32,0.48) {\small $1$};
	\node[draw,shape=circle, inner sep=0.2] at (0.66,0.59) {\small $2$};
	\node[draw,shape=circle, inner sep=0.2] at (0.32,-0.48) {\small $3$};
	\node[draw,shape=circle, inner sep=0.2] at (-0.66,-0.59) {\small $4$};
		\end{tikzpicture}\hspace{30pt}
		\begin{tikzpicture}[>=latex,scale=0.75]
				\draw(0,1)--(1.14,0.96)--(1.11,-0.43)--(0,-1)--(-1.14,-0.96)--(-1.11,0.43)--(0,1)--(0,-1) (-1.11,0.43)--(1.11,-0.43)
                (1.11,-0.43)--(1.01,-1.69)--(0.45,-2.61)--(-0.61,-2.01)--(-1.14,-0.96)--(-2.35,-0.39)
                (1.01,-1.69)--(0,-1)--(-0.61,-2.01)--(-1.96,-1.98)--(-2.35,-0.39)--(-1.11,0.43);
		\node at (0.29,0.19) {\small $\bbb$};
		\node at (-0.29,-0.19) {\small $\bbb$};
		\node at (0.08,-1.4) {\small $\ccc$};
		\node at (-1.35,-1.24) {\small $\bbb$};
		\node at (-0.61,-2.01) {\small $\bullet$};
	\node[draw,shape=circle, inner sep=0.2] at (-0.32,0.48) {\small $1$};
	\node[draw,shape=circle, inner sep=0.2] at (0.66,0.59) {\small $2$};
	\node[draw,shape=circle, inner sep=0.2] at (0.32,-0.48) {\small $3$};
	\node[draw,shape=circle, inner sep=0.2] at (-0.66,-0.59) {\small $4$};
	\node[draw,shape=circle, inner sep=0.2] at (0.69,-1.05) {\small $5$};
	\node[draw,shape=circle, inner sep=0.2] at (0.24,-1.93) {\small $6$};
	\node[draw,shape=circle, inner sep=0.2] at (-0.61,-1.34) {\small $7$};
	\node[draw,shape=circle, inner sep=0.2] at (-1.53,-0.34) {\small $8$};
	\node[draw,shape=circle, inner sep=0.2] at (-1.65,-1.56) {\small $9$};
		\end{tikzpicture}\hspace{30pt}
		\begin{tikzpicture}[>=latex,scale=0.75]
				\draw(0,1)--(1.14,0.96)--(1.11,-0.43)--(0,-1)--(-1.14,-0.96)--(-1.11,0.43)--(0,1)--(0,-1) (-1.11,0.43)--(1.11,-0.43)
                (-2.14,-1.79)--(-0.4,-3.3)--(1.16,-1.82)--(4.02,-0.87)--(2.14,1.79)
                (-0.42,-2.08)--(-0.4,-3.3)  (2.46,-0.39)--(4.02,-0.87);
				\foreach \a in {0,1}
				\draw[rotate=180*\a] 
                (-2.46,0.39)--(-1.14,-0.96)--(-0.42,-2.08)--(0,-1)--(1.16,-1.82)--(1.11,-0.43)--(2.46,-0.39)
                (-2.46,0.39)--(-2.14,-1.79)--(-0.42,-2.08)--(1.16,-1.82)--(2.46,-0.39);           
		\node at (0.29,0.19) {\small $\bbb$};
		\node at (-0.29,-0.19) {\small $\bbb$};
		\node at (1.35,1.24) {\small $\bbb$};
		\node at (-1.35,-1.24) {\small $\bbb$};
		\node at (1.16,-1.82) {\small $\bullet$};
	\node[draw,shape=circle, inner sep=0.2] at (-0.32,0.48) {\small $1$};
	\node[draw,shape=circle, inner sep=0.2] at (0.66,0.59) {\small $2$};
	\node[draw,shape=circle, inner sep=0.2] at (0.32,-0.48) {\small $3$};
	\node[draw,shape=circle, inner sep=0.2] at (-0.66,-0.59) {\small $4$};
	\node[draw,shape=circle, inner sep=0.2] at (0.69,-1.05) {\small $5$};
	\node[draw,shape=circle, inner sep=0.2] at (0.24,-1.63) {\small $6$};
	\node[draw,shape=circle, inner sep=0.2] at (-0.5,-1.34) {\small $7$};
	\node[draw,shape=circle, inner sep=0.2] at (1.5,0) {\small $9$};
	\node[draw,shape=circle, inner sep=0.2] at (1.5,-0.8) {\small $8$};
	\node[draw,shape=circle, inner sep=0.2] at (-0.69,1.05) {\footnotesize $10$};
	\node[draw,shape=circle, inner sep=0.2] at (-0.24,1.63) {\footnotesize $11$};
	\node[draw,shape=circle, inner sep=0.2] at (0.5,1.34) {\footnotesize $12$};
	\node[draw,shape=circle, inner sep=0.2] at (-1.5,0) {\footnotesize $14$};
	\node[draw,shape=circle, inner sep=0.2] at (-1.5,0.8) {\footnotesize $13$};
	\node[draw,shape=circle, inner sep=0.2] at (-1.8,-1.1) {\footnotesize $15$};
	\node[draw,shape=circle, inner sep=0.2] at (1.8,1.1) {\footnotesize $16$};
	\node[draw,shape=circle, inner sep=0.2] at (-0.85,-2.48) {\footnotesize $17$};
	\node[draw,shape=circle, inner sep=0.2] at (-0,-2.46) {\footnotesize $18$};
	\node[draw,shape=circle, inner sep=0.2] at (2.59,-0.92) {\footnotesize $19$};
	\node[draw,shape=circle, inner sep=0.2] at (2.88,0.11) {\footnotesize $20$};
		\end{tikzpicture}
		\caption{$\aaa^2\bbb^2$ can only be $\thin\aaa\thin\bbb\thin\aaa\thin\bbb\thin$.}
		\label{a2b26}
	\end{figure}

\nd Now the vertex $\aaa^2\bbb^2$ can only be $\thin\aaa\thin\bbb\thin\aaa\thin\bbb\thin$ which determines $T_1,T_2,T_3,T_4$ in Figure \ref{a2b26}. By similar proof of Lemma \ref{deg345}, we have $\bbb^2\cdots=\aaa^2\bbb^2$ . So $\thin\bbb\thin\bbb\thin\cdots$ is not a vertex. Then $\thin\ccc\thin\ccc\thin\cdots$ can never be a vertex and $\thin\aaa\thin\ccc\thin\cdots=\thin\aaa\thin\ccc\thin\bbb\thin\ccc\thin,\thin\aaa\thin\aaa\thin\ccc\thin\aaa\thin\ccc\thin,\aaa^4\ccc$.

When $\thin\aaa\thin\ccc\thin\cdots=\thin\aaa\thin\ccc\thin\bbb\thin\ccc\thin$, we have the AVC $\subset\{\thin\aaa\thin\bbb\thin\aaa\thin\bbb\thin,\thin\aaa\thin\ccc\thin\bbb\thin\ccc\thin\}$. Then one of $\thin\aaa_3\thin\ccc_4\thin\cdots$ and $\thin\aaa_3\thin\ccc_2\thin\cdots$ must be $\thin\aaa\thin\ccc\thin\ccc\thin\bbb\thin$, contradicting the AVC.  

When $\thin\aaa\thin\ccc\thin\cdots=\thin\aaa\thin\aaa\thin\ccc\thin\aaa\thin\ccc\thin$, we have the AVC $\subset\{\thin\aaa\thin\bbb\thin\aaa\thin\bbb\thin,\thin\aaa\thin\aaa\thin\ccc\thin\aaa\thin\ccc\thin\}$. Then we can determine $T_5,T_6,\cdots,T_9$ in the second of Figure \ref{a2b26}. But $\aaa\bbb\ccc\cdots$ appears, contradicting the AVC. 

When $\thin\aaa\thin\ccc\thin\cdots=\aaa^4\ccc$, we have the AVC $\subset\{\thin\aaa\thin\bbb\thin\aaa\thin\bbb\thin,\aaa^4\ccc\}$. Then we can determine $T_5,\cdots,T_{20}$ in the third of Figure \ref{a2b26}. But $\aaa^5\cdots$ appears, contradicting the AVC.

\subsubsection*{Case $\aaa\ccc^3$}
We have $\ccc>\frac{1}{2}\pi>\aaa$. The AAD $\thin\aaa\thin\ccc^\bbb\thin^\bbb\ccc\thin\ccc\thin$ of $\aaa\ccc^3$ gives $\thin\bbb\thin\bbb\thin\cdots$. By similar proof of Lemma \ref{deg345}, $\bbb^2\cdots$ is not a vertex.

\subsubsection*{Case $\aaa^2\ccc^2$}
We have $\ccc>\frac{1}{2}\pi>\aaa$. By similar proof of Lemma \ref{deg345}, $\bbb^2\cdots$ is not a vertex. So $\thin\ccc\thin\ccc\thin\cdots$ is not a vertex by Lemma \ref{AAD}. This implies $\aaa^2\ccc^2=\thin\aaa\thin\ccc\thin\aaa\thin\ccc\thin$ which determines $T_1,T_2,T_3,T_4$ in Figure \ref{a2r2}. Then we can get $\thin\aaa\thin\bbb\thin\cdots=\aaa^3\bbb$ and the AVC $\subset\{\aaa^2\ccc^2,\aaa^3\bbb\}$. Then $T_5,\cdots,T_{11}$ are determined and the 3D picture is the 14th of Figure \ref{Fig3}.\label{8,3}

	\begin{figure}[htp]
		\centering
		\begin{tikzpicture}[>=latex,scale=0.7]
				\foreach \a in {0,1}
				\draw[rotate=180*\a]      
				(0,0)--(-0.5,0.87)--(0,1.73)--(0.5,0.87)--(0,0) (0.5,0.87)--(-0.5,0.87)--(-2,0)--(-0.5,-0.87)
                 (2,0) arc (-10.95:109.16:1.53) (0,1.73) arc (70.84:190.95:1.53);
		\node at (0.25,0) {\small $\ccc$};
		\node at (-0.25,0) {\small $\ccc$};
		\node at (2.25,0) {\small $\ccc$};
	\node[draw,shape=circle, inner sep=0.2] at (0,0.5) {\small $1$};
	\node[draw,shape=circle, inner sep=0.2] at (1,0) {\small $2$};
	\node[draw,shape=circle, inner sep=0.2] at (0,-0.5) {\small $3$};
	\node[draw,shape=circle, inner sep=0.2] at (-1,0) {\small $4$};
	\node[draw,shape=circle, inner sep=0.2] at (0,1.2) {\small $5$};
	\node[draw,shape=circle, inner sep=0.2] at (1,1) {\small $6$};
	\node[draw,shape=circle, inner sep=0.2] at (-1,1) {\small $7$};
	\node[draw,shape=circle, inner sep=0.2] at (0,-1.2) {\small $8$};
	\node[draw,shape=circle, inner sep=0.2] at (1,-1) {\small $9$};
	\node[draw,shape=circle, inner sep=0.2] at (-1,-1) {\footnotesize $10$};
	\node[draw,shape=circle, inner sep=0.2] at (2.8,0) {\footnotesize $11$};
		\end{tikzpicture}
		\caption{$T(8a^3,3a^4;6\aaa^3\bbb,3\aaa^2\ccc^2)$.}
		\label{a2r2}
	\end{figure}

\subsubsection*{Case $\aaa\bbb\ccc^2$}
Similar to the proof of Lemma \ref{deg345}, $\bbb^2\cdots$ is not a vertex. So $\thin\ccc\thin\ccc\thin\cdots$ can never be a vertex by Lemma \ref{AAD}. Then $\aaa\bbb\ccc^2=\thin\aaa\thin\ccc\thin\bbb\thin\ccc\thin$ determines $T_1,T_2,T_3,T_4$ in Figure \ref{abr2}. We can get $\thin\bbb\thin\ccc\thin\cdots=\thin\bbb\thin\ccc\thin\aaa\thin\ccc\thin$ which determines $T_5,\cdots,T_{12}$. Then $\aaa^2\bbb\cdots=\aaa^4\bbb$ determines $T_{13},\cdots,T_{19}$ and the AVC $\subset\{\aaa\bbb\ccc^2,\aaa^4\bbb\}$. But $\aaa^5\cdots$ appears, contradicting the AVC.
	\begin{figure}[H]
		\centering
		\begin{tikzpicture}[>=latex,scale=0.7]
				\draw(-0.24,3.18)--(2.33,2.65)--(3.78,1.69)--(3.52,-0.32)--(2.33,-1.46)--(1.06,-2.3)--(-0.29,-2.54)--(-1.91,-1.96)--(-2.7,0.4)--(-0.85,0.95)--(-0.24,3.18)--(0.53,2.14)--(0.5,0.87)--(0,0)--(-0.53,-1.16)--(-0.29,-2.54)
                (-1.91,-1.96)--(-1.43,-0.11)--(-0.85,0.95)--(0.53,2.14)--(2.33,2.65)
                (3.78,1.69)--(2.22,1.38)--(0.5,0.87)--(-0.85,0.95)--(0,0)--(-1.43,-0.11)
                (3.52,-0.32)--(2.14,0)--(2.33,-1.46)--(0.77,-1.01)--(-0.53,-1.16)--(-1.43,-0.11)--(-2.7,0.4)
                (0.53,2.14)--(2.22,1.38)--(2.14,0)--(1,0)--(0,0) (0.5,0.87)--(1,0)--(0.77,-1.01)--(1.06,-2.3);
		\node at (1.14,0.26) {\small $\ccc$};
		\node at (1.24,-0.34) {\small $\bbb$};
		\node at (0.69,-0.24) {\small $\ccc$};
		\node at (0.56,-1.35) {\small $\ccc$};
		\node at (-0.85,-1.32) {\small $\ccc$};
		\node at (2.43,0.24) {\small $\ccc$};
		\node at (2.28,1.6) {\small $\ccc$};
		\node at (-0.85,0.95) {\small $\bullet$};
	\node[draw,shape=circle, inner sep=0.2] at (0.5,0.3) {\small $1$};
	\node[draw,shape=circle, inner sep=0.2] at (1.5,0.7) {\small $2$};
	\node[draw,shape=circle, inner sep=0.2] at (1.5,-0.7) {\small $3$};
	\node[draw,shape=circle, inner sep=0.2] at (0.2,-0.6) {\small $4$};
	\node[draw,shape=circle, inner sep=0.2] at (1.3,-1.55) {\small $5$};
	\node[draw,shape=circle, inner sep=0.2] at (0.1,-1.8) {\small $6$};
	\node[draw,shape=circle, inner sep=0.2] at (-0.7,-0.4) {\small $7$};
	\node[draw,shape=circle, inner sep=0.2] at (-1.2,-1.3) {\small $8$};
	\node[draw,shape=circle, inner sep=0.2] at (2.8,-0.5) {\small $9$};
	\node[draw,shape=circle, inner sep=0.2] at (3.1,0.8) {\footnotesize $10$};
	\node[draw,shape=circle, inner sep=0.2] at (1,1.5) {\footnotesize $11$};
	\node[draw,shape=circle, inner sep=0.2] at (1.9,2.1) {\footnotesize $12$};
	\node[draw,shape=circle, inner sep=0.2] at (-0.77,0.34) {\footnotesize $13$};
	\node[draw,shape=circle, inner sep=0.2] at (-0.08,0.56) {\footnotesize $14$};
	\node[draw,shape=circle, inner sep=0.2] at (0.13,1.27) {\footnotesize $15$};
	\node[draw,shape=circle, inner sep=0.2] at (-0.08,2.17) {\footnotesize $16$};
	\node[draw,shape=circle, inner sep=0.2] at (0.69,2.62) {\footnotesize $17$};
	\node[draw,shape=circle, inner sep=0.2] at (-1.59,0.34) {\footnotesize $18$};
	\node[draw,shape=circle, inner sep=0.2] at (-2.01,-0.42) {\footnotesize $19$};
		\end{tikzpicture}
		\caption{$\aaa\bbb\ccc^2$ is a vertex.}
		\label{abr2}
	\end{figure}

\subsubsection*{Case $\aaa^3\ccc$}
We have $\ccc>\frac{1}{2}\pi>\aaa$. Similar to the proof of Lemma \ref{deg345}, $\bbb^2\cdots$ is not a vertex. So $\thin\ccc\thin\ccc\thin\cdots$ can never be a vertex by Lemma \ref{AAD}. By $\aaa^3\ccc=\thin\ccc^\bbb\thin^\aaa\aaa\thin\cdots$, we get $\thin\aaa\thin\bbb\thin\cdots$ which is not a vertex.

\subsubsection*{Case $\bbb\ccc^3$}

We have $\ccc>\frac{1}{3}\pi$. Then the AAD $\thin\ccc^\bbb\thin^\bbb\ccc\thin\cdots$ of $\bbb\ccc^3$ gives $\thin\bbb\thin\bbb\thin\cdots=\aaa\bbb^3,\aaa\bbb^2\ccc$. 

\vspace{9pt}
\noindent{Subcase $\thin\bbb\thin\bbb\thin\cdots=\aaa\bbb^3$, AVC $\subset\{\bbb\ccc^3,\aaa\bbb^3\}$}

There is no the vertex $\aaa\ccc\cdots$, contradicting Lemma \ref{AAD}.

\vspace{9pt}
\noindent{Subcase $\thin\bbb\thin\bbb\thin\cdots=\aaa\bbb^2\ccc$, AVC $\subset\{\bbb\ccc^3,\aaa\bbb^2\ccc\}$}

	\begin{figure}[htp]
		\centering
		\begin{tikzpicture}[>=latex,scale=0.8]       
		\draw (0,1.41) -- (-1.41,0) (0.71,0.71) -- (-0.71,-0.71) (1.41,0) -- (0,-1.41)
              (-1.41,0) -- (0,-1.41) (-0.71,0.71) -- (0.71,-0.71) (0,1.41) -- (1.41,0)
              (-1.41,0) -- (-2.12,0.71) -- (-1.41,1.41) -- (-0.71,0.71)
              (-1.41,1.41) -- (0,1.41) -- (0.71,2.12) -- (1.41,1.41) -- (0.71,0.71) (1.41,1.41) -- (1.41,0);
        \draw (-2.12,0.71) -- (-2.83,1.41) -- (-2.12,2.12) -- (-1.41,1.41) 
                      (-2.12,2.12) -- (0.71,2.12);
        \draw (-2.83,1.41)--(-3.54,2.12)--(-2.83,2.83)--(-2.12,2.12) (-2.12,0.71)--(-2.83,0)--(-1.41,-1.41)--(0,-1.41) (-1.41,0)--(-2.12,-0.71)--(-0.71,-0.71);
		\node at (0,-0.3) {\small $\bbb$};
		\node at (0.2,0) {\small $\ccc$};
		\node at (-0.2,0) {\small $\ccc$};
		\node at (0,0.3) {\small $\ccc$};
		\node at (-1,0.71) {\small $\ccc$};
		\node at (0.71,1) {\small $\ccc$};
		\node at (-0.1,1.6) {\small $\bbb$};
		\node at (-1.7,1.41) {\small $\bbb$};
		\node at (-2.45,2.12) {\small $\ccc$};
		\node at (-2.09,0.35) {\small $\ccc$};
		\node at (-0.77,-1.01) {\small $\bbb$};
		\node at (-2.12,-0.71) {\small $\bullet$};
	\node[draw,shape=circle, inner sep=0.2] at (0,-0.7) {\small $1$};
	\node[draw,shape=circle, inner sep=0.2] at (0.7,0) {\small $2$};
	\node[draw,shape=circle, inner sep=0.2] at (0,0.7) {\small $3$};
	\node[draw,shape=circle, inner sep=0.2] at (-0.7,0) {\small $4$};
	\node[draw,shape=circle, inner sep=0.2] at (-0.7,1.1) {\small $5$};
	\node[draw,shape=circle, inner sep=0.2] at (-1.4,0.7) {\small $6$};
	\node[draw,shape=circle, inner sep=0.2] at (0.7,1.4) {\small $7$};
	\node[draw,shape=circle, inner sep=0.2] at (1.2,0.7) {\small $9$};
	\node[draw,shape=circle, inner sep=0.2] at (-0.7,1.7) {\small $8$};
	\node[draw,shape=circle, inner sep=0.2] at (-2.2,1.5) {\footnotesize $10$};
	\node[draw,shape=circle, inner sep=0.2] at (-2.9,2) {\footnotesize $11$};
	\node[draw,shape=circle, inner sep=0.2] at (-2.1,-0.1) {\footnotesize $12$};
	\node[draw,shape=circle, inner sep=0.2] at (-1.4,-0.4) {\footnotesize $13$};
	\node[draw,shape=circle, inner sep=0.2] at (-1.25,-1.05) {\footnotesize $14$};
		\end{tikzpicture}\hspace{20pt}
		\begin{tikzpicture}[>=latex,scale=0.8]       
		\draw (0,1.41) -- (-1.41,0) (0.71,0.71) -- (-0.71,-0.71) (1.41,0) -- (0,-1.41)
              (-1.41,0) -- (0,-1.41) (-0.71,0.71) -- (0.71,-0.71) (0,1.41) -- (1.41,0)
              (-1.41,0) -- (-2.12,0.71) -- (-1.41,1.41) -- (-0.71,0.71)
              (-1.41,1.41) -- (0,1.41) -- (0.71,2.12) -- (1.41,1.41) -- (0.71,0.71) (1.41,1.41) -- (1.41,0);
        \draw (-2.12,0.71) -- (-2.83,1.41) -- (-2.12,2.12) -- (-1.41,1.41) 
                      (-2.12,2.12) -- (0.71,2.12);
        \draw (-2.12,2.12)--(-2.83,2.83)--(1.41,2.83)--(2.12,2.12)--(2.12,0.71)--(1.41,0)
              (0.71,2.12)--(1.41,2.83) (1.41,1.41)--(2.12,2.12);
		\node at (0,-0.3) {\small $\bbb$};
		\node at (0.2,0) {\small $\ccc$};
		\node at (-0.2,0) {\small $\ccc$};
		\node at (0,0.3) {\small $\ccc$};
		\node at (-1,0.71) {\small $\ccc$};
		\node at (0.71,1) {\small $\ccc$};
		\node at (-0.1,1.6) {\small $\bbb$};
		\node at (-1.7,1.41) {\small $\bbb$};
		\node at (0.65,2.35) {\small $\ccc$};
		\node at (1.67,1.3) {\small $\bbb$};
		\node at (1.1,2.14) {\small $\bbb$};
		\node at (1.41,0) {\small $\bullet$};
	\node[draw,shape=circle, inner sep=0.2] at (0,-0.7) {\small $1$};
	\node[draw,shape=circle, inner sep=0.2] at (0.7,0) {\small $2$};
	\node[draw,shape=circle, inner sep=0.2] at (0,0.7) {\small $3$};
	\node[draw,shape=circle, inner sep=0.2] at (-0.7,0) {\small $4$};
	\node[draw,shape=circle, inner sep=0.2] at (-0.7,1.1) {\small $5$};
	\node[draw,shape=circle, inner sep=0.2] at (-1.4,0.7) {\small $6$};
	\node[draw,shape=circle, inner sep=0.2] at (0.7,1.4) {\small $7$};
	\node[draw,shape=circle, inner sep=0.2] at (1.2,0.7) {\small $9$};
	\node[draw,shape=circle, inner sep=0.2] at (-0.7,1.7) {\small $8$};
	\node[draw,shape=circle, inner sep=0.2] at (-2.2,1.5) {\footnotesize $10$};
	\node[draw,shape=circle, inner sep=0.2] at (-0.71,2.5) {\footnotesize $11$};
	\node[draw,shape=circle, inner sep=0.2] at (1.55,2.12) {\footnotesize $12$};
	\node[draw,shape=circle, inner sep=0.2] at (1.8,0.9) {\footnotesize $13$};
		\end{tikzpicture}\hspace{20pt}
		\begin{tikzpicture}[>=latex,scale=0.8]       
		\draw (0,1.41) -- (-1.41,0) (0.71,0.71) -- (-0.71,-0.71) (1.41,0) -- (0,-1.41)
              (-1.41,0) -- (0,-1.41) (-0.71,0.71) -- (0.71,-0.71) (0,1.41) -- (1.41,0)
              (-1.41,0) -- (-1.41,1.41)  (-0.71,0.71) -- (-1.41,1.41) -- (-0.71,2.12) -- (0,1.41)
              (0,1.41) -- (0.71,2.12) -- (1.41,1.41) -- (0.71,0.71) (1.41,1.41) -- (1.41,0);
        \draw (-0.71,2.12) -- (0.71,2.12);
		\node at (0,-0.3) {\small $\bbb$};
		\node at (0.2,0) {\small $\ccc$};
		\node at (-0.2,0) {\small $\ccc$};
		\node at (0,0.3) {\small $\ccc$};
		\node at (-0.71,1) {\small $\ccc$};
		\node at (0.71,1) {\small $\ccc$};
		\node at (0.71,2.12) {\small $\bullet$};
		\node at (-0.71,2.12) {\small $\bullet$};
	\node[draw,shape=circle, inner sep=0.2] at (0,-0.7) {\small $1$};
	\node[draw,shape=circle, inner sep=0.2] at (0.7,0) {\small $2$};
	\node[draw,shape=circle, inner sep=0.2] at (0,0.7) {\small $3$};
	\node[draw,shape=circle, inner sep=0.2] at (-0.7,0) {\small $4$};
	\node[draw,shape=circle, inner sep=0.2] at (-0.7,1.4) {\small $5$};
	\node[draw,shape=circle, inner sep=0.2] at (-1.2,0.7) {\small $6$};
	\node[draw,shape=circle, inner sep=0.2] at (0.7,1.4) {\small $7$};
	\node[draw,shape=circle, inner sep=0.2] at (1.2,0.7) {\small $8$};
	\node[draw,shape=circle, inner sep=0.2] at (0,1.8) {\small $9$};
		\end{tikzpicture}
		\caption{$\thin\bbb\thin\bbb\thin\cdots=\aaa\bbb^2\ccc$.}
		\label{br31}
	\end{figure}

In Figure \ref{br31}, $\bbb\ccc^3$ determines $T_1,T_2,T_3,T_4$. By the symmetry of the partial tiling and $\bbb^2\cdots=\aaa\bbb^2\ccc$, we might as well take $\thin\bbb_3\thin\bbb_4\thin\cdots=\thin\ccc\thin\aaa\thin\bbb_3\thin\bbb_4\thin$ which determines $T_5,T_6$ in the first and second; $\aaa\ccc\cdots=\aaa\bbb^2\ccc$ determines $T_7,T_8$; $\bbb^2\ccc\cdots=\aaa\bbb^2\ccc$ determines $T_9$; $\aaa\bbb\ccc$ $\cdots=\aaa\bbb^2\ccc$ determines $T_{10}$. When $\thin\ccc_7\thin\ccc_8\thin\cdots=\thin\ccc_7\thin\ccc_8\thin\bbb\thin\ccc\thin$, $\thin\bbb_8\thin\ccc_{10}\thin\cdots=\bbb\ccc^3$ determines $T_{11}$ in the first; $\thin\ccc_{10}\thin\ccc_6\thin\cdots=\thin\ccc_{10}\thin\ccc_6\thin\ccc\thin\cdots$ determines $T_{12}$; $\bbb^2\ccc\cdots=\aaa\bbb^2\ccc$ determines $T_{13}$; $\aaa\bbb\ccc\cdots=\aaa\bbb^2\ccc$ determines $T_{14}$. But $\aaa\ccc^2\cdots$ appears, contradicting the AVC. So $\thin\ccc_7\thin\ccc_8\thin\cdots=\thin\ccc_7\thin\ccc_8\thin\ccc\thin\bbb\thin$ determines $T_{11},T_{12}$ in the second; $\aaa\bbb\ccc\cdots=\aaa\bbb^2\ccc$ determines $T_{13}$. But $\aaa\ccc^2\cdots$ appears, contradicting the AVC.

Therefore, $\thin\bbb_3\thin\bbb_4\thin\cdots=\thin\aaa\thin\ccc\thin\bbb_3\thin\bbb_4\thin$ and $\thin\bbb_3\thin\bbb_2\thin\cdots=\thin\aaa\thin\ccc\thin\bbb_3\thin\bbb_2\thin$ which determine $T_5,\cdots,T_8$ in the third. Then $T_9$ is determined. But one of $\thin\aaa_9\thin\ccc_7\thin\cdots$ and $\thin\aaa_9\thin\ccc_5\thin\cdots$ can never be a vertex.

\subsubsection*{Case $\aaa\bbb^3$}
We have $\bbb<\frac{5}{9}\pi$ and $\ccc>\frac{4}{9}\pi$. Then the AAD $\thin\aaa\thin\bbb^\ccc\thin^\ccc\bbb\thin\bbb\thin$ of $\aaa\bbb^3$ gives $\thin\ccc\thin\ccc\thin\cdots=\aaa^3\ccc^2$ and the AVC $\subset\{\aaa\bbb^3,\aaa^3\ccc^2\}$. Similar to Figure \ref{rar}, $\aaa^2\bbb\cdots$ appears, contradicting the AVC.

\subsubsection*{Case $\aaa^2\bbb\ccc$}
By $2\aaa+2\bbb>2\aaa+\bbb+\ccc=2\pi$, we get $\aaa+\bbb>\pi$. Similar to the proof of Lemma \ref{deg345}, $\bbb^2\cdots$ is not a vertex. So $\thin\ccc\thin\ccc\thin\cdots$ can never be a vertex by Lemma \ref{AAD}. Then we have $\bbb\ccc\cdots=\aaa^2\bbb\ccc$, $\aaa\bbb\cdots=\aaa^2\bbb\ccc$. 

	\begin{figure}[htp]
		\centering
        \begin{tikzpicture}[>=latex,scale=0.8]       
		\draw (0,0) -- (2,0)  (0,1) -- (2,1) (1,-1) -- (1,2)
              (0,0) -- (1,-1) -- (2,0) -- (2,1) -- (1,2) -- (0,1) -- (0,0);
        \draw (-1,2)--(3,2)--(3,-1)--(-1,-1)--(-1,2)--(0,1) (2,1)--(3,2)
              (2,0)--(3,-1) (0,0)--(-1,-1);
        \draw (3,-1)--(2,-2)--(1,-1) (2,-2)--(5,-2)--(3,2);
		\node at (0.8,0.7) {\small $\ccc$};
		\node at (1.2,0.7) {\small $\bbb$};
		\node at (-0.2,0.7) {\small $\ccc$};
		\node at (2.2,0) {\small $\ccc$};
		\node at (3.2,-1) {\small $\ccc$};
		\node at (1,-1) {\small $\bullet$};
	\node[draw,shape=circle, inner sep=0.2] at (0.7,1.3) {\small $1$};
	\node[draw,shape=circle, inner sep=0.2] at (1.3,1.3) {\small $2$};
	\node[draw,shape=circle, inner sep=0.2] at (1.5,0.5) {\small $3$};
	\node[draw,shape=circle, inner sep=0.2] at (0.5,0.5) {\small $4$};
	\node[draw,shape=circle, inner sep=0.2] at (0.7,-0.3) {\small $5$};
	\node[draw,shape=circle, inner sep=0.2] at (1.3,-0.3) {\small $6$};
	\node[draw,shape=circle, inner sep=0.2] at (0,1.5) {\small $7$};
	\node[draw,shape=circle, inner sep=0.2] at (-0.5,0.5) {\small $8$};
	\node[draw,shape=circle, inner sep=0.2] at (0,-0.5) {\small $9$};
	\node[draw,shape=circle, inner sep=0.2] at (2,-0.5) {\footnotesize $10$};
	\node[draw,shape=circle, inner sep=0.2] at (2.5,0.5) {\footnotesize $11$};
	\node[draw,shape=circle, inner sep=0.2] at (2,1.5) {\footnotesize $12$};
	\node[draw,shape=circle, inner sep=0.2] at (2,-1.5) {\footnotesize $13$};
	\node[draw,shape=circle, inner sep=0.2] at (3.8,-1) {\footnotesize $14$};
		\end{tikzpicture}\hspace{40pt}
        \begin{tikzpicture}[>=latex,scale=0.8]       
		\draw (0,0) -- (2,0)  (0,1) -- (2,1) (1,-1) -- (1,2)
              (0,0) -- (1,-1) -- (2,0) -- (2,1) -- (1,2) -- (0,1) -- (0,0);
        \draw (-1,2)--(3,2)--(3,-1)--(-1,-1)--(-1,2)--(0,1) (2,1)--(3,2)
              (2,0)--(3,-1) (0,0)--(-1,-1);
        \draw (3,-1) arc (-40.05:-139.95:2.61) 
              (3,2) arc (40.05:139.95:2.61);
		\node at (0.8,0.7) {\small $\ccc$};
		\node at (1.2,0.7) {\small $\bbb$};
		\node at (-0.2,0.7) {\small $\ccc$};
		\node at (2.2,0) {\small $\ccc$};
		\node at (3.2,-1) {\small $\ccc$};
	\node[draw,shape=circle, inner sep=0.2] at (0.7,1.3) {\small $1$};
	\node[draw,shape=circle, inner sep=0.2] at (1.3,1.3) {\small $2$};
	\node[draw,shape=circle, inner sep=0.2] at (1.5,0.5) {\small $3$};
	\node[draw,shape=circle, inner sep=0.2] at (0.5,0.5) {\small $4$};
	\node[draw,shape=circle, inner sep=0.2] at (0.7,-0.3) {\small $5$};
	\node[draw,shape=circle, inner sep=0.2] at (1.3,-0.3) {\small $6$};
	\node[draw,shape=circle, inner sep=0.2] at (0,1.5) {\small $7$};
	\node[draw,shape=circle, inner sep=0.2] at (-0.5,0.5) {\small $8$};
	\node[draw,shape=circle, inner sep=0.2] at (0,-0.5) {\small $9$};
	\node[draw,shape=circle, inner sep=0.2] at (2,-0.5) {\footnotesize $10$};
	\node[draw,shape=circle, inner sep=0.2] at (2.5,0.5) {\footnotesize $11$};
	\node[draw,shape=circle, inner sep=0.2] at (2,1.5) {\footnotesize $12$};
	\node[draw,shape=circle, inner sep=0.2] at (1,-1.5) {\footnotesize $13$};
	\node[draw,shape=circle, inner sep=0.2] at (1,2.5) {\footnotesize $15$};
	\node[draw,shape=circle, inner sep=0.2] at (3.5,0.5) {\footnotesize $14$};
		\end{tikzpicture}
		\caption{$\thin\aaa\thin\aaa\thin\bbb\thin\ccc\thin$ appears.}
		\label{a2br1}
	\end{figure}
If there is a vertex $\aaa^2\bbb\ccc=\thin\aaa\thin\aaa\thin\bbb\thin\ccc\thin$, then $T_1,\cdots,T_6$ are determined in Figure \ref{a2br1}. By the symmetry of the partial tiling and $\aaa\bbb\cdots=\aaa^2\bbb\ccc$, we might as well take $\thin\bbb_4\thin\aaa_1\thin$ $\cdots=\thin\bbb_4\thin\aaa_1\thin\aaa\thin\ccc\thin$ which determines $T_7,T_8$. Then we can determine $T_{9},\cdots,T_{12}$. Since $4\aaa<2\pi$, $T_{13}$ and $T_{14}$ are determined when $5\aaa\le2\pi$. When $5\aaa<2\pi$, $\thin\aaa_9\thin\aaa_5\thin\aaa_6\thin\aaa_{10}\thin\aaa_{13}\thin\cdots$ is not a vertex. When $5\aaa=2\pi$, we get a tiling in the second. By calculation, we have $\bbb=2\aaa=2\ccc=\frac{4}{5}\pi$. Then the tiling is of the icosahedral type by Lemma \ref{ico}. So $\thin\bbb_4\thin\aaa_1\thin\cdots=\thin\bbb_4\thin\aaa_1\thin\ccc\thin\aaa\thin$ and $\thin\bbb_3\thin\aaa_6\thin\cdots=\thin\bbb_3\thin\aaa_6\thin\ccc\thin\aaa\thin$. Then we get a tiling as the first of Figure \ref{a2br.2}, which is the same as the third of Figure \ref{ab=r}. 

Now $\aaa^2\bbb\ccc$ can only be $\thin\aaa\thin\bbb\thin\aaa\thin\ccc\thin$, and we get a tiling as the second of Figure \ref{a2br.2} which is the same as the fourth of Figure \ref{ab=r}. We notice that the red line divides the tilings into two same modules and the first tiling becomes the second by flipping one of the two modules.

	\begin{figure}[htp]
		\centering
        \begin{tikzpicture}[>=latex,scale=0.7]       
		\draw (0,0) -- (2,0)  (0,1) -- (2,1) (1,-1) -- (1,2)
              (0,0) -- (1,-1) -- (2,0) -- (2,1) -- (1,2) -- (0,1) -- (0,0)
              (2,0) -- (3,0.5) -- (2,1) (0,0) -- (-1,0.5) -- (0,1) (1,2) -- (1,3.12) (1,-1) -- (1,-2.12)
              (3,0.5) -- (1,3.12) -- (-1,0.5) -- (1,-2.12) -- (3,0.5);
        \draw[color=red] (1,-2.12)--(1,3.12);
        \draw[color=red] (1,3.12) arc (97.74:-97.74:2.64);
		\node at (0.8,0.7) {\small $\ccc$};
		\node at (1.2,0.7) {\small $\bbb$};
		\node at (-0.1,1.2) {\small $\ccc$};
		\node at (1.1,2.1) {\small $\ccc$};
		\node at (2.1,-0.2) {\small $\ccc$};
		\node at (0.8,-1.2) {\small $\ccc$};
		\end{tikzpicture}\hspace{50pt}
		\begin{tikzpicture}[>=latex,scale=0.8]       
		\draw (0,1) -- (1,1) -- (1,0) -- (0,0);
        \draw (1,0) -- (2,0.5) -- (1,1) -- (0.5,2)  (0,1) -- (-1,0.5) -- (0,0)  (0.5,-1) -- (1,0)
              (1.79,1.79) -- (2,0.5) -- (1.79,-0.79)  (0.5,-1) -- (-0.79,-0.79) -- (-1,0.5) -- (-0.79,1.79) -- (0.5,2);
        \draw[color=red] (1.79,1.79)--(0.5,2)--(0,1)--(0,0)--(0.5,-1)--(1.79,-0.79);
        \draw (1.79,1.79) arc (45:135:1.84);
        \draw[color=red] (1.79,1.79) arc (45:-45:1.84);
        \draw (-0.79,-0.79) arc (-135:-45:1.84);
        \draw (-0.79,-0.79) arc (-135:-225:1.84);
		\node at (0.2,0.75) {\small $\bbb$};
		\node at (-0.1,1.1) {\small $\ccc$};
		\node at (0.8,1.7) {\small $\ccc$};
		\node at (1.2,-0.2) {\small $\ccc$};
		\node at (0,-0.7) {\small $\ccc$};
		\node at (-1,-1) {\small $\ccc$};
		\end{tikzpicture}
		\caption{$T(8a^3,6a^4;12\aaa^2\bbb\ccc)$.}
		\label{a2br.2}
	\end{figure}

\subsubsection*{Case $\aaa\bbb^2\ccc$}

By $2\bbb+2\ccc>2\pi=\aaa+2\bbb+\ccc$, we get $\ccc>\aaa$. Similar to the proof of Lemma \ref{deg345}, we can get $\bbb^2\cdots=\aaa\bbb^2\ccc$ and $\bbb\ccc\cdots=\aaa\bbb^2\ccc$. If there is a vertex $\aaa\bbb^2\ccc=\thin\aaa\thin\bbb^\ccc\thin^\ccc\bbb\thin\ccc\thin$, then we get $\thin\ccc\thin\ccc\thin\cdots=\ccc^5,\aaa\ccc^4,\aaa^2\ccc^3,\aaa^3\ccc^2$.

\vspace{9pt}
\noindent{Subcase $\thin\ccc\thin\ccc\thin\cdots=\ccc^5$, AVC $\subset\{\aaa\bbb^2\ccc,\ccc^5\}$}

	\begin{figure}[htp]
		\centering
		\begin{tikzpicture}[>=latex,scale=0.7]
				\foreach \a in {0,1,2,3,4}
				\draw[rotate=72*\a]
				(2.67,1.62)--(1.49,1.85)--(0.82,2.99)--(-0.71,3.04)--(0,2.02)--(0.59,0.81)--(1.92,0.62)
				(0,2.02)--(1.49,1.85)--(0.59,0.81)--(0,0);
		\node at (0,0.25) {\small $\ccc$};
		\node at (0.25,0.1) {\small $\ccc$};
		\node at (0.2,-0.3) {\small $\ccc$};
		\node at (-0.2,-0.3) {\small $\ccc$};
		\node at (-0.3,0.1) {\small $\ccc$};
		\node at (-0.6,1.2) {\small $\ccc$};
		\node at (-1.4,-0.2) {\small $\ccc$};
		\node at (-0.2,-1.4) {\small $\ccc$};
		\node at (1.2,-0.65) {\small $\ccc$};
		\node at (1,1) {\small $\ccc$};
		\node at (0.15,2.25) {\small $\bbb$};
		\node at (-2.05,0.8) {\small $\bbb$};
		\node at (-1.35,-1.9) {\small $\bbb$};
		\node at (1.15,-1.95) {\small $\bbb$};
		\node at (2.2,0.6) {\small $\bbb$};
	\node[draw,shape=circle, inner sep=0.2] at (0,1) {\small $1$};
	\node[draw,shape=circle, inner sep=0.2] at (0.8,0.35) {\small $2$};
	\node[draw,shape=circle, inner sep=0.2] at (0.5,-0.8) {\small $3$};
	\node[draw,shape=circle, inner sep=0.2] at (-0.5,-0.8) {\small $4$};
	\node[draw,shape=circle, inner sep=0.2] at (-0.8,0.35) {\small $5$};
	\node[draw,shape=circle, inner sep=0.2] at (-1.2,1.1) {\small $6$};
	\node[draw,shape=circle, inner sep=0.2] at (-0.5,2) {\small $7$};
	\node[draw,shape=circle, inner sep=0.2] at (-1.5,-0.9) {\small $8$};
	\node[draw,shape=circle, inner sep=0.2] at (-2,0) {\small $9$};
	\node[draw,shape=circle, inner sep=0.2] at (0.3,-1.65) {\footnotesize $10$};
	\node[draw,shape=circle, inner sep=0.2] at (-0.5,-2) {\footnotesize $11$};
	\node[draw,shape=circle, inner sep=0.2] at (1.6,-0.2) {\footnotesize $12$};
	\node[draw,shape=circle, inner sep=0.2] at (1.5,-1.2) {\footnotesize $13$};
	\node[draw,shape=circle, inner sep=0.2] at (0.7,1.5) {\footnotesize $14$};
	\node[draw,shape=circle, inner sep=0.2] at (1.55,1.2) {\footnotesize $15$};
	\node[draw,shape=circle, inner sep=0.2] at (0.65,2.5) {\footnotesize $16$};
	\node[draw,shape=circle, inner sep=0.2] at (-2.3,1.25) {\footnotesize $17$};
	\node[draw,shape=circle, inner sep=0.2] at (-1.8,-1.65) {\footnotesize $18$};
	\node[draw,shape=circle, inner sep=0.2] at (1.05,-2.45) {\footnotesize $19$};
	\node[draw,shape=circle, inner sep=0.2] at (2.55,0.2) {\footnotesize $20$};
		\end{tikzpicture}
		\caption{$T(20a^3,60a^4;60\aaa\bbb^2\ccc,12\ccc^5)$.}
		\label{ab2r1}
	\end{figure}

In Figure \ref{ab2r1}, $\ccc^5$ determines $T_1,\cdots,T_5$.  By the symmetry of the partial tiling and $\bbb^2\cdots=\aaa\bbb^2\ccc$, we might as well take $\thin\bbb_1\thin\bbb_5\thin\cdots=\thin\bbb_1\thin\bbb_5\thin\aaa\thin\ccc\thin$ which determines $T_6,T_7$. Then we can determine $T_8,\cdots,T_{20},\cdots$ and we get a tiling whose 3D picture is the 23rd of Figure \ref{Fig3}.\label{20,60}

\vspace{9pt}
\noindent{Subcase $\thin\ccc\thin\ccc\thin\cdots=\aaa\ccc^4$, AVC $\subset\{\aaa\bbb^2\ccc,\aaa\ccc^4\}$}

The vertex $\aaa\ccc^4$ determines $T_1,\cdots,T_5$ in Figure \ref{ab2r2}. By the symmetry of the partial tiling and $\aaa\bbb\cdots=\aaa\bbb^2\ccc$, we might as well take $\thin\aaa_1\thin\bbb_5\thin\cdots=\thin\aaa_1\thin\bbb_5\thin\bbb\thin\ccc\thin$ which determines $T_6,T_7$. Then $T_{8},\cdots,T_{39}$ are determined. But $\thin\bbb_{35}\thin\ccc_{36}\thin\ccc_{39}\thin$ $\cdots$ appears, contradicting the AVC.  

	\begin{figure}[htp]
		\centering
		\begin{tikzpicture}[>=latex,scale=0.7]       
		\draw (-1.19,1.55)--(0,1.2)--(1.19,1.55)--(1.16,0.31)--(1.80,-0.75)--(2.31,-1.97)--(1.57,-2.92)--(1.47,-3.95)--(0.06,-4.40)--(-1.45,-3.82)--(-1.68,-2.92)--(-2.37,-1.94)--(-1.8,-0.75)--(-1.16,0.31)--(-1.19,1.55)
        (0,1.2)--(0,0)--(0.6,-1.04)--(0.67,-1.97)--(0.01,-2.92)--(0.06,-4.40)
        (0.01,-2.92)--(-0.6,-1.94)--(-0.6,-1.04)--(0,0)
        (-1.16,0.31)--(0,0)--(1.16,0.31) (-1.8,-0.75)--(-0.6,-1.04)--(0.6,-1.04)--(1.8,-0.75)
        (-2.37,-1.94)--(-0.6,-1.94)--(0.67,-1.97)--(2.31,-1.97) (-1.68,-2.92)--(0.01,-2.92)--(1.57,-2.92)
        (0,1.2)--(-0.6,2.53)--(-0.44,3.78)--(-1.92,3.54)--(-3.33,3.78)--(-3.25,2.56)--(-3.93,1.53)--(-3.91,0.34)--(-3.43,-0.77)--(-3.54,-2.18)--(-1.68,-2.92)
        (-4.67,-2.6)--(-3.54,-2.18)--(-2.37,-1.94)
        (-1.92,3.54)--(-1.9,2.53)--(-2.45,1.53)--(-2.43,0.28)--(-1.80,-0.75)--(-3.43,-0.77)--(-4.73,-0.7)
        (-1.19,1.55)--(-1.9,2.53)
        (-6.13,1.95)--(-4.67,2.59)--(-3.25,2.56)--(-1.9,2.53)--(-0.6,2.53)
        (-7.27,-0.19)--(-6.29,0.36)--(-5.23,0.95)--(-3.93,1.53)--(-2.45,1.53)--(-1.19,1.55)
        (-3.33,3.78)--(-4.97,3.64)--(-6.87,2.96)--(-6.13,1.95)--(-5.23,0.95)--(-4.73,-0.7)--(-4.67,-2.6)--(-6.29,0.36)--(-6.13,1.95)
        (-4.97,3.64)--(-4.67,2.59)--(-3.93,1.53) (-4.73,-0.7)--(-3.91,0.34)--(-2.43,0.28)--(-1.16,0.31)
        (0.06,-4.4)--(-4.46,-5.11)--(-6.63,-3.16)--(-7.27,-0.19)--(-6.87,2.96)
        (-4.67,-2.6)--(-1.45,-3.82)--(-4.15,-4.27)--(-4.67,-2.6)
        (-1.45,-3.82)--(-4.15,-4.27)--(-7.27,-0.19) (-4.15,-4.27)--(-4.46,-5.11);   
		\node at (0.2,-0.2) {\small $\ccc$};
		\node at (0.2,0.3) {\small $\ccc$};
		\node at (-0.2,0.3) {\small $\ccc$};
		\node at (-0.25,-0.2) {\small $\ccc$};
		\node at (-0.85,-1.25) {\small $\bbb$};
		\node at (-0.4,-1.25) {\small $\ccc$};
		\node at (0.8,-1.25) {\small $\ccc$};
		\node at (-0.75,-2.2) {\small $\bbb$};
		\node at (0.75,-2.25) {\small $\bbb$};
		\node at (-0.2,-3.2) {\small $\ccc$};
		\node at (0.2,-3.2) {\small $\ccc$};
		\node at (-2.2,-0.6) {\small $\ccc$};
		\node at (-2.1,-1) {\small $\ccc$};
		\node at (-1.3,0.5) {\small $\ccc$};
		\node at (-2.6,0.5) {\small $\ccc$};
		\node at (-2.55,1.75) {\small $\bbb$};
		\node at (-0.85,1.7) {\small $\bbb$};
		\node at (-2.1,2.7) {\small $\ccc$};
		\node at (-3.5,2.8) {\small $\ccc$};
		\node at (-1.7,2.7) {\small $\ccc$};
		\node at (-4.3,1.5) {\small $\ccc$};
		\node at (-4.1,1.2) {\small $\ccc$};
		\node at (-5,2.6) {\small $\ccc$};
		\node at (-3.7,-1) {\small $\ccc$};
		\node at (-3.5,-2.5) {\small $\bbb$};
		\node at (-5,-0.7) {\small $\bbb$};
		\node at (-6.35,1.8) {\small $\bbb$};
		\node at (-6.35,-0.05) {\small $\bbb$};
		\node at (-1.5,-4.2) {\small $\bbb$};
		\node at (-4.4,-4.4) {\small $\bbb$};
		\node at (-7.27,-0.19) {\small $\bullet$};
	\node[draw,shape=circle, inner sep=0.2] at (0,-0.6) {\small $1$};
	\node[draw,shape=circle, inner sep=0.2] at (0.8,-0.5) {\small $2$};
	\node[draw,shape=circle, inner sep=0.2] at (0.5,0.8) {\small $3$};
	\node[draw,shape=circle, inner sep=0.2] at (-0.5,0.8) {\small $4$};
	\node[draw,shape=circle, inner sep=0.2] at (-0.8,-0.5) {\small $5$};
	\node[draw,shape=circle, inner sep=0.2] at (-1.3,-1.5) {\small $6$};
	\node[draw,shape=circle, inner sep=0.2] at (0,-1.5) {\small $7$};
	\node[draw,shape=circle, inner sep=0.2] at (1.3,-1.5) {\small $8$};
	\node[draw,shape=circle, inner sep=0.2] at (0,-2.3) {\small $9$};
	\node[draw,shape=circle, inner sep=0.2] at (-1.25,-2.5) {\footnotesize $10$};
	\node[draw,shape=circle, inner sep=0.2] at (1.2,-2.5) {\footnotesize $11$};
	\node[draw,shape=circle, inner sep=0.2] at (-0.8,-3.5) {\footnotesize $12$};
	\node[draw,shape=circle, inner sep=0.2] at (0.8,-3.5) {\footnotesize $13$};
	\node[draw,shape=circle, inner sep=0.2] at (-1.75,-0.1) {\footnotesize $14$};
	\node[draw,shape=circle, inner sep=0.2] at (-3,-0.3) {\footnotesize $15$};
	\node[draw,shape=circle, inner sep=0.2] at (-2.8,-1.5) {\footnotesize $16$};
	\node[draw,shape=circle, inner sep=0.2] at (-2.5,-2.25) {\footnotesize $18$};
	\node[draw,shape=circle, inner sep=0.2] at (-1.8,1) {\footnotesize $17$};
	\node[draw,shape=circle, inner sep=0.2] at (-3.1,1) {\footnotesize $19$};
	\node[draw,shape=circle, inner sep=0.2] at (-1.85,1.85) {\footnotesize $20$};
	\node[draw,shape=circle, inner sep=0.2] at (-2.9,2.15) {\footnotesize $21$};
	\node[draw,shape=circle, inner sep=0.2] at (-1.15,2.15) {\footnotesize $22$};
	\node[draw,shape=circle, inner sep=0.2] at (-1.2,3) {\footnotesize $23$};
	\node[draw,shape=circle, inner sep=0.2] at (-2.55,3) {\footnotesize $24$};
	\node[draw,shape=circle, inner sep=0.2] at (-3.9,2.25) {\footnotesize $25$};
	\node[draw,shape=circle, inner sep=0.2] at (-4.8,1.8) {\footnotesize $26$};
	\node[draw,shape=circle, inner sep=0.2] at (-4.5,0.5) {\footnotesize $27$};
	\node[draw,shape=circle, inner sep=0.2] at (-4,-0.4) {\footnotesize $29$};
	\node[draw,shape=circle, inner sep=0.2] at (-4,3.1) {\footnotesize $28$};
	\node[draw,shape=circle, inner sep=0.2] at (-5.6,2.8) {\footnotesize $30$};
	\node[draw,shape=circle, inner sep=0.2] at (-4,-1.5) {\footnotesize $31$};
	\node[draw,shape=circle, inner sep=0.2] at (-2.5,-3) {\footnotesize $32$};
	\node[draw,shape=circle, inner sep=0.2] at (-5.5,-0.1) {\footnotesize $33$};
	\node[draw,shape=circle, inner sep=0.2] at (-5.8,1.1) {\footnotesize $34$};
	\node[draw,shape=circle, inner sep=0.2] at (-6.6,1.2) {\footnotesize $35$};
	\node[draw,shape=circle, inner sep=0.2] at (-6,-1) {\footnotesize $36$};
	\node[draw,shape=circle, inner sep=0.2] at (-3.5,-3.5) {\footnotesize $37$};
	\node[draw,shape=circle, inner sep=0.2] at (-2.5,-4.35) {\footnotesize $38$};
	\node[draw,shape=circle, inner sep=0.2] at (-6,-3) {\footnotesize $39$};
		\end{tikzpicture}
		\caption{$\thin\ccc\thin\ccc\thin\cdots=\aaa\ccc^4$.}
		\label{ab2r2}
	\end{figure}

\vspace{9pt}
\noindent{Subcase $\thin\ccc\thin\ccc\thin\cdots=\aaa^2\ccc^3$, AVC $\subset\{\aaa\bbb^2\ccc,\aaa^2\ccc^3\}$}

	\begin{figure}[htp]
		\centering
		\begin{tikzpicture}[>=latex,scale=0.7]       
		\draw (0,1.91)--(1.1,0)--(0.55,-0.95)--(0,-1.91)--(-1.1,-1.91)--(-2.63,-1.52)--(-2.2,0)--(-1.1,1.91)--(0,1.91)
        (-1.1,1.91)--(-0.55,0.95)--(0,0)--(-0.55,-0.95)--(0,-1.91)
        (-2.2,0)--(-1.1,-1.91)
        (-2.63,-1.52)--(-1.65,-0.95)--(0.55,-0.95) (-1.65,-0.95)--(-1.1,0)--(1.1,0)
        (-1.65,0.95)--(-1.1,0)--(-0.55,0.95)--(0.55,0.95)
        (0,-1.91)--(0.55,-2.86)--(1.1,-1.91)--(1.65,-0.95)--(2.06,0.61)--(0.55,0.95)
        (1.1,-1.91)--(0.55,-0.95) (1.65,-0.95)--(1.11,0)--(2.06,0.61)
        (0,1.91)--(1.25,1.83)--(2.06,0.61) (-2.63,-1.52)--(-1.71,-2.73)--(-0.55,-2.86) (-1.1,-1.91)--(-0.55,-2.86)--(0.55,-2.86); 
        \draw (-1.71,-2.73)--(-2.32,-3.84)--(-3.67,-2.12)--(-3.03,1.08)--(-2.48,2.04)--(-1.1,1.91)
              (-3.67,-2.12)--(-2.63,-1.52) (-3.03,1.08)--(-2.2,0) (-2.48,2.04)--(-1.65,0.95)
              (-2.32,-3.84)--(-1.7,-5.08)--(-1.14,-4)--(0.95,-3.97)--(1.93,-2.99)
              (-1.71,-2.73)--(-1.14,-4)--(-0.55,-2.86) (0.55,-2.86)--(0.95,-3.97); 
        \draw (1.25,1.83)--(1.77,2.88)--(3.12,1.16)--(2.48,-2.04)--(1.93,-2.99)--(0.55,-2.86)
              (2.06,0.61)--(3.12,1.16) (1.65,-0.95)--(2.48,-2.04) (1.1,-1.91)--(1.93,-2.99)
              (-2.48,2.04)--(-1.5,3.02)--(0.59,3.04)--(1.15,4.13)--(1.77,2.88) 
              (-1.5,3.02)--(-1.1,1.91) (0,1.91)--(0.59,3.04)--(1.25,1.83); 
		\node at (0.1,0.25) {\small $\bbb$};
		\node at (0,-0.25) {\small $\bbb$};
		\node at (-0.3,-0.2) {\small $\ccc$};
		\node at (-1.3,0) {\small $\bbb$};
		\node at (-1.1,0.2) {\small $\ccc$};
		\node at (-1.4,-1.2) {\small $\ccc$};
		\node at (-0.45,1.2) {\small $\bbb$};
		\node at (1.1,-0.35) {\small $\ccc$};
		\node at (0.55,-1.3) {\small $\ccc$};
		\node at (0.6,1.2) {\small $\bbb$};
		\node at (-0.1,-2.2) {\small $\bbb$};
		\node at (-1.2,-2.2) {\small $\bbb$};
		\node at (-2.66,-1.8) {\small $\ccc$};
		\node at (-2.75,-1.35) {\small $\ccc$};
		\node at (2.1,0.9) {\small $\ccc$};
		\node at (2.23,0.34) {\small $\ccc$};
		\node at (-2.21,0.35) {\small $\ccc$};
		\node at (1.68,-1.46) {\small $\ccc$};
		\node at (-1,2.15) {\small $\ccc$};
		\node at (0.38,-3.2) {\small $\ccc$};
		\node at (1.2,2.3) {\small $\ccc$};
		\node at (-1.71,-3.28) {\small $\ccc$};
	\node[draw,shape=circle, inner sep=0.2] at (-0.5,0.3) {\small $1$};
	\node[draw,shape=circle, inner sep=0.2] at (0.4,0.5) {\small $2$};
	\node[draw,shape=circle, inner sep=0.2] at (0.3,-0.5) {\small $3$};
	\node[draw,shape=circle, inner sep=0.2] at (-0.8,-0.5) {\small $4$};
	\node[draw,shape=circle, inner sep=0.2] at (-1.7,0) {\small $5$};
	\node[draw,shape=circle, inner sep=0.2] at (-1.1,0.9) {\small $6$};
	\node[draw,shape=circle, inner sep=0.2] at (-0.15,1.54) {\small $7$};
	\node[draw,shape=circle, inner sep=0.2] at (-0.9,-1.4) {\small $8$};
	\node[draw,shape=circle, inner sep=0.2] at (-1.7,-1.5) {\small $9$};
	\node[draw,shape=circle, inner sep=0.2] at (-2.1,-0.8) {\footnotesize $10$};
	\node[draw,shape=circle, inner sep=0.2] at (0,-1.3) {\footnotesize $11$};
	\node[draw,shape=circle, inner sep=0.2] at (1.1,-0.9) {\footnotesize $12$};
	\node[draw,shape=circle, inner sep=0.2] at (1.5,-0.2) {\footnotesize $13$};
	\node[draw,shape=circle, inner sep=0.2] at (1.2,0.5) {\footnotesize $14$};
	\node[draw,shape=circle, inner sep=0.2] at (0.55,-1.9) {\footnotesize $15$};
	\node[draw,shape=circle, inner sep=0.2] at (1.2,1.3) {\footnotesize $16$};
	\node[draw,shape=circle, inner sep=0.2] at (-0.5,-2.3) {\footnotesize $17$};
	\node[draw,shape=circle, inner sep=0.2] at (-1.6,-2.3) {\footnotesize $18$};
	\node[draw,shape=circle, inner sep=0.2] at (-2.53,-2.59) {\footnotesize $19$};
	\node[draw,shape=circle, inner sep=0.2] at (-2.87,-0.34) {\footnotesize $20$};
	\node[draw,shape=circle, inner sep=0.2] at (2,1.69) {\footnotesize $21$};
	\node[draw,shape=circle, inner sep=0.2] at (2.31,-0.53) {\footnotesize $22$};
	\node[draw,shape=circle, inner sep=0.2] at (-2.34,1.11) {\footnotesize $23$};
	\node[draw,shape=circle, inner sep=0.2] at (-1.66,1.59) {\footnotesize $24$};
	\node[draw,shape=circle, inner sep=0.2] at (1.81,-2.12) {\footnotesize $25$};
	\node[draw,shape=circle, inner sep=0.2] at (1.2,-2.57) {\footnotesize $26$};
	\node[draw,shape=circle, inner sep=0.2] at (-1.63,2.35) {\footnotesize $27$};
	\node[draw,shape=circle, inner sep=0.2] at (-0.36,2.51) {\footnotesize $28$};
	\node[draw,shape=circle, inner sep=0.2] at (1.1,-3.31) {\footnotesize $29$};
	\node[draw,shape=circle, inner sep=0.2] at (-0.17,-3.52) {\footnotesize $30$};
	\node[draw,shape=circle, inner sep=0.2] at (0.59,2.3) {\footnotesize $31$};
	\node[draw,shape=circle, inner sep=0.2] at (1.18,3.07) {\footnotesize $32$};
	\node[draw,shape=circle, inner sep=0.2] at (-1.13,-3.28) {\footnotesize $33$};
	\node[draw,shape=circle, inner sep=0.2] at (-1.71,-4.1) {\footnotesize $34$};
		\end{tikzpicture}\hspace{10pt}
		\begin{tikzpicture}[>=latex,scale=0.7]   
				\foreach \a in {0,1,2}
				\draw[rotate=120*\a]
                (0.55,-0.32)--(0,0.64)--(-0.55,1.59)--(0,2.54)--(0.47,3.48)--(1.98,3.83)--(2.39,2.57)--(3.75,0.65)--(2.78,-2.15)
                (2.20,-1.27)--(2.77,0.33)--(1.76,1.63)--(1.1,2.54)--(0,2.54)
                (1.65,-0.32)--(1.1,0.64)--(0.55,1.59)--(1.1,2.54)--(1.98,3.83)
                (-0.55,1.59)--(0.55,1.59)--(1.76,1.63)--(2.39,2.57)
                (0,0.64)--(1.1,0.64)--(2.77,0.33)--(3.75,0.65) (1.65,-0.32)--(2.77,0.33);
		\node at (0.6,0) {\small $\bbb$};
		\node at (0.6,-0.65) {\small $\bbb$};
		\node at (0.2,-0.6) {\small $\ccc$};
		\node at (-0.8,-0.3) {\small $\bbb$};
		\node at (-0.55,0.1) {\small $\ccc$};
		\node at (-1.5,-1.1) {\small $\ccc$};
		\node at (0.1,0.9) {\small $\bbb$};
		\node at (-2.55,-1.1) {\small $\ccc$};
		\node at (-1.3,-2.8) {\small $\ccc$};
		\node at (0.08,-1.5) {\small $\bbb$};
		\node at (-1.3,0.6) {\small $\bbb$};
		\node at (1.2,0.85) {\small $\bbb$};
		\node at (-1,-2.8) {\small $\ccc$};
		\node at (-2.8,0) {\small $\bbb$};
		\node at (-2.4,0.9) {\small $\bbb$};
		\node at (-1.7,2.5) {\small $\ccc$};
		\node at (-0.1,1.8) {\small $\ccc$};
		\node at (0.3,2.7) {\small $\ccc$};
		\node at (1.4,2.6) {\small $\bbb$};
		\node at (2.1,1.65) {\small $\bbb$};
		\node at (1.65,-0.75) {\small $\ccc$};
		\node at (2.8,0) {\small $\ccc$};
		\node at (2.2,-1.8) {\small $\ccc$};
		\node at (0.9,-2.5) {\small $\ccc$};
	\node[draw,shape=circle, inner sep=0.2] at (0,0) {\small $1$};
	\node[draw,shape=circle, inner sep=0.2] at (1,0.2) {\small $2$};
	\node[draw,shape=circle, inner sep=0.2] at (1,-0.7) {\small $3$};
	\node[draw,shape=circle, inner sep=0.2] at (-0.2,-0.8) {\small $4$};
	\node[draw,shape=circle, inner sep=0.2] at (-1.2,-0.4) {\small $5$};
	\node[draw,shape=circle, inner sep=0.2] at (-0.55,0.8) {\small $6$};
	\node[draw,shape=circle, inner sep=0.2] at (0.45,1.1) {\small $7$};
	\node[draw,shape=circle, inner sep=0.2] at (-0.7,-1.65){\small $8$};
	\node[draw,shape=circle, inner sep=0.2] at (-1.1,1.5) {\small $9$};
	\node[draw,shape=circle, inner sep=0.2] at (1.75,0.2) {\footnotesize $10$};
	\node[draw,shape=circle, inner sep=0.2] at (-1.4,-1.65) {\footnotesize $11$};
	\node[draw,shape=circle, inner sep=0.2] at (-2,-0.65) {\footnotesize $12$};
	\node[draw,shape=circle, inner sep=0.2] at (0.4,-1.8) {\footnotesize $13$};
	\node[draw,shape=circle, inner sep=0.2] at (-1.65,0.9) {\footnotesize $14$};
	\node[draw,shape=circle, inner sep=0.2] at (1.7,1.2) {\footnotesize $15$};
	\node[draw,shape=circle, inner sep=0.2] at (-2.2,0) {\footnotesize $16$};
	\node[draw,shape=circle, inner sep=0.2] at (-0.2,-2.9) {\footnotesize $17$};
	\node[draw,shape=circle, inner sep=0.2] at (-2,-2) {\footnotesize $18$};
	\node[draw,shape=circle, inner sep=0.2] at (-3.2,-0.8) {\footnotesize $19$};
	\node[draw,shape=circle, inner sep=0.2] at (-3.15,0.3) {\footnotesize $20$};
	\node[draw,shape=circle, inner sep=0.2] at (-2.5,1.5) {\footnotesize $21$};
	\node[draw,shape=circle, inner sep=0.2] at (-0.6,2.05) {\footnotesize $22$};
	\node[draw,shape=circle, inner sep=0.2] at (-1,2.8) {\footnotesize $23$};
	\node[draw,shape=circle, inner sep=0.2] at (0.3,2.15) {\footnotesize $24$};
	\node[draw,shape=circle, inner sep=0.2] at (1.1,1.9) {\footnotesize $25$};
	\node[draw,shape=circle, inner sep=0.2] at (0.8,3.1) {\footnotesize $26$};
	\node[draw,shape=circle, inner sep=0.2] at (1.8,2.9) {\footnotesize $27$};
	\node[draw,shape=circle, inner sep=0.2] at (2.9,1) {\footnotesize $28$};
	\node[draw,shape=circle, inner sep=0.2] at (2.1,-0.5) {\footnotesize $29$};
	\node[draw,shape=circle, inner sep=0.2] at (1.7,-1.4) {\footnotesize $30$};
	\node[draw,shape=circle, inner sep=0.2] at (2.7,-1) {\footnotesize $31$};
	\node[draw,shape=circle, inner sep=0.2] at (1.1,-2) {\footnotesize $32$};
	\node[draw,shape=circle, inner sep=0.2] at (2.2,-2.5) {\footnotesize $33$};
	\node[draw,shape=circle, inner sep=0.2] at (1.5,-2.9) {\footnotesize $34$};
		\end{tikzpicture}
		\caption{$\{20a^3,36a^4;36\aaa\bbb^2\ccc,12\aaa^2\ccc^3:2\}$.}
		\label{ab2r3.1}
	\end{figure}

In Figure \ref{ab2r3.1}, $\thin\aaa\thin\bbb\thin\bbb\thin\ccc\thin$ determines $T_1,\cdots,T_4$. When $\thin\aaa_1\thin\bbb_4\thin\cdots=\thin\aaa_1\thin\bbb_4\thin\bbb\thin\ccc\thin$, $T_5,T_6,T_7$ are determined. We see that the partial tiling is rotationally symmetric by rotating $120^\circ$ around the center of $T_1$. Then we have the following situations: 

			\begin{itemize}
		\item $\thin\ccc_5\thin\ccc_4\thin\cdots=\thin\ccc_5\thin\ccc_4\thin\ccc\thin\aaa\thin\aaa\thin$;
		\item $\thin\ccc_5\thin\ccc_4\thin\cdots=\thin\ccc_5\thin\ccc_4\thin\aaa\thin\aaa\thin\ccc\thin$,$\thin\ccc_7\thin\ccc_6\thin\cdots=\thin\ccc_7\thin\ccc_6\thin\aaa\thin\cdots$,$\thin\ccc_3\thin\ccc_2\thin\cdots=\thin\ccc_3\thin\ccc_2\thin\aaa\thin\cdots$;
        \item $\thin\ccc_5\thin\ccc_4\thin\cdots=\thin\ccc_5\thin\ccc_4\thin\aaa\thin\ccc\thin\aaa\thin$,$\thin\ccc_7\thin\ccc_6\thin\cdots=\thin\ccc_7\thin\ccc_6\thin\aaa\thin\ccc\thin\aaa\thin$,$\thin\ccc_3\thin\ccc_2\thin\cdots=\thin\ccc_3\thin\ccc_2\thin\aaa\thin\ccc\thin\aaa\thin$.
	       \end{itemize}

In situation 1, we can determine $T_8,T_9,T_{10}$ in the first of Figure \ref{ab2r3.1}. Then $T_{11},\cdots,T_{26}$ are determined. When $\thin\ccc_7\thin\ccc_6\thin\aaa_{24}\thin\cdots=\thin\ccc_7\thin\ccc_6\thin\aaa_{24}\thin\ccc\thin\cdots$, we get $\thin\bbb\thin\aaa_{24}\thin\ccc_{23}\thin\cdots=\thin\bbb\thin\aaa_{24}\thin\ccc_{23}\thin\bbb\thin$. But $\thin\ccc\thin\bbb_{23}\thin\ccc_{20}\thin\cdots$ appears, contradicting the AVC. So $\thin\ccc_7\thin\ccc_6\thin\aaa_{24}\thin\cdots=\thin\ccc_7\thin\ccc_6\thin\aaa_{24}\thin\aaa\thin\ccc\thin$ determines $T_{27},$ $T_{28}$. Similarly, $T_{29},T_{30}$ are determined. Then we can get a tiling and the 3D picture is the 27th of Figure \ref{Fig3}.\label{20,36}

In situation 2, we can determine $T_8,T_9,\cdots,T_{12}$ in the second of Figure \ref{ab2r3.1}. Similarly, we get a different tiling and the 3D picture is the 28th of Figure \ref{Fig3}.

In situation 3, we can determine $T_8,\cdots,T_{13}$ in the first of Figure \ref{ab2r3.3}. Then $T_{14},T_{15},T_{16}$ are determined. But $\bbb\ccc^2\cdots$ appears, contradicting the AVC.

 	\begin{figure}[htp]
		\centering
		\begin{tikzpicture}[>=latex,scale=0.7]   
				\foreach \a in {0,1,2}
				\draw[rotate=120*\a]
                (0.55,-0.32)--(-0.55,1.59)--(0.55,1.59)--(1.65,-0.32) (0,0.64)--(1.1,0.64);  
				\foreach \a in {0,1}
				\draw[rotate=120*\a]
                (-1.1,0.64)--(-1.5,2.14)--(-1.5,3.24)--(-0.55,2.69)--(0.55,1.59) (-1.5,2.14)--(-0.55,1.59)--(-0.55,2.69);
                \draw (-1.5,2.14)--(-2.22,0.87)--(-1.65,-0.32) (-1.5,3.24)--(-2.45,2.33)--(-2.22,0.87)--(-2.64,-0.4)--(-2.05,-1.82);
		\node at (0.6,0) {\small $\bbb$};
		\node at (0.6,-0.65) {\small $\bbb$};
		\node at (0.2,-0.6) {\small $\ccc$};
		\node at (-0.8,-0.3) {\small $\bbb$};
		\node at (-0.55,0.1) {\small $\ccc$};
		\node at (0.1,0.9) {\small $\bbb$};
		\node at (-1.3,-1.6) {\small $\ccc$};
		\node at (-0.8,2) {\small $\ccc$};
		\node at (-1.3,0.7) {\small $\bbb$};
		\node at (-2.15,1.3) {\small $\ccc$};
		\node at (-2.2,0.4) {\small $\ccc$};
		\node at (-2.22,0.87) {\small $\bullet$};
	\node[draw,shape=circle, inner sep=0.2] at (0,0) {\small $1$};
	\node[draw,shape=circle, inner sep=0.2] at (1,0.2) {\small $2$};
	\node[draw,shape=circle, inner sep=0.2] at (1,-0.7) {\small $3$};
	\node[draw,shape=circle, inner sep=0.2] at (-0.2,-0.8) {\small $4$};
	\node[draw,shape=circle, inner sep=0.2] at (-1.2,-0.4) {\small $5$};
	\node[draw,shape=circle, inner sep=0.2] at (-0.55,0.8) {\small $6$};
	\node[draw,shape=circle, inner sep=0.2] at (0.45,1.1) {\small $7$};
	\node[draw,shape=circle, inner sep=0.2] at (-0.75,-1.65) {\small $8$};
	\node[draw,shape=circle, inner sep=0.2] at (-1.7,-2.1) {\small $9$};
	\node[draw,shape=circle, inner sep=0.2] at (-1.55,-1.05) {\footnotesize $10$};
	\node[draw,shape=circle, inner sep=0.2] at (-1,1.5) {\footnotesize $11$};
	\node[draw,shape=circle, inner sep=0.2] at (-1.1,2.5) {\footnotesize $12$};
	\node[draw,shape=circle, inner sep=0.2] at (-0.2,1.85) {\footnotesize $13$};
	\node[draw,shape=circle, inner sep=0.2] at (-1.7,0.9) {\footnotesize $14$};
	\node[draw,shape=circle, inner sep=0.2] at (-2,2) {\footnotesize $15$};
	\node[draw,shape=circle, inner sep=0.2] at (-2.1,-0.5) {\footnotesize $16$};
		\end{tikzpicture}\hspace{10pt}
		\begin{tikzpicture}[>=latex,scale=0.7]       
		\draw (0.55,2.05)--(1.67,1.83)--(3.44,2.09)--(4.66,0.98)--(4.1,0.03)--(3.97,-1.83)--(2.2,-1.91)--(1.1,-1.91)--(-1.1,-1.91)--(-1.65,-0.95)--(-0.55,0.95)--(0.55,2.05)
        (0.55,2.05)--(0.55,0.95)--(-0.55,-0.95)--(0,-1.91)--(1.1,0)--(0.55,0.95)
        (3.44,2.09)--(2.88,0.98)--(2.2,0)--(1.1,-1.91)
        (-0.55,0.95)--(1.65,0.95)--(2.88,0.98)--(4.1,0.03)--(3.04,-0.87)--(1.65,-0.95)--(-1.65,-0.95)
        (2.2,-1.91)--(3.04,-0.87)--(2.2,0)--(1.65,0.95)--(1.67,1.83) (-1.1,0)--(2.2,0);   
		\node at (1.2,0.25) {\small $\bbb$};
		\node at (1.1,-0.35) {\small $\bbb$};
		\node at (0.8,-0.2) {\small $\ccc$};
		\node at (-0.1,0.2) {\small $\bbb$};
		\node at (-0.25,-0.2) {\small $\ccc$};
		\node at (-0.6,-1.25) {\small $\bbb$};
		\node at (0.55,-1.2) {\small $\bbb$};
		\node at (0.8,1.2) {\small $\ccc$};
		\node at (1.9,1.2) {\small $\ccc$};
		\node at (3.2,1) {\small $\ccc$};
		\node at (2.9,0.65) {\small $\bbb$};
		\node at (1.75,-1.2) {\small $\bbb$};
		\node at (3.1,-1.1) {\small $\bbb$};
		\node at (4.1,0.03) {\small $\bullet$};
	\node[draw,shape=circle, inner sep=0.2] at (0.55,0.3) {\small $1$};
	\node[draw,shape=circle, inner sep=0.2] at (1.5,0.65) {\small $2$};
	\node[draw,shape=circle, inner sep=0.2] at (1.5,-0.65) {\small $3$};
	\node[draw,shape=circle, inner sep=0.2] at (0.2,-0.5) {\small $4$};
	\node[draw,shape=circle, inner sep=0.2] at (-0.8,-0.5) {\small $5$};
	\node[draw,shape=circle, inner sep=0.2] at (-0.35,0.5) {\small $6$};
	\node[draw,shape=circle, inner sep=0.2] at (0,-1.2) {\small $7$};
	\node[draw,shape=circle, inner sep=0.2] at (1,-1.5) {\small $8$};
	\node[draw,shape=circle, inner sep=0.2] at (-1,-1.4) {\small $9$};
	\node[draw,shape=circle, inner sep=0.2] at (0.2,1.25) {\footnotesize $10$};
	\node[draw,shape=circle, inner sep=0.2] at (1.2,1.5) {\footnotesize $11$};
	\node[draw,shape=circle, inner sep=0.2] at (2.25,0.55) {\footnotesize $12$};
	\node[draw,shape=circle, inner sep=0.2] at (2.5,1.5) {\footnotesize $13$};
	\node[draw,shape=circle, inner sep=0.2] at (3.7,1) {\footnotesize $14$};
	\node[draw,shape=circle, inner sep=0.2] at (3,0) {\footnotesize $15$};
	\node[draw,shape=circle, inner sep=0.2] at (2.24,-0.55) {\footnotesize $16$};
	\node[draw,shape=circle, inner sep=0.2] at (2.1,-1.5) {\footnotesize $17$};
	\node[draw,shape=circle, inner sep=0.2] at (3.5,-1.35) {\footnotesize $18$};
		\end{tikzpicture}\hspace{10pt}
		\begin{tikzpicture}[>=latex,scale=0.7]       
		\draw (-1.1,1.91)--(1.1,1.91)--(1.65,2.86)--(2.75,0.95)--(3.15,-0.55)--(4.47,-0.56)--(3.73,-1.85)--(3.04,-2.86)--(1.65,-2.86)--(-0.55,-2.86)--(-1.1,-1.91)--(-1.65,-0.95)--(-1.1,0)--(-0.55,0.95)--(-1.1,1.91)
        (1.65,2.86)--(1.1,1.91)--(0.55,0.95)--(-0.55,-0.95)--(0,-1.91)--(1.1,0)--(0.55,0.95)
        (-0.55,0.95)--(1.65,0.95)--(2.2,1.91) (1.1,1.91)--(2.2,0)--(-1.1,0)
        (2.2,0)--(3.15,-0.55)--(1.65,-0.95)--(-1.65,-0.95)
        (-1.1,-1.91)--(1.1,-1.91)--(2.33,-1.83)--(3.73,-1.85) (2.75,0.95)--(2.2,0)--(1.65,-0.95)--(1.1,-1.91)--(1.65,-2.86)--(2.33,-1.83)--(3.15,-0.55) (-0.55,-2.86)--(0,-1.91);   
		\node at (1.2,0.25) {\small $\bbb$};
		\node at (1.1,-0.35) {\small $\bbb$};
		\node at (0.8,-0.2) {\small $\ccc$};
		\node at (-0.1,0.2) {\small $\bbb$};
		\node at (-0.25,-0.2) {\small $\ccc$};
		\node at (-0.6,-1.25) {\small $\bbb$};
		\node at (0.55,-1.2) {\small $\bbb$};
		\node at (0.5,1.2) {\small $\ccc$};
		\node at (0.1,-2.15) {\small $\ccc$};
		\node at (1.65,1.15) {\small $\ccc$};
		\node at (1.9,1) {\small $\bbb$};
		\node at (1.8,-1.15) {\small $\bbb$};
		\node at (2.7,-1.55) {\small $\ccc$};
		\node at (2.4,-2.1) {\small $\bbb$};
		\node at (3.15,-0.55) {\small $\bullet$};
	\node[draw,shape=circle, inner sep=0.2] at (0.55,0.3) {\small $1$};
	\node[draw,shape=circle, inner sep=0.2] at (1.5,0.65) {\small $2$};
	\node[draw,shape=circle, inner sep=0.2] at (1.5,-0.6) {\small $3$};
	\node[draw,shape=circle, inner sep=0.2] at (0.2,-0.5) {\small $4$};
	\node[draw,shape=circle, inner sep=0.2] at (-0.8,-0.5) {\small $5$};
	\node[draw,shape=circle, inner sep=0.2] at (-0.35,0.5) {\small $6$};
	\node[draw,shape=circle, inner sep=0.2] at (0,-1.2) {\small $7$};
	\node[draw,shape=circle, inner sep=0.2] at (1,-1.5) {\small $8$};
	\node[draw,shape=circle, inner sep=0.2] at (-1,-1.4) {\small $9$};
	\node[draw,shape=circle, inner sep=0.2] at (-0.1,1.3) {\footnotesize $10$};
	\node[draw,shape=circle, inner sep=0.2] at (1.1,1.25) {\footnotesize $11$};
	\node[draw,shape=circle, inner sep=0.2] at (0.7,-2.5) {\footnotesize $12$};
	\node[draw,shape=circle, inner sep=0.2] at (-0.5,-2.2) {\footnotesize $13$};
	\node[draw,shape=circle, inner sep=0.2] at (1.65,2) {\footnotesize $14$};
	\node[draw,shape=circle, inner sep=0.2] at (2.3,1) {\footnotesize $15$};
	\node[draw,shape=circle, inner sep=0.2] at (2.6,0.1) {\footnotesize $16$};
	\node[draw,shape=circle, inner sep=0.2] at (2.3,-0.4) {\footnotesize $17$};
	\node[draw,shape=circle, inner sep=0.2] at (2.1,-1.5) {\footnotesize $18$};
	\node[draw,shape=circle, inner sep=0.2] at (1.75,-2.2) {\footnotesize $19$};
	\node[draw,shape=circle, inner sep=0.2] at (3.4,-1.4) {\footnotesize $20$};
	\node[draw,shape=circle, inner sep=0.2] at (2.85,-2.4) {\footnotesize $21$};
		\end{tikzpicture}
		\caption{$\thin\ccc\thin\ccc\thin\cdots=\aaa^2\ccc^3$.}
		\label{ab2r3.3}
	\end{figure}

When $\thin\aaa_1\thin\bbb_4\thin\cdots$ can only be $\thin\aaa_1\thin\bbb_4\thin\ccc\thin\bbb\thin$, $T_5,T_6,\cdots,T_9$ are determined in the second and third of Figure \ref{ab2r3.3}. By the symmetry of the partial tiling, we might as well take $\thin\ccc_2\thin\aaa_1\thin\ccc_6\thin\cdots=\thin\ccc_2\thin\aaa_1\thin\ccc_6\thin\aaa\thin\ccc\thin$. Then we can determine $T_{10},T_{11},\cdots,T_{18}$ in the second. But $\bbb\ccc^2\cdots$ appears, contradicting the AVC. Hence $\thin\ccc_2\thin\aaa_1\thin\ccc_6\thin\cdots=\thin\ccc_2\thin\aaa_1\thin\ccc_6\thin\ccc\thin\aaa\thin$ and $\thin\ccc_9\thin\aaa_7\thin\ccc_8\thin\cdots=\thin\ccc_9\thin\aaa_7\thin\ccc_8\thin\ccc\thin\aaa\thin$. Then we can determine $T_{10},\cdots,T_{19}$ in the third. Since the partial tiling $\{T_{19},T_{12},T_8,T_{18}\}$ appears, $\thin\aaa_{19}\thin\bbb_{18}\thin\cdots=\thin\aaa_{19}\thin\bbb_{18}\thin\ccc\thin\bbb\thin$ determines $T_{20},T_{21}$. But $\aaa^2\bbb\ccc\cdots$ appears, contradicting the AVC.

\vspace{9pt}
\noindent{Subcase $\thin\ccc\thin\ccc\thin\cdots=\aaa^3\ccc^2$, AVC $\subset\{\aaa\bbb^2\ccc,\aaa^3\ccc^2\}$}
	\begin{figure}[H]
		\centering
		\begin{tikzpicture}[>=latex,scale=0.7]
				\foreach \a in {0,1,2}
				\draw[rotate=120*\a]
				(1.65,-0.32)--(0.55,-0.32)--(0,0.64)--(1.1,0.64)--(1.65,-0.32)--(1.1,-1.27)--(0,-1.27)
				(-0.36,2.67)--(1.03,2.54) (3.4,1.38)--(4.4,-0.24)--(3.45,-2.04)--(1.68,-2.16)
                (1.10,0.64)--(2.49,0.39)--(3.34,-0.32)--(2.49,-1.02)--(1.1,-1.27)
                (1.03,2.54)--(1.96,1.43)--(2.49,0.39)--(2.49,-1.02)--(1.65,-0.32)--(2.49,0.39)
                (1.1,0.64)--(2.49,0.39)
                (0.55,1.59)--(1.96,1.43)--(3.4,1.38)--(3.34,-0.32)--(3.45,-2.04);
		\node at (0.6,-0.05) {\small $\bbb$};
		\node at (0.55,-0.6) {\small $\bbb$};
		\node at (0.2,-0.5) {\small $\ccc$};
		\node at (1.2,0.8) {\small $\bbb$};
		\node at (0.75,0.8) {\small $\ccc$};
		\node at (1.15,-1.5) {\small $\bbb$};
		\node at (0.7,-1.5) {\small $\ccc$};
		\node at (2.6,0.7) {\small $\ccc$};
		\node at (2.6,-1.2) {\small $\ccc$};
		\node at (3.5,-0.3) {\small $\ccc$};
		\node at (-0.25,0.55) {\small $\bbb$};
		\node at (-0.8,-0.3) {\small $\bbb$};
		\node at (-1.3,0.64) {\small $\bbb$};
		\node at (0.6,1.8) {\small $\bbb$};
		\node at (-1.7,2.05) {\small $\ccc$};
		\node at (-0.36,2.85) {\small $\ccc$};
		\node at (-1.4,3.25) {\small $\ccc$};
		\node at (-1.9,-0.3) {\small $\bbb$};
		\node at (-2.3,-1.75) {\small $\ccc$};
		\node at (-0.91,-2.5) {\small $\ccc$};
		\node at (-1.9,-3) {\small $\ccc$};
	\node[draw,shape=circle, inner sep=0.2] at (0,0) {\small $1$};
	\node[draw,shape=circle, inner sep=0.2] at (1,0.2) {\small $2$};
	\node[draw,shape=circle, inner sep=0.2] at (1,-0.9) {\small $3$};
	\node[draw,shape=circle, inner sep=0.2] at (-0.3,-0.9) {\small $4$};
	\node[draw,shape=circle, inner sep=0.2] at (1.8,0.15) {\small $5$};
	\node[draw,shape=circle, inner sep=0.2] at (2.2,-0.3) {\small $6$};
	\node[draw,shape=circle, inner sep=0.2] at (1.8,-0.75) {\small $7$};
	\node[draw,shape=circle, inner sep=0.2] at (1.6,1.1) {\small $8$};
	\node[draw,shape=circle, inner sep=0.2] at (0.2,1.2) {\small $9$};
	\node[draw,shape=circle, inner sep=0.2] at (1.5,-1.7) {\footnotesize $10$};
	\node[draw,shape=circle, inner sep=0.2] at (0,-2) {\footnotesize $11$};
	\node[draw,shape=circle, inner sep=0.2] at (2.8,-0.3) {\footnotesize $12$};
	\node[draw,shape=circle, inner sep=0.2] at (3,1) {\footnotesize $13$};
	\node[draw,shape=circle, inner sep=0.2] at (3,-1.5) {\footnotesize $14$};
	\node[draw,shape=circle, inner sep=0.2] at (4,-0.3) {\footnotesize $15$};
	\node[draw,shape=circle, inner sep=0.2] at (-0.7,0.55) {\footnotesize $16$};
	\node[draw,shape=circle, inner sep=0.2] at (-1.25,-0.4) {\footnotesize $17$};
	\node[draw,shape=circle, inner sep=0.2] at (-1,1.4) {\footnotesize $18$};
	\node[draw,shape=circle, inner sep=0.2] at (-0.8,2) {\footnotesize $19$};
	\node[draw,shape=circle, inner sep=0.2] at (-0.1,1.9) {\footnotesize $20$};
	\node[draw,shape=circle, inner sep=0.2] at (-1.65,1) {\footnotesize $21$};
	\node[draw,shape=circle, inner sep=0.2] at (1,2) {\footnotesize $22$};
	\node[draw,shape=circle, inner sep=0.2] at (-1.15,2.65) {\footnotesize $23$};
	\node[draw,shape=circle, inner sep=0.2] at (-2.2,2) {\footnotesize $24$};
	\node[draw,shape=circle, inner sep=0.2] at (0,3.2) {\footnotesize $25$};
	\node[draw,shape=circle, inner sep=0.2] at (-1.8,3.6) {\footnotesize $26$};
	\node[draw,shape=circle, inner sep=0.2] at (-1.6,-1.1) {\footnotesize $27$};
	\node[draw,shape=circle, inner sep=0.2] at (-1.3,-1.7) {\footnotesize $28$};
	\node[draw,shape=circle, inner sep=0.2] at (-0.7,-1.65) {\footnotesize $29$};
	\node[draw,shape=circle, inner sep=0.2] at (-2.3,-0.3) {\footnotesize $30$};
	\node[draw,shape=circle, inner sep=0.2] at (-1.7,-2.2) {\footnotesize $31$};
	\node[draw,shape=circle, inner sep=0.2] at (-2.8,-1.5) {\footnotesize $32$};
	\node[draw,shape=circle, inner sep=0.2] at (-0.8,-2.95) {\footnotesize $33$};
	\node[draw,shape=circle, inner sep=0.2] at (-2.3,-3.3) {\footnotesize $34$};
		\end{tikzpicture}
		\caption{$T(20a^3,24a^4;24\aaa\bbb^2\ccc,12\aaa^3\ccc^2)$.}
		\label{ab2r4}
	\end{figure}

In Figure \ref{ab2r4}, $\thin\aaa\thin\bbb\thin\bbb\thin\ccc\thin$ determines $T_1,\cdots,T_4$, $\ccc^2\cdots=\aaa^3\ccc^2$ determines $T_5,T_6,T_7$, etc. Then we get a tiling and the 3D picture is the 21st of Figure \ref{Fig3}.\label{20,24.1}

\vspace{9pt}

\nd Now $\aaa\bbb^2\ccc$ can only be $\thin\aaa\thin\bbb\thin\ccc\thin\bbb\thin$ which determines $T_1,\cdots,T_{13}$ in Figure \ref{ab2r5}. Similar to the proof of Lemma \ref{deg345}, we can get $\aaa^3\ccc\cdots=\aaa^3\ccc^2,\aaa^4\ccc$. When $\thin\aaa_{12}\thin\ccc_7\thin\aaa_8\thin\aaa_{13}\thin\cdots=\aaa^3\ccc^2$, we get $\thin\ccc_{11}\thin\aaa_{12}\thin\bbb\thin\cdots=\thin\ccc_{11}\thin\aaa_{12}\thin\bbb\thin\bbb\thin$, a contradiction. So $\thin\aaa_{12}\thin\ccc_7\thin\aaa_8\thin\aaa_{13}\thin\cdots=\aaa^4\ccc$ and $\thin\aaa_{6}\thin\ccc_3\thin\aaa_8\thin\aaa_{13}\thin$ $\cdots=\aaa^4\ccc$ determine $T_{14},T_{15}$, and we get the AVC $\subset\{\aaa\bbb^2\ccc,\aaa^4\ccc\}$. But $\thin\aaa_{14}\thin\aaa_{13}\thin\aaa_{15}\thin\cdots$ is not a vertex.

	\begin{figure}[htp]
		\centering
        \begin{tikzpicture}[>=latex,scale=0.75]       
		\draw (-0.55,2.86) -- (1.65,2.86) -- (2.75,0.95) -- (2.2,0) -- (-1.1,0) -- (-1.65,0.95) -- (-0.55,2.86) -- (1.1,0)--(2.2,1.91)
              (1.65,2.86) -- (0,0) -- (-1.1,1.91) (-1.1,0)--(0.55,2.86)--(2.2,0)
              (0,0)--(0.55,-0.95)--(1.1,0)
              (-1.1,0) arc (175.63:304.5:1.06)
              (2.2,0) arc (4.37:-124.5:1.06);
		\node at (1.35,1.95) {\small $\bbb$};
		\node at (1.1,1.55) {\small $\ccc$};
		\node at (0.8,1.91) {\small $\bbb$};
		\node at (1.9,1) {\small $\bbb$};
		\node at (0.3,1) {\small $\bbb$};
		\node at (-0.25,1.9) {\small $\bbb$};
		\node at (-0.85,0.95) {\small $\bbb$};
		\node at (0.55,-0.95) {\small $\bullet$};
	\node[draw,shape=circle, inner sep=0.2] at (1.1,2.55) {\small $1$};
	\node[draw,shape=circle, inner sep=0.2] at (1.75,2) {\small $2$};
	\node[draw,shape=circle, inner sep=0.2] at (1.1,1) {\small $3$};
	\node[draw,shape=circle, inner sep=0.2] at (0.4,2) {\small $4$};
	\node[draw,shape=circle, inner sep=0.2] at (2.35,1) {\small $5$};
	\node[draw,shape=circle, inner sep=0.2] at (1.65,0.3) {\small $6$};
	\node[draw,shape=circle, inner sep=0.2] at (-0.15,1) {\small $7$};
	\node[draw,shape=circle, inner sep=0.2] at (0.55,0.3) {\small $8$};
	\node[draw,shape=circle, inner sep=0.2] at (-0.7,2) {\small $9$};
	\node[draw,shape=circle, inner sep=0.2] at (0,2.6) {\footnotesize $10$};
	\node[draw,shape=circle, inner sep=0.2] at (-1.3,1) {\footnotesize $11$};
	\node[draw,shape=circle, inner sep=0.2] at (-0.55,0.3) {\footnotesize $12$};
	\node[draw,shape=circle, inner sep=0.2] at (0.55,-0.3) {\footnotesize $13$};
	\node[draw,shape=circle, inner sep=0.2] at (-0.5,-0.4) {\footnotesize $14$};
	\node[draw,shape=circle, inner sep=0.2] at (1.5,-0.4) {\footnotesize $15$};
		\end{tikzpicture}
		\caption{$\aaa\bbb^2\ccc$ can only be $\thin\aaa\thin\bbb\thin\ccc\thin\bbb\thin$.}
		\label{ab2r5}
	\end{figure}

\subsubsection*{Case $\aaa^3\bbb$}

By $\aaa<\frac{\pi}{2}$, we have $\aaa+\bbb>\pi$. Similar to the proof of Lemma \ref{deg345}, $\bbb^2\cdots$ is not a vertex. So $\thin\ccc\thin\ccc\thin\cdots$ can never be a vertex by Lemma \ref{AAD}, and $\aaa^2\cdots=\aaa^3\bbb,\aaa^2\bbb\ccc^2,\aaa^2\bbb\ccc^3,\aaa^5\ccc^k(1\le k\le5),\aaa^4\ccc^k(2\le k\le4),\aaa^3\ccc^3$. Obviously $\thin\aaa\thin\aaa\thin\bbb\thin\cdots$ can only be $\thin\aaa\thin\aaa\thin\aaa\thin\bbb\thin$. As shown in figure \ref{rar1}, when $\thin\ccc\thin\aaa\thin\ccc\thin\cdots$ appears, we get two $\thin\aaa\thin\bbb\thin\cdots$ which are $\thin\aaa\thin\aaa\thin\aaa\thin\bbb\thin$.

In the left of Figure \ref{a3bf}, $\aaa^3\bbb$ determines $T_1,\cdots,T_4$. When $\thin\aaa_3\thin\aaa_2\thin\cdots=\thin\aaa_3\thin\aaa_2\thin\ccc\thin\cdots$, $T_5$ is determined, and $\thin\bbb_5\thin\aaa_2\thin\ccc_1\thin\cdots=\thin\bbb_5\thin\aaa_2\thin\ccc_1\thin\aaa\thin\cdots=\thin\bbb_5\thin\aaa_2\thin\ccc_1\thin\aaa\thin\ccc\thin$ determines $T_6,T_7$. Then we get the AVC $\subset\{\aaa^3\bbb,\aaa^2\bbb\ccc^2\}$ which determines $T_8,\cdots,T_{13}$. But $\aaa^3\ccc\cdots$ appears, contradicting the AVC. So $\thin\aaa_3\thin\aaa_2\thin\cdots\not=\thin\aaa_3\thin\aaa_2\thin\ccc\thin\cdots$, $\thin\aaa_3\thin\aaa_4\thin\cdots\not=\thin\aaa_3\thin\aaa_4\thin\ccc\thin\cdots$ and they are $\aaa^4\ccc^2,\aaa^5\ccc,\aaa^5\ccc^2,\aaa^5\ccc^3,\aaa^3\bbb$.

	\begin{figure}[H]
		\centering
        \begin{tikzpicture}[>=latex,scale=0.7]       
		\draw (0,0) -- (0,1.2) -- (2,1.2) -- (2,0) -- (1,-0.5) -- (0,0) -- (1,0.5) -- (2,1.2)
              (0,1.2) -- (1,0.5) -- (2,0) 
              (0,1.2)--(-0.93,0.37)--(-0.9,-0.82)--(-0.79,-2.04)--(1.11,-1.67)--(2,0)
              (-0.9,-0.82)--(0,0)--(0,-1.16)--(1.11,-1.67)--(1,-0.5)--(0,-1.16)--(-0.79,-2.04)
              (-0.93,0.37)--(-0.79,1.93)--(0.61,2.3)--(3.39,1.46)--(2,0)
              (0,1.2)--(0.61,2.3) (2,1.2)--(3.39,1.46);
		\node at (1,0.2) {\small $\bbb$};
		\node at (-0.2,0.7) {\small $\ccc$};
		\node at (-0.2,-0.48) {\small $\ccc$};
		\node at (-0.1,1.3) {\small $\bbb$};
		\node at (0.3,1.4) {\small $\ccc$};
		\node at (2,0) {\small $\bullet$};
	\node[draw,shape=circle, inner sep=0.2] at (1,-0.2) {\small $1$};
	\node[draw,shape=circle, inner sep=0.2] at (0.3,0.5) {\small $2$};
	\node[draw,shape=circle, inner sep=0.2] at (1,0.9) {\small $3$};
	\node[draw,shape=circle, inner sep=0.2] at (1.7,0.5) {\small $4$};
	\node[draw,shape=circle, inner sep=0.2] at (-0.5,0) {\small $5$};
	\node[draw,shape=circle, inner sep=0.2] at (0.3,-0.5) {\small $6$};
	\node[draw,shape=circle, inner sep=0.2] at (-0.5,-1) {\small $7$};
	\node[draw,shape=circle, inner sep=0.2] at (0.7,-1.1) {\small $8$};
	\node[draw,shape=circle, inner sep=0.2] at (0,-1.6) {\small $9$};
	\node[draw,shape=circle, inner sep=0.2] at (1.3,-0.8) {\footnotesize $10$};
	\node[draw,shape=circle, inner sep=0.2] at (-0.5,1.5) {\footnotesize $11$};
	\node[draw,shape=circle, inner sep=0.2] at (1.2,1.5) {\footnotesize $12$};
	\node[draw,shape=circle, inner sep=0.2] at (2.4,0.9) {\footnotesize $13$};
		\end{tikzpicture}\hspace{30pt}
        \begin{tikzpicture}[>=latex,scale=0.7]      
		\draw (0,0) -- (0,1.2) -- (2,1.2) -- (2,0) -- (1,-0.5) -- (0,0) -- (1,0.5) -- (2,1.2)
              (0,1.2) -- (1,0.5) -- (2,0) 
              (0,1.2)--(-0.93,0.37)--(-0.9,-0.82)--(0,0);
		\node at (-0.2,0.7) {\small $\ccc$};
		\node at (0.35,0) {\small $\ccc$};
		\node at (0,1.2) {\small $\bullet$};
	\node[draw,shape=circle, inner sep=0.2] at (-0.5,0) {\small $1$};
	\node[draw,shape=circle, inner sep=0.2] at (0.3,0.5) {\small $2$};
	\node[draw,shape=circle, inner sep=0.2] at (1,0.9) {\small $3$};
	\node[draw,shape=circle, inner sep=0.2] at (1,-0) {\small $4$};
	\node[draw,shape=circle, inner sep=0.2] at (1.7,0.5) {\small $5$};
		\end{tikzpicture}
		\caption{$\thin\aaa_3\thin\aaa_2\thin\cdots=\thin\aaa_3\thin\aaa_2\thin\ccc\thin\cdots$ or $\thin\ccc\thin\aaa\thin\aaa\thin\cdots$ appears.}
		\label{a3bf}
	\end{figure}

		If $\thin\ccc\thin\aaa\thin\aaa\thin\cdots$ appears, then $T_1,T_2,T_3$ are determined in the right of Figure \ref{a3bf}. When $\thin\bbb_1\thin\aaa_2\thin\cdots=\thin\bbb_1\thin\aaa_2\thin\ccc\thin\cdots$, $T_4,T_5$ are determined, a contradiction. So $\thin\bbb_1\thin\aaa_2\thin\cdots=\thin\bbb_1\thin\aaa_2\thin\aaa\thin\cdots=\thin\bbb_1\thin\aaa_2\thin\aaa\thin\aaa\thin$.

\vspace{9pt}
\noindent{Subcase $\thin\aaa_3\thin\aaa_2\thin\cdots=\aaa^4\ccc^2$, AVC $\subset \{\aaa^3\bbb,\aaa^4\ccc^2\}$}

The vertex $\thin\aaa_3\thin\aaa_2\thin\cdots=\thin\aaa_3\thin\aaa_2\thin\aaa\thin\ccc\thin\aaa\thin\ccc\thin$ determines $T_5,\cdots,$ $T_8$ in Figure \ref{a3b1}. Then we can determine $T_9,\cdots,T_{26}$ and the 3D picture is the 17th of Figure \ref{Fig3}.\label{20,6}

	\begin{figure}[htp]
		\centering
		\begin{tikzpicture}[>=latex,scale=0.7]
				\foreach \a in {0,1}
				\draw[rotate=180*\a]
				(-0.9,0.52)--(0.3,0.52)--(1.5,0.52)--(0.9,-0.52)--(0.3,0.52)--(-0.3,1.56)--(0.3,2.6)--(0.9,1.56)--(0.3,0.52)
				(-2.21,1.97)--(-0.3,1.56)--(0.9,1.56)--(2.81,1.97)--(0.3,2.6)--(-2.21,1.97)--(0.3,3.63)--(2.81,1.97)--(2.21,-1.97)--(1.5,0.52)
                (0.9,-0.52)--(2.21,-1.97) (1.5,0.52)--(2.81,1.97) (0,0)--(0.3,0.52) (-0.3,-3.63)--(-1.17,-4.22)
                (-0.3,-3.63) arc (-98.73:40.68:3.42);
		\node at (0.85,-0.8) {\small $\bbb$};
		\node at (-0.1,0.7) {\small $\ccc$};
		\node at (0.6,0.7) {\small $\ccc$};
		\node at (0.5,2.8) {\small $\bbb$};
		\node at (-0.5,-2.8) {\small $\bbb$};
		\node at (-0.6,-0.8) {\small $\ccc$};
	\node[draw,shape=circle, inner sep=0.2] at (0.5,-1.2) {\small $1$};
	\node[draw,shape=circle, inner sep=0.2] at (0.3,-0.2) {\small $2$};
	\node[draw,shape=circle, inner sep=0.2] at (0.9,0.2) {\small $3$};
	\node[draw,shape=circle, inner sep=0.2] at (1.4,-0.5) {\small $4$};
	\node[draw,shape=circle, inner sep=0.2] at (-0.3,0.2) {\small $5$};
	\node[draw,shape=circle, inner sep=0.2] at (-0.9,1.2) {\small $6$};
	\node[draw,shape=circle, inner sep=0.2] at (0.3,1.1) {\small $7$};
	\node[draw,shape=circle, inner sep=0.2] at (1.5,1.2) {\small $8$};
	\node[draw,shape=circle, inner sep=0.2] at (2.1,0) {\small $9$};
	\node[draw,shape=circle, inner sep=0.2] at (-0.5,2) {\footnotesize $10$};
	\node[draw,shape=circle, inner sep=0.2] at (0.3,2) {\footnotesize $11$};
	\node[draw,shape=circle, inner sep=0.2] at (1.1,2) {\footnotesize $12$};
	\node[draw,shape=circle, inner sep=0.2] at (-0.9,-0.2) {\footnotesize $13$};
	\node[draw,shape=circle, inner sep=0.2] at (-1.5,-1.2) {\footnotesize $16$};
	\node[draw,shape=circle, inner sep=0.2] at (-0.3,-1.2) {\footnotesize $15$};
	\node[draw,shape=circle, inner sep=0.2] at (-1.4,0.5) {\footnotesize $14$};
	\node[draw,shape=circle, inner sep=0.2] at (-2,-0.5) {\footnotesize $17$};
	\node[draw,shape=circle, inner sep=0.2] at (-1.1,-2) {\footnotesize $18$};
	\node[draw,shape=circle, inner sep=0.2] at (-0.3,-2) {\footnotesize $19$};
	\node[draw,shape=circle, inner sep=0.2] at (0.4,-2) {\footnotesize $20$};
	\node[draw,shape=circle, inner sep=0.2] at (0,-3) {\footnotesize $21$};
	\node[draw,shape=circle, inner sep=0.2] at (0,3) {\footnotesize $23$};
	\node[draw,shape=circle, inner sep=0.2] at (2.8,-1.5) {\footnotesize $22$};
	\node[draw,shape=circle, inner sep=0.2] at (-2.8,1.5) {\footnotesize $24$};
	\node[draw,shape=circle, inner sep=0.2] at (3,3) {\footnotesize $25$};
	\node[draw,shape=circle, inner sep=0.2] at (-3,-3) {\footnotesize $26$};
		\end{tikzpicture}
		\caption{$T(20a^3,6a^4;12\aaa^3\bbb,6\aaa^4\ccc^2)$.}
		\label{a3b1}
	\end{figure}

\vspace{9pt}
\noindent{Subcase $\thin\aaa_3\thin\aaa_2\thin\cdots=\aaa^5\ccc$,  AVC $\subset \{\aaa^3\bbb,\aaa^5\ccc\}$}

In the left of Figure \ref{a3ba5r}, $\aaa^5\ccc$ determines $T_1,\cdots,T_6$; $\aaa\bbb\cdots=\aaa^3\bbb$ determines $T_7,\cdots,T_{10}$; $\aaa^2\ccc\cdots=\aaa^5\ccc$ determines $T_{11},T_{12},T_{13}$. By the symmetry of the partial tiling, we have the following situations:

			\begin{itemize}
		\item $\thin\aaa_7\thin\aaa_8\thin\aaa_{13}\thin\cdots=\thin\aaa_7\thin\aaa_8\thin\aaa_{13}\thin\bbb\thin$ and $\thin\aaa_9\thin\aaa_{10}\thin\aaa_{11}\thin\cdots=\thin\aaa_9\thin\aaa_{10}\thin\aaa_{11}\thin\bbb\thin$;
		\item $\thin\aaa_7\thin\aaa_8\thin\aaa_{13}\thin\cdots=\thin\aaa_7\thin\aaa_8\thin\aaa_{13}\thin\bbb\thin$ and $\thin\aaa_9\thin\aaa_{10}\thin\aaa_{11}\thin\cdots=\thin\aaa_9\thin\aaa_{10}\thin\aaa_{11}\thin\ccc\thin\aaa\thin\aaa\thin$;
        \item $\thin\aaa_7\thin\aaa_8\thin\aaa_{13}\thin\cdots=\thin\aaa_7\thin\aaa_8\thin\aaa_{13}\thin\bbb\thin$ and $\thin\aaa_9\thin\aaa_{10}\thin\aaa_{11}\thin\cdots=\thin\aaa_9\thin\aaa_{10}\thin\aaa_{11}\thin\aaa\thin\ccc\thin\aaa\thin$;
        \item $\thin\aaa_7\thin\aaa_8\thin\aaa_{13}\thin\cdots=\thin\aaa_7\thin\aaa_8\thin\aaa_{13}\thin\bbb\thin$ and $\thin\aaa_9\thin\aaa_{10}\thin\aaa_{11}\thin\cdots=\thin\aaa_9\thin\aaa_{10}\thin\aaa_{11}\thin\aaa\thin\aaa\thin\ccc\thin$;
		\item $\thin\aaa_7\thin\aaa_8\thin\aaa_{13}\thin\cdots,\thin\aaa_7\thin\aaa_2\thin\aaa_3\thin\cdots,\thin\aaa_9\thin\aaa_{10}\thin\aaa_{11}\thin\cdots,\thin\aaa_9\thin\aaa_6\thin\aaa_5\thin\cdots$ are not $\aaa^3\bbb$.
	       \end{itemize}

	\begin{figure}[htp]
		\centering
        \begin{tikzpicture}[>=latex,scale=0.7]       
		\draw (1.02,-2.35) -- (1.37,-0.71) -- (1.37,0.71) -- (1.02,2.35) -- (0,1.53) -- (0.6,0) -- (0,-1.53) -- (1.02,-2.35)
              (0,-1.53) -- (1.37,-0.71) -- (0.6,0)-- (1.37,0.71)--(0,1.53)
              (-1.02,-2.35) -- (-1.37,-0.71) -- (-1.37,0.71) -- (-1.02,2.35) -- (0,1.53) -- (-0.6,0) -- (0,-1.53) -- (-1.02,-2.35)
              (0,-1.53) -- (-1.37,-0.71) -- (-0.6,0)-- (-1.37,0.71)--(0,1.53)
              (1.02,2.35)--(-1.02,2.35) (1.02,-2.35)--(-1.02,-2.35);
        \draw(1.02,2.35)--(2.61,1.88)--(1.37,0.71)--(2.35,-0.71)--(1.02,-2.35)--(2.61,-3.25)--(0,-3.36)--(1.02,-2.35) 
             (2.61,1.88)--(2.35,-0.71)--(2.61,-3.25)
             (-1.02,2.35)--(-2.61,1.88)--(-1.37,0.71)--(-2.35,-0.71)--(-1.02,-2.35)--(-2.61,-3.25)--(0,-3.36)--(-1.02,-2.35) 
             (-2.61,1.88)--(-2.35,-0.71)--(-2.61,-3.25)
             (2.61,-3.25) arc (-25.62:25.62:5.94)
             (-2.61,-3.25) arc (-154.38:-205.62:5.94); 
        \draw[dashed] (0,3)--(0,-4);
		\node at (0,1) {\small $\ccc$};
		\node at (1.5,-0.71) {\small $\bbb$};
		\node at (-1.5,-0.71) {\small $\bbb$};
		\node at (0.8,3) {$L_1$};
	\node[draw,shape=circle, inner sep=0.2] at (0,0) {\small $1$};
	\node[draw,shape=circle, inner sep=0.2] at (0.7,0.5) {\small $2$};
	\node[draw,shape=circle, inner sep=0.2] at (0.7,1.5) {\small $3$};
	\node[draw,shape=circle, inner sep=0.2] at (0,2) {\small $4$};
	\node[draw,shape=circle, inner sep=0.2] at (-0.7,1.5) {\small $5$};
	\node[draw,shape=circle, inner sep=0.2] at (-0.7,0.5) {\small $6$};
	\node[draw,shape=circle, inner sep=0.2] at (1,0) {\small $7$};
	\node[draw,shape=circle, inner sep=0.2] at (0.7,-0.5) {\small $8$};
	\node[draw,shape=circle, inner sep=0.2] at (-1,0) {\small $9$};
	\node[draw,shape=circle, inner sep=0.2] at (-0.7,-0.6) {\small $10$};
	\node[draw,shape=circle, inner sep=0.2] at (0.7,-1.5) {\footnotesize $13$};
	\node[draw,shape=circle, inner sep=0.2] at (0,-2) {\footnotesize $12$};
	\node[draw,shape=circle, inner sep=0.2] at (-0.7,-1.5) {\footnotesize $11$};
	\node[draw,shape=circle, inner sep=0.2] at (1.9,-0.7) {\footnotesize $14$};
	\node[draw,shape=circle, inner sep=0.2] at (-1.9,-0.7) {\footnotesize $15$};
	\node[draw,shape=circle, inner sep=0.2] at (0,-2.8) {\footnotesize $16$};
	\node[draw,shape=circle, inner sep=0.2] at (1,-2.9) {\footnotesize $17$};
	\node[draw,shape=circle, inner sep=0.2] at (2,-2) {\footnotesize $18$};
	\node[draw,shape=circle, inner sep=0.2] at (1.5,1.8) {\footnotesize $19$};
	\node[draw,shape=circle, inner sep=0.2] at (2,0.8) {\footnotesize $20$};
	\node[draw,shape=circle, inner sep=0.2] at (2.7,-0.7) {\footnotesize $21$};
	\node[draw,shape=circle, inner sep=0.2] at (-1,-2.9) {\footnotesize $22$};
	\node[draw,shape=circle, inner sep=0.2] at (-2,-2) {\footnotesize $23$};
	\node[draw,shape=circle, inner sep=0.2] at (-1.5,1.8) {\footnotesize $24$};
	\node[draw,shape=circle, inner sep=0.2] at (-2,0.8) {\footnotesize $25$};
	\node[draw,shape=circle, inner sep=0.2] at (-2.7,-0.7) {\footnotesize $26$};
		\end{tikzpicture}\hspace{30pt}
        \begin{tikzpicture}[>=latex,scale=0.75]       
		\draw (-2.09,0.11) -- (0,1.16) -- (2.09,0.11) -- (2.57,-1.9) -- (1.35,-3.84) -- (-1.35,-3.84) -- (-2.57,-1.91) -- (-2.09,0.11)--(0,0)--(2.09,0.11)
              (-2.09,0.11) -- (-1.08,-1.11) -- (0,0)-- (1.08,-1.11)--(2.09,0.11)
              (-2.57,-1.91) -- (-1.08,-1.11) -- (1.08,-1.11) -- (2.57,-1.9) -- (0.69,-2.62) -- (1.08,-1.11) -- (-0.87,-2.22) -- (-2.57,-1.91)
              (-1.35,-3.84) -- (-0.87,-2.22) -- (0.69,-2.62)-- (-1.35,-3.84);
        \draw[dashed] (-3.42,-2.09)--(3.15,-0.66);
		\node at (0,0.85) {\small $\bbb$};
		\node at (-1.8,-1.8) {\small $\ccc$};
		\node at (1.9,-2.4) {\small $\ccc$};
		\node at (3.35,-0.66) {$L_2$};
	\node[draw,shape=circle, inner sep=0.2] at (0,0.35) {\small $27$};
	\node[draw,shape=circle, inner sep=0.2] at (-1.9,-1) {\small $28$};
	\node[draw,shape=circle, inner sep=0.2] at (-1.05,-0.5) {\small $29$};
	\node[draw,shape=circle, inner sep=0.2] at (0,-0.6) {\small $30$};
	\node[draw,shape=circle, inner sep=0.2] at (1.05,-0.5) {\small $31$};
	\node[draw,shape=circle, inner sep=0.2] at (1.9,-1) {\small $32$};
	\node[draw,shape=circle, inner sep=0.2] at (-1,-1.7) {\small $33$};
	\node[draw,shape=circle, inner sep=0.2] at (0.3,-2) {\small $34$};
	\node[draw,shape=circle, inner sep=0.2] at (1.4,-1.8) {\small $35$};
	\node[draw,shape=circle, inner sep=0.2] at (-1.5,-2.65) {\small $36$};
	\node[draw,shape=circle, inner sep=0.2] at (-0.5,-2.8) {\footnotesize $37$};
	\node[draw,shape=circle, inner sep=0.2] at (1,-3.2) {\footnotesize $38$};
		\end{tikzpicture}
		\caption{$\thin\aaa_7\thin\aaa_8\thin\aaa_{13}\thin\cdots=\aaa^3\bbb$ and $\thin\aaa_9\thin\aaa_{10}\thin\aaa_{11}\thin\cdots=\aaa^3\bbb$.}
		\label{a3ba5r}
	\end{figure}

In situation 1, we can determine $T_{14},T_{15},\cdots,T_{26}$ and the partial tiling is symmetric with respect to line $L_1$. The rest of the region is a regular heptagon of angle $\bbb$. By $\bbb=2\aaa+\ccc$, it is not difficult to find that the tiling of the regular heptagon is unique as the right of Figure \ref{a3ba5r} which is symmetric with respect to line $L_2$. Therefore, we get four different tilings and the 3D pictures are the 29th through the 32nd in Figure \ref{Fig3}.\label{32,6}

In situation 2, we can determine $T_{14},\cdots,T_{27}$ in Figure \ref{a3ba5r1}. By $\thin\aaa_{23}\thin\aaa_{22}\thin\aaa_{21}\thin\cdots\not=\thin\aaa_{23}\thin\aaa_{22}\thin$ $\aaa_{21}\thin\bbb\thin$, we can get a new tiling whose 3D picture is the 33rd of Figure \ref{Fig3}.

	\begin{figure}[H]
		\centering
        \begin{tikzpicture}[>=latex,scale=0.7]       
		\draw (1.02,-2.35) -- (1.37,-0.71) -- (1.37,0.71) -- (1.02,2.35) -- (0,1.53) -- (0.6,0) -- (0,-1.53) -- (1.02,-2.35)
              (0,-1.53) -- (1.37,-0.71) -- (0.6,0)-- (1.37,0.71)--(0,1.53)
              (-1.02,-2.35) -- (-1.37,-0.71) -- (-1.37,0.71) -- (-1.02,2.35) -- (0,1.53) -- (-0.6,0) -- (0,-1.53) -- (-1.02,-2.35)
              (0,-1.53) -- (-1.37,-0.71) -- (-0.6,0)-- (-1.37,0.71)--(0,1.53)
              (1.02,2.35)--(-1.02,2.35) (1.02,-2.35)--(-1.02,-2.35);
        \draw (1.02,2.35)--(2.74,1.24)--(2.53,-1.01)--(2.02,-3.52)--(1.02,-2.35)--(2.53,-1.01)--(1.37,0.71)--(2.74,1.24)
              (-1.37,0.71)--(-2.79,-0.45)--(-2.21,-2.12)--(-1.84,-3.73)--(-1.02,-2.35)
              (-2.79,-0.45)--(-1.37,-0.71)--(-2.21,-2.12)--(-3.72,-2.57) (-2.79,-0.45)--(-3.72,-2.57)--(-1.84,-3.73)--(1.02,-2.35) (-0.83,-5.03)--(2.02,-3.52)
              (-1.84,-3.73)--(2.02,-3.52) (-0.83,-5.03)--(-1.84,-3.73);
        \draw (2.74,1.24) arc (40.18:-57.24:3.2)
              (-3.72,-2.57) arc (-177.67:-83.27:2.58);
		\node at (0,1) {\small $\ccc$};
		\node at (1.5,-0.71) {\small $\bbb$};
		\node at (-1.4,-1.2) {\small $\ccc$};
	\node[draw,shape=circle, inner sep=0.2] at (0,0) {\small $1$};
	\node[draw,shape=circle, inner sep=0.2] at (0.7,0.5) {\small $2$};
	\node[draw,shape=circle, inner sep=0.2] at (0.7,1.5) {\small $3$};
	\node[draw,shape=circle, inner sep=0.2] at (0,2) {\small $4$};
	\node[draw,shape=circle, inner sep=0.2] at (-0.7,1.5) {\small $5$};
	\node[draw,shape=circle, inner sep=0.2] at (-0.7,0.5) {\small $6$};
	\node[draw,shape=circle, inner sep=0.2] at (1,0) {\small $7$};
	\node[draw,shape=circle, inner sep=0.2] at (0.7,-0.5) {\small $8$};
	\node[draw,shape=circle, inner sep=0.2] at (-1,0) {\small $9$};
	\node[draw,shape=circle, inner sep=0.2] at (-0.7,-0.6) {\small $10$};
	\node[draw,shape=circle, inner sep=0.2] at (0.7,-1.5) {\footnotesize $13$};
	\node[draw,shape=circle, inner sep=0.2] at (0,-2) {\footnotesize $12$};
	\node[draw,shape=circle, inner sep=0.2] at (-0.7,-1.5) {\footnotesize $11$};
	\node[draw,shape=circle, inner sep=0.2] at (1.95,-0.9) {\footnotesize $14$};
	\node[draw,shape=circle, inner sep=0.2] at (-1.65,-2) {\footnotesize $15$};
	\node[draw,shape=circle, inner sep=0.2] at (-2.05,-1.11) {\footnotesize $16$};
	\node[draw,shape=circle, inner sep=0.2] at (-1.81,-0.19) {\footnotesize $17$};
	\node[draw,shape=circle, inner sep=0.2] at (-0.65,-2.78) {\footnotesize $18$};
	\node[draw,shape=circle, inner sep=0.2] at (0.78,-3.15) {\footnotesize $19$};
	\node[draw,shape=circle, inner sep=0.2] at (1.84,-2.28) {\footnotesize $20$};
	\node[draw,shape=circle, inner sep=0.2] at (1.71,1.38) {\footnotesize $21$};
	\node[draw,shape=circle, inner sep=0.2] at (2.13,0.42) {\footnotesize $22$};
	\node[draw,shape=circle, inner sep=0.2] at (2.95,-1.22) {\footnotesize $23$};
	\node[draw,shape=circle, inner sep=0.2] at (-2.84,-1.72) {\footnotesize $24$};
	\node[draw,shape=circle, inner sep=0.2] at (-2.5,-2.75) {\footnotesize $25$};
	\node[draw,shape=circle, inner sep=0.2] at (-2.34,-4.13) {\footnotesize $26$};
	\node[draw,shape=circle, inner sep=0.2] at (-0.57,-4.21) {\footnotesize $27$};
		\end{tikzpicture}
		\caption{$\thin\aaa_7\thin\aaa_8\thin\aaa_{13}\thin\cdots=\aaa^3\bbb$ and $\thin\aaa_9\thin\aaa_{10}\thin\aaa_{11}\thin\cdots=\thin\aaa_9\thin\aaa_{10}\thin\aaa_{11}\thin\ccc\thin\aaa\thin\aaa\thin$.}
		\label{a3ba5r1}
	\end{figure}

In situation 3, we can determine $T_{14},\cdots,T_{29}$ as the first of Figure \ref{f}. But it is the same as the left of Figure \ref{a3ba5r}. 

	\begin{figure}[htp]
		\centering
        \begin{tikzpicture}[>=latex,scale=0.65]       
		\draw (1.02,-2.35) -- (1.37,-0.71) -- (1.37,0.71) -- (1.02,2.35) -- (0,1.53) -- (0.6,0) -- (0,-1.53) -- (1.02,-2.35)
              (0,-1.53) -- (1.37,-0.71) -- (0.6,0)-- (1.37,0.71)--(0,1.53)
              (-1.02,-2.35) -- (-1.37,-0.71) -- (-1.37,0.71) -- (-1.02,2.35) -- (0,1.53) -- (-0.6,0) -- (0,-1.53) -- (-1.02,-2.35)
              (0,-1.53) -- (-1.37,-0.71) -- (-0.6,0)-- (-1.37,0.71)--(0,1.53)
              (1.02,2.35)--(-1.02,2.35) (1.02,-2.35)--(-1.02,-2.35);
        \draw (-5.04,3.44)--(-4.91,-0.74)--(-4.93,-4.58)--(0.04,-3.57)--(-1.02,-2.35)--(-3.61,-2.86)--(-4.91,-0.74)--(-3.74,1.51)--(-5.04,3.44)--(-1.02,2.35)
              (-4.93,-4.58)--(-3.61,-2.86)--(-3,-1.43)--(-1.37,-0.71)--(-3.03,0.11)--(-1.37,0.71)--(-3.74,1.51)--(-3.03,0.11)--(-4.91,-0.74)--(-3,-1.43) (-3,-1.43)--(-1.02,-2.35);
        \draw (1.37,0.71)--(2.69,-0.95)--(2.18,-3.52)--(0.04,-3.57)--(1.02,-2.35)--(2.18,-3.52)
              (1.02,-2.35)--(2.69,-0.95);
        \draw (-5.03,3.44) arc (148.89:212.63:7.59);
		\node at (0,1) {\small $\ccc$};
		\node at (1.5,-0.71) {\small $\bbb$};
		\node at (-2,-0.71) {\small $\ccc$};
		\node at (-1.18,-2.75) {\small $\ccc$};
		\node at (-1.6,1.11) {\small $\ccc$};
	\node[draw,shape=circle, inner sep=0.2] at (0,0) {\small $1$};
	\node[draw,shape=circle, inner sep=0.2] at (0.7,0.5) {\small $2$};
	\node[draw,shape=circle, inner sep=0.2] at (0.7,1.5) {\small $3$};
	\node[draw,shape=circle, inner sep=0.2] at (0,2) {\small $4$};
	\node[draw,shape=circle, inner sep=0.2] at (-0.7,1.5) {\small $5$};
	\node[draw,shape=circle, inner sep=0.2] at (-0.7,0.5) {\small $6$};
	\node[draw,shape=circle, inner sep=0.2] at (1,0) {\small $7$};
	\node[draw,shape=circle, inner sep=0.2] at (0.7,-0.5) {\small $8$};
	\node[draw,shape=circle, inner sep=0.2] at (-1,0) {\small $9$};
	\node[draw,shape=circle, inner sep=0.2] at (-0.7,-0.6) {\small $10$};
	\node[draw,shape=circle, inner sep=0.2] at (0.7,-1.5) {\footnotesize $13$};
	\node[draw,shape=circle, inner sep=0.2] at (0,-2) {\footnotesize $12$};
	\node[draw,shape=circle, inner sep=0.2] at (-0.7,-1.5) {\footnotesize $11$};
	\node[draw,shape=circle, inner sep=0.2] at (1.95,-0.9) {\footnotesize $14$};
	\node[draw,shape=circle, inner sep=0.2] at (-1.73,-1.51) {\footnotesize $15$};
	\node[draw,shape=circle, inner sep=0.2] at (-3,-0.71) {\footnotesize $16$};
	\node[draw,shape=circle, inner sep=0.2] at (-1.87,0) {\footnotesize $17$};
	\node[draw,shape=circle, inner sep=0.2] at (-3.72,-1.72) {\footnotesize $18$};
	\node[draw,shape=circle, inner sep=0.2] at (-2.55,-2.17) {\footnotesize $19$};
	\node[draw,shape=circle, inner sep=0.2] at (-3.77,0.42) {\footnotesize $20$};
	\node[draw,shape=circle, inner sep=0.2] at (-2.79,0.69) {\footnotesize $21$};
	\node[draw,shape=circle, inner sep=0.2] at (-2.5,1.9) {\footnotesize $22$};
	\node[draw,shape=circle, inner sep=0.2] at (-4.46,1.38) {\footnotesize $23$};
	\node[draw,shape=circle, inner sep=0.2] at (-5.57,-0.48) {\footnotesize $24$};
	\node[draw,shape=circle, inner sep=0.2] at (-4.43,-2.54) {\footnotesize $25$};
	\node[draw,shape=circle, inner sep=0.2] at (0,-2.83) {\footnotesize $26$};
	\node[draw,shape=circle, inner sep=0.2] at (1.02,-3.02) {\footnotesize $27$};
	\node[draw,shape=circle, inner sep=0.2] at (1.94,-2.35) {\footnotesize $28$};
	\node[draw,shape=circle, inner sep=0.2] at (-2,-3.41) {\footnotesize $29$};
		\end{tikzpicture}\hspace{30pt}
        \begin{tikzpicture}[>=latex,scale=0.65]       
		\draw (1.02,-2.35) -- (1.37,-0.71) -- (1.37,0.71) -- (1.02,2.35) -- (0,1.53) -- (0.6,0) -- (0,-1.53) -- (1.02,-2.35)
              (0,-1.53) -- (1.37,-0.71) -- (0.6,0)-- (1.37,0.71)--(0,1.53)
              (-1.02,-2.35) -- (-1.37,-0.71) -- (-1.37,0.71) -- (-1.02,2.35) -- (0,1.53) -- (-0.6,0) -- (0,-1.53) -- (-1.02,-2.35)
              (0,-1.53) -- (-1.37,-0.71) -- (-0.6,0)-- (-1.37,0.71)--(0,1.53)
              (1.02,2.35)--(-1.02,2.35) (1.02,-2.35)--(-1.02,-2.35);
        \draw (1.37,0.71)--(2.58,0)--(2.69,-1.8)--(1.37,-0.71)--(2.58,0)--(4.14,-1.32)
              (1.37,0.71)--(2.82,1.38)--(2.58,0)--(4.46,0)--(4.14,-1.32)--(2.69,-1.8)--(2.39,-3.55)--(1.02,-2.35);
        \draw (1.02,2.35)--(2.54,2.96)--(2.82,1.38)--(4.68,1.46)--(4.46,0);
		\node at (0,1) {\small $\ccc$};
		\node at (1.52,-1.14) {\small $\ccc$};
		\node at (2.87,0.26) {\small $\ccc$};
		\node at (1.57,1.08) {\small $\ccc$};
		\node at (2.82,1.38) {\small $\bullet$};
	\node[draw,shape=circle, inner sep=0.2] at (0,0) {\small $1$};
	\node[draw,shape=circle, inner sep=0.2] at (0.7,0.5) {\small $2$};
	\node[draw,shape=circle, inner sep=0.2] at (0.7,1.5) {\small $3$};
	\node[draw,shape=circle, inner sep=0.2] at (0,2) {\small $4$};
	\node[draw,shape=circle, inner sep=0.2] at (-0.7,1.5) {\small $5$};
	\node[draw,shape=circle, inner sep=0.2] at (-0.7,0.5) {\small $6$};
	\node[draw,shape=circle, inner sep=0.2] at (1,0) {\small $7$};
	\node[draw,shape=circle, inner sep=0.2] at (0.7,-0.5) {\small $8$};
	\node[draw,shape=circle, inner sep=0.2] at (-1,0) {\small $9$};
	\node[draw,shape=circle, inner sep=0.2] at (-0.7,-0.6) {\small $10$};
	\node[draw,shape=circle, inner sep=0.2] at (0.7,-1.5) {\footnotesize $13$};
	\node[draw,shape=circle, inner sep=0.2] at (0,-2) {\footnotesize $12$};
	\node[draw,shape=circle, inner sep=0.2] at (-0.7,-1.5) {\footnotesize $11$};
	\node[draw,shape=circle, inner sep=0.2] at (1.79,0) {\footnotesize $14$};
	\node[draw,shape=circle, inner sep=0.2] at (2.16,-0.77) {\footnotesize $15$};
	\node[draw,shape=circle, inner sep=0.2] at (1.79,-2.04) {\footnotesize $16$};
	\node[draw,shape=circle, inner sep=0.2] at (3.11,-1.06) {\footnotesize $17$};
	\node[draw,shape=circle, inner sep=0.2] at (3.8,-0.45) {\footnotesize $18$};
	\node[draw,shape=circle, inner sep=0.2] at (2.24,0.66) {\footnotesize $19$};
	\node[draw,shape=circle, inner sep=0.2] at (1.90,1.88) {\footnotesize $20$};
	\node[draw,shape=circle, inner sep=0.2] at (3.70,0.69) {\footnotesize $21$};
		\end{tikzpicture}
        \begin{tikzpicture}[>=latex,scale=0.65] 
				\foreach \a in {0,1}
				\draw[rotate=180*\a]      
        (3.03,3.41)--(4.54,2.09)--(2.63,1.75)--(3.03,3.41)--(1.02,2.35)--(2.63,1.75)--(3.11,-0.98)--(2.58,-3.31)--(0,-3.73)--(-1.02,-2.35)
        (2.63,1.75)--(1.37,0.71)--(3.11,-0.98)--(1.02,-2.35)--(0,-3.73) 
        (1.02,-2.35)--(2.58,-3.31)
        (2.63,1.75) arc (63.71:-64.91:2.8)
        (4.54,2.09) arc (39.8:-115.67:3.78);
		\draw (1.02,-2.35) -- (1.37,-0.71) -- (1.37,0.71) -- (1.02,2.35) -- (0,1.53) -- (0.6,0) -- (0,-1.53) -- (1.02,-2.35)
              (0,-1.53) -- (1.37,-0.71) -- (0.6,0)-- (1.37,0.71)--(0,1.53)
              (-1.02,-2.35) -- (-1.37,-0.71) -- (-1.37,0.71) -- (-1.02,2.35) -- (0,1.53) -- (-0.6,0) -- (0,-1.53) -- (-1.02,-2.35)
              (0,-1.53) -- (-1.37,-0.71) -- (-0.6,0)-- (-1.37,0.71)--(0,1.53)
              (1.02,2.35)--(-1.02,2.35) (1.02,-2.35)--(-1.02,-2.35);
		\node at (0,1) {\small $\ccc$};
		\node at (1.5,-0.71) {\small $\bbb$};
		\node at (-1.57,-0.16) {\small $\ccc$};
		\node at (-3.37,-1.61) {\small $\ccc$};
		\node at (3.37,1.61) {\small $\ccc$};
		\node at (0,-3.73) {\small $\bullet$};
	\node[draw,shape=circle, inner sep=0.2] at (0,0) {\small $1$};
	\node[draw,shape=circle, inner sep=0.2] at (0.7,0.5) {\small $2$};
	\node[draw,shape=circle, inner sep=0.2] at (0.7,1.5) {\small $3$};
	\node[draw,shape=circle, inner sep=0.2] at (0,2) {\small $4$};
	\node[draw,shape=circle, inner sep=0.2] at (-0.7,1.5) {\small $5$};
	\node[draw,shape=circle, inner sep=0.2] at (-0.7,0.5) {\small $6$};
	\node[draw,shape=circle, inner sep=0.2] at (1,0) {\small $7$};
	\node[draw,shape=circle, inner sep=0.2] at (0.7,-0.5) {\small $8$};
	\node[draw,shape=circle, inner sep=0.2] at (-1,0) {\small $9$};
	\node[draw,shape=circle, inner sep=0.2] at (-0.7,-0.6) {\small $10$};
	\node[draw,shape=circle, inner sep=0.2] at (0.7,-1.5) {\footnotesize $13$};
	\node[draw,shape=circle, inner sep=0.2] at (0,-2) {\footnotesize $12$};
	\node[draw,shape=circle, inner sep=0.2] at (-0.7,-1.5) {\footnotesize $11$};
	\node[draw,shape=circle, inner sep=0.2] at (1.95,-0.9) {\footnotesize $14$};
	\node[draw,shape=circle, inner sep=0.2] at (-1.65,-1.56) {\footnotesize $15$};
	\node[draw,shape=circle, inner sep=0.2] at (-2.24,-0.61) {\footnotesize $16$};
	\node[draw,shape=circle, inner sep=0.2] at (-2.16,0.85) {\footnotesize $17$};
	\node[draw,shape=circle, inner sep=0.2] at (0,2.83) {\footnotesize $18$};
	\node[draw,shape=circle, inner sep=0.2] at (-1.18,3.1) {\footnotesize $19$};
	\node[draw,shape=circle, inner sep=0.2] at (-2.16,2.38) {\footnotesize $20$};
	\node[draw,shape=circle, inner sep=0.2] at (-3.64,0.71) {\footnotesize $21$};
	\node[draw,shape=circle, inner sep=0.2] at (1.68,1.56) {\footnotesize $22$};
	\node[draw,shape=circle, inner sep=0.2] at (2.26,0.64) {\footnotesize $23$};
	\node[draw,shape=circle, inner sep=0.2] at (0,-2.83) {\footnotesize $24$};
	\node[draw,shape=circle, inner sep=0.2] at (1.02,-3.02) {\footnotesize $25$};
	\node[draw,shape=circle, inner sep=0.2] at (1.94,-2.35) {\footnotesize $26$};
	\node[draw,shape=circle, inner sep=0.2] at (3.64,-0.71) {\footnotesize $27$};
	\node[draw,shape=circle, inner sep=0.2] at (2.26,2.43) {\footnotesize $28$};
	\node[draw,shape=circle, inner sep=0.2] at (3.43,2.41) {\footnotesize $29$};
	\node[draw,shape=circle, inner sep=0.2] at (4.83,0) {\footnotesize $30$};
	\node[draw,shape=circle, inner sep=0.2] at (-2.26,-2.43) {\footnotesize $31$};
	\node[draw,shape=circle, inner sep=0.2] at (-3.43,-2.41) {\footnotesize $32$};
	\node[draw,shape=circle, inner sep=0.2] at (-4.83,0) {\footnotesize $33$};
		\end{tikzpicture}
		\caption{The other situations of $\thin\aaa_3\thin\aaa_2\thin\cdots=\aaa^5\ccc$.}
		\label{f}
	\end{figure}

In situation 4, we can determine $T_{14},\cdots,T_{27}$ as the third of Figure \ref{f}. Then $\thin\aaa_{27}\thin\aaa_{23}\thin\aaa_{22}\thin\cdots$ $=\thin\aaa_{27}\thin\aaa_{23}\thin\aaa_{22}\thin\aaa\thin\aaa\thin\ccc\thin$ and $\thin\aaa_{21}\thin\aaa_{16}\thin\aaa_{15}\thin\cdots=\thin\aaa_{21}\thin\aaa_{16}\thin\aaa_{15}\thin\aaa\thin\aaa\thin\ccc\thin$ determine $T_{28},\cdots,T_{33}$. But $\thin\aaa_{24}\thin\aaa_{25}\thin\ccc_{30}\thin\cdots$ is not a vertex.

In situation 5, we can determine $T_{14}$ as the second of Figure \ref{f}. If $\thin\aaa_{14}\thin\aaa_7\thin\aaa_8\thin\aaa_{13}\thin\cdots=\thin\aaa\thin\ccc\thin\aaa_{14}\thin\aaa_7\thin\aaa_8\thin\aaa_{13}\thin$ and $\thin\aaa_{14}\thin\aaa_7\thin\aaa_2\thin\aaa_3\thin\cdots=\thin\aaa\thin\ccc\thin\aaa_{14}\thin\aaa_7\thin\aaa_2\thin\aaa_3\thin$, then $\aaa_{14}\bbb^2\cdots$ appears, a contradiction. So we might as well take $\thin\aaa_{14}\thin\aaa_7\thin\aaa_8\thin\aaa_{13}\thin\cdots=\thin\ccc\thin\aaa\thin\aaa_{14}\thin\aaa_7\thin\aaa_8\thin\aaa_{13}\thin$ which determines $T_{15},\cdots,T_{21}$. But $\aaa\bbb^2\cdots$ appears, a contradiction.

\vspace{9pt}
\noindent{Subcase $\thin\aaa_3\thin\aaa_2\thin\cdots=\aaa^5\ccc^2$, AVC $\subset \{\aaa^3\bbb,\aaa^5\ccc^2\}$}

The vertex $\thin\aaa_3\thin\aaa_2\thin\cdots$ is $\thin\aaa_3\thin\aaa_2\thin\aaa\thin\ccc\thin\aaa\thin\ccc\thin\aaa\thin$, $\thin\aaa_3\thin\aaa_2\thin\aaa\thin\ccc\thin\aaa\thin\aaa\thin\ccc\thin$ or $\thin\aaa_3\thin\aaa_2\thin\aaa\thin\aaa\thin\ccc\thin\aaa\thin\ccc\thin$. When $\thin\aaa_3\thin\aaa_2\thin\cdots=\thin\aaa_3\thin\aaa_2\thin\aaa\thin\ccc\thin\aaa\thin\ccc\thin\aaa\thin$, $T_5,\cdots,T_{14}$ are determined as the left of Figure \ref{f1}, contradicting $\thin\aaa_3\thin\aaa_4\thin\cdots\not=\thin\aaa_3\thin\aaa_4\thin\ccc\thin\cdots$. When $\thin\aaa_3\thin\aaa_2\thin\cdots=\thin\aaa_3\thin\aaa_2\thin\aaa\thin\ccc\thin\aaa\thin\aaa\thin\ccc\thin$, $T_5,\cdots,T_{23}$ are determined as the right of Figure \ref{f1}. But $\thin\aaa_{10}\thin\aaa_4\thin\ccc_1\thin\aaa_{13}\thin\aaa_{23}\thin\cdots$ is not a vertex. 

	\begin{figure}[htp]
		\centering
        \begin{tikzpicture}[>=latex,scale=0.8]       
		\draw (0,0) -- (0,1.2) -- (2,1.2) -- (2,0) -- (1,-0.5) -- (0,0) -- (1,0.5) -- (2,1.2)
              (0,1.2) -- (1,0.5) -- (2,0) 
              (0,0)--(-1,0.6)--(-2.35,0.95)--(-1.19,1.56)--(-0.16,2.35)--(0,1.2)--(-1.19,1.56)
              (-0.16,2.35)--(1.03,3.31)--(1.06,2.12)--(0,1.2)
              (1.03,3.31)--(2.17,2.57)--(1.06,2.12)--(2,1.2)--(2.17,2.57)--(3.41,2.99)--(3.2,1.77)--(2,1.2)--(3.1,0.5)--(3.2,1.77)
              (3.1,0.5)--(3.1,-0.69)--(2,0) (-1,0.6)--(0,1.2);
		\node at (1,0.2) {\small $\bbb$};
		\node at (-0.56,1.14) {\small $\ccc$};
		\node at (0.19,1.69) {\small $\ccc$};
		\node at (2.3,1.59) {\small $\ccc$};
		\node at (2.22,0.79) {\small $\ccc$};
		\node at (2,1.2) {\small $\bullet$};
	\node[draw,shape=circle, inner sep=0.2] at (1,-0.2) {\small $1$};
	\node[draw,shape=circle, inner sep=0.2] at (0.3,0.5) {\small $2$};
	\node[draw,shape=circle, inner sep=0.2] at (1,0.9) {\small $3$};
	\node[draw,shape=circle, inner sep=0.2] at (1.7,0.5) {\small $4$};
	\node[draw,shape=circle, inner sep=0.2] at (-0.32,0.61) {\small $5$};
	\node[draw,shape=circle, inner sep=0.2] at (-1.24,1.06) {\small $6$};
	\node[draw,shape=circle, inner sep=0.2] at (-0.42,1.69) {\small $7$};
	\node[draw,shape=circle, inner sep=0.2] at (0.42,2.2) {\small $8$};
	\node[draw,shape=circle, inner sep=0.2] at (1.01,1.61) {\small $9$};
	\node[draw,shape=circle, inner sep=0.2] at (1.75,1.93) {\footnotesize $10$};
	\node[draw,shape=circle, inner sep=0.2] at (1.4,2.62) {\footnotesize $11$};
	\node[draw,shape=circle, inner sep=0.2] at (2.59,2.12) {\footnotesize $12$};
	\node[draw,shape=circle, inner sep=0.2] at (2.7,1.14) {\footnotesize $13$};
	\node[draw,shape=circle, inner sep=0.2] at (2.54,0.34) {\footnotesize $14$};
		\end{tikzpicture}\hspace{30pt}
        \begin{tikzpicture}[>=latex,scale=0.8]          
		\draw (0,0) -- (0,1.2) -- (2,1.2) -- (2,0) -- (1,-0.5) -- (0,0) -- (1,0.5) -- (2,1.2)
              (0,1.2) -- (1,0.5) -- (2,0) ;
		\draw (0,0)--(0,-1.2)--(2,-1.2)--(2,0) (0,-1.2)--(1,-0.5)--(2,-1.2);
        \draw (0,0)--(-1,0.6)--(-2.43,0.93)--(-1.14,1.61)--(0,1.2)--(-1,0.6)
              (-1.14,1.61)--(-0.08,2.43)--(1.24,2.2)--(3.28,2.2)--(3.28,2.2)--(2,1.2)
               (-0.08,2.43)--(0,1.2)--(1.24,2.2);
        \draw (-2.43,0.93)--(-1.75,-0.29)--(-2.8,-1.51)--(-1.35,-1.27)--(-1.46,-2.38)--(-2.8,-1.51)
              (-1.75,-0.29)--(-1,0.6)--(0,0)--(-1.75,-0.29)
              (0,0)--(-1.35,-1.27)--(0,-1.2)--(-1.46,-2.38)--(-1.48,-3.78)--(-0.19,-2.73)--(0,-1.2)--(1.24,-2.33)--(-0.19,-2.73) (1.24,-2.33)--(3.12,-2.17)--(2,-1.2);
        \draw (2,0) arc (-66.43:5.92:2.15)
              (2,0) arc (62.25:-7.55:2.13);
		\node at (1,0.2) {\small $\bbb$};
		\node at (-0.48,1.14) {\small $\ccc$};
		\node at (0.61,1.4) {\small $\ccc$};
		\node at (-0.71,-0.37) {\small $\ccc$};
		\node at (-0.32,-1.77) {\small $\ccc$};
		\node at (0.66,-1.46) {\small $\ccc$};
		\node at (2,0) {\small $\bullet$};
	\node[draw,shape=circle, inner sep=0.2] at (1,-0.2) {\small $1$};
	\node[draw,shape=circle, inner sep=0.2] at (0.3,0.5) {\small $2$};
	\node[draw,shape=circle, inner sep=0.2] at (1,0.9) {\small $3$};
	\node[draw,shape=circle, inner sep=0.2] at (1.7,0.5) {\small $4$};
	\node[draw,shape=circle, inner sep=0.2] at (-0.37,0.58) {\small $5$};
	\node[draw,shape=circle, inner sep=0.2] at (-1.19,1.06) {\small $6$};
	\node[draw,shape=circle, inner sep=0.2] at (-0.45,1.72) {\small $7$};
	\node[draw,shape=circle, inner sep=0.2] at (0.32,1.96) {\small $8$};
	\node[draw,shape=circle, inner sep=0.2] at (1.56,1.61) {\small $9$};
	\node[draw,shape=circle, inner sep=0.2] at (2.54,0.95) {\footnotesize $10$};
	\node[draw,shape=circle, inner sep=0.2] at (-0.98,0.19) {\footnotesize $11$};
	\node[draw,shape=circle, inner sep=0.2] at (-1.67,0.37) {\footnotesize $12$};
	\node[draw,shape=circle, inner sep=0.2] at (1.65,-0.6) {\footnotesize $13$};
	\node[draw,shape=circle, inner sep=0.2] at (0.35,-0.61) {\footnotesize $14$};
	\node[draw,shape=circle, inner sep=0.2] at (1.03,-0.9){\footnotesize $15$};
	\node[draw,shape=circle, inner sep=0.2] at (-1.32,-0.71) {\footnotesize $16$};
	\node[draw,shape=circle, inner sep=0.2] at (-0.48,-0.82) {\footnotesize $17$};
	\node[draw,shape=circle, inner sep=0.2] at (-1.83,-1.77) {\footnotesize $18$};
	\node[draw,shape=circle, inner sep=0.2] at (-1.01,-1.61) {\footnotesize $19$};
	\node[draw,shape=circle, inner sep=0.2] at (-0.66,-2.38) {\footnotesize $20$};
	\node[draw,shape=circle, inner sep=0.2] at (0.34,-2.06) {\footnotesize $21$};
	\node[draw,shape=circle, inner sep=0.2] at (1.59,-1.77) {\footnotesize $22$};
	\node[draw,shape=circle, inner sep=0.2] at (2.54,-1.19) {\footnotesize $23$};
		\end{tikzpicture}
		\caption{$\thin\aaa_3\thin\aaa_2\thin\cdots$ is $\thin\aaa_3\thin\aaa_2\thin\aaa\thin\ccc\thin\aaa\thin\ccc\thin\aaa\thin$ or $\thin\aaa_3\thin\aaa_2\thin\aaa\thin\ccc\thin\aaa\thin\aaa\thin\ccc\thin$.}
		\label{f1}
	\end{figure}

Hence $\thin\aaa_3\thin\aaa_2\thin\cdots=\thin\aaa_3\thin\aaa_2\thin\aaa\thin\aaa\thin\ccc\thin\aaa\thin\ccc\thin$ determines $T_{5},\cdots,T_9$ in Figure \ref{a3b3}. Then we can get a tiling and the 3D picture is the 22nd of Figure \ref{Fig3}.\label{44,12}

	\begin{figure}[htp]
		\centering
		\begin{tikzpicture}[>=latex,scale=0.9] 
				\foreach \a in {0,1}
				\draw[rotate=180*\a]      
        (0,0)--(0,0.5)--(0.77,1)--(2.3,0)--(2.3,-1)--(-0.77,-1)
        (0,-0.5)--(0.77,0)--(2.3,-1) (0,0.5)--(0.77,0)--(1.53,0.5)--(1.53,-0.5)--(2.3,0)
        (0.77,1)--(0.77,0)--(0.4,-1)--(0.77,-2.01)--(1.14,-1)--(0.77,0)
        (-1.56,-1.35)--(-0.77,-1)--(0.77,-2.01)--(2.3,-1)
        (-1.53,-0.5)--(-1.56,-1.35)--(-0.79,-1.9)--(-0.77,-1)
        (2.3,-1) arc (11.65:-159.15:1.62);
		\node at (-0.19,-0.87) {\small $\ccc$};
		\node at (-1.11,0.43) {\small $\ccc$};
		\node at (-0.42,0.4) {\small $\ccc$};
		\node at (1.14,-0.47) {\small $\ccc$};
		\node at (-0.53,-1.43) {\small $\ccc$};
		\node at (0.53,1.43) {\small $\ccc$};
	\node[draw,shape=circle, inner sep=0.2] at (0.21,-0.74) {\small $1$};
	\node[draw,shape=circle, inner sep=0.2] at (-0.5,-0.5) {\small $2$};
	\node[draw,shape=circle, inner sep=0.2] at (-0.29,0) {\small $3$};
	\node[draw,shape=circle, inner sep=0.2] at (0.29,0) {\small $4$};
	\node[draw,shape=circle, inner sep=0.2] at (-1.03,-0.53) {\small $5$};
	\node[draw,shape=circle, inner sep=0.2] at (-1.27,0) {\small $6$};
	\node[draw,shape=circle, inner sep=0.2] at (-1.32,0.69) {\small $7$};
	\node[draw,shape=circle, inner sep=0.2] at (-0.77,0.66) {\small $8$};
	\node[draw,shape=circle, inner sep=0.2] at (-0.19,0.69) {\small $9$};
	\node[draw,shape=circle, inner sep=0.2] at (0.5,0.51) {\footnotesize $10$};
	\node[draw,shape=circle, inner sep=0.2] at (1.01,0.45) {\footnotesize $14$};
	\node[draw,shape=circle, inner sep=0.2] at (1.27,0) {\footnotesize $15$};
	\node[draw,shape=circle, inner sep=0.2] at (1.38,-0.76) {\footnotesize $16$};
	\node[draw,shape=circle, inner sep=0.2] at (0.77,-0.66) {\footnotesize $11$};
	\node[draw,shape=circle, inner sep=0.2] at (-1.8,0){\footnotesize $17$};
	\node[draw,shape=circle, inner sep=0.2] at (-2.06,0.51) {\footnotesize $18$};
	\node[draw,shape=circle, inner sep=0.2] at (-1.32,1.3) {\footnotesize $19$};
	\node[draw,shape=circle, inner sep=0.2] at (-0.74,1.27) {\footnotesize $20$};
	\node[draw,shape=circle, inner sep=0.2] at (-0.11,1.27) {\footnotesize $21$};
	\node[draw,shape=circle, inner sep=0.2] at (1.8,0) {\footnotesize $22$};
	\node[draw,shape=circle, inner sep=0.2] at (2.01,-0.5) {\footnotesize $23$};
	\node[draw,shape=circle, inner sep=0.2] at (1.4,-1.33) {\footnotesize $24$};
	\node[draw,shape=circle, inner sep=0.2] at (0.74,-1.27) {\small $12$};
	\node[draw,shape=circle, inner sep=0.2] at (0.16,-1.32) {\small $13$};
	\node[draw,shape=circle, inner sep=0.2] at (1.11,-2.35) {\footnotesize $25$};
	\node[draw,shape=circle, inner sep=0.2] at (-1.03,-1.43) {\footnotesize $26$};
	\node[draw,shape=circle, inner sep=0.2] at (-1.3,-1) {\footnotesize $27$};
	\node[draw,shape=circle, inner sep=0.2] at (-1.11,2.35) {\footnotesize $28$};
	\node[draw,shape=circle, inner sep=0.2] at (1.03,1.43) {\footnotesize $29$};
	\node[draw,shape=circle, inner sep=0.2] at (1.3,1) {\footnotesize $30$};
		\end{tikzpicture}
		\caption{$T(44a^3,12a^4;24\aaa^3\bbb,12\aaa^5\ccc^2)$.}
		\label{a3b3}
	\end{figure}

\vspace{9pt}
\noindent{Subcase $\thin\aaa_3\thin\aaa_2\thin\cdots=\aaa^5\ccc^3$, AVC $\subset \{\aaa^3\bbb,\aaa^5\ccc^3\}$}

The vertex $\thin\aaa_3\thin\aaa_2\thin\cdots=\thin\aaa_3\thin\aaa_2\thin\aaa\thin\ccc\thin\aaa\thin\ccc\thin\aaa\thin\ccc\thin$  determines $T_5,\cdots,T_{10}$ in Figure \ref{a5r3}. Then $T_{11},\cdots,T_{30}$ are determined. When $\thin\aaa_{29}\thin\aaa_{28}\thin\cdots=\aaa^5\ccc^3$, $\thin\aaa_{29}\thin\aaa_{28}\thin\aaa\thin\ccc\thin\aaa\thin\ccc\thin\aaa\thin\ccc\thin$ will appear, contradicting $\thin\aaa_{29}\thin\aaa_{28}\thin\cdots=\thin\aaa_{30}\thin\aaa_{29}\thin\aaa_{28}\thin\cdots$. So $\thin\aaa_{29}\thin\aaa_{28}\thin\cdots=\aaa^3\bbb$. Similarly we have $\thin\aaa_{21}\thin\aaa_{22}\thin\cdots=\aaa^3\bbb$. But $\thin\ccc\thin\aaa_{30}\thin\ccc_{27}\thin\aaa_{19}\thin\ccc_{16}\thin\aaa_{20}\thin\ccc\thin\cdots$ appears, contradicting the AVC.
	\begin{figure}[H]
		\centering
        \begin{tikzpicture}[>=latex,scale=0.8]       
		\draw (0,2.48)--(0.45,1.09)--(1.22,0.54)--(2.41,0.14)--(0.93,-0.34)--(0.5,-1.21)--(0,0)--(0.45,1.09)
              (1.22,0.54)--(0,0)--(0.93,-0.34) (0.5,-1.21)--(-0.5,-1.21);      
		\draw (0,2.48)--(-0.45,1.09)--(-1.22,0.54)--(-2.41,0.14)--(-0.93,-0.34)--(-0.5,-1.21)--(0,0)--(-0.45,1.09)
              (-1.22,0.54)--(0,0)--(-0.93,-0.34) (0.5,-1.21)--(-0.5,-1.21);
        \draw (0.5,-1.21)--(1.67,-1.47)--(0.45,-2.24)--(-0.5,-1.21)
              (0.93,-0.34)--(1.67,-1.47)--(2.41,0.14)--(2.75,-1.05)--(1.67,-1.47)--(3.02,-2)--(4.07,-1.45)--(2.75,-1.05) (-0.19,-3.14)--(0.45,-2.24)--(1.4,-2.69)--(1.67,-1.47)--(2.25,-2.64)--(3.02,-2)--(3.25,-3.19)--(2.25,-2.64)--(1.98,-3.78)--(1.4,-2.69)--(-0.19,-3.14)--(1.98,-3.78)--(3.25,-3.19)--(4.07,-1.45);
        \draw (-2.41,0.14)--(-1.83,-0.92)--(-2.99,-1.69) (-1.32,-2.74)--(-2.09,-2.11)--(-0.5,-1.21)--(-1.83,-0.92)--(-0.93,-0.34)--(-0.5,-1.21)
        (-0.19,-3.14)--(-0.5,-1.21)--(-1.32,-2.74)--(-2.46,-3.17)--(-2.99,-1.69)--(-2.09,-2.11)--(-2.46,-3.17);
        \draw (1.98,-3.78) arc (-63.13:-151.52:2.48)
              (1.98,-3.78) arc (-40.78:-154.81:2.67);
		\node at (0,-1.45) {\small $\ccc$};
		\node at (-0.61,0) {\small $\ccc$};
		\node at (0,0.48) {\small $\ccc$};
		\node at (0.69,0) {\small $\ccc$};
		\node at (2.28,-1.56) {\small $\ccc$};
		\node at (1.75,-2.16) {\small $\ccc$};
		\node at (-1.14,-1.32) {\small $\ccc$};
		\node at (-0.61,-1.9) {\small $\ccc$};
		\node at (1.98,-3.78) {\small $\bullet$};
	\node[draw,shape=circle, inner sep=0.2] at (0.48,-1.63) {\small $1$};
	\node[draw,shape=circle, inner sep=0.2] at (0,-0.84) {\small $2$};
	\node[draw,shape=circle, inner sep=0.2] at (0.5,-0.52) {\small $3$};
	\node[draw,shape=circle, inner sep=0.2] at (1.01,-0.94) {\small $4$};
	\node[draw,shape=circle, inner sep=0.2] at (-0.5,-0.49) {\small $5$};
	\node[draw,shape=circle, inner sep=0.2] at (-1.27,0.11) {\small $6$};
	\node[draw,shape=circle, inner sep=0.2] at (-0.56,0.62) {\small $7$};
	\node[draw,shape=circle, inner sep=0.2] at (0,1.2) {\small $8$};
	\node[draw,shape=circle, inner sep=0.2] at (0.58,0.54) {\small $9$};
	\node[draw,shape=circle, inner sep=0.2] at (1.32,0.09) {\footnotesize $10$};
	\node[draw,shape=circle, inner sep=0.2] at (1.59,-0.57) {\footnotesize $11$};
	\node[draw,shape=circle, inner sep=0.2] at (2.3,-0.84) {\footnotesize $15$};
	\node[draw,shape=circle, inner sep=0.2] at (2.88,-1.56) {\footnotesize $17$};
	\node[draw,shape=circle, inner sep=0.2] at (2.35,-2.16) {\footnotesize $18$};
	\node[draw,shape=circle, inner sep=0.2] at (1.8,-2.72) {\footnotesize $16$};
	\node[draw,shape=circle, inner sep=0.2] at (1.17,-2.16) {\footnotesize $12$};
	\node[draw,shape=circle, inner sep=0.2] at (0,-2.29) {\footnotesize $13$};
	\node[draw,shape=circle, inner sep=0.2] at (0.55,-2.65) {\footnotesize $14$};
	\node[draw,shape=circle, inner sep=0.2] at (1.24,-3.13) {\footnotesize $19$};
	\node[draw,shape=circle, inner sep=0.2] at (2.54,-3.11) {\footnotesize $20$};
	\node[draw,shape=circle, inner sep=0.2] at (2.8,-2.59) {\footnotesize $21$};
	\node[draw,shape=circle, inner sep=0.2] at (3.36,-2.24) {\footnotesize $22$};
	\node[draw,shape=circle, inner sep=0.2] at (-1.03,-0.79) {\footnotesize $23$};
	\node[draw,shape=circle, inner sep=0.2] at (-1.72,-0.42) {\footnotesize $24$};
	\node[draw,shape=circle, inner sep=0.2] at (-1.8,-1.47) {\footnotesize $25$};
	\node[draw,shape=circle, inner sep=0.2] at (-1.35,-2.06) {\footnotesize $26$};
	\node[draw,shape=circle, inner sep=0.2] at (-0.69,-2.8) {\footnotesize $27$};
	\node[draw,shape=circle, inner sep=0.2] at (-2.49,-2.29) {\footnotesize $28$};
	\node[draw,shape=circle, inner sep=0.2] at (-1.96,-2.64) {\footnotesize $29$};
	\node[draw,shape=circle, inner sep=0.2] at (-1.22,-3.78) {\footnotesize $30$};
		\end{tikzpicture}
		\caption{$\thin\aaa_3\thin\aaa_2\thin\cdots=\aaa^5\ccc^3$.}
		\label{a5r3}
	\end{figure}

\nd Now $\thin\aaa_3\thin\aaa_2\thin\cdots$ and $\thin\aaa_3\thin\aaa_4\thin\cdots$ can only be $\aaa^3\bbb$, which determines $T_5,T_6,T_7$ in Figure \ref{a3ba3b}. Then we get $\thin\ccc\thin\aaa\thin\ccc\thin\cdots=\aaa^2\bbb\ccc^2,\aaa^2\bbb\ccc^3,$ $\aaa^5\ccc^k(2\le k\le5),\aaa^4\ccc^k(2\le k\le4),\aaa^3\ccc^3$.
	\begin{figure}[H]
		\centering
        \begin{tikzpicture}[>=latex,scale=0.7]
				\foreach \a in {0,1,2}
				\draw[rotate=120*\a]      
                (-0.6,0.35)--(0,1.39)--(0.6,0.35)--(-0.6,0.35)
                (0,1.39)--(1.56,0.9)--(1.2,-0.69);
		\node at (-0.65,-0.9) {\small $\ccc$};
		\node at (-0.6,0.93) {\small $\ccc$};
		\node at (0.6,0.93) {\small $\ccc$};
	\node[draw,shape=circle, inner sep=0.2] at (0,-1.24) {\small $1$};
	\node[draw,shape=circle, inner sep=0.2] at (-0.6,-0.32) {\small $2$};
	\node[draw,shape=circle, inner sep=0.2] at (0,0) {\small $3$};
	\node[draw,shape=circle, inner sep=0.2] at (0.6,-0.32) {\small $4$};
	\node[draw,shape=circle, inner sep=0.2] at (0,0.72) {\small $5$};
	\node[draw,shape=circle, inner sep=0.2] at (-1.05,0.61) {\small $6$};
	\node[draw,shape=circle, inner sep=0.2] at (1.05,0.61) {\small $7$};
		\end{tikzpicture}
		\caption{$\thin\aaa_3\thin\aaa_2\thin\cdots=\aaa^3\bbb$ and $\thin\aaa_3\thin\aaa_4\thin\cdots=\aaa^3\bbb$.}
		\label{a3ba3b}
	\end{figure}

\vspace{9pt}
\noindent{Subcase $\thin\ccc\thin\aaa\thin\ccc\thin\cdots=\aaa^2\bbb\ccc^2$, AVC $\subset\{\aaa^3\bbb,\aaa^2\bbb\ccc^2\}$}

We might as well take $\thin\ccc_6\thin\aaa_5\thin\ccc_7\thin\cdots=\thin\ccc_6\thin\aaa_5\thin\ccc_7\thin\bbb\thin\aaa\thin$ which determines $T_8,T_9$ in the first of Figure \ref{a3b4}. If $\thin\aaa_8\thin\bbb_6\thin\cdots=\thin\aaa_8\thin\bbb_6\thin\aaa\thin\cdots$, then $\thin\aaa_8\thin\bbb_6\thin\cdots=\thin\aaa_8\thin\bbb_6\thin\aaa\thin\aaa\thin$. But $\thin\aaa\thin\aaa_8\thin\ccc_9\thin\cdots$ appears, a contradiction. So $\thin\ccc_1\thin\aaa_2\thin\ccc_6\thin\cdots=\thin\ccc_1\thin\aaa_2\thin\ccc_6\thin\bbb\thin\aaa\thin$ determines $T_{10},T_{11}$. Similarly, $T_{12},T_{13}$ are determined. Then we can determine $T_{14},\cdots,T_{32}$ and the 3D picture is the 20th of Figure \ref{Fig3}.\label{20,12.2}
	\begin{figure}[htp]
		\centering
		\begin{tikzpicture}[>=latex,scale=0.65]
				\foreach \a in {0,1,2}
				\draw[rotate=120*\a]
				(-1.26,2.41)--(0.01,3.89)--(1.83,4.53)--(2.07,2.62)--(1.56,0.9)--(0,1.39)--(-1.26,2.41)--(0.56,2.46)--(2.07,2.62)--(0.01,3.89)--(0.56,2.46)--(2.07,2.62)--(0.01,3.89)--(0.56,2.46)--(1.56,0.9)--(2.72,-0.11)--(1.2,-0.69)--(1.56,0.9)
                (-0.6,0.35)--(0,1.39)--(0.6,0.35)--(-0.6,0.35)
                (1.83,4.53) arc (80.47:175.42:5.74)
                (1.83,4.53) arc (51.48:-29.91:3.62);
		\node at (-0.65,-0.9) {\small $\ccc$};
		\node at (-0.6,0.93) {\small $\ccc$};
		\node at (0.6,0.93) {\small $\ccc$};
		\node at (-0.63,2.2) {\small $\ccc$};
		\node at (-1.58,-1.66) {\small $\ccc$};
		\node at (2.23,-0.55) {\small $\ccc$};
		\node at (-1.74,1.17) {\small $\ccc$};
		\node at (-0.1,-2.17) {\small $\ccc$};
		\node at (1.91,0.95) {\small $\ccc$};
		\node at (-1.39,2.7) {\small $\ccc$};
		\node at (-1.55,-2.55) {\small $\ccc$};
		\node at (3,-0.18) {\small $\ccc$};
	\node[draw,shape=circle, inner sep=0.2] at (0,-1.24) {\small $1$};
	\node[draw,shape=circle, inner sep=0.2] at (-0.6,-0.32) {\small $2$};
	\node[draw,shape=circle, inner sep=0.2] at (0,0) {\small $3$};
	\node[draw,shape=circle, inner sep=0.2] at (0.6,-0.32) {\small $4$};
	\node[draw,shape=circle, inner sep=0.2] at (0,0.72) {\small $5$};
	\node[draw,shape=circle, inner sep=0.2] at (-1.05,0.61) {\small $6$};
	\node[draw,shape=circle, inner sep=0.2] at (1.05,0.61) {\small $7$};
	\node[draw,shape=circle, inner sep=0.2] at (-0.94,1.56) {\small $8$};
	\node[draw,shape=circle, inner sep=0.2] at (0.3,1.93) {\small $9$};
	\node[draw,shape=circle, inner sep=0.2] at (-0.89,-1.56) {\footnotesize $10$};
	\node[draw,shape=circle, inner sep=0.2] at (-1.82,-0.76) {\footnotesize $11$};
	\node[draw,shape=circle, inner sep=0.2] at (1.81,0.08) {\footnotesize $12$};
	\node[draw,shape=circle, inner sep=0.2] at (1.54,-1.13) {\footnotesize $13$};
	\node[draw,shape=circle, inner sep=0.2] at (-2.4,0.21) {\footnotesize $14$};
	\node[draw,shape=circle, inner sep=0.2] at (-2.53,1.4) {\footnotesize $15$};
	\node[draw,shape=circle, inner sep=0.2] at (1.01,-2.22) {\footnotesize $16$};
	\node[draw,shape=circle, inner sep=0.2] at (0,-2.8) {\footnotesize $17$};
	\node[draw,shape=circle, inner sep=0.2] at (1.41,1.99) {\footnotesize $18$};
	\node[draw,shape=circle, inner sep=0.2] at (2.5,1.51) {\footnotesize $19$};
	\node[draw,shape=circle, inner sep=0.2] at (-0.26,2.97) {\footnotesize $20$};
	\node[draw,shape=circle, inner sep=0.2] at (-1.63,3.23) {\footnotesize $21$};
	\node[draw,shape=circle, inner sep=0.2] at (-2.35,-1.61) {\footnotesize $22$};
	\node[draw,shape=circle, inner sep=0.2] at (-1.84,-3.09) {\footnotesize $23$};
	\node[draw,shape=circle, inner sep=0.2] at (2.65,-1.35) {\footnotesize $24$};
	\node[draw,shape=circle, inner sep=0.2] at (3.69,-0.18) {\footnotesize $25$};
	\node[draw,shape=circle, inner sep=0.2] at (0.83,2.97) {\footnotesize $26$};
	\node[draw,shape=circle, inner sep=0.2] at (1.41,3.76) {\footnotesize $27$};
	\node[draw,shape=circle, inner sep=0.2] at (-3.01,-0.66) {\footnotesize $28$};
	\node[draw,shape=circle, inner sep=0.2] at (-3.88,-0.63) {\footnotesize $29$};
	\node[draw,shape=circle, inner sep=0.2] at (2.1,-2.17) {\footnotesize $30$};
	\node[draw,shape=circle, inner sep=0.2] at (2.58,-3.04) {\footnotesize $31$};
	\node[draw,shape=circle, inner sep=0.2] at (4.35,2.91) {\footnotesize $32$};
		\end{tikzpicture}\hspace{10pt}
		\begin{tikzpicture}[>=latex,scale=0.7]
				\foreach \a in {0,1,2}
				\draw[rotate=120*\a]
				(1.29,3.47)--(0.84,2.41)--(0,1.39)--(0.6,0.35)--(1.2,-0.69)--(2.51,-0.48)--(3.65,-0.62)
                (0,1.39)--(1.56,0.9)--(1.2,-0.69)
                (0.84,2.41)--(1.56,0.9)--(2.51,-0.48)
                (-0.6,0.35)--(0.6,0.35)
                (0.84,2.41) arc (75.58:-15.58:2.34);
		\node at (-0.65,-0.9) {\small $\ccc$};
		\node at (-0.6,0.93) {\small $\ccc$};
		\node at (0.6,0.93) {\small $\ccc$};
		\node at (0,1.88) {\small $\ccc$};
		\node at (-1.58,-0.9) {\small $\ccc$};
		\node at (1.58,-0.9) {\small $\ccc$};
	\node[draw,shape=circle, inner sep=0.2] at (0,-1.24) {\small $1$};
	\node[draw,shape=circle, inner sep=0.2] at (-0.6,-0.32) {\small $2$};
	\node[draw,shape=circle, inner sep=0.2] at (0,0) {\small $3$};
	\node[draw,shape=circle, inner sep=0.2] at (0.6,-0.32) {\small $4$};
	\node[draw,shape=circle, inner sep=0.2] at (0,0.72) {\small $5$};
	\node[draw,shape=circle, inner sep=0.2] at (-1.05,0.61) {\small $6$};
	\node[draw,shape=circle, inner sep=0.2] at (1.05,0.61) {\small $7$};
	\node[draw,shape=circle, inner sep=0.2] at (-0.76,1.59) {\small $8$};
	\node[draw,shape=circle, inner sep=0.2] at (0,2.99) {\small $9$};
	\node[draw,shape=circle, inner sep=0.2] at (0.76,1.59) {\footnotesize $10$};
	\node[draw,shape=circle, inner sep=0.2] at (-0.94,-1.4) {\footnotesize $11$};
	\node[draw,shape=circle, inner sep=0.2] at (-2.37,-1.35) {\footnotesize $12$};
	\node[draw,shape=circle, inner sep=0.2] at (-1.79,-0.13) {\footnotesize $13$};
	\node[draw,shape=circle, inner sep=0.2] at (1.79,-0.13) {\footnotesize $14$};
	\node[draw,shape=circle, inner sep=0.2] at (2.37,-1.35) {\footnotesize $15$};
	\node[draw,shape=circle, inner sep=0.2] at (0.94,-1.4) {\footnotesize $16$};
	\node[draw,shape=circle, inner sep=0.2] at (-1.9,1.19) {\footnotesize $17$};
	\node[draw,shape=circle, inner sep=0.2] at (-2.66,1.72) {\footnotesize $18$};
	\node[draw,shape=circle, inner sep=0.2] at (0,-2.22) {\footnotesize $19$};
	\node[draw,shape=circle, inner sep=0.2] at (0,-3.2) {\footnotesize $20$};
	\node[draw,shape=circle, inner sep=0.2] at (1.9,1.19) {\footnotesize $21$};
	\node[draw,shape=circle, inner sep=0.2] at (2.66,1.72) {\footnotesize $22$};
		\end{tikzpicture}
		\caption{$T(20a^3,12a^4;12\aaa^3\bbb,12\aaa^2\bbb\ccc^2)$ and $T(16a^3,6a^4;12\aaa^3\bbb,4\aaa^3\ccc^3)$.}
		\label{a3b4}
	\end{figure}

\vspace{9pt}
\noindent{Subcase $\thin\ccc\thin\aaa\thin\ccc\thin\cdots=\aaa^k\ccc^k(k=3,4,5)$, AVC $\subset\{\aaa^3\bbb,\aaa^k\ccc^k\}$}

The vertex $\aaa^k\ccc^k$ can only be $\thin\ccc\thin\aaa\thin\ccc\thin\aaa\thin\cdots\thin\ccc\thin\aaa\thin$. When $k=3$, we can determine $T_8,\cdots,T_{22}$ in the second of Figure \ref{a3b4} and the 3D picture is the 24th of Figure \ref{Fig3}. When $k=4,5$, we have the same tilings and the 3D pictures respectively are the 25th and 26th of Figure \ref{Fig3}.\label{16,6}

\vspace{9pt}
\noindent{Subcase $\thin\ccc\thin\aaa\thin\ccc\thin\cdots=\aaa^2\bbb\ccc^3$, AVC $\subset\{\aaa^3\bbb,\aaa^2\bbb\ccc^3\}$}

We might as well take $\thin\ccc_6\thin\aaa_5\thin\ccc_7\thin\cdots=\thin\ccc_6\thin\aaa_5\thin\ccc_7\thin\bbb\thin\ccc\thin\aaa\thin$ which determines $T_8,T_9,T_{10}$ in the first of Figure \ref{a2br3}. Then $T_{11},\cdots,T_{22}$ are determined. But $\thin\aaa_{19}\thin\ccc_{6}\thin\aaa_2\thin\ccc_1\thin\aaa_{21}\thin\cdots$ appears, contradicting the AVC.

	\begin{figure}[htp]
		\centering
        \begin{tikzpicture}[>=latex,scale=0.7]
				\foreach \a in {0,1,2}
				\draw[rotate=120*\a]      
                (-0.6,0.35)--(0,1.39)--(0.6,0.35)--(-0.6,0.35)
                (0,1.39)--(1.56,0.9)--(1.2,-0.69); 
        \draw (-1.56,0.9)--(-2.93,2.15)--(-1.29,3.84)--(-0.02,2.97)--(0,1.39)--(-1.29,2.44)--(-1.56,0.9)
              (-2.93,2.15)--(-1.29,2.44)--(-1.29,3.84) (1.2,2.44)--(1.56,0.9)--(2.55,2.15)
              (-0.02,2.97)--(1.2,2.44)--(2.55,2.15)--(3.98,2.17)--(3.13,1.03)--(1.56,0.9)--(2.65,-0.39)--(3.13,1.03) (2.65,-0.39)--(2.34,-1.82)--(1.2,-0.69)--(1.07,-2.3)--(2.36,-3.25)--(2.34,-1.82)
              (2.36,-3.25)--(-0.15,-3.28)--(1.07,-2.3)--(0,-1.8)--(-0.15,-3.28);
        \draw (-1.2,-0.69) arc (-107.29:-190.04:2.51)
               (-1.2,-0.69) arc (171.55:232.65:2.75);
		\node at (-0.65,-0.9) {\small $\ccc$};
		\node at (-0.6,0.93) {\small $\ccc$};
		\node at (0.6,0.93) {\small $\ccc$};
		\node at (-0.2,1.85) {\small $\ccc$};
		\node at (0.3,1.64) {\small $\bbb$};
		\node at (2.1,1.19) {\small $\ccc$};
		\node at (1.68,0.32) {\small $\ccc$};
		\node at (1.38,-1.28) {\small $\ccc$};
		\node at (-1.2,-0.69) {\small $\bullet$};
	\node[draw,shape=circle, inner sep=0.2] at (0,-1.24) {\small $1$};
	\node[draw,shape=circle, inner sep=0.2] at (-0.6,-0.32) {\small $2$};
	\node[draw,shape=circle, inner sep=0.2] at (0,0) {\small $3$};
	\node[draw,shape=circle, inner sep=0.2] at (0.6,-0.32) {\small $4$};
	\node[draw,shape=circle, inner sep=0.2] at (0,0.72) {\small $5$};
	\node[draw,shape=circle, inner sep=0.2] at (-1.05,0.61) {\small $6$};
	\node[draw,shape=circle, inner sep=0.2] at (1.05,0.61) {\small $7$};
	\node[draw,shape=circle, inner sep=0.2] at (-1,1.59) {\small $8$};
	\node[draw,shape=circle, inner sep=0.2] at (-0.57,2.7) {\small $9$};
	\node[draw,shape=circle, inner sep=0.2] at (0.72,1.96) {\footnotesize $10$};
	\node[draw,shape=circle, inner sep=0.2] at (1.78,1.8) {\footnotesize $11$};
	\node[draw,shape=circle, inner sep=0.2] at (2.79,1.51) {\footnotesize $12$};
	\node[draw,shape=circle, inner sep=0.2] at (2.47,0.53) {\footnotesize $13$};
	\node[draw,shape=circle, inner sep=0.2] at (1.94,-0.42) {\footnotesize $14$};
	\node[draw,shape=circle, inner sep=0.2] at (1.75,-1.88) {\footnotesize $15$};
	\node[draw,shape=circle, inner sep=0.2] at (0.72,-1.64) {\footnotesize $16$};
	\node[draw,shape=circle, inner sep=0.2] at (-1.9,1.8) {\footnotesize $17$};
	\node[draw,shape=circle, inner sep=0.2] at (-1.71,2.78) {\footnotesize $18$};
	\node[draw,shape=circle, inner sep=0.2] at (-1.98,0.56) {\footnotesize $19$};
	\node[draw,shape=circle, inner sep=0.2] at (0.33,-2.41) {\footnotesize $20$};
	\node[draw,shape=circle, inner sep=0.2] at (-0.63,-1.88) {\footnotesize $21$};
	\node[draw,shape=circle, inner sep=0.2] at (1.12,-2.83) {\footnotesize $22$};
		\end{tikzpicture}\hspace{60pt}
        \begin{tikzpicture}[>=latex,scale=0.75]
        \draw (-0.5,-0.87)--(0.5,0.87)--(1,0)--(1.37,-1.37)--(0,-1.73)--(-0.5,-0.87)
              (0,-1.73)--(1,0)--(0,0)--(0.5,-0.87)--(1.37,-1.37);
		\node at (0,-0.45) {\small $\ccc$};
		\node at (1,0) {\small $\bullet$};
	\node[draw,shape=circle, inner sep=0.2] at (0,-1) {\small $1$};
	\node[draw,shape=circle, inner sep=0.2] at (0.5,-0.26) {\small $2$};
	\node[draw,shape=circle, inner sep=0.2] at (0.5,0.26) {\small $3$};
	\node[draw,shape=circle, inner sep=0.2] at (0.9,-0.74) {\small $4$};
	\node[draw,shape=circle, inner sep=0.2] at (0.58,-1.3) {\small $5$};
		\end{tikzpicture}
		\caption{$\thin\ccc\thin\aaa\thin\ccc\thin\cdots$ is others.}
		\label{a2br3}
	\end{figure}

\vspace{9pt}
\noindent{Subcase $\thin\ccc\thin\aaa\thin\ccc\thin\cdots$ is others}

There must be $\thin\ccc\thin\aaa\thin\aaa\thin\cdots$, which determines $T_{1},T_2,T_3$ in the second of Figure \ref{a2br3}. Then we can determine $T_4,T_5$. But $\thin\aaa_4\thin\aaa_2\thin\cdots\not=\thin\aaa_4\thin\aaa_2\thin\bbb\thin\aaa\thin$, a contradiction.

	\subsection{The vertex types of degree 5}

Similar to the proof of Lemma \ref{deg345}, $\bbb^2\cdots$ can never be a vertex. So $\thin\ccc\thin\ccc\thin\cdots$ can never be a vertex by Lemma \ref{AAD}. This implies that $\ccc^5,\bbb\ccc^4,\aaa\ccc^4,\aaa^2\ccc^3$ or $\aaa\bbb\ccc^3$ can never be a vertex.

\subsubsection*{Case $\aaa^5$}

In Figure \ref{a5}, $\aaa^5$ determines $T_1,\cdots,T_5$. By similar proof of Lemma \ref{deg345}, $\thin\aaa\thin\aaa\thin\bbb\thin\cdots$ is not a vertex. Then $\thin\aaa_1\thin\aaa_2\thin\cdots=\thin\aaa\thin\aaa_1\thin\aaa_2\thin\aaa\thin\cdots$ $=\aaa^5$ determines $T_6,T_7,T_8$. Similarly, $T_9,\cdots,T_{20}$ are determined. But it is a monohedral tiling, a contradiction.
	\begin{figure}[H]
		\centering
		\begin{tikzpicture}[>=latex,scale=0.75]
				\foreach \a in {0,1,2,3,4}
				\draw[rotate=72*\a]
				(0,1)--(1.01,1.41)--(0.95,0.31)--(0,0) (0,1)--(0.95,0.31) (1.6,2.2)--(1.01,1.41)
                (1.01,1.41) arc (46:134.7:1.46);
	\node[draw,shape=circle, inner sep=0.2] at (0.35,0.4) {\small $1$};
	\node[draw,shape=circle, inner sep=0.2] at (0.5,-0.2) {\small $2$};
	\node[draw,shape=circle, inner sep=0.2] at (0,-0.5) {\small $3$};
	\node[draw,shape=circle, inner sep=0.2] at (-0.5,-0.2) {\small $4$};
	\node[draw,shape=circle, inner sep=0.2] at (-0.35,0.4) {\small $5$};
	\node[draw,shape=circle, inner sep=0.2] at (0.65,0.86) {\small $6$};
	\node[draw,shape=circle, inner sep=0.2] at (1.06,-0.4) {\small $7$};
	\node[draw,shape=circle, inner sep=0.2] at (1.4,0.36) {\small $8$};
	\node[draw,shape=circle, inner sep=0.2] at (0,1.43) {\small $9$};
	\node[draw,shape=circle, inner sep=0.2] at (-0.66,0.89) {\footnotesize $10$};
	\node[draw,shape=circle, inner sep=0.2] at (-1.38,0.42) {\footnotesize $11$};
	\node[draw,shape=circle, inner sep=0.2] at (-1.06,-0.37) {\footnotesize $12$};
	\node[draw,shape=circle, inner sep=0.2] at (-0.82,-1.22) {\footnotesize $13$};
	\node[draw,shape=circle, inner sep=0.2] at (0,-1.22) {\footnotesize $14$};
	\node[draw,shape=circle, inner sep=0.2] at (0.85,-1.19) {\footnotesize $15$};
	\node[draw,shape=circle, inner sep=0.2] at (0,2.35) {\footnotesize $16$};
	\node[draw,shape=circle, inner sep=0.2] at (-2.24,0.73) {\footnotesize $17$};
	\node[draw,shape=circle, inner sep=0.2] at (-1.38,-1.91) {\footnotesize $18$};
	\node[draw,shape=circle, inner sep=0.2] at (1.38,-1.91) {\footnotesize $19$};
	\node[draw,shape=circle, inner sep=0.2] at (2.24,0.73) {\footnotesize $20$};
		\end{tikzpicture}
		\caption{$T(20a^3;12\aaa^5).$}
		\label{a5}
	\end{figure}
\subsubsection*{Case $\aaa^3\ccc^2$}

The vertex $\aaa^3\ccc^2=\thin\aaa\thin\aaa\thin\ccc\thin\aaa\thin\ccc\thin$ determines $T_1,\cdots,T_5$ in Figure \ref{a3r2}. Similar to the proof of Lemma \ref{deg345}, $\aaa\bbb\cdots=\aaa^4\bbb$ determines $T_6,\cdots,T_{16}$. Then we have the AVC $\subset \{\aaa^3\ccc^2,\aaa^4\bbb\}$, and $T_{17},\cdots,T_{34}$ are determined. But $\aaa^2\bbb\ccc\cdots$ appears, contradicting the AVC.

	\begin{figure}[H]
		\centering
		\begin{tikzpicture}[>=latex,scale=0.6]
				\draw(-1.69,4.47)--(0.98,3.33)--(2.04,1.77)--(2.66,0)--(2.04,-1.77)--(0.11,-2.28) (-2.80,-2.28)--(-5.03,-2.33)--(-5.75,0)--(-5.03,2.33)--(-3.76,3.78)--(-1.69,4.47)
                (0.98,3.33)--(-1.11,3.33)--(-2.8,2.28)--(-3.59,1.16)--(-4.3,0)--(-5.03,-2.33)
                (-4.3,0)--(-5.03,2.33)--(-2.8,2.28)--(-1.69,4.47)--(-1.11,3.33)
				(-3.59,1.16)--(-5.03,2.33)
                (-1.11,3.33)--(-1.32,2.28)--(-2.17,1.27)--(-2.77,0)--(-3.59,-1.16)--(-5.03,-2.33)
                (-2.77,0)--(-3.59,1.16)--(-2.17,1.27)
                (0.11,2.28)--(-0.71,1.16)--(-1.53,0)--(-2.17,-1.27)--(-3.59,-1.16)--(-4.3,0)--(-2.77,0)--(-2.17,-1.27)
                (-2.8,2.28)--(-1.32,2.28)--(0.11,2.28)--(2.04,1.77)--(1.08,0.93)--(0,0)--(-0.71,-1.16)--(-1.32,-2.28)
                (-1.11,3.33)--(0.11,2.28)--(1.08,0.93)--(2.66,0) (1.08,0.93)--(1.08,-0.93)--(2.66,0)
                (-1.32,2.28)--(-0.71,1.16)--(0,0)--(1.08,-0.93)--(2.04,-1.77) (0.11,-2.28)--(1.08,-0.93)
                (0.11,-2.28)--(-0.71,-1.16)--(-1.53,0)--(-2.17,1.27)--(-0.71,1.16)
                (-3.59,-1.16)--(-2.8,-2.28)--(-1.32,-2.28)--(-2.17,-1.27)--(-0.71,-1.16) (-1.53,0)--(0,0)
                (-1.32,-2.28)--(0.11,-2.28);
		\node at (0,0.3) {\small $\ccc$};
		\node at (0,-0.4) {\small $\ccc$};
		\node at (-1.8,0) {\small $\bbb$};
		\node at (-2.2,1.5) {\small $\ccc$};
		\node at (-2.2,-1.55) {\small $\ccc$};
		\node at (0.2,2.5) {\small $\ccc$};
		\node at (-2.95,2.5) {\small $\ccc$};
		\node at (-4.6,0) {\small $\bbb$};
		\node at (-5.03,2.33) {\small $\bullet$};
	\node[draw,shape=circle, inner sep=0.2] at (0.7,0) {\small $1$};
	\node[draw,shape=circle, inner sep=0.2] at (0,-1) {\small $2$};
	\node[draw,shape=circle, inner sep=0.2] at (-0.7,-0.5) {\small $3$};
	\node[draw,shape=circle, inner sep=0.2] at (-0.7,0.5) {\small $4$};
	\node[draw,shape=circle, inner sep=0.2] at (0,1) {\small $5$};
	\node[draw,shape=circle, inner sep=0.2] at (-1.5,0.7) {\small $6$};
	\node[draw,shape=circle, inner sep=0.2] at (-1.4,1.6) {\small $7$};
	\node[draw,shape=circle, inner sep=0.2] at (-0.7,1.8) {\small $8$};
	\node[draw,shape=circle, inner sep=0.2] at (-1.5,-0.7) {\small $9$};
	\node[draw,shape=circle, inner sep=0.2] at (-1.4,-1.6) {\footnotesize $10$};
	\node[draw,shape=circle, inner sep=0.2] at (-0.7,-1.85) {\footnotesize $11$};
	\node[draw,shape=circle, inner sep=0.2] at (1.5,0) {\footnotesize $12$};
	\node[draw,shape=circle, inner sep=0.2] at (1.8,-1) {\footnotesize $13$};
	\node[draw,shape=circle, inner sep=0.2] at (1.2,-1.5) {\footnotesize $14$};;
	\node[draw,shape=circle, inner sep=0.2] at (1.8,1) {\footnotesize $15$};
	\node[draw,shape=circle, inner sep=0.2] at (1.2,1.5) {\footnotesize $16$};
	\node[draw,shape=circle, inner sep=0.2] at (-2.3,0) {\footnotesize $17$};
	\node[draw,shape=circle, inner sep=0.2] at (-2.8,0.8) {\footnotesize $18$};
	\node[draw,shape=circle, inner sep=0.2] at (-2.65,1.8) {\footnotesize $19$};
	\node[draw,shape=circle, inner sep=0.2] at (-2.8,-0.8) {\footnotesize $20$};
	\node[draw,shape=circle, inner sep=0.2] at (-2.65,-1.8) {\footnotesize $21$};
	\node[draw,shape=circle, inner sep=0.2] at (-1.6,2.6) {\footnotesize $22$};
	\node[draw,shape=circle, inner sep=0.2] at (-0.8,2.6) {\footnotesize $23$};
	\node[draw,shape=circle, inner sep=0.2] at (0.8,2.7) {\footnotesize $24$};
	\node[draw,shape=circle, inner sep=0.2] at (-0.8,3.7) {\footnotesize $25$};
	\node[draw,shape=circle, inner sep=0.2] at (-1.8,3.5) {\footnotesize $26$};
	\node[draw,shape=circle, inner sep=0.2] at (-3.55,0.35) {\footnotesize $27$};
	\node[draw,shape=circle, inner sep=0.2] at (-3.55,-0.35) {\footnotesize $28$};
	\node[draw,shape=circle, inner sep=0.2] at (-4.2,1) {\footnotesize $29$};
	\node[draw,shape=circle, inner sep=0.2] at (-3.75,1.8) {\footnotesize $30$};
	\node[draw,shape=circle, inner sep=0.2] at (-3.5,3.1) {\footnotesize $31$};
	\node[draw,shape=circle, inner sep=0.2] at (-4.2,-1) {\footnotesize $32$};
	\node[draw,shape=circle, inner sep=0.2] at (-3.75,-1.8) {\footnotesize $33$};
	\node[draw,shape=circle, inner sep=0.2] at (-5.2,0) {\footnotesize $34$};
		\end{tikzpicture}
		\caption{$\aaa^3\ccc^2$ is a vertex.}
		\label{a3r2}
	\end{figure}

\subsubsection*{Case $\aaa^2\bbb\ccc^2$}
	\begin{figure}[H]
		\centering
		\begin{tikzpicture}[>=latex,scale=0.8]
				\draw (0,-1)--(-1.85,0)--(0,1)--(1.85,0)--(0,-1)--(0,0)
                      (-0.93,0.5)--(0.93,-0.5) (0.93,0.5)--(-0.93,-0.5);
		\node at (-0.58,0) {\small $\ccc$};
		\node at (0.58,0) {\small $\ccc$};
		\node at (-0.42,0.5) {\small $\ccc$};
		\node at (0.93,-0.5) {\small $\bullet$};
	\node[draw,shape=circle, inner sep=0.2] at (0.32,-0.5) {\small $1$};
	\node[draw,shape=circle, inner sep=0.2] at (-0.32,-0.5) {\small $2$};
	\node[draw,shape=circle, inner sep=0.2] at (-1,0) {\small $3$};
	\node[draw,shape=circle, inner sep=0.2] at (0,0.5) {\small $4$};
	\node[draw,shape=circle, inner sep=0.2] at (1,0) {\small $5$};
		\end{tikzpicture}
		\caption{$\aaa^2\bbb\ccc^2=\thin\aaa\thin\aaa\thin\ccc\thin\bbb\thin\ccc\thin$.}
		\label{a2br2}
	\end{figure}

By $2\bbb+2\ccc>2\pi$ and $3\aaa+\bbb+\ccc>2\pi$, we have $\bbb>2\aaa$ and $\aaa>\ccc$. Similar to the proof of Lemma \ref{deg345}, we get $\aaa^2\bbb\cdots=\aaa^2\bbb\ccc^2$. Then $\thin\aaa\thin\aaa\thin\bbb\thin\cdots$ is not a vertex. When $\thin\ccc\thin\aaa\thin\ccc\thin\cdots$ appears, we get two $\thin\aaa\thin\bbb\thin\cdots$ which are not vertices. So $\aaa^2\bbb\ccc^2=\thin\aaa\thin\aaa\thin\ccc\thin\bbb\thin\ccc\thin$ determines $T_1,\cdots,T_5$ in Figure \ref{a2br2}. But $\thin\aaa_1\thin\bbb_5\thin\cdots$ is not a vertex.

\subsubsection*{Case $\aaa^4\ccc$}
By $2\aaa<\pi$, we have $2\aaa+\ccc>\pi$. The AAD $\thin\ccc^\bbb\thin^\aaa\aaa\thin\cdots$ of $\aaa^4\ccc$ gives $\thin\aaa\thin\bbb\thin\cdots$. Similar to the proof of Lemma \ref{deg345}, $\thin\aaa\thin\bbb\thin\cdots$ is not a vertex.

\subsubsection*{Case $\aaa^4\bbb$}
By $3\aaa+\bbb+\ccc>2\pi$, we have $\ccc>\aaa$. Similar to the proof of Lemma \ref{deg345}, $\aaa\bbb\ccc\cdots$ and $\aaa^2\ccc\cdots$ can never be vertices. The AAD $\thin\bbb^\ccc\thin^\aaa\aaa\thin\cdots$ of $\aaa^4\bbb$ gives $\thin\aaa\thin\ccc\thin\cdots$, which is not a vertex.

\newpage
		\section*{Appendix: Exact geometric data}

\begin{table}[htp]
	\begin{center}
	
		\bgroup
		\def\arraystretch{1.5}
		  \begin{tabular}[htp]{  c|c|c }

				$(f_\tri ,f_\dia )$ & $\aaa,\bbb,\ccc;\,\cos a = \cot \alpha \cot \frac{\alpha}{2} = \cot \frac{\beta}{2}\cot \frac{\gamma}{2}$ & all vertices  \\ \hline
				
				$(2,3)$ &  \parbox[c][2cm][c]{8cm}{$\aaa\in[\arccos(\frac{1}{8}),\pi),\ccc=2\pi-\aaa-\bbb, \\ \bbb=2\arccot \left(-\frac{\sqrt{-8\cos \aaa+1}+2\cos \aaa-1}{2\sin \aaa}\right)$. \\ When $\aaa=\arccos(\frac{1}{8})$, we have $\bbb=\ccc=\arccos(-\frac{3}{4})$.}
				 & $6\aaa\bbb\ccc$  \\
				\hline
				$(8,6)$& \parbox[c][2cm][c]{8cm}{$\aaa\in[\arccos(\frac{1}{3}),\frac{1}{2}\pi), \ccc=2\pi-2\aaa-\bbb, \\ \bbb=2\arccot \left(-\frac{\sin \aaa\sqrt{\frac{3\cos \aaa-1}{\cos \aaa-1}}+\cos(2\aaa)+\cos \aaa}{\sin(2\aaa)}\right)$. \\ When $\aaa=\arccos(\frac{1}{3})$, we have $\bbb=\ccc=\arccos(-\frac{1}{3})$.} 
                 &  $12\aaa^2\bbb\ccc$  \\
				\hline
				
		\end{tabular}
		\egroup
			\caption{Two 1-parameter families of protosets.}\label{tab-1}
	\end{center}
\end{table}

\begin{table}[htp]
	\begin{center}
	
		\bgroup
		\def\arraystretch{1.5}
		  \begin{tabular}[htp]{  c|c|c }

				$(f_\tri ,f_\dia )$ & $\aaa,\bbb,\ccc;\,\cos a = \cot \alpha \cot \frac{\alpha}{2} = \cot \frac{\beta}{2}\cot \frac{\gamma}{2}$ & all vertices  \\ \hline
				
				$(2,6n-3),n\ge3$ &  \parbox[c][2.2cm][c]{6cm}{$\ccc=2\pi-2\bbb,\aaa=(1-n)2\pi+(2n-1)\bbb, \\ 
\cos((2n-1)\bbb) = \frac{\cos \bbb}{2 \cos \bbb -1}$. \\ $\lim\limits_{n\to+\infty}\beta=\pi$, $\lim\limits_{n\to+\infty}\ccc=0$, \\ $\lim\limits_{n\to+\infty}a=\frac{\pi}{3}$,$\lim\limits_{n\to+\infty}\aaa=\arccos\frac{1}{3}$.}
				 & $6\aaa\bbb\ccc^n,(6n-6)\bbb^2\ccc$  \\
				\hline
				
		\end{tabular}
		\egroup
			\caption{A sequence of protosets.}\label{tab-2}
	\end{center}
\end{table}

\begin{table}[htp]
	\begin{center}
	
		\bgroup
		\def\arraystretch{1.5}
		  \begin{tabular}[htp]{  c|c|c }

				$(f_\tri ,f_\dia )$ & $\aaa,\bbb,\ccc;\,\cos a = \cot \alpha \cot \frac{\alpha}{2} = \cot \frac{\beta}{2}\cot \frac{\gamma}{2}$ & all vertices  \\ \hline
				
				$(8,2)$ &  \parbox[c][1.6cm][c]{8cm}{$\aaa=\frac{2\pi-\bbb}{3},\ccc=\bbb,\\ \bbb=3\arccos\left(\frac{1-\sqrt{2}+\sqrt{3+6\sqrt{2}}}{4}\right)$}
				 & $8\aaa^3\bbb$  \\
				\hline
				$(32,6)$ & \parbox[c][2.0cm][c]{8cm}{$\aaa=\frac{2\pi-\bbb}{4},\ccc=\bbb,\\ \bbb=4\arcsin\left(\frac{\left(19+3\sqrt{33}\right)^{\frac{2}{3}}-2\left(19+3\sqrt{33}\right)^{\frac{1}{3}}+4}{6\left(19+3\sqrt{33}\right)^{\frac{1}{3}}}\right)$} 
				&  $24\aaa^4\bbb$\\
				\hline
				$(8,18)$& \parbox[c][1.3cm][c]{8cm}{$\aaa=2\pi-3\bbb,\ccc=\bbb,\\ \bbb=2\arctan\sqrt{7-4\sqrt{2}}$} &  $24\aaa\bbb^3$  \\
				\hline
(4,4) & \parbox[c][1.3cm][c]{8cm}{$\aaa=2\pi-2\bbb,\ccc=2\bbb-\pi, \\ \bbb=\arccos(\frac{1-\sqrt{17}}{4})$} & $4\aaa\bbb^2,4\aaa^2\ccc^2$\\
\hline  
(8,3) & \parbox[c][1.6cm][c]{8cm}{$\aaa=\frac{2\pi-\bbb}{3},\ccc=\frac{\bbb+\pi}{3}, \\ \bbb=6\arctan\Big(\sqrt{3(9-4\sqrt{5})}\Big)$} & $6\aaa^3\bbb,3\aaa^2\ccc^2$\\
\hline  
(4,12) & \parbox[c][2.5cm][c]{8cm}{$\aaa=\frac{4}{3}\pi-2\ccc,\bbb=\frac{2}{3}\pi, \\ \ccc=\arccos\left(\frac{c^{\frac{1}{6}}b^{\frac{1}{4}}-b^{\frac{3}{4}}+\sqrt{\left(-b+9c^{\frac{1}{3}}\right)\sqrt{b}+8\sqrt{c}}}{4c^{\frac{1}{6}}b^{\frac{1}{4}}}\right) \\   c=17+2\sqrt{41}, b=c^{\frac{2}{3}}+3c^{\frac{1}{3}}+5$.} & $12\aaa\bbb\ccc^2,4\bbb^3$\\
\hline
(8,12) & \parbox[c][2cm][c]{8cm}{$\aaa=\pi-\bbb,\ccc=\frac{\pi}{2}, \\ \bbb=\arccos\left(-\frac{c^{2}+2c-14}{6c}\right) \\  c=\left(8+6\sqrt{78}\right)^{\frac{1}{3}}$.} & $6\aaa^2\bbb^2,12\aaa\bbb\ccc^2$\\
\hline 
(20,6) & \parbox[c][1.5cm][c]{8cm}{$\aaa=\frac{2\pi-\bbb}{3},\ccc=\frac{2\bbb-\pi}{3}, \\ \bbb=-6\arctan\Big(3+2\sqrt{3}-\sqrt{24+14\sqrt{3}}\Big)$} & $12\aaa^3\bbb,6\aaa^4\ccc^2$\\
\hline   
(8,24) & \parbox[c][2.2cm][c]{8cm}{$\ccc=\frac{\pi}{2},\bbb=\frac{3}{4}\pi-\frac{\aaa}{2}, \\ \aaa=\arcsin\left(\frac{\sqrt{6}}{12}\left(b+\sqrt{\frac{-bc^2+12\sqrt{6}c-12b}{bc}}\right)\right)
\\ c=\left(108+12\sqrt{69}\right)^{\frac{1}{3}},b=\sqrt{\frac{c^2+12}{c}}.$} & $24\aaa\bbb^2\ccc,6\ccc^4$\\
\hline
(20,12) & \parbox[c][0.8cm][c]{8cm}{$\aaa=\pi-\bbb,\ccc=\frac{3\bbb-\pi}{2}, \bbb=2\arcsin\left(\frac{1+\sqrt{33}}{8}\right)$} & $12\aaa^2\bbb^2,12\aaa^3\ccc^2$\\
\hline 

		\end{tabular}
		\egroup
	\end{center}
\end{table}

\begin{table}[htp]
	\begin{center}
	
		\bgroup
		\def\arraystretch{1.5}
		  \begin{tabular}[htp]{  c|c|c }

				$(f_\tri ,f_\dia )$ & $\aaa,\bbb,\ccc;\,\cos a = \cot \alpha \cot \frac{\alpha}{2} = \cot \frac{\beta}{2}\cot \frac{\gamma}{2}$ & all vertices  \\ \hline
(20,12) & \parbox[c][4cm][c]{8cm}{$\aaa=\frac{2\pi-\bbb}{3},\ccc=\frac{2\pi-\bbb}{6}, \\ \bbb=\arcsin\left(\frac{\left(-3^{\frac{4}{3}}2^{\frac{7}{12}}d^{\frac{1}{2}}-3\cdot12^{\frac{1}{3}}b^{\frac{3}{2}}+252b^{\frac{1}{2}}c^{\frac{1}{6}}\right)^{\frac{1}{2}}}{12b^{\frac{1}{4}}c^{\frac{1}{12}}}\right) \\ 
c=1597545+5291\sqrt{1689}, \\
b=\sqrt{\left(18c^2\right)^{\frac{1}{3}}+120\left(12c\right)^{\frac{1}{3}}+35592}, \\
d=-35592\cdot2^{\frac{1}{6}}b+240\cdot2^{\frac{5}{6}}b\left(3c\right)^{\frac{1}{3}}-\sqrt{2}b\left(9c^2\right)^{\frac{1}{3}} \\ {\quad\quad}+2340\cdot2^{\frac{1}{6}}\sqrt{c}.$} & $12\aaa^3\bbb,12\aaa^2\bbb\ccc^2$\\
\hline   
(20,24) & \parbox[c][1.5cm][c]{8cm}{$\aaa=4\bbb-2\pi,\ccc=4\pi-6\bbb,\cos\bbb\in(-\frac{1}{2},0), \\ 
\cos\bbb=\frac{RootOf(Z^6+Z^5-4Z^4-3Z^3+3Z^2-1)}{2}.$} & $24\aaa\bbb^2\ccc,12\aaa^3\ccc^2$\\
\hline 
(44,12) & \parbox[c][4cm][c]{8cm}{$\aaa=\frac{2\pi-\bbb}{3},\ccc=\frac{5\bbb-4\pi}{6}, \\ \bbb=\arcsin\left(\frac{\left(-6\sqrt{3}b^{\frac{3}{4}}-6d^{\frac{1}{2}}+306b^{\frac{1}{4}}c^{\frac{1}{2}}\right)^{\frac{1}{2}}}{12b^{\frac{1}{8}}c^{\frac{1}{4}}}\right) \\
c=\left(38969189+564000\sqrt{1473}\right)^{\frac{1}{3}},  \\
b=c^2-115c+101641, \\
d=\left(-3b-1035c\right)\sqrt{b}+36576\sqrt{3c^3}.$} & $24\aaa^3\bbb,12\aaa^5\ccc^2$\\
\hline  
(20,60) & \parbox[c][4cm][c]{8cm}{$\aaa=\frac{8}{5}\pi-2\bbb,\ccc=\frac{2}{5}\pi, \\ \bbb=\arccos\left(\frac{\sqrt{6}\left(bd+33\sqrt{3}c^{\frac{3}{2}}\right)^{\frac{1}{2}}+3c^{\frac{1}{2}}d^{\frac{1}{2}}-\sqrt{3}d^{\frac{3}{2}}}{24c^{\frac{1}{2}}d^{\frac{1}{2}}}\right) \\
c=\left(-774\sqrt{5}+10070+6\sqrt{85830-7890\sqrt{5}}\right)^{\frac{1}{3}}, \\
b=-c^2+24\sqrt{5}+43c-460,  \\
d=\sqrt{2c^2+43c-48\sqrt{5}+920}.$} & $60\aaa\bbb^2\ccc,12\ccc^5$\\
\hline 
(16,6) & \parbox[c][2cm][c]{8cm}{$\aaa=\frac{2\pi-\bbb}{3},\ccc=\frac{\bbb}{3},\bbb=\arcsin\left(\frac{\sqrt{bc}}{18c}\right) \\
c=\left(18559+3321\sqrt{-47}\right)^{\frac{1}{3}}, \\ b=\left(\sqrt{-3}-1\right)c^2+164c-952-952\sqrt{-3}.$} & $12\aaa^3\bbb,4\aaa^3\ccc^3$\\
\hline  
(32,12) & \parbox[c][1.5cm][c]{8cm}{$\aaa=\frac{2\pi-\bbb}{3},\ccc=\frac{2\bbb-\pi}{6}, \\ \bbb=\arcsin\left(\frac{5-\sqrt{7}}{8}\right)$} & $24\aaa^3\bbb,6\aaa^4\ccc^4$\\
\hline  
(80,30) & \parbox[c][2cm][c]{8cm}{$\aaa=\frac{2\pi-\bbb}{3},\ccc=\frac{5\bbb-4\pi}{15}, \\ \bbb=\arcsin\left(\frac{\left(-2\sqrt{5}-6\right)\sqrt{5-2\sqrt{5}}+\sqrt{130-10\sqrt{5}}}{20}\right)$} & $60\aaa^3\bbb,12\aaa^5\ccc^5$\\
\hline 

		\end{tabular}
		\egroup
	\end{center}
\end{table}

\begin{table}[htp]
	\begin{center}
	
		\bgroup
		\def\arraystretch{1.5}
		  \begin{tabular}[htp]{  c|c|c }

				$(f_\tri ,f_\dia )$ & $\aaa,\bbb,\ccc;\,\cos a = \cot \alpha \cot \frac{\alpha}{2} = \cot \frac{\beta}{2}\cot \frac{\gamma}{2}$ & all vertices  \\ \hline
(20,36) & \parbox[c][6cm][c]{8cm}{$\aaa=4\pi-6\bbb,\ccc=4\bbb-2\pi, \\ \bbb=\arcsin\left(\frac{\left(19516464-24642\sqrt{3}b+111\sqrt{6}\left(bd+e\right)^{\frac{1}{2}}\right)^{\frac{1}{2}}}{5328}\right)  \\ 
c=\left(108+12\sqrt{849}\right)^{\frac{1}{3}},  \\
b=\sqrt{\left(-6\sqrt{849}+54\right)c^2+1152c+15552},  \\
d=\left(5217\sqrt{3}+2997\sqrt{283}\right)c^2+3271392\sqrt{3}+ \\{\quad\quad} \left(-223776\sqrt{3}+10656\sqrt{283}\right)c, \\
e=\left(147852\sqrt{849}-1330668\right)c^2-28387584c \\ {\quad\quad} +766464768.$} & $36\aaa\bbb^2\ccc,12\aaa^2\ccc^3$\\
\hline 
(32,6) & \parbox[c][1.85cm][c]{8cm}{$\aaa=\frac{2\pi-\bbb}{3},\ccc=\frac{5\bbb-4\pi}{3},\bbb=\arcsin\left(\frac{\sqrt{-3bc}}{6c}\right) \\
c=\left(-756+84\sqrt{-3}\right)^{\frac{1}{3}}, \\ b=\left(1+\sqrt{-3}\right)c^2+15c+84-84\sqrt{-3}.$} & $12\aaa^3\bbb,12\aaa^5\ccc$\\
\hline      
(20,24) & \parbox[c][2.2cm][c]{8cm}{$\aaa=2\pi-2\bbb,\ccc=\frac{3}{2}\bbb-\pi, \\ \bbb=4\arccos\left(\frac{\sqrt{-3bc}}{12c}\right) \\ c=\left(28+84\sqrt{-3}\right)^{\frac{1}{3}}, \\ b=\left(1+\sqrt{-3}\right)c^2-32c+28-28\sqrt{-3}.$} & \tabincell{c}{$(12+k)\aaa\bbb^2,k\aaa^3\ccc^4$,\\ $(24-2k)\aaa^2\bbb\ccc^2$ \\ $0\le k <12$}\\
\hline  
				
		\end{tabular}
		\egroup
			\caption{20 sporadic protosets with some irrational angle.}\label{tab-3}
	\end{center}
\end{table}

\end{document}